\newif\ifpupbook 
\pupbookfalse      

\newif\ifcolor 
 \colorfalse 

\ifpupbook
\documentclass[7x10,11pt]{pupbook2015}
\else
\documentclass[11pt]{amsbook}
\usepackage{amsthm}
\fi

\usepackage{amsmath}
\usepackage{amsxtra}
\usepackage{amscd}
\usepackage[T1]{fontenc}
\usepackage{amsfonts}
\usepackage{amssymb}
\usepackage{eucal}
\usepackage[lowtilde]{url} 
\usepackage[all]{xypic}

\makeatletter
\renewenvironment{proof}[1][\proofname]{\par
  \pushQED{\qed}%
  \normalfont \topsep6\p@\@plus6\p@\relax
  \trivlist
  \itemindent0pt 
  \item[\hskip\labelsep
        \itshape
    #1\@addpunct{.}]\ignorespaces
}{%
  \popQED\endtrivlist\@endpefalse
}

\def\th@plain{%
  \let\thm@indent\noindent 
  \itshape 
}
\def\th@definition{%
  \let\thm@indent\noindent 
  \normalfont 
}
\def\th@remark{%
  \let\thm@indent\noindent 
  \normalfont 
}

\makeatother
\theoremstyle{plain}
\newtheorem{thm}[subsection]{Theorem}
\newtheorem{cor}[subsection]{Corollary}
\newtheorem{lem}[subsection]{Lemma}
\newtheorem*{claim}{Claim}
\newtheorem*{claim1}{Claim 1}
\newtheorem*{claim2}{Claim 2}
\newtheorem{prop}[subsection]{Proposition}

\newtheorem{lemdef}[subsection]{Lemma-Definition}

\theoremstyle{definition}

\newtheorem{defn}[subsection]{Definition}

\newtheorem{rem}[subsection]{Remark}

\ifpupbook\fi
\newtheorem{remark}[subsection]{Remark}
\newtheorem{question}[subsection]{Question}
\newtheorem{example}[subsection]{Example}

\newtheorem{warning}[subsection]{Warning}

\usepackage{yhmath}

\def\origwidehat{\mathaccent "0362\relax}
\newcommand{\doublewidehat}[2][9pt]{\origwidehat{\vphantom{\rule{1pt}{#1}}\smash{\origwidehat {#2}}}}
\DeclareSymbolFont{largesymbols}{OMX}{yhex}{m}{n}
\DeclareMathAccent{\widetilde}{\mathord}{largesymbols}{"65}

\newcommand{\thmref}[1]{Theorem~\ref{#1}}

\newcommand{\lemref}[1]{Lemma~\ref{#1}}
\newcommand{\defref}[1]{Definition~\ref{#1}}
\newcommand{\propref}[1]{Proposition~\ref{#1}}
\newcommand{\corref}[1]{Corollary~\ref{#1}}

\newcommand{\remref}[1]{Remark~\ref{#1}}
\newcommand{\exref}[1]{Example~\ref{#1}}

\newcommand{\nc}{\newcommand}
\nc{\renc}{\renewcommand}

\renc{\d}{{\delta}}
\nc{\Aa}{{\mathbb{A}}}
 \nc{\Gg}{{\mathbb{G}}}  
\nc{\Hh}{{\mathbb{H}}}
 \nc{\Nn}{{\mathbb{N}}}
\nc{\Pp}{{\mathbb{P}}}
\nc{\Rr}{{\mathbb{R}}}
\nc{\BV}{{\mathbb{V}}}
\nc{\BW}{{\mathbb{W}}}
\nc{\Zz}{{\mathbb{Z}}}
\nc{\Qq}{{\mathbb{Q}}}
\nc{\Ss}{{\mathbb{S}}}
\nc{\Cc}{{\mathbb{C}}}
\nc{\Ff}{{\mathbb{F}}}

\nc{\CA}{{\mathcal{A}}}
\nc{\CB}{{\mathcal{B}}}

\nc{\CE}{{\mathcal{E}}}
\nc{\CF}{{\mathcal{F}}}
\nc{\CG}{{\mathcal{G}}}
\nc{\CL}{{\mathcal{L}}}
\nc{\CC}{{\mathcal{C}}}
\nc{\CM}{{\mathcal{M}}}
\def\Mm{\CM}
\nc{\CN}{{\mathcal{N}}}
\nc{\Oo}{{\mathcal{O}}}
\nc{\CP}{{\mathcal{P}}}
\nc{\CQ}{{\mathcal{Q}}}
\nc{\CR}{{\mathcal{R}}}
\nc{\CS}{{\mathcal{S}}}
\nc{\CT}{{\mathcal{T}}}
\nc{\CU}{{\mathcal{U}}}
\nc{\CV}{{\mathcal{V}}}
\nc{\CK}{{\mathcal{K}}}
\nc{\CW}{{\mathcal{W}}}
\nc{\CZ}{{\mathcal{Z}}}

\nc{\cM}{{\check{\mathcal M}}{}}
\nc{\csM}{{\check{\mathcal A}}{}}
\nc{\oM}{{\overset{\circ}{\mathcal M}}{}}
\nc{\obM}{{\overset{\circ}{\mathbf M}}{}}
\nc{\oCA}{{\overset{\circ}{\mathcal A}}{}}
\nc{\obA}{{\overset{\circ}{\mathbf A}}{}}
\nc{\ooM}{{\overset{\circ}{M}}{}}
\nc{\osM}{{\overset{\circ}{\mathsf M}}{}}
\nc{\vM}{{\overset{\bullet}{\mathcal M}}{}}
\nc{\nM}{{\underset{\bullet}{\mathcal M}}{}}
\nc{\oD}{{\overset{\circ}{\mathcal D}}{}}
\nc{\obD}{{\overset{\circ}{\mathbf D}}{}}
\nc{\oA}{{\overset{\circ}{\mathbb A}}{}}
\nc{\op}{{\overset{\bullet}{\mathbf p}}{}}
\nc{\cp}{{\overset{\circ}{\mathbf p}}{}}
\nc{\oU}{{\overset{\bullet}{\mathcal U}}{}}
\nc{\oZ}{{\overset{\circ}{\mathcal Z}}{}}
\nc{\ofZ}{{\overset{\circ}{\mathfrak Z}}{}}
\nc{\oF}{{\overset{\circ}{\fF}}}

\nc{\fa}{{\mathfrak{a}}}
\nc{\fb}{{\mathfrak{b}}}
\nc{\fg}{{\mathfrak{g}}}
\nc{\fgl}{{\mathfrak{gl}}}
\nc{\fh}{{\mathfrak{h}}}
\nc{\fj}{{\mathfrak{j}}}
\nc{\fm}{{\mathfrak{m}}}
\nc{\fn}{{\mathfrak{n}}}
\nc{\fu}{{\mathfrak{u}}}
\nc{\fp}{{\mathfrak{p}}}
\nc{\fr}{{\mathfrak{r}}}
\nc{\fs}{{\mathfrak{s}}}
\nc{\fsl}{{\mathfrak{sl}}}
\nc{\hsl}{{\widehat{\mathfrak{sl}}}}
\nc{\hgl}{{\widehat{\mathfrak{gl}}}}
\nc{\hg}{{\widehat{\mathfrak{g}}}}
\nc{\chg}{{\widehat{\mathfrak{g}}}{}^\vee}
\nc{\hn}{{\widehat{\mathfrak{n}}}}
\nc{\chn}{{\widehat{\mathfrak{n}}}{}^\vee}

\nc{\fA}{{\mathfrak{A}}}
\nc{\fB}{{\mathfrak{B}}}
\nc{\fD}{{\mathfrak{D}}}
\nc{\fE}{{\mathfrak{E}}}
\nc{\fF}{{\mathfrak{F}}}
\nc{\fG}{{\mathfrak{G}}}
\nc{\fK}{{\mathfrak{K}}}
\nc{\fL}{{\mathfrak{L}}}
\nc{\fM}{{\mathfrak{M}}}
\nc{\fN}{{\mathfrak{N}}}
\nc{\fP}{{\mathfrak{P}}}
\nc{\fU}{{\mathfrak{U}}}
\nc{\fV}{{\mathfrak{V}}}
\nc{\fZ}{{\mathfrak{Z}}}

\nc{\bb}{{\mathbf{b}}}
\nc{\bc}{{\mathbf{c}}}
\nc{\bd}{{\mathbf{d}}}
\nc{\be}{{\mathbf{e}}}
\nc{\bj}{{\mathbf{j}}}
\nc{\bn}{{\mathbf{n}}}
\nc{\bp}{{\mathbf{p}}}
\nc{\bq}{{\mathbf{q}}}
\nc{\bF}{{\mathbf{F}}}
\nc{\bu}{{\mathbf{u}}}
\nc{\bv}{{\mathbf{v}}}
\nc{\bx}{{\mathbf{x}}}
\nc{\bs}{{\mathbf{s}}}
\nc{\by}{{\mathbf{y}}}
\nc{\bw}{{\mathbf{w}}}
\nc{\bA}{{\mathbf{A}}}
\nc{\bK}{{\mathbf{K}}}
\nc{\bI}{{\mathbf{I}}}
\nc{\bB}{{\mathbf{B}}}
\nc{\bG}{{\mathbf{G}}}
\nc{\bC}{{\mathbf{C}}}
\nc{\bD}{{\mathbf{D}}}
\nc{\bP}{{\mathbf{P}}}
\nc{\bH}{{\mathbf{H}}}
\nc{\bM}{{\mathbf{M}}}
\nc{\bN}{{\mathbf{N}}}
\nc{\bV}{{\mathbf{V}}}
\nc{\bU}{{\mathbf{U}}}
\nc{\bL}{{\mathbf{L}}}
\nc{\bT}{{\mathbf{T}}}
\nc{\bW}{{\mathbf{W}}}
\nc{\bX}{{\mathbf{X}}}
\nc{\bY}{{\mathbf{Y}}}
\nc{\bZ}{{\mathbf{Z}}}
\nc{\bS}{{\mathbf{S}}}

\nc{\sA}{{\mathsf{A}}}
\nc{\sB}{{\mathsf{B}}}
\nc{\sC}{{\mathsf{C}}}
\nc{\sD}{{\mathsf{D}}}
\nc{\sF}{{\mathsf{F}}}
\nc{\sG}{{\mathsf{G}}}
\nc{\sK}{{\mathsf{K}}}
\nc{\sM}{{\mathsf{M}}}
\nc{\sO}{{\mathsf{O}}}
\nc{\sQ}{{\mathsf{Q}}}
\nc{\sP}{{\mathsf{P}}}
\nc{\sZ}{{\mathsf{Z}}}
\nc{\sfp}{{\mathsf{p}}}
\nc{\sr}{{\mathsf{r}}}
\nc{\sg}{{\mathsf{g}}}
\nc{\sff}{{\mathsf{f}}}
\nc{\sfb}{{\mathsf{b}}}
\nc{\sfc}{{\mathsf{c}}}
\nc{\sd}{{\ltimes}}

\nc{\tA}{{\widetilde{\mathbf{A}}}}
\nc{\tB}{{\widetilde{\mathcal{B}}}}
\nc{\tg}{{\widetilde{\mathfrak{g}}}}
\nc{\tG}{{\widetilde{G}}}
\nc{\TM}{{\widetilde{\mathbb{M}}}{}}
\nc{\tO}{{\widetilde{\mathsf{O}}}{}} 
\nc{\tU}{\widetilde{U}}
\nc{\TZ}{{\tilde{Z}}}
\nc{\tx}{{\tilde{x}}}
\nc{\tq}{{\tilde{q}}}

\nc{\tfP}{{\widetilde{\mathfrak{P}}}{}}
\nc{\tz}{{\tilde{\zeta}}}
\nc{\tmu}{{\tilde{\mu}}}

 \def\RES{\mathrm{RES}}
  \def\tV{{\widetilde{V}}}

 \def\VF{\mathrm{VF}}

  \def\Sym{\operatorname{Sym}}
   \def\Aut{\operatorname{Aut}}
 
  \def\spec{\operatorname{Spec}}

 \def\dom{\operatorname{dom}}

 \def\GL{\operatorname{GL}}
 \def\Id{\operatorname{Id}}
\def\St{\operatorname{St}}

\def\val{\operatorname{val}}
\def\res{\operatorname{res}}
  
\def\Ind#1#2#3{{#1} {\downarrow}_{#3} {#2} }  

\let\tensor\otimes
\let\meet\cap
\let\union\cup

\def\si{\sigma}
\def\g{\gamma}
\def\G{\Gamma}

\def\<{\begin}
\def\>{\end}

\def\m{\smallsetminus}

\nc{\seq}[1]{\stackrel{#1}{\sim}}

\def\inv{^{-1}}

\def\step#1{{\noindent \bfseries Step #1.\ }}
\def\beq#1{{\begin{equation} \label{#1}}  }

\def\Uu{\mathbb U}

\def\prf{\begin{proof}}
\def\pv{\end{proof}}
\def\eprf{\end{proof}}

\def\acl{\operatorname{acl}}
\def\dcl{\operatorname{dcl}}

\def\a{\alpha}

\def\kk{{\mathrm k}}
 \renc{\b}{{\beta}}

\def\std#1{{\widehat{#1}}}
 \def\tK{{\widetilde{K}}}

    \def\RES{{\mathrm{RES}}}
 \def\tV{{\widetilde{V}}}
 
 \def\tv{\tilde{v}}

\def\Def{\mathrm{Def}}
\def\ProDef{\mathrm{ProDef}}
\def\IndDef{\mathrm{IndDef}}
\def\tp{\mathrm{tp}}

\let\limproj\varprojlim
\let\limind\varinjlim
\def\Ind#1#2{#1\setbox0=\hbox{$#1x$}\kern\wd0\hbox to 0pt{\hss$#1\mid$\hss}
\lower.9\ht0\hbox to 0pt{\hss$#1\smile$\hss}\kern\wd0}
\def\dnf{\mathop{\mathpalette\Ind{}}}
 
 \def\Lam{\Lambda}

\def\Lam{\Lambda}

\ifpupbook
\usepackage{makeidx}
\fi
\makeindex
\makeatletter
\renewcommand{\@idxitem}  {\par\hangindent 2em\rightskip0pt plus 4em}
\makeatother
\usepackage[intoc,refpage]{nomencl}
\usepackage{multicol}

\makeatletter
\def\nommark{\markboth{\MakeUppercase{\nomname}}{\MakeUppercase{\nomname}}}%
\def\thenomenclature{%
  \@ifundefined{chapter}%
  {
    \section*{\nomname}
\ifpupbook
    \if@intoc\addcontentsline{toc}{section}{\nomname}\fi%
\fi%
  }%
  {
    \chapter*{\nomname}
\ifpupbook
      \nommark 
    \if@intoc\addcontentsline{toc}{chapter}{\nomname}\fi%
\fi%
  }%
  \nompreamble
  \list{}{%
\begin{multicols}{2}
    \labelwidth\nom@tempdim
    \leftmargin\labelwidth
    \advance\leftmargin\labelsep
    \rightskip0pt plus 1fil
    \itemsep\nomitemsep\parsep.5\parsep
    \let\makelabel\nomlabel}}
\def\endthenomenclature{%
\end{multicols}
  \endlist
  \nompostamble}
\makeatother
\makenomenclature

  \def\nomname{List of notations}
  
\let\tcb\relax\let\tg\relax
  \usepackage{color} 
\ifcolor
  \definecolor{BLUE}{rgb}{0,0,1}
  \def\tcb{\textcolor{blue}}
  \definecolor{RED}{rgb}{1,0,0}
  \def\tg#1{\textcolor{red}{#1}}
\fi

\def\kk{\mathbf{k}}
\def\alg{\mathrm{alg}}

 \def\Ups{\Upsilon} 
 
 \def\redu{\mathrm{red}}
 \def\nsubset{\subset}
 \def\nsubseteq{\subset}
 \def\llp{\mathopen{(\!(}}

\def\rrp{\mathopen{)\!)}}

\def\vol{\operatorname{vol}}

\def\mD{\mathcal D}
\def\Fn{\mathrm{Fn}}
\def\stda#1{ {#1}^{\#}}

\def\uF{\overline{F}}

\def\Sym{\mathrm{Sym}}
 
 \def\rto{R_{21}^{20}}

 \def\h{\underline{h}}  
\def\Th{\mathrm{Th}}
 \def\v{\underline{v}}
\def\GL{\mathrm{GL}}
  
\def\ACVF{\mathrm{ACVF}}
\def\ACF{\mathrm{ACF}}
\def\ACV2F{\mathrm{ACV}^2\mathrm{F}}
\def\NIP{\mathrm{NIP}}
\def\DOAG{\mathrm{DOAG}} 
\def\RCF{\mathrm{RCF}} 
\def\ZC{\mathrm{ZC}} 
\def\IZC{\mathrm{IZC}}
\def\RV{\mathrm{RV}}
\def\pt{\mathrm{point}}

\def\rv{\operatorname{rv}}

\makeindex

\def\tcr{}

\begin{document}
\frontmatter

\ifpupbook
\else
 
 \title[Non-archimedean tame topology]{Non-archimedean tame topology and stably dominated types}

\author{Ehud Hrushovski}
\address{Department of Mathematics, The Hebrew University, Jerusalem, Israel} \email{ehud@math.huji.ac.il}

\author{Fran\c cois Loeser}
\address{Sorbonne Universit\'es, UPMC Univ Paris 06, UMR 7586 CNRS, Institut Math\'ematique de Jussieu, F-75005 Paris, France}
\email{Francois.Loeser@upmc.fr} 

\fi

\ifpupbook
\else
 \begin{abstract}Let $V$ be a quasi-projective algebraic variety over a non-archimedean valued field.  We introduce topological methods into the model theory of valued fields,
define an analogue $\widehat {V}$ of the Berkovich analytification $V^{an}$ of $V$, and deduce several new results on Berkovich spaces from it.  In particular
 we show that  $V^{an}$ retracts to a finite simplicial complex and is locally contractible, without any smoothness assumption on $V$.   When
 $V$ varies in an algebraic family, we show that the homotopy type of $V^{an}$ takes only a finite number of values.   The space $\widehat {V}$ is obtained
 by defining a topology on the pro-definable set of stably dominated types on $V$.  The key result is the construction of  a pro-definable
  strong retraction of $\widehat {V}$ to an o-minimal subspace, the skeleton, definably homeomorphic to a space definable over the value group with its piecewise linear structure.   
  \end{abstract}
  
\fi

\ifpupbook
\else
\maketitle
\fi


\setcounter{tocdepth}{1}
\tableofcontents

\mainmatter

\chapter{Introduction}

 Model theory rarely deals directly with topology; the great exception is the   theory of o-minimal structures,
 where the topology arises naturally from an ordered   structure, especially 
 in the setting of ordered fields.    See \cite{vddries-tame} for a basic introduction.  
  Our goal in this work is to create a framework of this kind for valued fields.   
 
 A fundamental tool, imported from stability theory, will be the notion of a definable type; it will play a number of roles, starting from the definition of a point of the fundamental
 spaces that will concern us.    
  A definable type on a definable set $V$ is a uniform decision, for each definable subset $U$ (possibly
 defined with parameters from larger base sets), of whether $x \in U$; here $x$ should be viewed as a kind of ideal element of $V$.  A good example is given by any semi-algebraic function $f$ from $\Rr$ to a real variety $V$.  Such
 a function has a unique limiting behavior at $\infty$:  for any semi-algebraic subset $U$ of $V$, either $f(t) \in U$ for all large enough $t$,
 or $f(t) \notin U$ for all large enough $t$.  In this way $f$ determines a definable type.  
 
 One of the roles  of definable types will be to be a substitute for the classical notion of a sequence,
 especially in situations where one is willing to refine to a subsequence.  The classical notion of
 the limit of a sequence makes little sense in a saturated setting.  In o-minimal situations it can often be replaced
 by the limit of a definable curve; notions such as definable compactness  are defined 
 using continuous definable maps from the field $R$ into a variety $V$.  Now to discuss the limiting behavior of $f$ at $\infty$ (and thus to define notions such
 as compactness), we really require only the answer to this dichotomy\tcr{\textemdash}is $f(t) \in U$ for large $t$?\tcr{\textemdash}uniformly, for all $U$; i.e. knowledge of the definable type associated with $f$.   
  For the spaces we consider, curves will not always be sufficiently plentiful to define compactness, 
but definable types will be, and our main notions will all be defined in these terms.   In particular the limit of a definable
type on a space with a definable topology is a point whose every neighborhood is large in the sense of the definable type.
  
  A different example of a definable type is the generic type of the valuation ring $\Oo$, or of  a closed
  ball $B$ of $K$, for $K$ a non-archimedean valued field, or of  $V(\Oo)$
  where $V$ is a smooth scheme over $\Oo$.     Here again, for any definable subset $U$ of $\Aa^1$, 
 we have $v \in U$ for all sufficiently generic $v \in V$, or else $v \notin U$ for all   sufficiently generic $v \in V$;
 where ``sufficiently generic'' means ``having residue outside $Z_U$'' for a certain    proper Zariski closed subset $Z_U$ of $V(\kk)$, depending only on $U$.  Here  $\kk$ is the residue field.  Note that the generic type of $\Oo$
 is invariant under multiplication by $\Oo^*$ and addition by $\Oo$, 
 and hence induces a definable type on any closed ball.  Such definable types are {\em stably dominated}, being determined by a function into objects over the residue field, in this case the residue map into $V(\kk)$.  They can 
 also be characterized as {\em generically stable}.  Their basic properties were developed in \cite{hhm};
 some results are now seen more easily using the general theory of $\NIP$,  \cite{hp}.

 Let $V$ be an algebraic variety over a field $K$.  A valuation or ordering on $K$ induces 
 a topology on $K$, hence on $K^n$, and finally on $V(K)$.    We view this topology as an object   of the definable world; for any model $M$, we obtain a topological space whose set of points is $V(M)$.  In this sense, the topology is on $V$. 
 
In the valuative case however, it has been recognized since the early days of the theory that this topology is inadequate for geometry.   The valuation topology is totally disconnected, and does not afford a useful  globalization of local questions.   Various 
 remedies have been proposed, by Krasner, Tate, Raynaud, Berkovich \tcb{and Huber}.  Our approach can be viewed as
 a lifting of Berkovich's to the definable category.  We will mention below a number of applications to classical Berkovich spaces, that indeed motivated the direction of our work. 

The fundamental topological spaces we will consider will not live on algebraic varieties.  Consider instead the set
of semi-lattices in $K^n$.  These are $\Oo^n$-submodules of $K^n$ isomorphic to $\Oo^k \oplus K^{n-k}$ for some $k$.  Intuitively, a sequence $\Lambda_n$ of semi-lattices approaches a semi-lattice $\Lambda$ if
for any $a$, if  $a \in \Lambda_n$ for infinitely many $n$ then $a \in \Lambda$; and if $a \notin \Mm \Lambda_n$
for infinitely many $n$, then $a \notin \Mm \Lambda$.  The actual definition is the same, but using definable types.
A definable set of semi-lattices is {\em closed} if it is closed under limits of definable types.  The set of
closed balls in the affine line $\Aa^1$ can be viewed as a closed subset of the set of semi-lattices in $K^2$.
In this case the limit of a decreasing sequence of balls is the intersection of these balls; the limit of the generic type of the valuation ring $\Oo$ (or of small closed balls around generic points of $\Oo$) is the closed ball $\Oo$. 
We also consider subspaces of these spaces of semi-lattices.  They tend to be definably connected and compact,
as tested by definable types.  For instance the set of all semi-lattices in $K^n$ cannot be split into two disjoint closed \tcb{nonempty} definable subsets.  

 To each algebraic variety $V$ over a valued field $K$ we will associate
in a canonical way a projective limit $\std{V}$ of spaces of the type described above.  A point of $\std{V}$
does not correspond to a point of $V$, but rather to a stably dominated definable type on $V$.  We call $\std{V}$ the {\em stable completion} of $V$. For instance  
when $V=\Aa^1$, $\std{V}$ is the 
 set of closed balls of $V$; the stably dominated type associated to a closed ball is just the generic type of that
 ball (which may be a point, or larger).  In this case, and in general for curves, $\std{V}$ is definable
 (more precisely, a definable set of some imaginary sort), and no projective limit is needed.


While $V$ admits no definable functions of interest from the value group $\G$, there do exist 
 definable functions from $\G$ to $\std{\Aa^1}$:  for any point $a$ of $\Aa^1$, one can consider
 the closed ball $B(a;\alpha) = \{x: \val(a-x) \geq \alpha \}$ as a definable function of $\alpha \in \G$. 
These functions will serve to connect the space $\std{\Aa^1}$.  In \cite{hhmcrelle} the imaginary sorts were classified,
and moreover the definable functions from $\G$ into them were classified; in the case of $\std{\Aa^1}$, essentially
the only definable functions are the ones mentioned above.  It is this kind of fact that is the basis of 
the geometry of imaginary sorts that we study here.

 At present we remain in a purely algebraic setting.  The applications to Berkovich spaces are thus  only  to 
 Berkovich spaces of algebraic varieties.  This limitation has the merit of showing  that Berkovich spaces can be
 developed purely algebraically;   historically, Krasner and Tate introduce analytic functions immediately  even when interested in  algebraic varieties, so that the name of the subject is rigid {\em analytic} geometry, but this is not necessary, a rigid algebraic geometry exists as well.

While we discussed o-minimality as an analogy, our real goal is a {\em reduction} of questions over valued fields
to the o-minimal setting.  The value group $\G$ of a valued field is o-minimal of a simple kind, where all definable
objects are piecewise $\Qq$-linear.  Our main result is that for any quasi-projective variety $V$ over $K$, $\std{V}$ admits a definable
 deformation  retraction to a subset $S$, \tcr{called a skeleton,} which is definably homeomorphic to a space defined over $\G$.
 There is a delicate point here:  the definable homeomorphism is valid semi-algebraically, but if one stays in the (tropical) locally semi-linear setting, one must take into account subspaces of $\G_\infty^n$, where 
$\G_\infty$ is a partial completion of $\G$ by the addition of a point at $\infty$.  The intersection of the space with the points at $\infty$ contains valuable additional information.
 In general, such a  skeleton is non-canonical.
At this point, o-minimal results such as triangulation can be quoted.  As a corollary we obtain an 
 equivalence of categories between the 
category of \tcb{definable subsets of quasi-projective} varieties over $K$, with homotopy classes of definable continuous maps $\std{U} \to \std{V}$  as morphisms $U \to V$, and a \tcb{homotopy} category of definable spaces over the o-minimal $\G$.

In case the value group is $\Rr$, our  results specialize  to similar tameness theorems  for Berkovich spaces.
In particular we obtain local contractibility for Berkovich spaces associated to algebraic varieties,
a result which was proved by Berkovich  under smoothness assumptions \cite{berk1}, \cite{berk2}.
We also show that for projective varieties, the corresponding 
Berkovich space is homeomorphic to a projective limit of finite-dimensional simplicial complexes that are deformation retracts of itself. 
 We further 
obtain finiteness statements that were not known classically; we refer to  Chapter \ref{berkovich}
for these applications.

We now present the contents of the chapters  and a sketch 
of the proof of the main theorem.

Chapter \ref{sec2}  includes some background material on definable sets, definable types, orthogonality and domination, especially in the valued field context.   
  In    \ref{ss2.10}  we present the main result
 of \cite{hhm} with a new insight regarding one point, that will be used in several critical points later in the paper.
 We know that every nonempty definable set over an algebraically closed substructure of a model of $\ACVF$
 extends to a definable type.  A definable type $p$ can be decomposed into  a
 definable type $q$ on $\G^n$, and a  map $f$ from this type to stably dominated definable types.  In previous
 definitions of metastability, this decomposition involved an uncontrolled base change that prevented any
 canonicity.  We note here that the $q$-germ of $f$ is defined with no additional parameters, and that it is this germ
 that really determines $p$.  Thus a general definable type is a function from a definable type on $\G^n$
 to stably dominated definable types.  
 
   In Chapter \ref{sec3} we introduce the space $\std{V}$ of stably dominated types 
 on a definable set $V$.  We show that $\std{V}$ is pro-definable; this is in fact true in any $\NIP$ theory,
 and not just in $\ACVF$.  We further show that $\std{V}$ is strict pro-definable, i.e. the image of $\std{V}$
 under any projection to a definable set is definable.  This uses metastability, and also a classical definability
 property of irreducibility in algebraically closed fields.  In the case of curves, we note later that $\std{V}$ is 
 in fact definable; for  many purposes strict pro-definable sets behave in the same way.  Still in Chapter \ref{sec3},
 we define a topology on $\std{V}$, and study the connection between this topology and $V$.  Roughly 
 speaking, the topology on $\std{V}$ is generated by $\std{U}$, where $U$ is a definable set cut out by
 {\em strict} valuation inequalities.  The space $V$ is a dense subset of $\std{V}$, so a continuous map
 $\std{V} \to \std{U}$ is determined by the restriction to $V$.  Conversely, given a definable map $V \to \std{U}$,
 we explain the conditions for extending it to $\std{V}$.  This uses the interpretation of $\std{V}$
 as a set of definable types.   We determine the Grothendieck topology on $V$ itself induced from 
 the topology on $\std{V}$; the closure or continuity of  definable subsets or of  functions on $V$ can be described
 in terms of this Grothendieck topology without reference to $\std{V}$, but we will see that this viewpoint is more limited.

 In Chapter \ref{definable-compactness} we define the central notion of definable compactness; we give a general definition that may
 be useful whenever one has definable topologies with enough definable types.   The o-minimal formulation
 regarding limits of curves is replaced by limits of definable types.  
 We relate definable compactness to 
 being closed and bounded.  We show the expected properties hold, in particular the image of a definably
 compact set under a continuous definable map is definably compact.  
 
 The definition of $\std{V}$ is a little abstract.  In Chapter \ref{sec5}  we give a   concrete representation of $\std{\Aa^n}$
 in terms of spaces of semi-lattices.  This was already alluded to in the first paragraphs of the introduction.
 
 A major issue in this paper is the frontier between the definable and the topological categories.  In o-minimality
 automatic continuity theorems play a role.  Here we did not find such results very useful.  At all events in  \ref{ss6.2}
 we characterize topologically those subspaces of $\std{V}$ that can be definably parameterized by $\G^n$.
 They turn out to be o-minimal in the topological sense too.  We use here in an essential way the
 construction of $\std{V}$ in terms of spaces of semi-lattices, and the characterization in \cite{hhmcrelle}
 of definable maps from $\G$ into such spaces. 
 We shall prove that our retractions provide skeleta lying in the subspace $\stda{V}$ of $\std{V}$ 
 \tcb{of {\em strongly stably dominated types}  introduced in \ref{ss6.5}.}   \tcb{This is another canonical
 space associated with $V$, ind-definable this time, admitting a natural continuous map into $\std{V}$
 which  restricts to a topological embedding on definable subsets.   We study it further in Chapter \ref{section:ssd}; our uniformity results for
 $\std{V}$ depend on it.}

Chapter \ref{sec7} is concerned with the case of curves.  We show that $\std{C}$ is definable (and not just pro-definable)
 when $C$ is a curve.  The case of $\Pp^1$ is elementary, and in equal characteristic zero it is possible to reduce everything to this case.  But in general we use model-theoretic methods.  
 We \tcb{construct} a definable  deformation  retraction from $\std{C}$ into a $\G$-internal subset. 
 We consider relative curves too, i.e. varieties $V$ with maps $f: V \to U$, whose fibers
 are of dimension one.  In this case we  \tcb{prove the existence of} a  deformation  retraction of all fibers that is 
 globally continuous and takes $\std{C}$ into  a $\G$-internal subset for almost all fibers $C$, i.e. all outside
 a proper subvariety of $U$.   On curves lying over this variety, the motions on nearby curves do not converge to any continuous motion.  
  
Chapter \ref{specializations} contains some algebraic criteria for the verification of continuity.  For the Zariski topology on algebraic varieties,
 the valuative criterion is useful:  a constructible set is closed if it is invariant under specializations.  Here we are led
 to doubly valued fields.  These can be obtained from valued fields {\em either} by adding a valued field structure
 to the residue field, {\em or} by enriching the value group with a new convex subgroup.  The functor
 $\std{X}$ is meaningful  for definable sets of this theory as well, and interacts well with the various specializations.  These criteria are used in Chapter \ref{sec9} to verify  the continuity of the relative homotopies of  Chapter \ref{sec7}.  

Chapter \ref{sec9} includes some additional \tcb{material} on homotopies. In particular, for a smooth variety $V$, there exists
 an ``inflation'' homotopy, taking a simple point to the generic type of a small neighborhood of that point.  This
 homotopy has an image that is properly a subset of $\std{V}$, and cannot be understood directly in terms 
 of definable subsets of $V$.  The image of this homotopy retraction has the merit of being contained in $\std{U}$
 for any \tcb{dense} Zariski open subset $U$ of $V$.  
 
Chapter \ref{sec10} contains the  statement and proof of the main theorem.   For any quasi-projective algebraic variety $V$, we  \tcb{prove the existence of} a definable homotopy retraction from $\std{V}$ to an o-minimal subspace of the type described in  \ref{ss6.2}.  After some preliminary reductions,  we  may assume 
$V$ fibers over a variety $U$ of lower dimension and the fibers are curves.  On each fiber, a homotopy retraction
 can be constructed with o-minimal image, as in Chapter \ref{sec7}; above a certain Zariski open subset $U_1$ of $U$,
 these homotopies can be viewed as the fibers of a single homotopy $h_1$.   We require however a global homotopy.
The homotopy $h_1$ itself does not extend
 to the complement of $U_1$; but in the smooth case, one can first apply an inflation homotopy whose image
 lies in $\std{V_1}$, where $V_1$ is the pullback of $U_1$.  If $V$ has singular points, a more delicate preparation
 is necessary.  Let $S_1$ be the image of the homotopy $h_1$.   Now
 a relative version of the results of  \ref{ss6.2} applies (\thmref{relative}); after pulling back the situation 
 to a finite covering $U'$ of $U$, we show that $S_1  $ embeds topologically into $U' \times \G_\infty^N$.  Now any homotopy  retraction of $\std{U}$, \tcb{lifting to $\std{U'}$ and fixing} certain functions into $\G^m$, can be extended to a homotopy retraction  of $S_1$ (\thmref{omin-finite}).    Using induction on dimension, we apply this to a homotopy
 retraction taking $U$ to an o-minimal set; we obtain a retraction of $V$ to a subset $S_2$ of $S_1$ lying over
 an o-minimal set, hence itself o-minimal.  
At this point   o-minimal topology as in \cite{coste} applies to $S_2$, and hence to the homotopy type of $\std{V}$.   
In  \ref{ss10.7} we give a uniform version of
\thmref{1} with respect to parameters.
\tcb{In Chapter \ref{secns} we examine the simplifications occuring in the proof of the main theorem in the smooth case and in 
Chapter \ref{secec} we deduce an
 equivalence of categories between a certain
 homotopy category
 of definable subsets of quasi-projective varieties over a given valued field and a suitable homotopy category of 
 definable spaces over the o-minimal $\G$.}

Chapter  \ref{berkovich} contains various applications to classical Berkovich spaces. 
 Let $V$ be a quasi-projective variety over a  field 
 $F$ endowed with a non-archimedean norm and let $V^{\textrm{an}}$ be the corresponding Berkovich space.   
 We
 deduce from our  main theorem several new results on the topology of
 $V^{\textrm{an}}$ which were not known previously in such a  level of generality.
 In particular we show that $V^{\textrm{an}}$ admits a strong deformation retraction to a subspace
homeomorphic to a finite simplicial complex and that
$V^{\textrm{an}}$
 is locally contractible.
We prove a finiteness statement  for the homotopy type of fibers in families. 
We also show that if $V$ is projective,  $V^{\textrm{an}}$ is homeomorphic to a projective limit of finite-dimensional simplicial complexes that are deformation retracts of $V^{\textrm{an}}$. 

We do not assume any previous knowledge of Berkovich spaces, but highly recommend the survey \cite{duc-b}, as well as \cite{duc-b2} for an introduction to the model-theoretic viewpoint, and a sketch of proof of \thmref{1}.

 \medskip
 We are grateful to Vladimir Berkovich, Antoine Chambert-Loir, Zo\'e Chatzidakis, Antoine Ducros,  Martin Hils, Dugald Macpherson, Kobi Peterzil, Anand Pillay, and Sergei Starchenko for their very useful comments.   
 \tcb{We address special thanks to Antoine Chambert-Loir for his invaluable help in preparing the final version of the text.}
 The paper has also benefited greatly from highly extensive and thorough comments by anonymous referees, and we are very grateful to them.

\medskip{During the preparation of this paper, the research of the authors has been partially supported by the following grants: 
E. H. by ISF  1048/07 and the European Research Council under the European Union's Seventh Framework Programme (FP7/2007-2013)/ERC Grant Agreement No. 291111;  F.L. by ANR-06-BLAN-0183 and the  European Research Council under the European Union's Seventh Framework Programme (FP7/2007-2013)/ERC Grant Agreement No. 246903/NMNAG.}
\[
\star \quad \star \quad \star
\]


\chapter{Preliminaries}\label{sec2}

{\small \noindent \textbf{Summary.}
In \ref{ss2.1}-\ref{ss2.5} we recall some  model theoretic notions we shall use in an essential way in this work: definable, pro-definable and ind-definable sets, definable types,
orthogonality to a definable set, stable domination. In 
\ref{ss2.6}-\ref{ss2.8} we consider  more specifically these concepts in the framework of the theory  $\ACVF$ of algebraically closed valued fields and recall in particular some results of
\cite{hhmcrelle} and \cite{hhm} we rely on. In \ref{ss2.9} we describe the
definable types concentrating  on a stable definable $V$ as an ind-definable set.
In \ref{ss2.10},  we prove a key \tcb{result} allowing us to view
definable types  as  integrals of stably dominated types along some definable type on the value group sort.
Finally, in \ref{pseudogalois} we discuss the notion of pseudo-Galois coverings that we shall use in Chapter \ref{secgammaint}. \par\bigskip
}

We will rapidly recall the basic model theoretic notions of which we make use, 
but we recommend to the non-model theoretic reader an introduction such as
\cite{pillay-mt} (readers seeking a more comprehensive  text on  stability may
also consult \cite{pillay}).

\section{Definable  sets}\label{ss2.1}

Let us fix a first order language $\CL$ and a complete  theory $T$ \nomenclature{$\CL$}{first order language}
over $\CL$. The language $\CL$ may be multisorted.
If ${\mathcal S}$ is a sort, and $A$ is an $\CL$-structure, we denote by
${\mathcal S} (A)$, the part of $A$ belonging to the sort ${\mathcal S}$. \nomenclature{${\mathcal S} (A)$}{part of $A$ belonging to the sort ${\mathcal S}$}
For $C$ a set of parameters in a model of $T$ and
$x$ any set of variables, we denote by $\CL_C$ the language $\CL$
with symbols of constants for element of $C$ added and
 by $S_x (C)$ the set of types over $C$  \nomenclature{$S_x (C)$}{type space}
in the variables $x$. Thus, 
$S_x (C)$ is the Stone space of the Boolean algebra of formulas with free variables contained in $x$
\tcb{and parameters from $C$}
up to equivalence over $T$.
\tcb{If $A$ is a tuple or a set of parameters and $B$ is a  set of parameters, we shall denote by
$\tp (A / B)$ the type of $A$ over $B$. We write $\tp (A / B) \vdash \tp (A / BC)$ to mean that $\tp (A / B)$ implies  $\tp (A / BC)$, i.e.
 $\tp (A / BC) = \tp (A' / BC)$
whenever
$\tp (A / B) = \tp (A' / B)$.}

We shall work in a large saturated model $\Uu$ (a universal domain for $T$). \nomenclature{$\Uu$}{universal domain}
More precisely, we shall  fix some  uncountable cardinal $\kappa $ larger than any cardinality of interest, and consider
a model $\Uu$ of cardinality $\kappa$ such that for 
every $A \nsubset  U$ of cardinality $< \kappa$,
every $p$ in 
$S_x(A)$ is realized in $\Uu$, for $x$ any  set of variables of \tcb{cardinality $< \kappa$}.
Such a $\Uu$ is unique up to isomorphism. 
Set theoretic issues
involved in the choice of $\kappa$ turn out to be unimportant and resolvable in numerous 
ways; cf. \cite{chang-keisler} or \cite{bedlewo}, Appendix A.

All sets of parameters $A$ we shall consider will be {\em small} \index{small} subsets of $\Uu$, that is
of cardinality $< \kappa$, 
and all
models $M$ of $T$ we shall consider will be elementary substructures
of $\Uu$ with cardinality $< \kappa$.
By a {\em substructure} \index{substructure} of $\Uu$ we shall generally mean a small definably closed subset
of $\Uu$.

If
$\varphi$ is a formula in
$\CL_C$, involving some sorts ${\mathcal S}_i$ with arity $n_i$, for every
small model $M$ containing $C$, one can consider
the set
$Z_{\varphi} (M)$ of tuples $a$ in the cartesian product of the 
 ${\mathcal S}_i (M)^{n_i}$
such that $M \models \varphi (a)$.
One can view $Z_{\varphi}$ \nomenclature{$Z_{\varphi}$}{definable set associated to $\varphi$} as a functor from the category of models and 
 elementary embeddings, to the category of sets.
Such functors will be called {\em definable sets} \index{definable set} over
$C$. 
Note  that a definable set $X$ is completely determined by
the (large) set
$X (\Uu)$, so we may identify
definable sets with subsets of cartesian products of sets  ${\mathcal S}_i (\Uu)^{n_i}$.
Definable sets over $C$ form a category $\Def_C$ in a natural way. \nomenclature{$\Def_C$}{definable sets over $C$}
Under the previous identification
a definable morphism between
definable sets
$X_1 (\Uu)$  and $X_2 (\Uu)$ is a function
$X_1 (\Uu) \rightarrow X_2 (\Uu)$
whose graph is definable.

By a definable set, we mean definable over some $C$. When $C$ is empty
one says $\varnothing$-definable or $0$-definable. \index{$\varnothing$-definable} \index{$0$-definable}
A subset of a given definable set $X$ which is an intersection
of $< \kappa$ definable subsets of $X$ is said to
be {\em $\infty$-definable}. \index{$\infty$-definable}

\tcb{When the theory $T$ has quantifier elimination, any definable set can be defined by a quantifier-free formula, and in any place where it matters we will always suppose
that it is so defined.}

Sets of $\Uu$-points of definable sets satisfy the following form of compactness:
if $X$ is a definable set such that
 $X (\Uu) = 
\tcb{\bigcup}_{i \in I}X_i (\Uu)$,  
with 
 $(X_i)_{i \in I}$  a small family of definable sets, then
 $X = \tcb{\bigcup}_{i \in A} X_i$ with $A$ a finite subset of $I$.

Recall that if $C$ is a subset of a model $M$ of $T$, by the {\em algebraic closure} of \index{algebraic closure}
$C$, denoted by $\acl (C)$, \nomenclature{$\acl$}{algebraic closure} one denotes the subset of those elements $c$ of $M$,
such that, for some formula $\varphi$ over $C$ with one free variable, $Z_{\varphi} (M)$ is finite and contains $c$.
The {\em definable closure}\index{definable closure} $\dcl (C)$ \nomenclature{$\dcl$}{definable closure}  of $C$ is 
the subset of those elements $c$ of $M$,
such that, for some formula $\varphi$ over $C$ with one free variable, $Z_{\varphi} (M)
= \{c\}$.

If $X$ is a $C$-definable set  and $C \nsubseteq B$,
we write $X (B)$ for $X (\Uu) \cap \dcl (B)$.

\section{Pro-definable  and ind-definable sets}\label{ss2.2}  We define the category $\ProDef_C$ of {\em pro-definable} sets over $C$  as the category of pro-objects  in the category $\Def_C$ indexed by a small directed partially  ordered set. \index{pro-definable} \nomenclature{$\ProDef_C$}{pro-definable sets over $C$}Thus, if
$X = (X_i)_{i \in I}$ and
$Y = (Y_j)_{i \in J}$ are two objects in 
$\ProDef_C$,
\[\mathrm{Hom}_{\ProDef_C} (X, Y) = 
\limproj_j \limind_i \mathrm{Hom}_{\Def_C} (X_i, Y_j).\]
Elements of 
$\mathrm{Hom}_{\ProDef_C} (X, Y)$
will be called $C$-pro-definable morphisms (or $C$-definable for short) between
$X$ and $Y$.

By a result of Kamensky
\cite{kamensky},
the functor of ``taking $\Uu$-points'' induces an equivalence of categories between the 
category $\ProDef_C$  and
the sub-category of the category of sets whose objects and morphisms are
inverse  limits of $\Uu$-points of definable sets  indexed by a small directed partially ordered set 
(here the word ``co-filtering'' is also used, synonymously with ``directed''). By pro-definable, we mean pro-definable over some $C$.
\tcb{In this paper we shall freely identify a pro-definable set $X = (X_i)_{i \in I}$ with the set $X (\Uu) = \limproj_i X_i (\Uu)$.
For any set $B$ with $C \nsubseteq B \nsubseteq \Uu$, we set $X (B) = X (\Uu) \cap \dcl (B)=  \limproj_i X_i (B)$.}
Pro-definable is thus the same as $\ast$-definable in the sense of Shelah, that is,  a small projective limit of definable subsets.  
One defines similarly the category $\IndDef_C$ of {\em ind-definable} sets over $C$ \index{ind-definable} \nomenclature{$\IndDef_C$}{ind-definable sets over $C$}
for which a similar equivalence holds.

 Let $X $ be a pro-definable set. We shall say  
 it is {\em strict} pro-definable \index{strict pro-definable} if it may  be represented as a pro-object
 $ (X_i)_{i \in I}$, with surjective
 transition morphisms $X_j (\Uu)  \rightarrow X_i (\Uu)$.  
 Equivalently, it is 
a $\ast$-definable set, such that the projection to any finite number
of coordinates is definable.    

Dual definitions apply to ind-definable  sets; thus ``strict''  \index{strict ind-definable}
 means that the maps are injective: 
in $\Uu$, a small union of definable sets is a strict ind-definable 
set.

By a morphism from an ind-definable set $X = \limind_i X_i $ to a pro-definable one $Y = \limproj_j Y_j$, we mean a compatible family of morphisms
$X_i \to Y_j$.  A morphism $Y \to X$ is defined dually; it is always  represented by  a morphism $Y_j \to X_i$, for some $j,i$.

 \begin{rem}  
   Any strict ind-definable set $X$ with a definable point   admits a bijective morphism 
to a   strict pro-definable set.
     On the other hand, if $Y$ is strict pro-definable and $X$ is strict ind-definable,
 a morphism $Y \to X$ always has definable image.   
 \end{rem}

\prf  Fix a definable point $p$ \tcb{in $X$}.
If $f:X_i\to X_j$ is
injective, define $g:X_j\to X_i$ by setting it equal to $f^{-1}$ on $\mathrm{Im}(f)$,
constant equal to $p$ outside $\mathrm{Im} (f)$.   The second statement is clear \tcb{by compactness}.  
 \eprf
 
\begin{defn}  \label{indproetc}
Let  $Y  = \limproj_i Y_i$ be pro-definable. Assume given, for each $i$,  $X_i \nsubseteq Y_i$ such that
the transition maps $Y_i \to Y_{i'}$, for $i \geq i'$, restrict to maps
$X_i \to X_{i'}$ and set $X=\limproj_i X_i$.
 
 \begin{enumerate} \item  If  each $X_i$ is definable and, for some $i_0$, the maps $X_i \to X_{i'}$ are bijections for all $i \geq i' \geq i_0$,
 we say $X$ is {\em iso-definable}. \index{iso-definable}
 \item If  each $X_i$ is $\infty$-definable and, for some $i_0$, the maps $X_i \to X_{i'}$ are bijections for all $i \geq i' \geq i_0$,
 we say $X$ is {\em  iso-$\infty$-definable}. \index{iso-$\infty$-definable}
 \item If there exists a definable set $W$ and a pro-definable morphism $g : W \to Y$  such that for each $i$,  the composition of $g$ and the  projection $Y \to Y_i$ has image $X_i$, we say $X$ is {\em definably parameterized}.  \index{definably parameterized}
 \end{enumerate} 
\end{defn}  

\tcb{In \exref{nonisodefinable} we shall give an  example, for the spaces we will consider,
of a definably parameterized subset which is not iso-definable.  In Question \ref{f7q} we formulate an open problem about inverse images of iso-definable subsets
under finite morphisms.
We now give  two conditions under which definably parameterized sets are iso-definable.  
}

\begin{lem} \label{b3} Let $W$ be a definable set, $Y$ a pro-definable set, and let $f: W \to Y$
be a pro-definable map.  Then the image of $W$ in $Y$ is strict pro-definable. 
If $f$ is injective, or more generally if the equivalence relation \tcb{on $W$
defined by $f(w)=f(w')$} is definable,  then $f(W)$ is iso-definable.   \end{lem}

\prf  
 Say $Y  = \limproj_i Y_i$.  Let $f_i$ be  the composition $W \to Y \to Y_i$.  Then $f_i$ is a 
function   whose graph is $\infty$-definable.  By compactness there exists a definable function
$F: W \to Y_i$ whose graph contains $f_i$; but then clearly $F=f_i$ and so the image $X_i=f_i(W)$
and $f_i$ itself are definable.  Now $f(W)$ is the projective limit of the system $(X_i)$, with maps
induced from $(Y_i)$; the maps $X_i \to X_j$ are surjective for $i>j$, since $W \to X_j$ is 
surjective.    Now assume there exists a definable equivalence relation $E$ on $\tcb{W}$ such that
$f(\tcb{w})=f(\tcb{w'})$ if and only if $(\tcb{w}, \tcb{w'}) \in E$.  If $(w,w') \in W^2 \m E$ then $w$ and $w'$
have distinct images in some $X_i$.  By compactness, for some $i_0$,  if
$(w,w') \in W^2 \m E$ then $f_{i_0}(w) \neq f_{i_0}(w')$.  So  for any $i$ mapping to $ i_0$ the map $  X_i \to X_{i_0}$   is injective.        \eprf

 \begin{cor} \label{b3c}Let $Y$ be pro-definable and let
$X \nsubseteq Y$ be a pro-definable subset. Then $X$ is iso-definable if and only if $X$ is in  \textup{(}pro-definable\textup{)} bijection with a 
definable set.     \qed \end{cor}

\begin{lem}\label{b3cc}  Let $Y$ be pro-definable, $X$ an iso-definable subset.  Let $G$ be a finite group acting on $Y$,
and leaving $X$ invariant.  Let $f: Y \to Y'$  be a map of pro-definable sets, whose fibers are exactly the orbits
of $G$.  Then $f(X)$ is iso-definable.  \end{lem}

\prf  Let $U$ be a definable set, and $h: U \to X$ a pro-definable bijection.   Define $g(u)=u'$ if $ gh(u) = h(u')$. 
This induces a definable action of $G$ on $U$.   We have $f(h(u))=f(h(u'))$ iff there exists $g$ such that $gu =u'$.  
Thus the equivalence relation $f(h(u)) = f(h(u'))$ is definable; by \lemref{b3}, the image is iso-definable.
\eprf

 We shall call a subset $X$ of a pro-definable set $Y$  {\em relatively definable} \index{relatively definable} in $Y$
  if   $X$ is cut out from $Y$ by a single formula.
More precisely, 
if $Y  = \limproj_i Y_i$ is pro-definable, $X$ will be relatively definable  in
$Y$ if there exists some index $i_0$ and a definable
subset $Z$ of $Y_{i_0}$, such that, denoting by $X_i$ the inverse image
of $Z$ in $Y_i$ for $i \geq i_0$, $X = \limproj_{i\geq i_0} X_i$.
\tcb{A subset  of a pro-definable set $Y$ is called {\em relatively $\infty$-definable} \index{relatively $\infty$-definable} in $Y$
if it is
the intersection of a small family of relatively definable subsets of $Y$.}

Iso-definability and relative definability  are   related somewhat as 
finite dimension is related to finite codimension; so they rarely hold together.  
In this terminology, a  semi-algebraic subset of $\std{V}$, that is,  a subset of the form $\std{X}$, where $X$ is a definable
subset of $V$,
will be  relatively definable, but most often  not iso-definable.

\begin{lem} \label{b5}\begin{enumerate}\item Let $X$ be pro-definable, and assume that the equality relation $\Delta_X$ is a
 relatively definable subset of $X^2$.  Then $X$ is iso-$\infty$-definable.
\item A pro-definable subset of an iso-$\infty$-definable set is iso-$\infty$-definable.
\end{enumerate}\end{lem}

\prf  (1)  $X$ is the projective limit of an inverse system $\{X_i\}$, with maps $f_i: X \to X_{i}$.  We have $(x,y) \in \Delta_X$ if and only if $f_i(x)=f_i(y)$ for each $i$.
It follows that for some $i$, $(x,y) \in \Delta_X$ if and only if $f_i(x)=f_i(y)$.  For otherwise,
for any finite set $I_0$ of indices, we may find $(x,y) \notin \Delta_X$ with $f_i(x)=f_i(y)$ for every $i \in I_0$. 
But then by compactness, and using the relative definability of (the complement of)
$\Delta_X$, there exist $(x,y) \in X^2 \m \Delta_X$ with $f_i(x)=f_i(y)$ for all $i$,
a contradiction.  Thus the map $f_i$ is injective.  (2) follows from (1), or can be proved directly.  \eprf

 \begin{lem}\tcb{Let $f : X \to Y$ be a morphism between pro-definable sets.
 If $Y$ is \textup{(}isomorphic to\textup{)} a definable set, then $\mathrm{Im} f$ is $\infty$-definable. In general
 $\mathrm{Im} f$ is pro-definable.}
 \end{lem}

\prf  \tcb{Follows easily from compactness.}
\eprf

\begin{lem}Let $f : X \to Y$ be a morphism between strict pro-definable sets.
Then $ \mathrm{Im} f$ is strict pro-definable, as is the graph of $f$.
\end{lem}

\prf  We can represent $X$ and $Y$ as respectively projective limit of definable sets $X_{i}$ and $ Y_{j}$ with surjective transition mappings
and $f$ by $f_{j}: X_{c(j)} \to Y_{j}$, for some function $c$ between the index sets.  The projection of 
\tcb{$\mathrm{Im} f$}
to $Y_{j}$ is the same as the image of $f_{j}$, using the surjectivity of the maps between the sets $X_{\tcb{c(j)}}$
and $f_{j}(X_{\tcb{c(j)}})$.  The graph of $f$ is the image of $\mathrm{Id} \times f: X \to (X \times Y)$. 
\eprf

\begin{remark}[On terminology] \label{yoneda}   
We often have a  set $D(A)$ depending functorially on a structure $A$.    We say that $D$ is pro-definable   if there exists a pro-definable set $D'$ such that $D'(A)$ and $D(A)$ are in canonical bijection; in other words
$D$ and $D'$ are isomorphic functors.

In practice we have in mind a choice of $D'$ arising naturally from the definition of $D$; usually various interpretations
are possible, but all are isomorphic as pro-definable sets.
Once $D'$ is specified, so is, for any pro-definable $W$ and any $A$, the set  of $A$-definable maps $W \to D'$.   If worried about the identity of $D'$, it suffices to  specify what we mean by an $A$-definable map $W \to D$.  Then 
Yoneda's lemma ensures the uniqueness of a pro-definable set $D'$ compatible with this notion.

The same applies for ind.  For instance, let
 $\Fn (V,V')(A)$ \nomenclature{$\Fn$}{set of definable functions} be the set of $A$-definable functions between two given  $\varnothing$-definable sets $V$ and $V'$.  Then $\Fn(V,V')$
  is an ind-definable set.  The representing   ind-definable  set is clearly determined by the description.

  To avoid all doubts,  we specify that 
  $\Fn(U, \Fn(V,V')) = \Fn(U \times V, V')$.
 \end{remark}

\subsection{Maps from ind-definable to pro-definable sets}
Let $X = \limind_i X_i$ be an ind-definable set, and let $Y = \limproj_j Y_j$ be a pro-definable set.    
Recall that
$\mathrm{Hom}(X,Y) = \limproj_{i,j} \mathrm{Hom}(X_i,Y_j)$, where one denotes by $\mathrm{Hom}(X_i,Y_j)$ the set of definable maps $X_i \to Y_j$.
Clearly, if $f \in \mathrm{Hom}(X,Y)$ then $f$ induces a map $f_M: X(M) \to Y(M)$, for any model $M$.  In case $X$ is strict ind-definable,
we call $f$ {\em injective} if in any model, $f_M$ is injective.   If $X$ is  strict ind-definable and $f$ is represented by $(f_{i,j})$, then $f$ is injective
iff for each $i$, for some $j$, $f_{ij}$ is injective;  since if for arbitrarily large $j$ there exist distinct $x,x' \in X_i$ with $f_{ij}(x) = f_{ij}(x')$,
then by compactness we can find   a pair $x \neq x' \in X_i$ such that for all $j$, $f_{ij}(x) = f_{ij}(x')$.

\begin{defn}  \index{strict ind-definable} \label{strict-ind-def}
 Let $X$ be a subset of a pro-definable set.  By a  {\em  strict ind-definable structure} on $X$ we shall mean a
strict ind-definable set $Z$ together with an injective morphism $g : Z \to Y$ with image $X$. Two such data  $g:Z \to Y$ and $g' : Z' \to Y$ will be considered
to  induce the same structure if there exists an isomorphism $h: Z \to Z'$ of ind-definable sets with $g = h \circ g'$.  
\end{defn}
 
 We will say that ``$X$ is strict ind-definable'' if a strict ind-definable structure is fixed.   In this situation  we will view $X$ itself as being ind-definable, and can apply any notion 
appropriate for ind-definable sets.  Notably we can speak of {\em definable} subsets of $X$; these are iso-definable, but in general 
an iso-definable subset of a strict ind-definable set need not be definable in the sense of the given structure.

\begin{lem} \label{strongind}  \tcb{Let $Y$ be pro-definable. Assume $W \nsubseteq Y$ admits a strict ind-definable structure   $f: X \to W$, such that  
  for each definable $X' \nsubset X$,   for some definable quotient $\pi: Y \to Y'$, the restriction $\pi | f(X)$     
is injective above $\pi(f(X'))$.   Then $W$ has a unique such  ind-definable structure, i.e. 
if $W=f'(X')$ with the same property, then there exists an isomorphism $g: X \to X'$ of ind-definable
sets with $f = g \circ f'$.}   \end{lem}

\prf   \tcb{Let $W$ be strict ind-definable via $f: X \to Y$ and via $f': X' \to Y$ having the above properties.  We need to show that $f \inv \circ f': X' \to X$ is
an isomorphism of ind-definable sets.   
As $f \inv \circ f'$ is a bijection on points, and since the restriction of the graph 
of this bijection to any product $U \times U'$ of definable subsets of $X$ and $X'$ respectively is $\infty$-definable, it suffices to show that $(f \inv \circ f' )(U')$ is contained
in a definable subset of $X'$, for any definable $U' \nsubset X'$ (and vice versa).  Let $\pi: Y \to D$ be a morphism to a definable
set $D$, such that $\pi \circ f'$ is injective above $\pi(f'(U'))$.  
  Now $U' \nsubset \bigcup_{U} ((f') \inv \circ f)(U)$, where $U$ ranges
over all definable subsets of $X$, defined over a given set of definition for $X$.  For $u \in U, u' \in U'$, we have 
$u' = ((f') \inv \circ f)(u) $ iff $f(u)=f'(u')$ iff $\pi \circ f (u) = \pi \circ f'(u')$; this is a definable condition.   So $((f') \inv \circ f) (U)$ is definable.
By compactness, $U'$ is contained in a finite union of sets $((f') \inv \circ f) (U)$; as the union of finitely many definable subsets of $X$ is definable, it is contained in such a set.}
 \eprf

\tcb{Let $Y$ be pro-definable, and consider an injective morphism $f$ from an ind-definable set $X = \limind_i X_i$ to $Y$.  Then
$f(X)$ is strict pro-definable assuming  
 that the equivalence relation $E_i$
on $X_i$ defined by $f(x)=f(x')$ be definable; for then $f_\Uu(X(\Uu)) = g_\Uu(\bX(\Uu))$, where $\bX = \limind_i (X_i/E_i)$, 
$\pi: X \to \bX$ is the natural quotient, 
and $g$ is the map such that $f=g \circ \pi$; note that $\bX$ is strict ind-definable.}

\tcb{From this, and the fact that strict ind-definable sets are closed under disjoint unions, we obtain:}

\begin{lem}\label{unionstrict}\tcb{Let $Y$ be pro-definable.
Let $S_k \nsubseteq Y$ admit a  strict ind-definable structure, via  ind-definable sets $X_k$ and injective morphisms $f_k$ with 
  $f_k(X_k)=S_k$.  
 Assume the pullback to $S_k \times S_{k'}$ of the diagonal 
$\Delta_Y \subset Y \times Y$ is   piecewise definable; i.e.   
  $(f_k \times f_{k'}) \inv  (\Delta_Y ) \meet (X \times X')$
is definable, for any definable $X \nsubset X_k, X' \nsubset X_{k'}$.  Then  $\union_k S_k$ admits a strict ind-definable structure.}
 \qed \end{lem}

\section{Definable types}\label{ss2.3}
 
\tcb{For any set $z$ of variables, we shall denote by $\CL_z$ \nomenclature{$\CL_z$}{set of $\CL$-formulas in variables  in $z$} the set of $\CL$-formulas in variables  in $z$ up to equivalence in the theory $T$.
A type $p(x)$ in variables $x = (x_1, \dots, x_n)$ can be viewed as a Boolean homomorphism from
$\CL_x$ to the two-element Boolean algebra. 
Now consider variables $y_i$  running through all finite
products of sorts. A {\em $0$-definable type} \index{definable type} $p(x)$ is defined to be a 
function $d_p  x : \CL_{x,y_1, \dots,} \rightarrow \CL_{y_1, \dots,}$, such that for any finite $y=(y_1,\ldots,y_n)$, 
$d_p x$ restricts to a 
Boolean retraction $ \CL_{x,y} \rightarrow \CL_{y}$.   An $A$-definable type $p$ is a $0$-definable type of the theory $T_A$ obtained
by naming constants for the elements of the substructure $A$. Sometimes we shall also say $p$ is {\em based on $A$}.  
By a {\em definable type} we mean a $\Uu$-definable type.
   The image of a formula $\phi(x,y)$ under $d_p x$ is called the $\phi$-definition of $p$.
Note that this definition makes sense for any, possibly infinite, set of variables $x$.
When there is no risk of confusion, we sometimes will write $d_p$ instead of $d_p x$.}

\tcb{Given such a retraction, and given any model $M$ of $T$, one obtains
a type over $M$, namely \[p \vert M := \{{\varphi (x, b_1, \dots, b_n) : M \models  (d_p x) (\varphi )(b_1, \dots, b_n)}\}.\]
This type over $M$ determines $p$; this   explains the use of the term definable type.  However viewed as above,
a definable type is really not a type but a different kind of object.
We will often identify $p$ with the type $p \vert \Uu$ which is  $\Aut(\Uu)$-invariant,
and determines $p$.   For any $B \nsubset \Uu$, we denote by
$p \vert B$ the restriction of $p \vert \Uu$ to $B$.
Similarly, for any $C \nsubset \Uu$, replacing $\CL$ by $\CL_C$ one gets the notion of
$C$-definable type.
If $p$ is $C$-definable, then
the type $p\vert \Uu$ is 
$\Aut(\Uu /C )$-invariant.}

\tcb{If $p$ is a definable type and
$X$ is a definable set, or a pro-definable set, one says $p$ is on $X$ if all realizations
of $p\vert \Uu$ lie in $X$. 
One denotes by 
$S_{def, X}$  \nomenclature{$S_{def, X}$}{definable types on $X$} the set of definable types on $X$.
Let  $f : X \rightarrow Z$ be a definable map
between
definable sets, or a pro-definable map between pro-definable sets.
For $p$ in 
$S_{def, X} $ one denotes by 
$f_{*} (p)$ \nomenclature{$f_{*} (p)$}{pushforward of the definable type $p$ by the function $f$} 
the definable type
defined
by
$(d_{f_{*}(p)} z) (\varphi (z, y)) = (d_p x)(\varphi (f (x), y))$.
This gives rise to a mapping
$f_{*}: S_{def, X} \to S_{def, Z}$.
}

\tcb{For a $\varnothing$-definable set $V$, let $L_V$ denote the Boolean algebra of 
$\varnothing$-definable subsets of $V$.  Then a type on $V$ corresponds to  an element of $\mathrm{Hom}(L_V,2)$ and a definable type
on $V$ is the same as a compatible family of elements of $\mathrm{Hom}_W( L_{V \times W}, L_W)$,
for $W$ running over the $\varnothing$-definable sets, where $\mathrm{Hom}_W$ denotes
the set of Boolean homomorphisms $h$ such that $h(V \times X) = X$ for $X \nsubseteq W$.  
Let $U$ be a pro-definable set.  
 By a {\em definable function} $U \to S_{def, V}$, \index{definable function $U \to S_{def, V}$} we mean a compatible family of elements of $\mathrm{Hom}_{W \times U}( L_{V \times W \times U},  L_{W \times U})$ for $W$ running over the $\varnothing$-definable sets.  Any  element $u \in U$
 gives a Boolean retraction $L_{W \times U} \to L_W(u)$ by $Z \mapsto Z(u) = \{z: (z,u) \in Z \}$, with $L_W (u)$ the
 Boolean algebra of $u$-definable subsets of $W$.  So 
 a definable function $U \to S_{def, V}$ gives indeed a $U$-parameterized family of definable types on $V$.
}

\tcb{Let $p$ be a partial $\Uu$-type.}
Let us say $p$ is {\em definably generated} \index{definably generated type} over $A$  if it is generated by a    partial 
type of the form $\bigcup_{(\phi,\theta) \in S} P(\phi,\theta)$, where $S$ is a set of pairs of formulas
$(\phi(x,y), \theta(y))$ over $A$, and $P(\phi,\theta) = \{\phi(x,b):  \tcb{\Uu \models} \theta(b) \}$.  

\begin{lem} \label{defgen}   Let $p$ be a type over $\Uu$.  If $p$ is definably generated over $A$, then $p$ is $A$-definable.
\end{lem} 

\prf  This follows from Beth's theorem:  if one adds a predicate for the $p$-definitions of all formulas $\phi(x,y)$, with the obvious axioms, there is a unique
   interpretation of these predicates in $\Uu$, hence they must be definable.

  Alternatively, let $\phi(x,y)$ be any formula.   From the
fact that $p$ is definably generated it follows easily that 
$\{b: \phi(x,b) \in p \}$ is a \tcb{strict} ind-definable set over $A$. Indeed,  $\phi(x,b) \in p$ if and only if for some $(\phi_1,\theta_1),\ldots,(\phi_m,\theta_m) \in S$, $(\exists c_1,\ldots,c_m)(\theta_i(c_i) \wedge (\forall x) (\bigwedge_i \phi_i(x,c) \implies \phi(x,b)))$.  
Applying this to $\neg \phi$, we see that the complement of $\{b: \phi(x,b) \in p \}$ is also \tcb{strict} ind-definable.
Hence $\{b: \phi(x,b) \in p \}$ is $A$-definable.
\eprf

\begin{cor}\label{dom1}  Let $f: X \to Y$ be an $A$-definable  \tcb{\textup{(}or pro-definable\textup{)}} function.  Let $q$ be an $A$-definable type on $Y$.
\tcb{Let $p_A$ be a type over $A$ such that
$f_* p_A = q|A$ and such that, for any $B$ containing $A$ there exists a unique type $p_B$  over $B$ such that
$p_B$ contains $p_A$, and $f_* p_B = q|B$.   
Here $f_* p_B$ denotes the type generated by  $\mathcal{L}  (B)$-formulas $\varphi (y)$ such that
$\varphi (f (x))$ belongs to $p_B$.
Then there exists a unique $A$-definable type $p$ such that for all $B$, $p|B=p_B$.} 
\qed
\end{cor}

\begin{defn}\label{domby}  In the situation of \corref{dom1}, $p$ is said to be {\em dominated by} \index{dominated by} $q$ via $f$.  \end{defn}

Let us recall that a theory $T$ is said to have 
 {\em elimination of imaginaries} \index{elimination of imaginaries} if, for any $M\models T$, any collection $\mathcal{S}_1,\ldots, \mathcal{S}_k$ of sorts in $T$, and any $\varnothing$-definable 
equivalence relation $E$ on $\mathcal{S}_1(M)\times \cdots \times \mathcal{S}_k(M)$, there is a $\varnothing$-definable function 
$f$ from $\mathcal{S}_1(M)\times \cdots \times \mathcal{S}_k(M)$ into 
a product of sorts of $M$, such that
for any $a,b\in \mathcal{S}_1(M) \times\cdots \times \mathcal{S}_k(M)$, we have $Eab$ if and only if $f(a)=f(b)$. 
Given a complete theory $T$, it is possible to extend it to a complete theory $T^{eq}$ \nomenclature{$T^{eq}$}{imaginary completion of $T$} over a 
language $\CL^{eq}$ by adjoining, 
for each collection $\mathcal{S}_1,\ldots,\mathcal{S}_k$ of sorts and $\varnothing$-definable equivalence relation $E$ on 
$\mathcal{S}_1\times \cdots \times \mathcal{S}_k$, a sort $(\mathcal{S}_1\times \cdots \times \mathcal{S}_k)/E$, together
 with a function symbol for the natural map $a\mapsto a/E$. Any model $M$ of $T$ can be canonically
 extended to a model of $T^{eq}$, denoted $M^{eq}$.
 We shall refer to the new sorts of $T^{eq}$ as imaginary sorts, and to 
elements of them as {\em imaginaries}. \index{imaginaries}

Suppose that $D$ is a definable set in $M\models T$, defined say by the formula $\phi(x,a)$. There is a
$\varnothing$-definable equivalence relation $E_\phi(y_1,y_2)$, where $E_\phi(y_1,y_2)$ holds if and only if
$\forall x(\phi(x,y_1)\leftrightarrow \phi(x,y_2))$. Now 
$a/E_\phi$ is identifiable with an element of an imaginary sort; it is determined uniquely
(up to interdefinability over $\varnothing$) by $D$, and will often be referred to as a {\em code} \index{code}
for $D$, and denoted $\lceil D\rceil$. \nomenclature{$\lceil D\rceil$}{code for $D$}
 We prefer to think of $\lceil D \rceil$ as a fixed object 
(e.g. as  a member of $\Uu^{eq}$) rather than
as an equivalence class of $M$; for  viewed as an equivalence class it is formally a different set
(as is $D$ itself) in elementary extensions of $M$.

\begin{lem}\label{findef}Assume  the theory $T$ has elimination of imaginaries.
Let $f : X \to Y$ 
be a $C$-definable mapping between $C$-definable  sets.
Assume  $f$ has finite fibers, say of cardinality bounded by some integer $m$.
Let $p$ be a $C$-definable type on $Y$.
Then, any global type 
$q$ on $X$ such that
$f_{*} (q) = 
p\vert \Uu$ is $\acl (C)$-definable.
\end{lem}

\prf Let $\bp= p \vert \Uu$.   The partial type
$\bp(f (x))$
admits at most $m$ distinct extensions $q_1$, \dots, $q_{\ell}$ to a complete type.
Choose $C' \supset C$ such that
all $q_i \vert C'$ are distinct.
Certainly the $q_i$ are $C'$-invariant.
It is enough to prove they are $C'$-definable,
since then, 
for every formula $\varphi$, 
the
$\Aut(\Uu /C )$-orbit of $d_{q_i} (\varphi)$
is finite, hence $d_{q_i} (\varphi)$ is equivalent to a formula
in $\CL (\acl (C))$.
To prove $q_i$ is
$C'$-definable
note that 
\[
p (f (x)) \cup (q_i \vert C' )(x) \vdash q_i (x).
\]
Thus,
there is a set
$A$ of formulas
$\varphi (x, y)$ in $\CL$,
a mapping
$\varphi (x, y) \to \vartheta_{\varphi} (y)$
assigning to formulas in 
$A$ formulas in 
$\CL (C')$ 
such that
$q_i $ is generated by
$\{ \varphi (x, b) :  \Uu \models \vartheta_{\varphi} (b)\}$.
It then follows from 
\lemref{defgen}  that 
$q_i$ is indeed 
$C'$-definable.
\eprf

\section{Stable embeddedness}  

A $C$-definable set $D$ in $\Uu$ is  {\em stably embedded}  \index{stably embedded} if, for any definable set $E$ (with parameters $a$ from $\Uu$)
and 
$r>0$, $E \cap D^r$ is definable over $C\cup D$.  To state a more explicit version that does not use $\Uu$:   
 for any formula $\phi(x,y)$ there is a formula $\psi(x,z)$ such that for all $a$ there is a sequence $d$ from $D$ such that
\[\{x\in D^r:\models\phi(x,a)\}=\{x\in D^r:\models \psi(x,d)\}.\]   



For more on stably embedded sets, we refer to the Appendix of \cite{ch}.


\begin{lem} \label{stablyemb} \tcb{Let $T$ be a complete theory in a language $L$ and $D$ a stably embedded sort.  
Let $L_D$ be the restriction of $L$ to $D$ and $L_D^*$ any enrichment of $L_D$.  
Let $T_D$ be the restriction of $T$ to $D$ and
let $T_D^*$ be any expansion of $T_D$ to a complete theory in $L_D^*$.  
Let  $T_D^\sharp$ be the relativization of $T_D^*$ to $D$, i.e the theory that states
that $D \models T_D^*$.  
Let $L^*=L \union L_D^*$ and let $T^*= T \union T_D^\sharp$. 
Then $T^*$ is complete, $D$ is stably embedded in $T^*$, and the $L_D^*$-type of a tuple $b$ of elements of
$D$  implies its $L^*$-type. Moreover, assume $T$ and $T_D^*$ admit quantifier elimination, and for any tuple $a$ in a model of $T$,
$\dcl(a) \meet \dcl(D) = \dcl((f_i(a)))$ where $(f_i)$ enumerates term functions with values in $D$.
Then $T^*$ admits quantifier elimination.}
\end{lem}

\prf  \tcb{Let $M^*,N^*$ be two saturated models of $T^*$ of the same cardinality.  To prove completeness, we must show that 
$M^* \cong N^*$.  To prove stable embeddedness, we must show that any isomorphism $f:D_M^* \to D_{N}^*$ extends to an isomorphism
$M^* \to N^*$.  But $D$ is stably embedded with respect to $L$, so $f$ extends to an $L$-isomorphism $M \to N$, which is by definition
also an $L^*$-isomorphism.  This proves both stable embeddedness and completeness;    completeness also follows since by completeness of $T_D^*$,
we  do have $D_M^* \cong D_N^*$.
The statement about the type of a tuple $b$   follows from the completeness result applied to $T$ and $T_D^*$,
each expanded  by constants for $b$.}

\tcb{To prove the ``moreover'' statement, we must show that if $a,b$ are tuples from $M^*$ respectively with the same
quantifier-free type, then there exists an automorphism of $M^*$ with $a \mapsto b$. 
Let $c=(f_i(a))$ and $d=(f_i(b))$ where $(f_i)$ enumerates term functions with values in $D$.  
Then $c$ and $d$ have the same quantifier-free type in $D(M^*)$ so there exists an automorphism
of $D$ as an $L^*$-structure with $a \mapsto b$.  As above this automorphism extends to $M^*$;
so we may assume it is the identity.  Now in the restriction $M$ of $M^*$ to $L$, we have 
$\tp(a/D)=\tp(b/D)$ so there exists an automorphism $\si$ of $M$ fixing $D$ pointwise with $a \mapsto b$;
and as it fixes $D$, $\si$ is also an $L^*$-automorphism.}
\eprf

\section{Orthogonality to a definable set}\label{ss2.4}

Let $Q$ be a fixed $\varnothing$-definable set.   We give definitions of orthogonality to $Q$ that are convenient for our purposes,
and are equivalent to the usual ones when  $Q$ is stably embedded and admits elimination of imaginaries; this is the only case we will need.

Let $A$ be a substructure of $\Uu$.  
A type $p = \tp (c/A)$ is said to be {\em almost orthogonal} \index{almost orthogonal} to $Q$ if $Q (A(c)) = Q(A)$. Here $A(c)$ is the
substructure generated by $c$ over $A$, and $Q(A) = Q \meet \dcl(A)$ is the set of points
of $Q$ definable over $A$.

An $A$-definable type $p$ is said to be {\em orthogonal} \index{orthogonal} to $Q$,  and one writes
$p \perp Q$, \nomenclature{$\perp$}{orthogonality relation} if $p \vert B$  is almost orthogonal to $Q$ for any
 substructure $B$ containing 
$A$.    Equivalently, for any $B$ and any $B$-definable function $f$ into $Q$
the pushforward $f_{*}(p)$ is a type concentrating on one point, i.e. including a formula of the form $y=\gamma$.

Let us recall that for  $F$ a structure containing $C$, 
$\Fn (W,Q) (F)$ denotes the family of $F$-definable functions $W \to Q$
and that
$\Fn(W,Q) = \Fn(W,Q)( \Uu)$  is an ind-definable set.

Let $V$ be a $C$-definable set.  Let $p$ be 
 a  definable type on $V$, orthogonal to $Q$.  
 Any $\Uu$-definable function $f: V \to Q$ is generically constant  on $p$.   Equivalently,
 any $C$-definable function $f: V \times W \to Q$ (where $W$ is some $C$-definable set)
 depends only on the $W$-argument, when the $V$-argument is a generic realization of $p$.  More precisely,
 we have a mapping
  \[p_{*}^W: \Fn(V \times W, Q) \longrightarrow  \Fn(W,Q)\]
 (denoted by $p_{*}$  when there is no possibility of confusion) \nomenclature{$p_*$, $p_{*}^W$}{pushforward of functions by the orthogonal type $p$}
given by
 $p_{*}(f) (w) = \gamma$ if $(d_p v)(f(v, w)=\gamma)$ holds in $\Uu$.

  Uniqueness of $\gamma$ is clear for any definable type.    Orthogonality to $Q$ is precisely the statement that for any $f$, $p_{*}(f)$ is a function on $W$, i.e.  for any $w$, such an element $\gamma$ exists.      The advantage
  of the presentation $f \mapsto p_{*}(f)$, rather than the two-valued $\phi \mapsto p_{*}(\phi)$, is that it makes
  orthogonality to $Q$ evident from the very data.

 Let $S_{def, V}^Q(A)$ \nomenclature{$S_{def, V}^Q(A)$}{$A$-definable types on $V$ orthogonal to $Q$} denote the set of $A$-definable types on $V$ orthogonal to $Q$.    
 It will be useful to note
 the (straightforward) conditions for pro-definability of $S_{def, V}^Q$.    Given a function $g: S \times W \to Q$,
 we let $g_s(w)=g(s,w)$, thus viewing it as a family of functions $g_s: W \to Q$.

 \begin{lem}\label{prodefcond}  Assume the theory $T$ eliminates imaginaries, and that for any
 formula $\phi(v,w)$ without parameters,  there exists a formula $\theta(w,s)$ without parameters
 such that for any $p \in S_{def, V}^Q$, for some $e$, 
 \[\phi(v,c) \in p \iff \theta(c,e).\]
  Then $S_{def, V}^Q$     is pro-definable,  i.e. there exists a canonical pro-definable
  $Z$ and a canonical bijection $Z(A) = S_{def, V}^Q(A)$ for every $A$.  
 \end{lem} 

\prf We first extend the hypothesis a little.  Let 
  $f: V \times W \to Q$ be $\varnothing$-definable.  Then 
   there exists a $\varnothing$-definable  $g: S \times W \to Q$   such that for any $p \in S_{def, V}^Q$, for some $s \in S$, $p_{*}(f)=g_s$.  Indeed, let $\phi(v,w,q)$ be the formula
   $f(v,w)=q$ and let $\theta(w,q,s)$ 
   \tcb{be}
   the corresponding formula
   provided by 
the hypothesis of the lemma.  Let $S$
   be the set of all $s$ such that for any $w \in W$ there exists a unique $q \in Q$
   with $\theta(w,q,s)$.  Now, by setting  $g(s,w)=q$ if and only if  $ \theta(w,q,s)$ holds, one gets the more general 
   statement. 

 Let $f_i: V \times W_i \to Q$ be an enumeration of all   
 $\varnothing$-definable functions
   $f: V \times W \to Q$, with $i$  running over some index set $I$. 
Let $g_i: S_i \times W_i \to Q$ be the corresponding functions provided by the previous paragraph.  Elimination of imaginaries allows us to 
assume that $s$ is a canonical parameter for the function  $g_{i,s}(w)=g_i(s,w)$, i.e. for no other $s'$
do we have $g_{i,s} = g_{i,s'}$.  We then have a natural map $\pi_i: S_{def, V}^Q \to  S_i$, namely 
$\pi_i(p)=s$ if $p_{*}(f_i) = g_{i,s}$.  Let $\pi= \Pi_i \pi_i: S_{def, V}^Q \to \Pi_i S_i$ be the product map.  
Now $\Pi_i S_i$ is canonically a pro-definable set, and the map $\pi$ is injective.  
So it suffices to show
that the image is relatively $\infty$-definable in $\Pi S_i$.  Indeed, $s=(s_i)_i$ lies in the image if and only if for each
finite tuple of indices $i_1,\ldots,i_n \in I$  (allowing repetitions),
\[(\forall w_{\tcb{i_1}} \in W_{\tcb{i_1}}) \cdots (\forall w_{\tcb{i_n}} \in W_{\tcb{i_n}}) (\exists v \in V) \bigwedge_{j=1}^n f_{\tcb{i_j}}(v,w_{\tcb{i_j}})=g_{\tcb{i_j}}(s_{\tcb{i_j}},w_{\tcb{i_j}}).\]
 For given this consistency condition,
there exists $a \in V(\Uu')$ for some $\Uu \prec \Uu'$ such that $f_i(a,w) =g_i(s,w)$ for all $w \in W_i$ and all $i$.
It follows immediately that $p= \tp(a/\Uu)$ is definable and orthogonal to $Q$, and $\pi(p)= s$.
Conversely if $p \in S_{def, V}^Q(\Uu)$ and $a \models p| \Uu$, for any $w_1 \in W_1(\Uu),\ldots,w_n \in W_n(\Uu)$,
the element $a$ witnesses the existence of $v$ as required.  So the image is cut out by
a set of formulas concerning the $s_i$.
\eprf

 If $Q$ is a two-element set, any definable type is orthogonal to $Q$, and $\Fn(V,Q)$ can 
 be identified with the algebra of formulas on $V$, via characteristic functions.  
 The presentation 
 of definable types as a Boolean retraction from formulas on $V \times W$ to formulas on $W$ can be generalized 
  to definable types orthogonal to $Q$, \tcb{for arbitrary $Q$}.  An element $p$ of  $S_{def, V}^Q(A)$ yields a compatible   system of retractions
  $p_{*}^W: \Fn(V \times W, Q) \longrightarrow  \Fn(W,Q)$.    These retractions are
  also compatible with definable functions $g:Q^m \to Q$, namely $p_{*}(g \circ( f_1,\ldots,f_m))=g \circ (p_{*}f_1,\ldots,p_{*}f_m)$.     One can restrict attention to $\varnothing$-definable functions
  $Q^m \to Q$
  along with compositions of the following form:  given $F: V \times W \times Q \to Q$
  and $f: V \times W \to Q$, let $F \circ '  f (v,w) = F(v,w,f(v,w))$.  Then $p_{*}(F \circ' f) = p_{*}(F) \circ' p_{*}(f)$.   
    It can be shown that any compatible system of retractions compatible with these compositions arises from a unique
  element $p$ of $S_{def, V}^Q(A)$.     This can be shown by the usual two way translation between sets and functions:
  a set can be coded by a function into a two-element set (in case two constants are not available, one can add
  variables $x$, $y$,  and consider functions whose values are among the variables).  On the other hand a function
  can be coded by a set, namely its graph.   This characterization will not be used, and we will leave the details to the reader.   It does give a slightly different way  to see the  $\infty$-definability of the image in \lemref{prodefcond}.

\section{Stable domination}\label{ss2.5}

We shall assume from now on that the theory $T$ has elimination of imaginaries.

\begin{defn} 
A 
$C$-definable set $D$ in $\Uu$ is said to be {\em stable} \index{stable subset} if the structure with domain $D$,
when equipped with all the $C$-definable relations, is stable.
\end{defn}

One considers the multisorted structure
$\St_C$ \nomenclature{$\St_C$}{stable part} whose sorts $D_i$ are the $C$-definable, stable and stably embedded subsets of $\Uu$.
For each finite set of sorts $D_i$, all the $C$-definable relations on their union are considered as
$\varnothing$-definable relations $R_j$. The structure $\St_C$ is stable by 
Lemma 3.2 of \cite{hhm}.

For any $A \nsubset \Uu$, one sets $\St_C (A) = \St_C \cap \dcl (CA)$.

\begin{lem}\label{eqvda}For any subsets $A$, $B$ and $C$ of  $\Uu$ the following conditions are equivalent:
\begin{enumerate}
\item \tcb{$\tp (B/ C\St_C (A)) \vdash \tp (B/ CA)$;}
\item \tcb{$\tp (A/ C\St_C (B)) \vdash \tp (A/ CB)$;}
\item \tcb{$\tp (A/ C\St_C (A)) \vdash \tp (A/ C \St_C (A)B)$.}
\end{enumerate}
\end{lem}

\prf \tcb{The equivalence of (1) and (2) is Lemma 3.8 (i) of \cite{hhm}.
The equivalence of (1) and (3) follows from the elementary fact that
$\tp (A / D) \vdash \tp (A / DB)$ is always equivalent to $\tp (B / D) \vdash \tp (B / DA)$,
cf.   \cite{hhm} p.~29.}
\eprf

\begin{rem}\label{std-ext-rem}\tcb{For any subsets $A$, $B$ and $C$ of  $\Uu$,
it is a consequence of stable embeddedness that
$\tp (A/ C\St_C (A)) \vdash \tp (A/ C \St_C (A)\St_C (B))$, as explained in Remark 3.7 of  \cite{hhm}.}
\end{rem}

\begin{defn}\label{dst}A type $\tp (A / C)$ is {\em stably dominated} \index{stably dominated type} if, for any $B$ such that
$\St_C (A)  \dnf_{\St_C (C)} \St_C (B)$, \tcb{the  conditions in \lemref{eqvda} are satisfied}.
\end{defn}

By \cite{hhm} 3.13, if 
$\tp(a / C)$ is stably dominated, then it has  
an $\acl (C)$-definable extension $p$ to $\Uu$; this definable type will also
be referred to as stably dominated.  In fact our focus is really on stably dominated {\em definable types},
and we will rarely refer to {\em types} as stably dominated.

The unique stably dominated extension of $\tp(a/\acl(C))$   will be denoted  
by $\tp(a/\acl(C)) | \Uu$; more generally,    for any $B$ with $\acl(C) \leq B \leq \Uu$, write
$p|B =\tp(a/\acl(C)) | B$.

We say that a stably dominated type $\tp(a/C)$ is {\em stationary} \index{stationary type}
if it has a $C$-definable extension $p$ to $\Uu$.    Equivalently, \[ \dcl(Ca) \meet \acl(C) = \dcl(C). \]

\tcb{One direction of the equivalence follows from the remark immediately following, 
applied to $N=\Uu$; for the other see, e.g.,  \cite{finiteim} Lemma 2.9.}

 For any $|C|^+$-saturated, $|C|^+$-homogenous extension $N$ of $C$, $p |N$ is
 the unique  $\Aut (N / \acl (C))$-invariant extension of
$\tp(a / \acl (C))$.  We will need a slight extension of this:

\begin{lem}  \label{std-ext-2}  Let $p = \tp(a/C)$ be a stably dominated $C$-definable type, $C=\acl(C)$.
Let $C \nsubseteq B=\dcl(B)$.
Assume that, for any $b \in \St_C(B) \m \tcb{C}$, there exists $b' \in B$, $b' \neq b$,
such that  $b$ and $b'$ are $\Aut (B / C)$-conjugate.  
   Then  $p |B$ is  the unique  $\Aut (B / C)$-invariant extension of 
$\tp(a / C)$. \end{lem}

\prf   By hypothesis, $p$ is stably dominated via some $C$-pro-definable function $h$ into $\St_C$. 
Let $q$ be an $\Aut (B / C)$-invariant extension of $\tp(a/C)$, say $q = \tp (d/B)$.  
Let $h_{*} q = \tp (h (d) /  \St_C(B))$ 
 be its pushforward.   
 Let $b$ enumerate the canonical base of  $h_{*}q$, so
 that $h_{*}q$ is the unique nonforking extension to
 $ \St_C(B)$ of $h_{*}q|C(b)$, and any automorphism fixing $q$ also fixes $b$.  As
 $q$ is $C$-invariant, any automorphism of $B/C$ fixes $b$.  
But by assumption, if $b \notin C$  there is an elementary permutation $\si$
 of  $ \St_C(B)$  over $C$ with $\si(b) \neq b$;  
 it follows that $b$ is contained in $C$.  
 Thus $h_{*} q$ does not fork over $C$,
 so   $h_{*} q = h_{*} p$.  
By definition of stable domination, it follows that $q=p$.  \eprf

  \begin{prop}[\cite{hhm} Proposition 6.11]\label{transitive}Assume the types $\tp(a / C)$
and 
$\tp(b / a C)$ are stably dominated,
then $\tp(a b / C)$ is stably dominated.
\end{prop}


\medskip

A formula  $\varphi (x, y)$ is said to {\em shatter} \index{to shatter} a subset $W$ of a model of $T$ if 
for any two  finite disjoint subsets $U,U'$ of $W$ there exists $b$ with $\phi(a,b)$ for $a \in U$, and $\neg \phi(a',b)$
for $a' \in U'$.  Shelah says that a formula $\varphi (x, y)$ has the {\em  independence property}   \index{independence property}
if it shatters arbitrarily large finite sets; otherwise, it has $\NIP$.    Finally, $T$ has $\NIP$ if every formula has $\NIP$. \index{$\NIP$}
 Stable and o-minimal theories are $\NIP$, as is  $\ACVF$.

 If   $\varphi (x, y)$ \tcb{has} $\NIP$ then there exists a positive integer  $k$, such that for any \tcb{finite (or infinite)} indiscernible sequence $(a_1,\ldots,a_n)$ and any $b$ in a model of $T$, $\{i: \phi(a_i,b) \}$ is the union of $\leq k$ convex segments.  
 If  $\{a_1,\ldots,a_n\}$ is an indiscernible 
 set, i.e. the type of $(a_{\si(1)}, \ldots,a_{\si(n)})$ does not depend on $\si \in \Sym(n)$, it follows that
 $\{i: \phi(a_i,b) \}$ has size $\leq k$, or else the complement has size $\leq k$.

\begin{defn}
If $T$ is a  $\NIP$-theory, 
 and  $ p$ is  an $\Aut(\Uu /C)$-invariant
type over $\Uu$, one says that $p$ is {\em generically stable} \index{generically stable type} over $C$ if
it is $ C$-definable and  finitely satisfiable in any model
containing $C$ (that is, for any formula $\varphi (x)$ in $p$ and any  model $M$ containing $C$,
there exists $c$ in $M$ such that $\Uu \models \varphi (c)$). 
\end{defn}

In general, when $p(x),q(y)$ are   $\Aut(\Uu/C)$-invariant types, there exists a unique $\Aut(\Uu/C)$-invariant 
type $r(x,y)$, such that for any $C' \supset C$, $(a,b) \models r(x, y)$ if and only if 
$a \models p|C$ and $b \models q| C(a)$.  This type is denoted $p(x) \tensor q(y)$.  In general $\tensor$ \nomenclature{$\tensor$}{tensor product of invariant types}
is associative but not necessarily symmetric.  We define
$p^n$ by $p^{n+1}=p^n \tensor p$.

  The following characterization of generically stable types in  $\NIP$
 theories from   \cite{hp} will be useful:

\begin{lem} [\cite{hp} Proposition 3.2]\label{equivgs} Assume $T$ has $\NIP$.  An  $\Aut(\Uu /C)$-invariant
type  $ p(x)$ is generically stable over $C$ 
 if and only if   $p^n$ is symmetric with respect to permutations of the variables $x_1,\ldots,x_n$.

  For any formula
 $\varphi (x, y)$, there exists a natural number
 $n$ such that whenever $p$ is generically stable and $(a_1,\ldots,a_N) \models p^N |C$ with $N>2n$, 
 for every $c$ in $\Uu$,
 $\varphi (x, c) \in p$ if and only if
 $\Uu \models  \bigvee_{i_0 < \dots < i_n} \varphi (a_{i_0}, c) \wedge \cdots \wedge \varphi (a_{i_n}, c)$.   \end{lem}
 
 The second part of the lemma is  an easy consequence of the definition of a $\NIP$ formula, or rather the remark
 on indiscernible sets just below the definition.

We remark  that \propref{transitive} also follows from   the characterization of generically stable definable
 types in $\NIP$ theories  as those with symmetric tensor powers in \lemref{equivgs}, cf. \cite{hp}. 

 We also recall the notion of a {\em strongly stably dominated type} from \cite{mst}.  
 These are the stably dominated types that are dominated within a single formula, rather than a type.
 The distinction is analogous to that between $\omega$-stability and stability, or regular and strongly regular types in stability theory.

\begin{defn}\label{ssd}
  \tcb{Say $\tp(a/C)$ is {\em strongly stably dominated}  \index{strongly stably dominated type}
if 
there exists  $\phi(x)  \in \tp(a/\St_C(a))$ such that for any tuple $b$ with
$\St_C(a) \dnf_{\St_C(C)} \St_C(b)$, $\phi$ implies $\tp(a/\St_C(a)b)$.   If $h$ is a $C$-definable function and $h(a)$ are the parameters of $\phi$,
we say $\tp(a/C)$ is strongly stably dominated via $\phi$ and $h$.}  \end{defn}

We say that a definable type $p$ is {\em strongly stably dominated} if for some $A=\acl(A)$ such that $p$ \tcb{is $A$-definable}, $p|A$ is strongly stably dominated.  Note that it follows that $p$ is stably dominated.

\begin{rem}\label{was2}
Assume $\tp(a/C)$ is stably dominated.  Then 
 $\tp(a/C)$ is strongly stably dominated iff $\tp(a/C')$ is isolated for some (or any) $C'$ with \tcb{$\St_C(a)  \nsubset C' \nsubset \St_C$.}
Indeed, by \remref{std-ext-rem},
 \tcb{$\tp(a/\St_C(a)) \vdash \tp(a/ C')$
for any
$C'$ with $\St_C(a)  \nsubset C' \nsubset \St_C$.}
 \end{rem}

  \tcb{For part (3) of the following proposition, we will need a special hypothesis (we refer to the beginning of \ref{gint} for the notion of internality):}

\smallskip
\noindent \tcb{(Sp)  There exists a  sort (or union of sorts) $\kk$ such that $\kk$  is $\omega$-stable, and for any $A$, $\St_A$ is $\kk$-internal, i.e.
any  $A$-definable stable and  stably embedded is $\kk$-internal.}
  
 \smallskip  Note that (Sp) holds in $\ACVF$ with $\kk$ the residue field sort. \tcr{This follows from  \lemref{3.4.11}} since, as recalled in \propref{2.1.3}, in this case the induced structure on $\kk$ is that of a pure 
 algebraically closed field, which is $\omega$-stable.

 \begin{lem}\label{Splem}
 \tcb{Assume \textup{(Sp)} holds. Then, for any $A$, any $c \in \St_A$, and any $A$-definable set $D$  containing $\kk$, or union of such sets, $\tp(c/ A \union D)$ is
 isolated.}
 \end{lem}

\prf \tcb{To see this let $P$ be the solution set of $\tp(c/A \union D)$.  Then $P$ is $\kk$-internal, so the automorphism group $G=\Aut(P)$ is an $\infty$-definable group internal to 
 $\kk$ by Theorem B.1$'$ (1) in \cite{bedlewo}. Since $\kk$  is $\omega$-stable, it follows then from Corollaire 5.19 in \cite{poizat} that $G$ is definable.
 Thus $P$, being a $G$-orbit,  is a definable set.}
 \eprf 
 
\begin{prop}\label{ssd1}\tcb{Let $p$  be a strongly stably dominated definable type, based on $A=\acl(A)$.
 \begin{enumerate} 
\item    $p|A$ is strongly stably dominated.
 \item  If $f$ is a definable function such that $p$ belongs to  its domain, then $f_* (p)$ is  strongly stably dominated.
 \item   Assume \textup{(Sp)} holds.  Let $b \models p|A$, and let $\tp(c/\acl(Ab))$ be strongly stably dominated.   Then so is $\tp(c/A)$. \end{enumerate}}
 \end{prop}
 
 \prf  (1)  If $p$ is based on $A$ and on $A'=\acl(A')$, we have to show that $p|A$   is strongly stably dominated iff $p|A'$ is strongly stably dominated.
We may assume here that $A \nsubset A'$.

\tcb{Let us show $p|A'$ is  strongly stably dominated, assuming the same for $p|A$.   Let $a \models p|A'$.
Now $p|A = \tp (a / A)$ is strongly stably dominated, say via $\phi(x,e)$ with $e \in \St_{A}(a)$; $e=h(a)$ for 
some $A$-definable function $h$.  Now if $b$ is such that 
     $e \dnf_{\St_{A'}(A')} \St_{A'}(b)$, 
     we have to show that  $\phi(x,e)$ implies $\tp(a/A'be)$.  Since
     $e \models h_*p | A'$, and    $e \dnf_{\St_{A'}(A')} \St_{A'}(b)$, we have      $e \models h_*p | A'b$.  
     In particular,  $e \models h_*p | \St_{A}(A'b)$ so 
     $e \dnf_{\St_A(A)} \St_{A}(A'b)$, i.e. 
       $e \dnf_{\St_A(A)}  \St_{A}(A'b)$.   By  stable domination, 
       it follows that   $\tp(a/Ae)$ implies  $\tp(a/A'be)$;  but $\phi(x,e)$ implies $\tp(a/Ae)$,
        so $\phi(x,e)$ implies $\tp(a/A'be)$.}
 
\tcb{Assume now that $p|A'$ is strongly stably dominated.   We have to show that $\tp(a/E)$ is isolated where $E = \St_A(a)$.  Let $E'=\dcl(A' \union E)$, so $\St_{A'}(a) \nsubseteq E'$ by \cite{hhm} 6.10 (iv).
   Then $\tp(a/E')$ is isolated, i.e. $\tp(a/EA')$ is isolated, say by $\psi(x,e,a')$.  But $\tp(a/E)$ implies $\tp(a/Ea')$.  So some $\theta(x,e) \in \tp(a/E)$ implies $\psi(x,e,a')$,
   and this $\theta(x,e)$ clearly isolates $\tp(a/E)$.}

\tcb{(2)  Say $p$ and $f$ are defined over $A$.  Let $c \models p |A$.  Then $\tp(c/ \St_A(c))$ is isolated, so $\tp(f(c)/\St_A(c))$ is isolated.
Since $\St_A(f(c)) \nsubset \St_A(c) \nsubset  \St_A$, and $\tp(f(c)/\St_A(f(c)) \vdash \tp(f(c) / \St_A))$, as noted above it follows that $\tp(f(c)/A)$ is strongly stably dominated.}
 
 
\tcb{(3)  We have $\tp(c / \St_{\acl(Ab)})$ isolated by some formula over $Ab'c'$, where $b' \in \acl(Ab)$ and $c' \in \St_{Abb'}(c)$.
 In particular $\tp(c/bb'c'\St_A)$ is isolated.  By \lemref{Splem} $\tp(c'/ \St_{Abb'})$ is also isolated. By transitivity of stable domination, \propref{transitive},
 $\tp(bb'c/A)$ is stably dominated.  
 Moreover $\tp(bb'c/ \St_A)$ is isolated, since $\tp(b/\St_A)$, $\tp(b'/b,\St_A)$, $\tp(c'/ \St_{Abb'})$  and $\tp(c/bb'c'\St_A)$ are all isolated.
 By \remref{was2}, $\tp(bb'c/A)$ is strongly stably dominated, and by (2) so is $\tp(c/A)$.} \eprf


\section{Review of $\ACVF$}\label{ss2.6}

A valued field consists of a field $K$ together
with a homomorphism $v$ from the multiplicative group
to an ordered  abelian group $\Gamma$, such that $v(x+y) \geq \min{( v(x),v(y))}$, for every $x$ and $y$ in $K^{\times}$.
In this paper we shall  write the law on $\Gamma$ additively.
We shall write $\Gamma_{\infty}$ for $\Gamma$ with an element $\infty$ added with usual conventions, \tcb{namely $\infty$ is larger than any element of  $\Gamma$ and is absorbing for the addition}. \nomenclature{$\Gamma_{\infty}$}{augmented value group sort}
In particular we extend $v$ to $K \rightarrow \G_{\infty}$ by setting $v  (0) = \infty$.
We denote by $\Oo$ or by $R$ the valuation ring, by $\Mm$ the maximal ideal and by $k$ the residue field. \nomenclature{$\Oo$}{valuation ring}\nomenclature{$R$}{valuation ring} \nomenclature{$\Mm$}{maximal ideal}

Now assume  $K$ is  algebraically closed and $v$ is surjective.
The value group $\Gamma$ is then divisible and the residue field $k$ 
is algebraically closed.
We shall denote by $\ACVF$ \nomenclature{$\ACVF$}{the theory of nontrivially valued algebraically
closed fields}
of
algebraically
closed valued fields with nontrivial valuation.
By a classical result of A. Robinson, the completions of $\ACVF$ are the theories $\ACVF_{p_1, p_2}$ of nontrivially valued algebraically 
closed fields of characteristic \tcb{$p_1$} and residue characteristic \tcb{$p_2$}. \nomenclature{$\ACVF_{p_1, p_2}$}{completion of $\ACVF$}
Several quantifier elimination results hold  for $\ACVF$. In particular $\ACVF$ admits quantifier elimination in the 3-sorted language
$\CL_{k, \Gamma}$, with  sorts $\VF$, $\Gamma$ and $\kk$ for the valued field, value group and residue field sorts, with respectively the ring, ordered abelian group and ring language, and additional symbols for the valuation $v$ and the map $\mathrm{Res} : \VF^2 \rightarrow k$
sending $(x, y)$  to the residue of
$x y^{-1}$  if $v (x) \geq v (y)$ and $y \not=0$ and to $0$ otherwise (cf. \cite{hhmcrelle} Theorem 2.1.1).\nomenclature{$\CL_{k, \Gamma}$}{3-sorted valued field language}\nomenclature{$\VF$}{valued field sort}\nomenclature{$\Gamma$}{value group sort}\nomenclature{$\kk$}{residue field sort} Sometimes we shall also write $\val$ instead of $v$ for the valuation. \nomenclature{$\val$}{valuation}
In this paper we shall use the extension $\CL_{\CG}$ of $\CL_{k, \Gamma}$ \nomenclature{$\CL_{\CG}$}{extended valued field language}
considered in section 3.1 of \cite{hhmcrelle}
for which elimination of imaginaries holds.
In addition to sorts $\VF$, $\Gamma$ and $\kk$, there are {\em geometric sorts} $S_n$ and $T_n$, $n \geq 1$. \nomenclature{$S_n$, $T_n$}{geometric sorts}
The sort $S_n$ is the collection of all codes for free rank $n$ $R$-submodules of $K^n$.
For $s \in S_n$, we denote by
$\redu (s)$ the reduction modulo the maximal ideal of the lattice $\Lambda (s)$ coded by $s$. \nomenclature{$\redu (s)$}{reduction of the lattice coded by $s$}
This has $\varnothing$-definably the structure of a rank $n$ $k$-vector space.
We denote by $T_n$ the set of codes for elements in
$\cup_{s \in S_n} \{\redu(s)\}$. Thus an element of $T_n$ is a code
for the coset of some element of $\Lambda (s)$ modulo
$\Mm \Lambda (s)$.
For each $n \geq 1$, we have symbols $\tau_n$
for the functions
$\tau_n : T_n \rightarrow S_n$
defined by $\tau_n (t) = s$ if $t$ codes an element of
$\redu(s)$.
We shall set $\CS = \cup_{n \geq 1}S_n$
and
$\CT = \cup_{n \geq 1}T_n$.
The main result of 
\cite{hhmcrelle} is that 
$\ACVF$ admits elimination of imaginaries in $\CL_{\CG}$.
\tcb{It is also proved in \cite{hhmcrelle} that 
$\ACVF$ admits elimination of quantifiers in $\CL_{\CG}$.}

\medskip

With our conventions, if $C \nsubset \Uu$, we write
$\Gamma (C)$ for \nomenclature{$\Gamma (C)$}{$\dcl (C) \cap \Gamma$}
$\dcl (C) \cap \Gamma$
and
$\kk (C)$ for \nomenclature{$\kk (C)$}{$\dcl (C) \cap \kk$}
$\dcl (C) \cap \kk$.
If $K$ is a subfield of $\Uu$, one denotes by $\Gamma_K$
the value group, thus $\G (K) = \mathbb{Q}\otimes \G_K$.
If the valuation induced on $K$ is nontrivial,
then the model theoretic algebraic closure $\acl (K)$ is a model
of 
$\ACVF$.
In particular the structure $\G (K)$ has definable choice, hence is  Skolemized,
being an expansion by constants of 
a divisible ordered
abelian group (cf. \propref{2.1.3}).
\medskip

We shall denote in the same way a finite cartesian product
of sorts and the
corresponding definable set. For instance, we shall denote
by $\Gamma$ the definable set which to
any model $K$ of $\ACVF$ assigns
$\G (K)$ and by $\kk$ the definable set which to $K$ assigns its residue field.
We shall also sometimes write $K$ for the sort $\VF$.
\medskip

For a field $F$, we denote by $F^{\alg}$ an algebraic closure of $F$. \nomenclature{$F^{\alg}$}{algebraic closure of $F$}
\tcr{By an algebraic variety over a field $F$ (or variety for short), we shall always mean a reduced and separated scheme of finite type over $F$.}
\index{algebraic variety over a field}\index{variety over a field}
\medskip

 The  following follows from the different versions of quantifier elimination (cf. \cite{hhmcrelle}
 Proposition 2.1.3):
\begin{prop}\label{2.1.3}\leavevmode
\begin{enumerate}
\item The definable set  $\Gamma$ is o-minimal in the sense that every definable
subset of $\Gamma$ is a finite union of intervals.
\item Any $K$-definable subset of $k$ is finite or cofinite \textup{(}uniformly in the parameters\textup{)}, i.e.  $k$ is strongly minimal.
\item The definable set  $\Gamma$ is stably embedded.
\item If $A \nsubseteq K$, then $\acl (A) \cap K$ is equal to the field algebraic closure of $A$ in $K$.
\item If $S \nsubseteq k$ and $\alpha \in k$ belongs to $\acl (S)$ in the $K^{\mathrm{eq}}$ sense,
then $\alpha$ belongs to the field algebraic closure of $S$.
\item The definable set $k$ is stably  embedded.
\end{enumerate}
In fact, $\G$ is endowed with the structure of a pure divisible ordered
abelian group and $k$ with the structure of a pure algebraically closed field.
\end{prop}

\begin{lem}[\cite{hhmcrelle} Lemma 2.1.7]\label{2.1.7}
Let $C$ be an algebraically closed valued  field, and let $s \in S_n(C)$, with $\Lambda=\Lambda_s \nsubseteq K^n$
the corresponding lattice.  Then $\Lambda$ is   $C$-definably isomorphic to $R^n$, and the torsor
$\redu(s)$ is $C$-definably isomorphic to $k^n$.   
\end{lem}

A $C$-definable set $D$ is called $k$-{\em internal} \index{$k$-internal} if there exists a finite $F \nsubset \Uu$ such that $D \nsubset \dcl (k \cup F)$
(this is a special case of the more general definition given at the beginning of section \ref{gint}).

We have the following characterizations of $k$-internal sets:

\begin{lem}[\cite{hhmcrelle} Lemma 2.6.2]\label{2.6.2}
Let $D$ be a $C$-definable set. Then the following conditions are equivalent:
\begin{enumerate}
\item $D$ is $k$-internal;
\item for any $m \geq 1$, there is no surjective definable map from $D^m$ to an infinite interval in $\G$;
\item $D$ is finite or, up to permutation of coordinates, is contained in a finite union of sets of the form
$\redu (s_1) \times \cdots \times \redu (s_m) \times F$, where $s_1$, \dots, $s_m$
are $\acl (C)$-definable elements of $\CS$ and $F$ is a $C$-definable finite set of tuples from
$\CG$.
\end{enumerate}
\end{lem}

For any parameter set $C$, let $\mathrm{VC}_{k, C}$ be the many-sorted structure whose sorts are $k$-vector spaces $\redu( s)$ \nomenclature{$\mathrm{VC}_{k, C}$}{many-sorted structure whose sorts are $k$-vector spaces}
with $s$ in $\dcl (C)\cap \CS$. Each sort 
$\redu( s)$ is endowed with a $k$-vector space structure. In addition, as its $\varnothing$-definable relations, $\mathrm{VC}_{k, C}$
has all $C$-definable relations on products of sorts.

By Proposition 3.4.11 of \cite{hhmcrelle}, we have:
\begin{lem}[\cite{hhmcrelle} Proposition 3.4.11]\label{3.4.11}
Let $D$ be a $C$-definable set of $K^{\mathrm{eq}}$. Then the following conditions are equivalent:
\begin{enumerate}
\item $D$ is $k$-internal;
\item $D$ is stable and stably embedded;
\item $D$ is contained in $\dcl (C \cup \mathrm{VC}_{k, C})$.
\end{enumerate}
\end{lem}

By combining \propref{2.1.3}, \lemref{2.1.7} and \lemref{3.4.11}, one sees that (over a model) the $\phi$-definition of a stably dominated type 
factors through some function into $k^n$, where $k$ is the 
 residue field.   

\begin{cor}\label{useful-cor}
Let $C$ be a model of $\ACVF$, let
$V$ be a $C$-definable set and let $a \in V$.
Assume $p = \tp (a / C)$
is a stably dominated type.  Let $\phi(x,y)$ be a formula over $C$.
  Then there exists  a $C$-definable map $g : V \rightarrow k^n$ and a formula $\theta$ over $C$
such that, if
$g (a) \dnf_{\kk (C)} \St_C (b)$, then $\phi(a,b)$ holds if and only if $\theta(g(a),b)$.   \end{cor}

\tcb{The following lemma from \cite{hhmcrelle} is also useful:}

\begin{lem}[\cite{hhmcrelle} Lemma 3.4.12]\label{3.4.12}
\tcb{If $B = \acl (B)$, then, for any $\alpha \in \G$,
$\acl (B \alpha) = \dcl (B \alpha)$.}
\end{lem}

\section{$\G$-internal sets}\label{gint}

Let $Q$ be an $F$-definable set.
An $F$-definable set $X$ is {\em $Q$-internal} \index{$Q$-internal} if there exists $F' \supset F$, and an $F'$-definable
surjection $h: Y \to X$, with $Y$ an $F'$-definable subset of $Q^n$ for some $n$.  When $Q$ is stably embedded
and eliminates imaginaries, as is the case of $\G$ in $\ACVF$, we can take $h$ to be a bijection, by factoring out the kernel.  If one can take $F'=F$ we say that $X$ is {\em directly $Q$-internal}. \index{directly $Q$-internal}
\tcb{We shall say an iso-definable subset of a pro-definable set is $Q$-internal if it is pro-definably isomorphic to some $Q$-internal set.}
\index{$Q$-internal iso-definable}

 In the case of $Q=\G$
in $\ACVF$, we mention some equivalent conditions. \index{$\G$-internal} \index{directly $\G$-internal}

\begin{lem}\label{omin0} Let $X$ be an $F$-definable set. The following conditions are equivalent:
\begin{enumerate} \item  $X$ is $\G$-internal;
\item  $X$ is internal to some o-minimal  definable linearly ordered set;
\item $X$ admits a definable linear ordering;
\item  every stably dominated type on $X$ \textup{(}over any base set\textup{)} is constant \textup{(}i.e. contains a formula $x=a$\textup{)};
\item there exists    an $\acl(F)$-definable injective map
  $h: X  \to \G^n$ for some $n \geq 0$.  
\end{enumerate}
\end{lem}

\prf The fact that  (2) implies (3) follows easily from  elimination of imaginaries in $\ACVF$: by inspection of the geometric
sorts,   the only o-minimal one is $\G$ itself.      Condition (3) clearly implies (4) by the symmetry 
property of generically stable types $p$:  $p(x) \tensor p(y)$ has $x \leq y$ if and only if $y \leq x$, hence $x=y$.  
We now prove that (4) implies (5) using elimination of imaginaries in $\ACVF$, and inspection of the geometric sorts. 
Namely, let $A=\acl(F)$ and let $c \in Y$.  Assuming (4), let us show that $c \in \dcl(A \cup \G)$. This reduces
to the case that $\tp(c/A)$ is unary in the sense of section 2.3 of  \cite{hhmcrelle};
for if  $c=(c_1,c_2)$ and the implication holds for $\tp(c_2/A)$ and
for $\tp(c_1/A(c_2))$  we obtain $c_2 \in \acl(A,\G,c_1)$; it follows that (4) holds for $\tp(c_1/A)$, so $c_1 \in \dcl(A,\g)$
and the result follows since $\acl(A,\g)=\dcl(A,\g)$ for $\g \in \G^m$ by \lemref{3.4.12}.   So assume $\tp(c/A)$ is unary, i.e. it is the 
type of a sub-ball $b$ of a free  $\Oo$-module $M$.   
The radius $\g$ of $b$ is well-defined.  Now $\tp(c/A,\g(b))$ is a type of balls of constant radius; if $c \notin \acl(A,\g(b))$ then there are infinitely many balls realizing this type, and 
their union fills out a set containing a larger closed sub-ball.  In this case the generic type of the closed sub-ball
induces a stably dominated type on a subset of $\tp(c/A,\g(b))$, contradicting (4).  Thus $c \in \acl(A,\g(b)) = \dcl(A,\g(b))$.  
This provides an $\acl(F)$-definable surjection
from a definable subset  of some $\G^n$ to $X$. Using definable Skolem functions, one obtains
a definable injection from  $X$ to some $\G^n$.

It remains to prove that (1) implies (2) and (5) implies (1), which is clear.     
  \eprf
 

Let $U$ and $V$ be definable sets.
A definable map $f: U \to V$ with all fibers $\G$-internal is called
a {\em $\G$-internal cover}. \index{$\G$-internal cover}
 If $f: U \to V$ is an $F$-definable map, such that
for every $v \in V$ the fiber is $F(v)$-definably isomorphic to a definable set in $\G^n$,
then by compactness and stable embeddedness of $\G$, $U$ is isomorphic over
$V$ to a fiber product $V \times_{g,h} Z$, where $g: V \to Y \nsubseteq \G^m$, 
and $Z \nsubseteq \G^n$, and $h: Z \to Y$.    We call such a cover {\em directly $\G$-internal.} \index{directly $\G$-internal cover}
 
Any finite cover of $V$ is  $\G$-internal, and so is any  directly $\G$-internal cover.

\begin{lem} \label{omin1} Let $V$ be a definable set in $\ACVF_F$.  Then any $\G$-internal cover $f: U \to V$
is isomorphic over $V$
to a finite  disjoint union of sets which are a 
 fiber product over $V$ of a finite cover and a directly
$\G$-internal cover.  \end{lem}

\prf  It suffices to prove this at a complete type $p=\tp(c/F)$ of $U$, since 
the statement will then be true (using compactness) above a 
(relatively) definable neighborhood of $f_{*}(p)$, and so (again by compactness, on $V$)
everywhere.    Let
$F'  = F(f(c))$.  
 By assumption, $f \inv (f(c))$ is $\G$-internal.  So over $F'$ there exists a finite
definable set $H$, for $t \in H$ an $F'(t)$-definable bijection $h_t: W_t \to U$, with
$W_t \nsubseteq \G^n$, and $c \in \mathrm{Im}(h_t)$.  We can assume $H$ is an orbit of
$G=\Aut(\acl(F')/F')$.  In this case, since $\G$ is linearly ordered,  $W_t$ cannot depend on $t$, 
so $W_t=W$.  Similarly let $G_c = \Aut(\acl(F)(c) / F(c)) \leq G$.  Then \tcb{the element}
$h_t \inv(c)$ \tcb{of $W$} depends only on the $G_c$-orbit of $h_t$.  Let $H_c$ be such
an orbit (defined over $F(c)$), and set $h\inv (c)=h_t \inv (c)$ for $t$ in this orbit and some $h \in H_c$.    Then $H_c$ has a canonical code $g_1(c) $, and
we have $g_1(c) \in \acl(F(f(c)))$, and $c \in \dcl(F(f(c),g_1(c), h \inv (c)))$.
   Let $g(c) = (f(c),g_1(c))$.   
Then $\tp(g(c)/F)$ is naturally a finite cover of $\tp(f(c)/F)$, and $\tp(f(c),h \inv (c)/F)$ 
is a directly $\G$-internal cover.  
\eprf

We write $\mathrm{VF}^*$ for 
$\mathrm{VF}^n$ when we do not need to specify $n$ (similarly for $\VF^* \times \G^*$).

\begin{lem} \label{omin2} Let  $F$ be a definably closed substructure of  $\VF^* \times \G^*$, let 
 $B \nsubseteq \VF^m$ be $\ACVF_F$-definable, and let
$B'$ be a definable set in any sorts  \textup{(}including possibly imaginaries\textup{)}.  
 Let $g: B' \to B$ be a 
definable map with finite fibers.  Then there exists a definable $B'' \nsubseteq \VF^{m+\ell}$
and a definable bijection $B' \to B''$ over $B$. \end{lem}


\prf  By compactness, working over $F(b)$ for $b \in B$,
 this reduces to the case that $B$ is a point.  So $B'$ is a finite   $\ACVF_F$-definable
 set, and we must show that $B'$ is definably isomorphic to a subset of $\VF^{\ell}$.
Now we can write
 $F=F_0(\g)$ for some $\g \in \G^*$ with
 $F_0 = F \meet \VF$.
   By  \lemref{3.4.12}, $\acl(F)= \tcb{\dcl}(\acl(F_0)(\g))$.
So $B' = \{f(\g): f \in B''\}$ where $B''$ is some finite $F_0$-definable set of functions on $\G$.
Replacing $F$ by $F_0$ and $B'$ by $B''$, we may assume $F$ is a field.

\medskip

\begin{claim}
$\acl(F) = \dcl(F^{\alg})$.
\end{claim}

\begin{proof}[Proof of the claim] This is clear if $F$ is not
 trivially valued since then $F^{\alg}$ is an elementary substructure of $\Uu$.
 
   When $F$ is trivially valued, suppose $e \in \acl(F)$; we wish to show that $e \in \dcl(F^{\alg})$;
   we may assume $F=F^{\alg}$.  The easiest proof is by inspection of the geometric imaginaries:  the only
   $F$-\tcb{algebraic} sublattice of $K^n$ is $\Oo^n$, and the elements of the sort $T_n$ above it are indexed by $k^n$.
 (Here is a sketch of  a more direct proof: let $t$ and $t'$ be elements with $0 < \val(t) \ll \val(t')$.  Then $e \in \dcl(F(t)^{\alg})$
 and $e \in \dcl(F(t')^{\alg})$ by the nontrivially valued case.  But by  the stationarity lemma (\cite{hhm} 8.11), 
 $\tp((e,t) /F) \union \tp((e,t')/F)$ generates $\tp ( (e,t),(e,t') /F)$, forcing $e \in \dcl(F)$.)   
  \end{proof}
  
  Now we have $B' \nsubseteq \acl(F) = \dcl(F^{\alg})$.  
 Using induction on $|B'|$ we may assume $B'$ is irreducible, and also admits no nonconstant
 $\ACVF_F$-definable 
  map to a smaller definable set.  If $B'$ admits a nonconstant definable
 map into $\VF$ then it must be 1-1 and we are done.    Let $b \in B'$ and let $F'= \mathrm{Fix} (\Aut(F^{\alg}/F(b)))$.
 Then $F'$ is a field, and if $d \in F' \m F$, then $d=h(b)$ for some definable map $h$,
 which must be nonconstant since $d \notin F$. 
   If $F'=F$ then by Galois theory, $b \in \dcl(F)$,
 so again the statement is clear.  \eprf
 
 Note that the last part of the argument is valid in any expansion of the theory of fields:  if $C$ is definably closed
 and $F \nsubset C \nsubset \dcl(F')$, with $F'$ an algebraic extension of $F$, then $C = \dcl(C \meet F')$.

\begin{cor} \label{omin3} The composition of two definable maps with $\G$-internal fibers also has $\G$-internal fibers.
In particular if $f$ has finite fibers and $g$ has $\G$-internal fibers then $g \circ f$
and $f \circ g$ have $\G$-internal fibers.\end{cor}

\prf    As pointed out by a referee this follows from characterization (4) in \lemref{omin0}, which is clearly closed under
towers.  Let us also give a direct proof.    We may work over a model $A$.   By \lemref{omin1} and the definition, the class of $\G$-internal covers is the same as compositions $g \circ f$ of definable maps $f$ with finite fibers, and $g$ with
directly $\G$-internal covers.  Hence to show that this class is closed under composition
it suffices to show that  if $f$ has finite fibers and $g$ has directly $\G$-internal covers,
then $f \circ g$ has $\G$-internal fibers; in other words that if $b \in \acl(A(a))$, $a \in \dcl(A \union \{\g\})$ 
with $\g$ a tuple from $\G$, then $(a,b) \in \dcl(A \union \G)$.  But 
$\acl(A,\g)=\dcl(A,\g)$ for $\g \in \G^m$ by \lemref{3.4.12}, so $(a,b) \in \dcl(A \union \G)$.      \eprf

\begin{warning}  The corollary applies to definable maps between definable sets, hence also to iso-definable sets.
However if $f:X \to Y$ is  a map between pro-definable sets and $U$ is a $\G$-internal, iso-definable subset of $Y$, 
we do not know if $f \inv(U)$ must be $\G$-internal, even if $f$ is $\leq 2$-to-one. 
\end{warning}

\section{Orthogonality to $\Gamma$}\label{ss2.8}



 %


Let $A$ be a substructure of $\Uu$.
\begin{prop}\label{equivstd}\textup{(a)} Let $p$ be an
$A$-definable type. The following conditions are equivalent:
\begin{enumerate}
\item $p$ is stably dominated;
\item $p$ is orthogonal to $\Gamma$;
\item $p$ is generically stable.
 \end{enumerate}
\textup{ (b)} A type $p$ over $A$ extends to at most one generically stable $A$-definable type.
\end{prop}
\prf 
The equivalence of (1) and (2) follows from \cite{hhm} 10.7 and 10.8.  Using Proposition 10.16 in
\cite{hhm}, and \cite{hp}, Proposition 3.2(v), we see that (2) implies (3).   (In fact (1) implies (3) is easily
seen to be true in any theory, in a similar way.)  
 To see that (3) implies (2) (again in any theory),  note that if $p$ is generically stable and $f$ is a definable function,
 then $f_{*}p$ is generically stable (by any of the criteria of \cite{hp} 3.2, say the symmetry of indiscernibles).
 Now a generically stable definable type on a linearly ordered set must concentrate on a single point:  a
 two-element Morley sequence $(a_{1},a_{2})$ based on $p$ will otherwise consist of distinct elements,
 so either   $a_{1}<a_{2}$ or $a_{1}>a_{2}$, neither of which
 can be an indiscernible set.  
 The statement on unique extensions follows from \cite{hp}, Proposition 3.2(v).
\eprf

We shall use the following statement, Theorem 12.18 from \cite{hhm}:
\begin{thm}\label{maxcomp}\leavevmode\begin{enumerate}
\item Suppose that $C \leq L$ are 
valued fields with 
$C$ maximally complete,
$\kk (L)$ is a regular extension of 
$\kk (C)$ and
$\Gamma_L / \Gamma_C$ is torsion free.
 Let $a$ be a sequence in $\Uu$, $a \in \mathrm{\dcl} (L)$.
Then
$\tp(a /C \cup \Gamma (Ca))$
is stably dominated.
\item Let $C$ be a maximally complete algebraically closed valued field, and
$a$ be a sequence in $\Uu$. Then $\tp(\mathrm{acl} (Ca)/C \cup \Gamma (Ca))$
is stably dominated.
\end{enumerate}
\end{thm}

We use this especially when $C$ is algebraically closed, so that the conditions on regularity and torsion-freeness are redundant.

In particular, if $C=\acl(C)$ and $\G(C)=\Rr$, every type of elements of $\G$ over $M$ is definable, so every type over $C$ is definable.   This is relevant to Berkovich spaces.   We note another instance of this, when the value group is extended only by infinite or infinitesimal elements.

\begin{lem}\label{lemhahn}  Let $A$ be a divisible abelian group.
Let $B$ be an extension of $A$ containing no proper  extension of $A$ in which $A$ is order-dense.
  Then every type realized in $B$ over $A$ is definable.
\end{lem} 

\prf
Indeed, let $B$ be a finitely generated
extension of $A$.  We show that $B/A$ is definable by induction on $\mathrm{rk}(B/A)$.  If there are any
positive elements  $b \in B$ with $b<a$ for any $0< a \in A$, one can find   such a $b$ with smallest
archimedean class; so any element $b'$ of $B$ with $0< b'<b$ has the form $\alpha b$, $\alpha \in \Qq$.
Let $B'= \{b' \in B:  b \ll |b'| \}$.  Let $B'' = \{b'' \in B: (\exists n \in \Nn) (|b''| < nb )\}$.  
Then $B \cong  B'' \oplus B'$, by induction $B'/A$ is definable,
and as $B''/B'$ is definable by inspection, the result follows.  Similarly, though slightly less canonically,
if there are any $b \in B$ with $b > A$,
find such a $b$ with maximal archimedean class.  Pick a maximal set of $\Qq$-linearly independent
elements $b_i$ in the same archimedean class as $b$.  
Let $B' = \{b' \in B: |b'| \ll b\}$.  
Then again $B = B' \oplus \oplus \Qq b_i$, $\tp(b_1,\ldots,b_m )/ B'$ is definable, and the result follows.  
Finally, if there are no infinitesimal nor any infinite elements in $B$ over $A$, then by assumption 
we have $A=B$, and certainly $B/A$ is definable.
\eprf

\section{$\std{V}$ for stable definable $V$}   \label{ss2.9}

We end   with a description of the set $\std{V}$ of definable types concentrating  on a stable definable $V$, as an ind-definable set.   
The notation $\std{V}$ is compatible with the one that will be introduced in greater generality  in  \ref{ss3.1}.
 Such
a representation will not be possible for algebraic varieties $V$ in $\ACVF$ and so the picture here is  not at all suggestive of the case  
that will mainly interest us, but it is simpler and will be lightly used at one point.

A family $X_a$ of definable sets is said to be {\em uniformly definable} \index{uniformly definable}  in the parameter
$a$ if there exists a definable $X$ such that $X_a= \{x: (a,x) \in X \}$.  
An ind-definable set $X_a$ depending on a parameter $a$ is said to be uniformly
definable in $a$ if it can be presented as the direct limit of a system $X_{a,i}$, with
each $X_{a,i}$ and the morphisms $X_{a,i} \to X_{a,j}$ definable uniformly in $a$.  
If $U$ is a definable set, and $X_u= \limind_i X_{u;i}$ is (strict) ind-definable uniformly in $u$, then the disjoint 
union of the $X_u$ is clearly (strict) ind-definable too.

Recall $\kk$ denotes the residue field sort.
Given a Zariski closed subset $W \nsubseteq \kk^n$, define $\deg(W)$ to be the degree of the Zariski
closure of $W$ in projective $n$-space.   Let $\ZC_d(\kk^n)$ be the family of 
Zariski closed subsets of degree $\leq d$ and let $\IZC_d(\kk^n)$ be the sub-family of
absolutely irreducible varieties.  It is well-known that  $\IZC_d(\kk^n)$ is definable (cf., for instance,
chapter17 of \cite{FJ}). \nomenclature{$\ZC_d(\kk^n)$, $\IZC_d(\kk^n)$}{families of Zariski closed subsets}
These families are invariant under $\GL_n(\kk)$, hence for any definable $\kk$-vector space $V$ of dimension $n$, we may consider  their pullbacks
 $\ZC_d ( V )$ and
 $\IZC_d( V )$  to families of subsets of $V$, under   a  $\kk$-linear isomorphism
$V \to \kk^n$.    Then $\ZC_d(V)$ and $\IZC_d(V)$ are definable, uniformly in any definition of $V$.

\begin{lem} \label{fst1}  If $V$ is a finite-dimensional $\kk$-space, then $\std{V}$ is 
strict ind-definable. 

The disjoint union $D_{st}$ of the $\std{V_{\Lambda}}$ 
with $V_{\Lambda} =\Lambda / \Mm \Lambda $ and where
$\Lambda$
ranges over the definable family
$S_n$ of lattices in $K^n$
 is also  strict ind-definable. \end{lem}

\prf  Since  $\std{V}$  can be identified
with the limit over all $d$ of $\IZC_d(V)$, it  is strict ind-definable uniformly 
in $V$.    The family of lattices $\Lambda$ in $K^n$ is a definable family,
so the disjoint union of $\std{V_{\Lambda}}$ over all such $\Lambda$ is strict
ind-definable.    \eprf

If $K$ is a valued field, we set
$\RV = K^{\times} / 1 + \Mm$ \nomenclature{$\RV$}{$K^{\times} / 1 + \Mm$}
and \tcb{denote by $\rv$ the canonical morphism
$K^{\times} \to \RV$}. \nomenclature{$\rv$}{canonical morphism
$K^{\times} \to \RV$}
So we have an exact sequence of abelian groups
$0 \to k^{\times} \to \RV \to \G \to 0$. For $\g \in \G$,
denote by  $V_\g^{\times}$ the preimage of $\g$ in $\RV$. It is a principal homogeneous
space for $k^{\times}$. It becomes a $k$-vector space $V_\g$ after adding an element $0$
and defining addition in the obvious way.
For $m \geq 0$, we denote by  $\std{\RV^m}$ the set of stably dominated types on $\RV^m$.

\begin{lem} \label{f2} For $m \geq 0$, $\std{\RV^m}$ is strict ind-definable.
\end{lem}

\prf  Note that $\RV$ is the union over $\g \in \G$ of   the principal homogeneous
spaces $V_\g^{\times}$.
For ${\overline{\g}}= (\g_1,\ldots,\g_n) \in \G^n$,  let $V_{\overline{\g}} = \Pi_{i=1}^n V_{\g_i}$.  
Since the  image of  a 
stably dominated type
on $\RV^m$ under the morphism $\RV^m \to \G^m$
is constant, 
any stably dominated type must concentrate on a finite product $V_{\overline{\g}}$. 
Thus
it suffices to show, uniformly in $\overline{\g}  \in \G^n$, that 
$\std{V_{\overline{\g}}}$ is strict ind-definable.   Indeed $\std{V_{\overline{\g}}}$ can be identified 
with the limit over all $d$ of $\IZC_d(V_{\overline{\g}})$.  \eprf 

 \begin{rem}By the above proof,
 the function $\dim$ on $\IZC_d(V_{\overline{\g}})$ induces a constructible function on 
  $\std{\RV^m}$, that is, having definable fibers on each definable piece of $\std{\RV^m}$.
\end{rem}

\section{Decomposition of definable types}\label{decompose}\label{ss2.10}

We seek to understand a definable type in terms of a definable type $q$ on $\G^n$, and the germ of a  definable
map from $q$ to stably dominated types.

\tcb{Let us start by recalling the notion of an $A$-definable germ, cf.  Definition 6.1 in \cite{hhm}. 
 Let $p$ be an $A$-definable type on some $A$-definable set $X$. Let $\varphi(x,y,b)$ be a formula defining a function
 $f_b(x)$ whose domain contains all realizations of $p$.  The {\em germ of $f_b $ on $p$}, \index{germ of a definable map} or 
{\em $p$-germ of $f_b$}, is the equivalence class of $b$ under the equivalence relation $\sim$,
where $b\sim b'$ if the formula $f_b(x)=f_{b'}(x)$ is in $p$. Equivalently, $b\sim b'$ if and only if for 
any $a\models p|Abb'$, $f_b(a)=f_{b'}(a)$. As $p$ is $A$-definable, $\sim$ is also $A$-definable, and the germ of 
$f_b$ on $p$ is a definable object.}

\tcb{Now  assume  $Y = \limproj_i Y_i$ is an $A$-pro-definable set and let $h$ and $h'$ be two 
pro-definable maps 
over $B \supset A$ taking  values in  $Y$  whose domain contains all realizations of $p$.
We say $h$ and $h'$ have the same {\em $p$-germ} \index{germ of a pro-definable map}
 if $h(e)=h'(e)$ when $e \models p | B$.  The $p$-germ of $h$ is the equivalence
 class of $h$.  Thus, $h$ and $h'$ have the same $p$-germ if and only if for every $i$ the maps $h_i$
 and $h'_i$ given by composing $h$ and $h'$ \tcb{with} the projection to $Y_i$ 
 have the same $p$-germ;
 and the $p$-germ of $h$ is determined by the sequence of $p$-germs of the $h_i$.}

Let $p$ be an $A$-definable type.  Define $\mathrm{rk}_\G(p) = \mathrm{rk}_\Qq \G(M(c)) / \G(M)$, \nomenclature{$\mathrm{rk}$}{$\G$-rank of a definable type}
where $A \leq M \models \ACVF$ and $c \models p|M$.  Since $p$ is definable, this rank does not depend on the choice of $M$,
but for the present discussion it suffices to take $M$ somewhat saturated, to make it easy to 
see that $\mathrm{rk}_\G(p)$ is well-defined.

If $p$ has rank $r$, then there exists a definable function to $\G^r$ whose image is not contained in a smaller
dimensional set.  We show first that at least the germ of such a function can be chosen $A$-definable.

\begin{lem}  \label{defty1}  Let $p$ be an $A$-definable type and set $r=\mathrm{rk}_\G(p)$.  Then there exists a nonempty $A$-definable set $Q''$ and for
$b \in Q''$ a function $\g_b= ((\g_b)_1,\ldots,(\g_b)_r)$ \tcb{from a definable set containing $p$} into $\G^r$, definable uniformly in $b$, such that 
\begin{enumerate}
\item If $b \in Q''$ and   $c \models p |A(b)$ then the image of $\g_b(c)$ in  $\G(A(b,c))/\G(A(b))$ is a  $\Qq$-linearly independent $r$-tuple;
\item
 if $b,b' \in Q''$ and $c \models p |A(b,b')$ then $\g_b(c)=\g_{b'}(c)$.
 \end{enumerate}
\end{lem}

\prf 
\tcb{Take $M$ sufficiently saturated} and consider an $M$-definable function $\g=(\g_1,\ldots,\g_r)$
into $\G^r$, such that if $c \models p |M$ then $\g_1(c),\ldots,\g_r(c)$ have
$\Qq$-linearly independent images in $\G(M(c))/\G(M)$.   Say $\g=\g_a$ \tcb{is defined over $A (a)$ with $a$ a finite tuple} and let $Q=\tp(a/A)$.
If $b \in Q$ there exist a unique $N(a,b) \in \GL_r(\Qq)$ and $\g'=\g'(a,b) \in \G\tcb{(M)}^r$ such that
for   $c \models p | M(b)$, $\g_b(c) = N (a, b) \g_a(c) + \g'$.  By compactness \tcb{and because $p$ is $A$-definable}, as $b$ varies
the matrices $N(a,b)$ vary among a finite number of possibilities $N_1,\ldots,N_k$
\tcb{and there exist 
finitely many $A$-definable functions $\gamma'_i:Q\times Q\to\Gamma^r$ 
such that whenever $a',b\in Q$, then for some $i$, and for any $c\models 
p|M(a',b)$, $\gamma_{a'}(c)=N_i\gamma_b(c)+\gamma'_i(a',b)$. By compactness 
again, the same is true for some $A$-definable \tcb{set} $Q'$ containing $Q$.}

\tcb{Consider the $A$-definable} equivalence relation $E$ on $Q'$ defined by   $b' E b$ if $(d_p x) (\g_{b'}(x) = \g_b(x))$.  Then
by the above discussion, $Q'/E \nsubseteq \dcl(A(a),\G)$
(in particular  $Q'/E$ is $\G$-internal, cf. \ref{gint}).
By \lemref{omin0}  it follows that $Q'/E \nsubseteq \acl(A,\G)$,
and there exists a definable map $g: Q'/E \to \G^{\ell}$ with finite fibers.   

We can consider the following  partial orderings on $Q'$:  $b' \leq_i b$ if and only if $(d_p x)((\g_{b'})_i(x) \leq (\g_b)_i(x))$.
These induce partial orderings on $Q'/E$, such that if $x \neq y$ then $x <_i y$ or $y<_i x$ for some $i$.
This permits a choice of an element from any given finite subset of $Q'/E$; thus
the map $g$ admits a definable section \tcb{$e$}.    

It follows in particular there exists a nonempty  $A$-definable subset $Y \nsubseteq \G^{\tcb{\ell}}$ and for $y \in Y$
 an  element $e(y) \in Q'/E$.  If $Y$ has an $A$-definable element then there exists an $A$-definable 
  $E$-class in $Q'/E$;  let $Q''$ be this class.  This is always the case   unless  
  $\G(A) = (0)$, and $Y \tcb{\subset} (0)^{\tcb{\ell}_1} \times  (0,\infty)^{\tcb{\ell}_2}  \times (-\infty,0)^{\tcb{\ell}_3}$, with $\tcb{\ell}_2+\tcb{\ell}_3 > 0$; but we give another argument that works in general.
  
  For $y \in Y$
we have a $p$-germ of a function $\g[y]$ into $\G^r$, and the germs of $y,y' \in Y$
differ by an element $(M(y,y'), d(y,y'))$ of $\GL_r(\Qq) \ltimes \G^r$.  It is 
easy to cut down $Y$ so that $M(y,y')=1$ for all $y,y'$.  Indeed, let $q$ be any definable
type on $Y$; then for some $M_0 \in \GL_r(\Qq)$, for $y \models q $ and $y' \models q |y$
we have $M(y,y')=M_0$. It follows that $M_0^2=M_0$ so that $M_0=1$, \tcb{hence we  may impose that $(d_qy')(M(y,y')=1)$
holds on $Y$}.  Now    $d(y,y'')=d(y,y')+d(y',y'')$.
Pick $a \in Y$, and let $d_a(y) = d(y,a)$; then 
we have $d(y,y') = d_a(y) -d_a(y')$.  Let $\tcb{\g_0} \in \G_\infty$ be some $A$-definable limit point of $Y$.  
(Such a point exists by induction on dimension; consider the boundary.)  Then $d_a$ has a finite number of limit
values at $\tcb{\g_0}$, being piecewise linear; let $c_a$ be the smallest of them.  So $d'_a = d_a  - c_a$ still satisfies
$d(y,y') = d'_a(y) - d'_a(y')$, and now $0$ is a limiting value of $d'_a$.  Any conjugate $d'_{a'}$ of $d'_a$ differs from $d'_{a}$
by a constant, and only finitely many constants are possible (since both functions have $0$ as a limit value at $\tcb{\g_0}$).  Thus $d'_{a}$
has only finitely many conjugates, so it is $\acl(\tcb{A})$-definable; as above it follows that it is \tcb{$A$}-definable. Set $d'=d'_a$ and replace
each germ $\g[y]$ by $\g[y] - d'(y)$.  The result is another family of germs 
 with $M(y,y')=1$ and $d(y,y')=0$.  This means that the germ does not depend on the choice of 
 $y \in Y$.   
\eprf

\begin{lem}\label{defty3}
Let $p$ be an $A$-definable type on some $A$-definable set $V$ and set
$r = \mathrm{rk}_{\G} (p)$.
There  exists an $A$-definable germ of maps $\delta : p \to \G^r$ of maximal rank. 
Furthermore for any such $\delta$ the definable type
$\delta_{*} (p)$ is 
 $A$-definable.
\end{lem}

\prf The existence of the germ $\delta$ follows from \lemref{defty1}.
It is clear
that any two such germs differ by composition with an element of $\GL_r(\Qq) \ltimes \G(A)^r$.
So,  if one fixes such a germ, it is represented by any element of the $A$-definable family $(\g_a: (a \in Q''))$ in \lemref{defty1}.  
The definable type $\delta_{*} (p)$ on $\G^r$ does not depend on the choice
of $\delta$ within this family, hence $\delta_{*}(p)$ is an $A$-definable type.
\eprf

In the remainder of this section, we will use the notation
 $\std{V}$ for the space of stably dominated types on $V$, for  $V$ an
$A$-definable set, introduced in  \ref{ss3.1}. In \thmref{prodef} we prove that 
$\std{V}$ can be canonically identified with a strict $A$-pro-definable set.
 \tcb{More generally, if $V = \limproj_i V_i$ is an $A$-pro-definable set, we denote by
  $\std{V}$ the set of stably dominated types on $V$. Note that
    $\std{V}$ is canonically isomorphic to $\limproj_i \std{V_i}$, hence is $A$-pro-definable.}

\begin{lemdef}\tcb{Let $V$ and $W$ be $A$-definable sets, or $A$-pro-definable sets.}
If $q$ is an $A$-definable type on  $V$
and $h : V \rightarrow \std{W}$ is an $A$-\tcb{pro-}definable map, 
there exists a unique
$A$-definable type $r$ on $W$ such that for any model  $M$ containing $A$, if $e \models q |M$
and $b \models h(e) | Me$ then $b \models r | M$.  We refer to $r$ as the integral $\int_q h$ \index{integral of a map along a definable type}
of 
$h$ along $q$.  As by definition $r$ depends only o\tcb{n} the $q$-germ $\underline{h}$ of $h$, we set
$\int_q \underline{ h} := \int_q h$. \nomenclature{$\int_q h$}{integral of a map $h$ along a definable type $q$}
\end{lemdef}

 Note that  for $h$ as above, if the $q$-germ 
  $\underline{h}$ is $A$-definable (equivalently $\Aut(\Uu/A)$-invariant),
  then so is $r$; again the definition of $r$ depends \tcb{only} on $\underline{h}$
  hence if $\underline{h}$ is $\Aut(\Uu/A)$-invariant then so is $r$ (even if $h$ is not).

\begin{rem}\label{doublest}One can consider 
the space $\doublewidehat V$ of stably dominated
types on the strict pro-definable
set $\std{V}$, for $V$ a definable set.
There is a canonical \tcb{pro-definable} map
$ \vartheta : \doublewidehat{V} \to \std{V}$
sending
a stably dominated
type $q$ on 
$\std{V}$ to $\vartheta (q) = \int_{q} \text{id}_{\std{V}}$.
So $\vartheta (q)$ is a definable type, and by  \propref{transitive} it is stably dominated. 
\end{rem}

The following key \thmref{defty2} states that any definable type may be viewed as an integral of stably dominated types along some definable type on $\G^r$.     The proposition states the existence of certain $A$-definable germs
of pro-definable functions; there may be no $A$-definable function with this germ.

\begin{thm} \label{defty2} Let $p$   be an $A$-definable type on some $A$-definable set $V$ 
and let 
$\delta : p \to \G^r$ be as in \textup{\lemref{defty3}}. Let $s=\delta_* p$.  
There exists an $A$-definable $s$-germ  $f: s \to \std{V}$ such that
$p = \int_{s} f$.
 \end{thm}

\prf  
Let $M$ be a maximally complete model, and let $c \models p|M$, $t=\d(c)$.  Then $\G(M(c))$ is generated over $\G(M)$ by $\d(c)$.  By \cite{hhmcrelle}, Corollary 3.4.3  
and Theorem 3.4.4,  
$M(t) := \dcl(M \union \{t)\})$ is algebraically closed.  
By \thmref{maxcomp}
$\tp(c / M(t) )$ is stably dominated, hence extends to a unique element $f_M(t)$ of
 $\std{V} (M(t))$.  We will show that $f_M$ does not depend on $M$.
 
   Let $M \leq N \models \ACVF$, with $N$ large and saturated.  We first show that $f_M=f_N$.  Let $c \models p|N$, $t=\d(c)$; so
   $t \models s|N$.    We will show that the homogeneity hypotheses of \lemref{std-ext-2} hold.
Consider an element $b$ of $N(t) \m M(t)$; it has the form $h(e,t)$ with  $e \in N$.   Let $\bar{e}$ be the class
of $e$ modulo the definable equivalence relation: $x \sim x'$ if $(d_st)(h(x,t)=h(x',t))$.
Since $b$ is not $M(t)$-definable, $\bar{e} \notin M$.  Hence there exists $e' \in N$ with $\tp(e'/M)= \tp(e/M)$, 
but $e' \not \sim e$; and there exists an automorphism of $N$ over $M$, taking $e$ to $e'$ \tcb{which} may be extended to an automorphism of $N(t)/M(t)$,
taking $b$ to $b'=h(e',t) \neq b$ \tcb{with} $\tp(b'/M(t))=\tp(b/M(t))$.   Since $\tp(c/N(t))$ is $\Aut(N(t)/M(t))$-invariant, 
by \lemref{std-ext-2},  $\tp(c/N(t)) = f_M(t) | N(t)$.   Hence $f_M(t) | N(t) = f_N(t) | N(t)$; but as above $N(t)$ is
algebraically closed, so two stably dominated types based on $N(t)$ and with the same restriction to $N(t)$ must be equal; hence
$f_M(t)=f_N(t)$, so $f_M=f_N$.

  Given two maximally complete fields
$M$ and $M'$  we see by choosing $N$ containing both that $f_M(t)|N(t)=f_{M'}(t)|N(t)$;
another use of \lemref{std-ext-2}, this time over $N(t)$ and with $\Uu$ as the homogeneous larger model,
gives $f_M(t)=f_{M'}(t)$.   So $f_M(t)$ does not depend on $M$ and  can be denoted
$f(t)$.    We obtain a \tcb{pro-}definable map $f: P \to \std{V}$, where $P=\tp(t /A)$. 
  The $\delta_{*}(p)$-germ of this function $f$ does not depend on the choice
of $\delta$.  It follows that the germ is $\Aut(\Uu/A)$-invariant, hence $A$-definable; and by construction
we have $p = \int_{\delta_{*} (p)} f$.
\eprf

\section{Pseudo-Galois \tcb{coverings}}\label{pseudogalois}    We finally recall a notion of Galois covering at the level of points;  it is essentially the notion of a Galois covering in the category of varieties in which 
radicial morphisms (EGA I, (3.5.4)) are viewed as invertible.  \tcb{Recall a morphism of schemes $V \to W$ is {\em radicial} \index{radicial morphism} if for every field $K$, the morphism
$V (K) \to W (K)$ is injective.}

Following \cite{VSF} p.~52, we call
  a finite surjective morphism $ Y \rightarrow X$ of integral separated noetherian schemes
 a {\em pseudo-Galois covering} \index{pseudo-Galois covering} if the 
 field extension $F(Y) / F(X)$ is normal and 
 the canonical group homomorphism
$\Aut_X (Y) \to  \mathrm{Gal} (F(Y), F(X))$
is an isomorphism, 
where $\mathrm{Gal} (F(Y), F(X))$ is by definition the group
$\Aut_{F (X)} (F (Y))$. 
  Injectivity follows from the irreducibility of $Y$ and the separateness assumption.  

Let $V$ be a normal irreducible  variety over a field $F$. \tcr{Recall that by a variety over $F$ we mean a reduced and separated scheme of finite type over $F$.
  We take normality to imply irreducibility, but sometimes will repeat the  irreducibility condition.
  Let $K'$ be a finite, normal field extension of $F(V)$. Then} the normalization $V'$ of $V$ in $K'$ is a pseudo-Galois covering since the canonical morphism
$\Aut_V (V') \to  G =  \mathrm{Gal} (K', F(V))$ is an isomorphism.   This is a special case of  the functoriality in $K'$ of the map taking $K'$ to the normalization of 
$V$ in $K'$.     The action of
$g \in G$ on $V'$ may be described as follows.  To $g$ corresponds  a rational map $V' \to V'$; let $W_g$ be the graph of this map,
a closed subvariety of $V' \times V'$.  Each of the projections $W_g \to V'$ is birational, and finite.  Since $V'$ is normal, these projections are isomorphisms, so $g$ is the graph of an isomorphism $V' \to V'$.

As observed in \cite{VSF} p.~53,  if  $ Y \rightarrow X$ is a pseudo-Galois covering and $X$ is normal, for any
morphism  $X' \to X$  with $X'$ an integral noetherian scheme, the Galois group
$G=\mathrm{Gal} (F(Y), F(X))$ acts transitively on the components of
$X' \times_X Y$.  
Here is a brief argument.  Note  that if $X'$ is the normalization of $X$ in a finite purely inseparable extension $K'$ of its function field $F (X)$,
the morphism $X' \to X$ is radicial. Indeed one may assume $X = \spec \, A$, $X' = \spec \, A'$ and the characteristic is $p$. For some integer $h$, $F(X)$ contains $K'^{p^h}$ and
an element $x$ of $K'$ lies in  $A'$ if and only if $x^{p^h} \in A$.
It follows that the morphism $Y/G \to X$ is radicial,
hence $G$ is transitive on fibers of $Y/X$.   So there are no proper $G$-invariant subvarieties of $Y$ mapping onto $X$.
It  is clear from Galois theory that $G$ acts transitively on the components of
$X' \times_X Y$ mapping dominantly to $X'$; it follows that the union of these components is 
a $\mathrm{Gal} (F(Y), F(X))$-invariant subset mapping onto $X'$, hence is all of $X' \times_X Y$.  So there are no other components.

If  $Y$ is a finite disjoint union of  nonempty integral noetherian schemes $Y_i$, we say
 a finite surjective morphism $ Y \rightarrow X$ is a pseudo-Galois covering if each
 restriction $Y_i \to X$ is a pseudo-Galois covering. Also,
 if $X$ 
is a finite disjoint union of  nonempty integral noetherian schemes $X_i$,
we shall say $Y \to X$ is a pseudo-Galois covering if its pullback over each $X_i$ is a pseudo-Galois covering.

\chapter{The space  $\std{V}$ of stably dominated types} \label{sec3}

\medskip

{\small \noindent \textbf{Summary.}
The core of this chapter is the study of the space $\std{V}$ of stably dominated types on a definable set $V$.
It is endowed  with a canonical structure of a (strict) pro-definable set
in \ref{ss3.1}. 
Some examples of stably dominated types are given in \ref{some_examples}.
Then, in \ref{topspace} we endow it with the structure of a definable topological space, a notion defined in \ref{deftop}.
 The properties of that definable topology are discussed  in \ref{ss3.4}. In \ref{sp} we study the canonical embedding of $V$ in $\std{V}$ as the set of simple points.
An essential feature in our approach is the existence of a  canonical extension for a definable function on $V$ to $\std{V}$. This is discussed 
in \ref{canon-ext} where continuity criteria are given. They  rely on the notion of v-, g-, and v+g-continuity
introduced in \ref{ggvv}. In \ref{ph} we introduce basic notions of (generalized) paths and homotopies.
In the  remaining  \ref{Smetrics}-\ref{schedis} we introduce notions of use in later chapters:  good metrics, 
 Zariski topology, schematic distance.
 \par\bigskip}

\medskip

\section{$\std{V}$ as a pro-definable set}\label{ss3.1}We shall now work in a big saturated model
$\Uu$ of
 $\ACVF$ in the language $\CL_{\CG}$.
  We fix a substructure $C$ of $\Uu$.
 If $X$ is an algebraic variety defined over the valued field 
part of $C$, we can view $X$ as embedded as a constructible in affine $n$-space,
via some affine chart.   
Alternatively we could make new sorts for $\Pp^n$,
and consider only quasi-projective varieties.  At all events we will treat $X$ as we treat
the basic sorts.   By a ``definable set'' we mean:  a definable subset of some product
of sorts (and varieties), unless otherwise specified.  \index{definable set}

For a $C$-definable set $V$, and any substructure 
$F$ containing
$C$, we denote by 
$\std{V} (F)$ the set of $F$-definable stably dominated types $p$ on 
$V$ (that is such that
$p \vert F$ contains the formulas defining $V$).

We will now construct the fundamental object of the present work,  initially as a pro-definable set.  We will later define a topology on $\std{V}$.  

We show that
there exists a canonical pro-definable set $E$ and a canonical identification
$\std{V}(F) = E(F)$ for any $F$.   We will later denote $E$ as $\std{V}$.  We call $\std{V}$ the {\em stable completion} of $V$.  \index{stable completion}
Here ``stable'' makes reference to  the notion of  stably dominated or generically stable type, and ``completion'' refers to the density of simple points, cf. \lemref{simple}. 
\nomenclature{$\std{V}$}{stable completion of $V$}

\begin{thm}\label{prodef}  Let $V$ be a $C$-definable set. Then there exists 
a canonical  strict $C$-pro-definable set $E$ and a canonical identification
$\std{V}(F) = E(F)$ for any $F$.   \tcb{Moreover, if $f : V \to W$ is a morphism
of $C$-definable sets, the induced map 
$\std{V} \to \std{W}$ is a morphism of 
$C$-pro-definable sets.}
  \end{thm}

\begin{remark}(This very formal remark can be skipped with no loss of understanding.)  
The canonical pro-definable set $E$ described in the proof will
  be denoted as $\std{V}$ throughout the rest of the paper.   
  
If one wishes to bring the choice of $E$ out of the proof and into a formal definition,  a Grothendieck-style
approach can be adopted.   The pro-definable structure of $E$ determines
in particular the notion of a pro-definable map $U \to E$, where $U$ is any 
pro-definable set.    
We thus have a functor from the category of pro-definable
sets to the category of sets, $U \mapsto E(U)$, where $E(U)$ is the set 
of (pro-)definable maps from $U$ to $\std{V}$.   This includes the functor
$F \mapsto E(F)$ considered above:  in case $U$ is a complete type associated with an enumeration of a structure $A$, then $\std{V}(U)$ can be identified with $\std{V}(A)$.  Now instead of describing $E$ we can
explicitly describe this functor.  Then the representing object $E$ is uniquely determined, by Yoneda's lemma, and can be called $\std{V}$.   
Yoneda's lemma also automatically yields the functoriality of the map  
$V \mapsto \std{V}$  from the category of 
$C$-definable  sets to the category
of   $C$-pro-definable sets.   

In the present case,  any reasonable choice of pro-definable 
structure satisfying the theorem will be pro-definably isomorphic to the $E$ we chose,
so the  more category-theoretic approach does not appear to us necessary.  
As usual in model theory, we will say ``$Z$ is pro-definable'' to mean \tcr{that $Z$ can be canonically identified with a pro-definable set $E$}, where no ambiguity regarding $E$ is possible.  
   \end{remark}

One more remark before beginning the proof.  
Suppose $Z$ is a   strict ind-definable set of pairs $(x,y)$, and let
$\pi (Z)$ be the projection of $Z$ to the $x$-coordinate.  If $Z = \union Z_n$
with each $Z_n$ definable, then $\pi (Z) = \union \pi (Z_n)$. Hence $\pi (Z)$ is naturally
represented as an ind-definable set (and is itself strict).

\begin{proof}[Proof of  \thmref{prodef}]
A definable type $p$ is stably dominated if and only if it is generically stable   (\propref{equivstd}).
The definition of $\phi(x,c) \in p$ stated in  \lemref{equivgs} clearly runs
over a uniformly definable family of formulas.  Hence by \lemref{prodefcond},
$\std{V}$ is pro-definable.  

To show strict pro-definability, let $f: V \times W \to \G_{\infty}$ be a definable function.  Write $f_w(v) = f(\tcb{v,w})$,  
and define $p_{*}(f): W \to \G_{\infty}$ by $p_{*}(f) (w) = p_{*}(f_w)$.   
Let
$Y_{W, f}$  be the subset
of $\Fn(W, \G_{\infty})$
consisting of all functions $p_{*} (f)$, for $p$ varying in $\std{V} (\Uu)$.  
By the proof of \lemref{prodefcond} it is enough to prove that $Y_{W,f}$
is definable.  
Since by  pro-definability of $\std{V}$, $Y_{W,f}$ is $\infty$-definable, 
it remains to show that it is ind-definable.

Set $Y = Y_{W, f}$ and consider
the set 
$Z$ of quadruples $(g,h,q, L)$ such that:
\begin{enumerate}
\item $L=k^n$ is a finite-dimensional $\kk$-vector space;
\item
$q \in \std{L}$;
\item $h$ is a definable function $V \to L$ (with parameters);
\item 
$g: W \to \G_{\infty}$ is a function satisfying:     $g(w)=\g$ if and only if  
\[(d_q \bar{v}) ((\exists v \in V)(h(v)=\bar{v}) \wedge (\forall v \in V)( h(v)= \bar{v} \implies f(v,w)=\g))\]
 i.e. 
for $ \bar{v} \models q$, $h \inv( \bar{v})$ is nonempty, and for any $v \in h \inv( \bar{v})$,
$g(w) = f(v,w)$. 
\end{enumerate}

Let $Z_1$ be the projection of $Z$ to the first coordinate.  
Note that   $Z$ is strict ind-definable by \lemref{fst1} 
and hence $Z_1$ is also strict ind-definable.

Let us prove $Y \nsubseteq Z_1$.  
Take $p$ in $\std{V} (\Uu)$, and let $g = p_{*} (f)$.  We have to show that $g \in Z_1$.
Say $p \in \std{V}(C')$, with   $C'$  a model
of
$\ACVF$ and let $a \models p |C'$.    
By \corref{useful-cor} there exists  
a ${C'}$-definable function $h : V \rightarrow L=k^n$ and a formula $\theta$ over $C'$ such that if $C' \nsubseteq B$
and $b, \gamma \in B$, 
if $h (a)  \dnf_{\kk ({C'})} \St_{C'} (B)$,
then $f(a,b)=\gamma$ if and only if $\theta(h(a),b,\gamma)$.  
Let  $q = \tp( h (a) /C')$.  Then (1-4) hold and   $(g,h,q, L)$ lies in $Z$.

Conversely, let $(g,h,q, L) \in Z$; say they are defined over some base set $M$; we may take $M$ to be a maximally complete model of $\ACVF$.  Let $\bar{v} \models q | M$, and pick $v \in V$ with $h(v)=\bar{v}$.
Let $\bar{\g}$  generate  $\G(M(v))$ over $\G(M)$.   
By  \thmref{maxcomp}  
 $\tp(v / M(\bar{\g}))$ is stably dominated. Let $M'  = \acl(M(\bar{\g}))$ (actually $\dcl(M(\bar{\g}))$ is algebraically closed).  Let $p$ be the unique element
of $\std{V} (M')$ such that $p | M' = \tp(v / M')$.  We need not have $p \in \std{V} (M)$,
i.e. $p$ may not be $M$-definable, but since $\kk$ and $\G$ are orthogonal and $k$ is stably embedded, 
$h_{*}(p) $ is $M$-definable.  Thus $h_{*}(p)$ is the unique $M$-definable type 
whose restriction to $M$ is $\tp(\bar{v}/M)$, i.e. $h_{*}(p)=q$.   By definition of $Z$ it follows
that $p_{*}(f) = g$.  Thus $Y=Z_1$ and
$Y_{W, f}$ is strict ind-definable, hence $C$-definable. 

\tcb{Now let $f : V \to W$ be a morphism of $C$-definable sets and denote by
$\std{f} : \std{V} \to \std{W}$ the corresponding map.
For any definable 
map
$g : W \times Z \to \G_{\infty}$, let $\tilde g := g \circ (f \times \Id_Z)$.
Since for any $p \in \std{V}$ we have $p_* (\tilde g) = \std{f} (p)_* (g)$, 
there is a definable inclusion $Y_{Z, \tilde g} \hookrightarrow Y_{Z, g}$.
The composition of $\std{f}$ with the projection $\std{W} \to Y_{Z, g}$
factors through that inclusion, and it follows that $\std{f}$ is a morphism of $C$-pro-definable sets.}
\end{proof}

If $f : V \rightarrow W$ is a morphism of definable sets, we shall denote
by 
$\std{f} : \std{V} \rightarrow \std{W}$
the corresponding morphism \tcb{of pro-definable sets}. Sometimes we shall write $f$ instead of
$\std{f}$. \nomenclature{$\std{f}$}{stable completion of the morphism $f$}

\section{Some examples}\label{some_examples}

\begin{example}\label{exp1}
If $b$ is a
 closed ball in $\Aa^1$, let $p_b \in \std{\Aa^1}$ be the generic type of $b$:\nomenclature{$p_b$}{generic type of the closed ball $b$}
it can be defined by $(p_b)_{*}(f) = \min \{\val (f(x)): x \in b \}$, for any polynomial $f$.  
This applies
even when $b$ has valuative radius $\infty$, i.e. consists of a single point.  
The generic type
of a finite product of closed balls is defined by exactly the same formula.
If $b$ and $b'$ are (finite products of) closed balls, in the notation of 
\remref{tensor}, $p_{b \times b'} = p_b \tensor p_{b'}$. 
By  \cite{hhmcrelle}, 2.3.6, 2.3.8, and 2.5.5, 
 $\std{\Aa^1}$ 
is equal to the set of all generic types of
closed balls of valuative radius in $\Gamma_{\infty}$.
As a set, $\std{\Pp^{1}}$ consists of the disjoint union of
$\std{\Aa^{1}}$ and  the definable type concentrating on the point $\infty$. 
\end{example}

\begin{example}\label{exotic}Let us give examples of a more exotic nature.
\tcb{Let $F$ be  a field and set $K = F (t)$ with valuation trivial on $F$ and $\val (t) = 1$.}
Let $\varphi = \sum_{i = 0}^{\infty} a_i x^i$ be a formal series with coefficients
$a_i \in \tcb{F}$. Assume $\varphi$ is not algebraic.
Let $p_0(x; y)$ 
consist of all formulas over $K$
of the form
\[
\val(y - \sum_{i = 0}^{n} a_i (tx)^i )\geq \tcb{n + 1}.
\]
Then $p_0(x; y) + (p_{\Oo}\vert \Uu)(x)$ generates a complete type $p_\varphi$
which  is a stably dominated type. \tcb{However,  $p_\varphi \in \std{\Aa^2}$ is not
strongly stably dominated  in the sense of  \defref{ssd}.}
\end{example}

 \prf Let $M$ be any   valued field   extension of $K^{\alg}$
such that $\Zz$ is cofinal in $\G(M)$.   
 For a series $\beta= \sum_{i = 0}^{\infty} b_i z^i$,   $b_i \in \Oo_M$, 
 define $p_{0, \beta}$ 
to consist of all formulas
\[ \val( y -  \sum_{i=0}^n  b_i (xt)^i ) \geq \tcb{n + 1}.
\]
 Let $c \models p_\Oo | M$.  First suppose $p_{0, \beta} (c;0)$ holds.  
 Then \[\min_{i \leq n} (\val(b_i) + i )= \val(\sum_{i=0}^n  b_i (ct)^i)  \geq \tcb{n    + 1},\]
 since $c \models p_\Oo | M$.
So $\val(b_i) \geq \tcb{n + 1} -i$.  Letting $n \to \infty$  we see that $b_i=0$; so $\beta =0$.

Next suppose just that $p_0(c;d)$ holds for some $d \in M(c)^{\alg}$.  So $Q(c,d)=0$ for some polynomial $Q \in \Oo_M[x,y]$.
Let $\varphi'=Q(x, \varphi (t x))$ be the power series obtained by substituting $\varphi (t x)$ for $y$. Then $p_{0, \varphi'} (c; 0)$ holds.
Hence by the previous paragraph, $\varphi'=0$, so $\varphi (tx)$ is algebraic, and $\varphi$ is also algebraic.  

Thus, $p_0(c; y)$ defines an infinite intersection $b$ of balls over $M(c)$, with no algebraic point.  Hence $b$ contains no nonempty $M(c)$-definable subset.  So $p_0 + \tp(c/M)$ generates a complete type $p_\varphi$ over $M(c)$. Now assume $M$ is maximally complete and let  $(c, d) \models p_\varphi | M$. \tcb{Since $\Gamma (M (c, d)) = \Gamma (M)$, it follows from} \thmref{maxcomp} that $\tp ((c,d) / M)$ is stably dominated.
\tcb{One has $\mathrm{trdeg}_M M(c, d) = 2$, while the corresponding residue field extension has transcendence degree $1$ by Lemma 2.5.5 of \cite{hhmcrelle}.
By \propref{st-st-dom} it follows that $p_\varphi$ is not strongly stably dominated.}
\eprf

\begin{example}\label{nonabhy}\tcb{By Example 13.1 in \cite{hhm}, which is rather similar to \exref{exotic}, over any valued field $K$, there exist points $p$ of $\std{\Aa^2}$  defined over some extension  $M$ of
$K$ such that if $c \models p|M$, then 
$\mathrm{trdeg}_M M(c) = 2$ while  the residue field extension has transcendence degree $1$.  By \propref{st-st-dom} such points are not 
strongly stably dominated  in the sense of  \defref{ssd}.}
\end{example}

\section{The notion of a definable topological space}\label{deftop}
We will consider topologies on definable and pro-definable sets $X$.  With the formalism 
of the universal domain $\Uu$, we can view these as certain topologies on $X(\Uu)$, in the usual sense. 

If $M$ is a model, the space $X(M)$ will not be a subspace of $X(\Uu)$;   indeed in the case of an order topology, or any Hausdorff Ziegler topology in the sense defined below, the induced topology from a saturated model on a small set is always discrete.  Instead we define
$X(M)$ to be  the topological space
whose underlying set is $X(M)$, and whose topology is generated by sets $U(M)$ with $U$ 
an   $M$-definable open set.   

We will not have occasion to consider $X(A)$ when $A$ is a substructure, which is not a model.  We remark however
that if $\acl(A)=M$ is a model, then the induced topology on $X(A)$ from $X(M)$ is induced by the $A$-definable open sets.
Indeed if $p \in X(A)$ and $p \in U$ with $U$ definable over $M$, let $\bU$ be the intersection of all $\Aut(M/A)$-conjugates
of $U$; then $\bU$ is open, $A$-definable, and $p \in \bU \nsubset U$.

 We will say that a topological space $X$ is {\em definable in the sense of Ziegler}  if the underlying set  $X$ is definable, and there exists a definable family
$B$ of definable subsets of $X$ forming a neighborhood basis at each point.   This allows for
a good topological logic, see \cite{ziegler}.  But it is too restrictive for our purposes.  An algebraic variety with the Zariski topology is not a definable space in this sense; nor is the topology even generated by a definable family.

Let $X$ be an $A$-definable, \tcb{resp.  pro-definable,} set.  Let $\mathcal{T}$ 
be a topology on $X(\Uu)$, and let $\mathcal{T}_d$ be the intersection of
$\mathcal{T}$ with the class of 
relatively
 $\Uu$-definable   subsets of $X$. 
  We will say that $\mathcal{T}$ is an
 {\em $A$-definable topology} \index{definable topology} if it is generated by $\mathcal{T}_d$, and  
 for any $A$-definable family $\mathcal{W}=(W_u:  u \in U)$ of \index{definable topological space}\index{pro-definable topological space}
 relatively
 definable subsets of $X$, $\mathcal{W} \meet \mathcal{T}$ is ind-definable over $A$.     
  The second condition is equivalent to the statement that
   $\{(x,W): x \in W, W \nsubseteq X, W \in \mathcal{W} \meet  {\mathcal{T}} \}$ is ind-definable over $A$.   
   An equivalent definition is that the topology
is generated by 
 an ind-definable family of 
relatively
definable sets over $A$.  We  
will also say that $(X,\mathcal{T})$
is a definable space over $A$, \tcb{resp. a pro-definable space over $A$,} or just that  $X$ is a definable, \tcb{resp. pro-definable,} space over $A$ when there is no ambiguity about $\mathcal{T}$.  
We say $X$   is a (pro-)definable space if it is a (pro-)$A$-definable  space for some small $A$.  As usual the smallest such $A$ may be recognized Galois theoretically.  \index{definable space} \index{pro-definable space}

If $\mathcal{T}_0$ is any ind-definable family of relatively definable subsets of $X$, the set
$\mathcal{T}_1$  of finite intersections of elements of $\mathcal{T}_0$ is also ind-definable.  
Let $\mathcal{T}$ be the family of   subsets of $X(\Uu)$ that are unions of sets $Z(\Uu)$, with $Z \in \mathcal{T}_1$.
Then $\mathcal{T}$ is a topology on $X(\Uu)$, generated by the
relatively
 definable sets within it.  
By compactness, a
relatively
 definable set $Y \nsubseteq X$ is in $\mathcal{T}$ if and only if  for some definable
$T' \nsubseteq \mathcal{T}_1$,  $Y$ is a union of sets $Z(\Uu)$ with $Z \in T'$.  It follows that the topology $\mathcal{T}$  generated by $\mathcal{T}_0$ is a definable
topology.  
  In the above situation, note also that if $Y$ is $A$-relatively definable, then $Y$ is an $A$-definable union of
relatively definable open sets from $T'$. Indeed,
let $Y' = \{Z \in T': Z \nsubseteq Y \}$,
then $Y = \union_{Z \in Y'} Z$.    In general $Y$ need not be a union of sets from $\mathcal{T}_1 (A)$, for any small $A$.

As is the case with groups, the notion of a pro-definable space is more general than that of
pro-(definable spaces).
However the spaces we will consider will be pro-(definable spaces).

When $Y$ is a definable   topological space, and $A$ a base substructure,
the set $Y(A)$ is topologized using the family of $A$-definable open subsets of $Y$.  
We do not use externally definable open subsets (i.e. $A'$-definable for larger $A$)
to define the  topology on $Y(A)$; if we did,
we would obtain the discrete topology on $Y(A)$ whenever $Y$ is Hausdorff.  
The same applies in the pro-definable case;  thus in the next section we shall topologize $\std{X}(K)$ 
using the $K$-definable open subsets of $\std{X}$, restricted to $\std{X}(K)$.  

When we speak of the topology of $Y$ 
without mention of $A$, we mean to take $A=\Uu$, the universal domain; often, any model will also do.

\section{$\std{V}$ as a topological space}\label{topspace}

Assume  that $V$ comes with a  definable topology $\mathcal{T}_V$, and an ind-definable sheaf $\Oo$ of definable functions into $\G_{\infty}$.  We define a topology on $\std{V}$ as follows. A pre-basic open set 
has the form:
$\{p \in \std{O}:  p_{*}(\phi) \in U \}$, where $O \in \mathcal{T}_V$, $U \nsubseteq \Gamma_{\infty}$ is a definable set, open for the order topology,
and $\phi \in \Oo(O)$.  A basic open set is by definition a   finite intersection  of pre-basic open sets. 

When $V$ is an algebraic variety, we take the topology to be the Zariski topology, and
the sheaf to be the sheaf of regular functions composed with $\val$.  

When $X$ is a definable subset of a given algebraic variety $V$, we give $\std{X}$ the subspace topology.

\section{The affine case}\label{ss3.4}
Assume $V$ is a definable subset of some affine variety.
Let $\Fn_r(V, \G_{\infty})$ denote the functions of the form $\val (F)$, where $F$ is a regular function 
on the Zariski closure of $V$.  \nomenclature{$\Fn_r(V, \G_{\infty})$}{functions of the form $\val (F)$,  $F$ a regular function}
By quantifier elimination any definable function \tcb{in $\Fn (V, \G_{\infty})$} is 
piecewise a difference of the form
$\frac{1}{n} f - \frac{1}{m}g$ with $f$ and $g$ in $\Fn_r$ and $n$ and $m$ positive integers.  Moreover, by piecewise we mean
sets cut out by Boolean combinations of sets of the form $ f \leq  g$, where
$f,g \in \Fn_r(V,\G_{\infty})$.  It follows that if $p$ is a definable type and  $p_{*}(f)$ is defined   for $f \in \Fn_r(V,\G_{\infty})$, then $p$ is stably dominated, and determined by $p_{*} | \Fn_r(V \times W,\G_{\infty})$
for all $W$.   
A  {\em basic open} set is defined by finitely many strict inequalities $p_{*}(f) < p_{*}(g)$,
with $f,g \in \Fn_r(V,\G_{\infty})$.  (In case $f = \val (F)$ and
$g = \val (G)$ with 
 $G=0$, this is the same as $F \neq 0$.)  
It is easy to verify that the topology generated by these basic open sets coincides with the definition of the topology on $\std{V}$ above, for
the Zariski topology and the sheaf of functions $\val ( f)$, $f$ regular.

Note that if $F_1,\ldots, F_n$  are regular functions on $V$, and each 
$p \mapsto p_{*}(f_i)$ is continuous, with $f_i = \val  (F_i)$, then 
$p \mapsto (p_{*}(f_1(x)),\ldots,p_{*}(f_n(x)))$ is continuous.   Thus
the topology on $\std{V}$ is the coarsest one such that all $p \mapsto p_{*}(f)$ are continuous,
for $f \in \Fn_r(V,\G_{\infty})$.  So the basic open sets with $f$ or $g$ constant suffice to generate
the topology.

The topology on $\std{V}$ is strict pro-definably generated in the following sense:
for each definable set $W$, one endows $\textrm{Fn} (W, \G_{\infty})$
with the product topology induced by the order topology on
$\G_{\infty}$. Now for a definable
function
$f : V \times W \rightarrow \G_{\infty}$
the topology induced on the  definable set  $Y_{W, f}$
is    generated by a definable family
of definable subsets 
of $Y_{W, f}$ (recall that $Y_{W, f}$  is the subset
of $\Fn(W, \G_{\infty})$
consisting of all functions $p_{*} (f)$, for $p$ varying in $\std{V} (\Uu)$). 
By definition, the pullbacks to $\std{V}$ of the definable open subsets of the
$\textrm{Fn} (W, \G_{\infty})$ generate the topology on $\std{V}$.
 In particular, $\std{V}$ is a 
 pro-definable space 
 in the sense of \ref{deftop}.

 %

When $V$ is a definable subset of an algebraic variety over $\VF$, the topology on $\std{V}$
can also be defined
by gluing the affine pieces.  It is easy to check that this is consistent (if $V'$ is an affine open of the affine
$V$, obtained say by inverting $g$, then any function $\val (f/g)$ can be written $\val (f)-\val (g)$,
hence is continuous on $\std{V'}$ in the topology induced from $\std{V}$).
Moreover, this coincides with the topology defined via the sheaf of regular functions.

\tcb{For any definable set $X$, we
 have an embedding $X \to \std{X}$, taking a point $x$ to the definable type $\tp (x / \Uu)$ concentrating on $x$.}
 \begin{lem}\label{simple0}
If $X$ is a definable subset of $\G^n_{\infty}$ then $X = \std{X}$ canonically.  More generally if $U$ is a definable
subset of $\VF^n$  or a definable subset of an algebraic variety over $\VF$
and $W$ is a definable subset of  $\G^m_{\infty}$, 
then the canonical map
$  \std{U} \times W \to 
\std{ U \times W}$ is a bijection.
\end{lem}

\prf  
Let $f: U \times W \to U$ and $g:  U \times W \to W$ be the projections. If $p \in \std{U \times W}$
we saw that $g_{*} (p) $ concentrates on some $a \in W$;  so $p= f_{*}(p) \times g_{*}(p)$
(i.e. $p(u,w)$ is generated by $f_{*}(p)(u) \union g_{*}(p)(w)$).
\eprf

If $U$ is a  definable subset of an algebraic variety over $\VF$,
we endow
$  \std{U} \times \G^m_{\infty} \simeq
\std{U \times \G^m_{\infty}}$
with the quotient topology
for the surjective mapping
$ \std{ U \times \Aa^m} \to \std{ U \times \G^m_{\infty}}$
induced by 
$\text{id} \times \val$.

We will see below (as a special case of \lemref{prodtop}) that the topology on $\G_\infty = \std{\G_\infty}$ is the order topology, and the topology on  $\std{\G^m_{\infty}} = \G^m_{\infty}$, is  the product topology.

For $\g = (\g_1,\ldots,\g_n) \in \G_\infty^n$, let $b(\g) = \{x=(x_1,\ldots,x_n) \in \Aa_1^n: \val(x_i) \geq \g_i, i=1,\ldots,n\}$.
\tcb{We set $p_\g = p_{b(\g)}$ with the notation from \exref{exp1}.}

 \begin{lem} \label{gtop0}  The map $j: \std{\Aa^n} \times \G_\infty \to \std{\Aa^{n+1}}$, 
$(q,\g) \mapsto q \tensor p_\g$ is continuous for the product topology of $\std{\Aa^n}$ with the order topology on $\G_\infty$.
 \end{lem}
 
 \prf  We have to show that for each polynomial $f(x_1, \ldots, x_n, y)$ with coefficients in $\VF$,
 the map $(p,\g) \mapsto j(p,\g)_{*} \tcb{(\val (f))}$ is continuous.  
 The functions $\min$ and $+$ extend naturally to continuous functions
$\G_\infty^2 \to \G_\infty$.  
Now if $f (x_1, \ldots, x_n, y)$ is a polynomial with coefficients in $\VF$,
there exists a function
$P (\gamma_1, \ldots, \gamma_n, \tau)$ obtained by
composition of  $\min$ and $+$,
and polynomials $h_i$
such that 
\[\min_{\val (y) = \alpha} \val (f (x_1, \ldots, x_n, y)) = 
P (\val (h_1 (x)), \ldots, \val (h_d (x)), \alpha),\] namely, $\min _{\val(y)=\alpha} \val (\sum h_i(x)y^i)
= \min_i (\val (h_i(x)) + i \alpha$).    So $P: \G_\infty^{n+1} \to \G_\infty$ is continuous.  
Hence $j(p,\g)_{*} f = P( p_{*}(h_1),\ldots,p_{*}(h_d), \g)$.  Continuity follows, by composition.
\eprf

\begin{lem}\label{prodtop}If
 $U$ is a definable
subset of $\Aa^n \times \G^{\ell}_\infty$ and $W$ is a definable subset of  $\G^m_{\infty}$, 
the induced topology on   $\std{ U \times W} = \std{U} \times W$
coincides with the product topology.
\end{lem}

\prf  We  have seen that the natural map $\std{ U \times W} \to \std{U} \times W$ is bijective; 
it is clearly continuous, where $\std{U} \times W$ is given the product topology.   To show that it is closed,
it suffices to show that the inverse map is continuous, and we may take
 $U = \Aa^n$ and  $W = \G^m_{\infty}$.    By factoring
$\std{U \times \G^m_{\infty}} \to \std{U \times \G^{m-1}_\infty} \times \G_\infty \to \std{U} \times \G_{\infty}^{m-1} \times \G_\infty$, we may assume $m=1$.  Having said this, by pulling back to $\Aa^{n+ \ell}$ we may assume $\ell=0$.  
The inverse map is equal to the composition of $j$ as in  \lemref{gtop0} with a projection, hence is continuous.    
\eprf


Let $U$ be a definable subset of $V$,  over a structure $A$.    
Say $\std{U}(A)$ is {\em explicitly $A$-open} \index{explicitly $A$-open}
if for any $p \in \std{U}(A)$, 
 there exists  a Zariski open $V'$ with $p \in \std{V'}$, and
  regular functions $G_1,\ldots,G_n$ on $V'$, $g_i = \val (G_i) :   V' \to \G_\infty$ and   open neighborhoods
$E_i$ of $p_*(g_i)$, all defined over $A$, such that $\meet_i g_i \inv(E_i) \nsubset \std{U}$.

\medskip

The following lemma will be used in  Chapter \ref{berkovich} for structures of the form  $\bF=(F,\Rr)$. 
\begin{lem}  \label{rattop}  
Let $\bF$  be any structure consisting of field points and $\G$-points
 including at least one positive element of $\G$.
Let
$V$ be a variety defined over $\bF$ and let $U$ be an  $\bF$-definable subset of $V$. 
  If $\std{U}$ is open  in $\std{V}$, then $\std{U}(\bF)$
is   explicitly $\bF$-open.   
\end{lem}

\prf  Covering $V$ by affine subsets, we may assume $V$ is affine.    Let $F=\VF(\bF)$ be the field points.

We first show that if the statement holds for $(F^{\alg},\G(F^{\alg}))$, then it holds for
$\bF = (F, \G (F))$.  Note that it is enough to show it holds for $(F,\G(F^{\alg}))$ since $g_i^{-1} ((\a,\b)) = (n g_i)^{-1} ((n \a, n \b))$.
Let $p \in \std{U}(F)$.  There exist regular functions $G_1,\ldots, G_n$ over $F^{\alg}$, 
and intervals $I_j$
of $\G_\infty$, defined over $\G(F)$, such that $p \in \meet_j \std{g_j} \inv (I_j) \nsubset \std{U}$, with 
$g_j = \val (G_j)$.    So it suffices to show, for each $j$, that the intersection
of  the Galois conjugates of $\std{g_j}\inv(I_j)$ contains an open neighborhood of $p$ in $\std{V}(F)$.  Let $G=G_j$, $g = g_j$ and $I=I_j$, and let $G^\nu$ be the Galois conjugates of $G$ over $F$,
$g^\nu = \val (G^\nu)$.

Let $b \models p$.  Then the $G^\nu$ are Galois conjugate over $F(b)$, 
$p$ being $F$-definable.  The elements $c_\nu = G^\nu(b)$ are Galois conjugate
over $F(b)$; they are the roots of a polynomial $H(b,y)= \Pi_\nu (y-G^\nu(b))=\sum_m h_\mu(b)y^m$.  For all $b'$ in some $F$-definable Zariski open set $U'$ containing $b$, 
 the set of roots of $H(b',y)$ is equal to $\{G^\nu(b') \} $.  Within $U'$,
the set of $b'$ such that, for all $\nu$,  $g^\nu(b') \in I$ can therefore be written in terms of the Newton polygon
of $H(b',y)$, i.e. in terms of certain inequalities between convex expressions in $\val (h_k(b'))$.
This shows that the intersection
of  Galois conjugates of ${\std{G}}\inv(I)$ contains an open neighborhood of $p$. 

This permits us to assume   $F$ is algebraically closed, as we will do from now on.

Assume first $\bF \nsubseteq \dcl(F)$.  In particular, by assumption, 
  $F$ is not trivially valued and since $\bF=\acl(F)$ is an elementary submodel, the statement is clear.   

We now have to deal with the case that $\bF$ is bigger than $F$; we may 
assume $\bF$ is generated over $F$ by finitely many elements of $\G$,
and indeed,  adding one element at a time, that
$\bF = F(\g)$ for some $\g \in \G$.  Let $c$ be a field element with $\val(c)=\g$;
it  suffices to show that if $U$ is open over $F(c)$, then it is over $F$ too.    Let
$G(x,c)= \sum G_i(x) c^i$ be a polynomial (where $x=(x_1,\ldots,x_n), V \leq \Aa^n$).   
Let $g(p,c)$ be the generic value of $\val (G (x, c))$ at $p$ and
$g_i (p)$ the one of $\val (G_i)$.
Then $g(p,c) = \min_i g_i (p) +i \gamma$.  From this the statement is clear.  
\eprf

When  $\G(\bF) = (0)$, \lemref{rattop} is not true as stated.
Here is a counterexample:
Let $V=\Aa^2$, $U = \{(x,y) \in V: \val (x) < \val (y) \}$, let $p$ be the generic type of $\Oo \times \{0\}$.
If $G(x,y) $ is any polynomial over $F=\Oo_F$, 
then $p_* (G) =0$ unless $y | G$ and then $p_* (G) = \infty$.  The only conditions about $G$  one can form around $p$  over $F$ are, in  case $G=yG_1$,
that $p_*(G) > 0$.    So no $F$-explicit open set can be contained in $U$, since one can always take $0 \ll \val(y) <\val(x)$ to satisfy $p_* (yG_1) >0$.
 But nevertheless, we still have:

\begin{cor}     \label{rattop-triv}  
Let $\bF$  be any structure consisting of field points and $\G$-points.   Let  
$V$ be a variety \tcb{defined over $\bF$} and let $U$ be an \tcb{$\bF$}-definable subset of $V$. 
\tcb{Let $p \in \std{U}(\bF)$.}
    If $\std{U}$ is open  in $\std{V}$, then   there exists   a  definable function  $\a :  V \to \G_\infty^n$, an   open neighborhood 
$E$ of $p_*(\a)$,  and a Zariski open $V'$ with $p \in \std{V'}$, all defined over $\bF$, such that $\a  \inv(E ) \nsubset \std{U}$
is   explicitly $\bF$-open and    $\a$ has the form $(\val (G_1),\ldots, \val (G_n))$ for some regular functions $G_i$ on  $V'$.
\end{cor}

\prf This follows from \lemref{rattop} unless $\G(\bF)=(0)$.  Assume therefore that
$\G(\bF)=0$, so that all elements of $\G$ of positive valuation have the same type over $\bF$.
Let $\gamma$ be such an element.  By
 \lemref{rattop}, there exist $G_1,\ldots,G_n,V',E_\gamma$ as required but over $\bF(\gamma)$.  So
 $G_1,\ldots,G_n,V'$ are defined over $\bF$;  $E_\gamma$ depends on $\gamma$.  Let $E = \union_{\gamma >0} E_\gamma$.
 Then clearly $E$ is open and the statement holds.   
\eprf


\section{Simple points}\label{sp}


Recall that for any definable set $V$, we
 have an embedding $V \to \std{V}$, taking a point $x$ to the definable type $\tp (x / \Uu)$ concentrating on $x$.  The points of the image are said to be {\em simple}.    \index{simple point}

\begin{lem}\label{simple}
Let $X$ be a definable subset of $\mathrm{VF}^n$.
  \begin{enumerate}
\item  The set of simple points  of $\std{X}$ \textup{(}which we identify with
$X$\textup{)} is an iso-definable and relatively definable dense subset of $\std{X}$. If $M$ is a model of $\ACVF$, then $X(M)$ is dense in $\std{X}(M)$.
\item The induced topology on $X$ agrees with the valuation topology on $X$.   
\end{enumerate}
\end{lem}

\prf (1) The fact that $X$ is iso-definable in $\std{X}$ is clear. For relative definability, note that a point
of
$\std{X}$ is simple if and only if each of its projections to $\std{\Aa^1}$ is simple
and that
on $\Aa^1$, the points are a definable subset of the \tcb{set of} closed balls (cf. \exref{exp1}).  
For density, consider (for instance) $p \in \std{X} (M)$ with $p_{*}(f) > \a$. Then
$\val (f( x)) > \alpha \wedge x \in X$ is satisfiable in $M$, hence there exists a simple point $q \in \std{X}(M)$
with $q_{*}(f) > \a$. 

(2)  Clear from the definitions.  The basic open subsets of the valuation topology are
of the form $\val (f(x)) > \a$ or $\val (f(x)) < \a$.
\eprf

\begin{lem}\label{simple2}
 Let  $f: U \to V$ be
a definable map between definable subsets of
$\mathrm{VF}^*$. If $f$
has finite fibers, then  the preimage of a simple point of $\std{V}$ under
$\std{f}$ is simple in
$\std{U}$.
\end{lem}

\prf
It is enough to prove that if $X$ is a finite definable subset of
$\mathrm{VF}^n$, then $X = \std{X}$, which is clear by
(1) of \lemref{simple}.
\eprf

\begin{rem}  \label{tensor}
The natural projection $S_{def, U \times V} \to S_{def, U} \times S_{def, V}$ 
induces a continuous map 
$\std{U \times V} \to \std{U} \times \std{V}$. On the other hand, it
admits a natural
 section,  namely
 $\tensor: S_{def, U} \times S_{def, V} \to S_{def, U \times V}$, which
 restricts to a section of   $\std{U \times V} \to \std{U} \times \std{V}$.
 This map is not continuous in the  logic topology, nor is its  restriction to
$\std{U} \times \std{V} \to \std{U \times V}$ continuous.  Indeed when $U=V$ the pullback of the diagonal 
$\std{\Delta_U}$ consists of 
simple points on the diagonal $\Delta_\std{U} $. But over a model,   the set of simple points
is dense, and hence not closed.  
 \end{rem}


\section{v-open and g-open subsets, v+g-continuity}  \label{ggvv}
\begin{defn}Let $V$ be 
an algebraic  variety over a valued field $F$.  
A definable subset of   $V$  is said to be {\em v-open} if \index{v-open}
it is open for the valuation topology.  
It is called {\em g-open} \index{g-open} if it is defined
by a positive Boolean combination of Zariski closed and open sets, 
and sets of the form $\{u \in U: \val (f(u)) > \val (g(u)) \}$, \tcb{for $f$ and $g$ regular functions on a Zariski open set $U  \subset V$}.
More generally,
if $V$ is a definable subset of an algebraic variety $W$,
a definable subset of $V$ is said to be v-open (resp. g-open) if it is of the form
$V \cap O$ with $O$ v-open (resp. g-open) in $W$.
A definable
subset of $V \times \G_\infty^m$ is called v- or g-open if its  pullback to $V \times \Aa^m$ via $\mathrm{id} \times \val$ is.
\tcb{The complement of a v-open (resp. g-open) subset is said to be {\em v-closed} (resp. {\em g-closed}).}
\index{v-closed}\index{g-closed}
\end{defn}

\begin{rem} If $X$ is $A$-definable, the  regular functions $f$ and $g$ in the definition of  
g-openness  are {\em not} assumed
to be $A$-definable; in general when $A$ consists of imaginaries, no such $f,g$ can be found.
    However when $A=\dcl(F)$ with $F$ a valued field, they may be taken to be $F$-definable,  by \lemref{gcriterion0}.
\end{rem}

\begin{prop} \label{vandg}
 \tcb{Let $V$ be an affine variety and  $X$ be a definable subset that is both v-closed and g-closed.   Then $X$ may be defined by a positive Boolean combination
of subvarieties and sets defined by weak valuation inequalities  $ \val (f(x)) \leq \val(g(x))$, where  $f,g$ are regular functions on $V$.  
 A similar statement may be made for $V$ projective, using homogeneous polynomials.}
 \end{prop}

\prf \tcb{We prove this by induction on $\dim(V)$; assume
the statement holds for varieties of lower dimension.  We may assume $V$ is irreducible.   As $X$ is g-closed, it is defined by 
  weak valuation inequalities along with algebraic equalities and inequalities; thus away from some proper subvariety $V'$ of $V$,
  $X$ coincides with a set  $X'$  cut out by the inequalities $\val (f_i) \leq \val (g_i)$,  $i=1,\ldots,n$.  Thus  $X' \m V' = X \m V'$; by induction,
  $X \meet V'$ has the right form;   if we also show that $X' \meet V' \subset X$, then $X = X' \union (X \meet V')$ will have the
  promised form.  Thus it suffices to show that   $X' \subset X$.   As $X'$ is v-closed, this follows from \lemref{blv}.}\eprf  

\begin{lem}\label{blv} \tcb{Let $V$ be an affine variety, let $f_i$ and $g_i$, $1 \leq i \leq n$,  be nonzero regular functions on $V$ and let $V'$ be a proper subvariety of $V$. Let $Y$ be the subset of $V$ defined by the inequalities $\val (f_i) \leq \val (g_i)$,  $i=1,\ldots,n$. Then any  point $b$ of $Y$ lies arbitrarily close to a point
  of $Y \m V'$ in the valuation topology.}
  \end{lem}
  
\prf \tcb{Let $p: \tV \to V$ be the result of blowing up the ideal $(f_1,g_1)$ on $V$; let $b'$ be a point of $\tV$ lying above $b$, and let $\tV'$
be an affine open of $\tV$ containing $b'$.   If we show the existence
of points of $p \inv(Y) \meet \tV'$   arbitrarily close to $b'$, avoiding the exceptional divisor as well as $p \inv(V')$, then by continuity of $p$ the claim will be proved.
Now on $\tV'$, there is a regular function $u_1$ such that $f_1=g_1u_1$ or $f_1u_1=g_1$; so the inequality $\val (f_1) \leq \val (g_1)$ can be replaced by $\val (u_1) \leq 0$,
or $\val (u_1) \geq 0$.  Iterating this construction, we may assume $Y$ is defined by a conjunction of inequalities  $\val (u_i) \leq 0 $ or $\val (u_i) \geq 0$ for some    regular functions $u_i$, $1 \leq i \leq n$.
 Now if we take any point of $V$ very close to $b$ in the valuation topology (but avoiding $V'$) these inequalities are preserved.}
 \eprf

\begin{defn}Let $V$ be 
an algebraic variety over a valued field $F$ or a definable subset of such a variety.  A definable function $h: V \to \G_\infty$
is called {\em v-continuous} (resp. {\em g-continuous}) if the pullback of any v-open (resp. g-open) set is v-open (resp. g-open).  \index{v-continuous}\index{g-continuous}
A function 
$h: V \to \std{W}$ with $W$ an affine $F$-variety is called v-continuous (resp. g-continuous) if, for any regular function
$f: W \to \Aa^1$,
$\val \circ f \circ h$ is v-continuous (resp. g-continuous).  \end{defn}

Note  that the topology generated by v-open subsets on $\G_\infty$ is discrete on $\G$, while the 
neighborhoods of $\infty$ in this topology are the same as in the order topology.  The   topology generated by g-open subsets
is the order topology on $\G$, with $\infty$
isolated.  We also have the topology on $\G_\infty$ coming from its canonical
identification with $\std{\G_\infty}$, or the v+g-topology; this is the intersection
of the two previous topologies, that is, the order topology on $\G_\infty$.  

\tcb{Let $V$ be 
an algebraic variety over a valued field $F$ and let $X$ be  a definable subset of $V \times \G_{\infty}^m$.} We say that \tcb{$X$} is {\em {v$+$g-}open} if it is both v-open and g-open.   \index{v+g-open}
\tcb{The complement of a v+g-open subset is said to be {\em v$+$g-closed}.}
\index{v+g-closed}
If $W$ has a definable topology,
a definable function
$X \to W$ is called {\em v$+$g-continuous} \index{v+g-continuous} if the pullback of a definable open subset of $W$ is both v- and g-open, and
similarly for functions to $V$.

\begin{rem}Note that v, g  and v+g-open sets are {\em definable} sets.   Over any given model  it is possible to extend v to a topology in the usual sense, the valuation topology, whose restriction to definable sets is the family of v-open sets.  But   this is not true of g and of v+g; in fact they are not closed
under definable unions, as the example $\Oo = \cup_{a \in \Oo} \tcb{a +} \Mm$ shows.
\end{rem}

Any g-closed subset $W$  of an algebraic variety is defined by a disjunction $\bigvee_{i=1}^m (\neg H_i  \wedge \phi_i)$,
with $\phi_i$ a  finite conjunction of   weak valuation inequalities $v(f) \leq v(g)$ and equalities, and $H_i$ defining a 
 Zariski closed subset.   If $W$ is also v-closed, $W$ is equal to the union of the v-closures of the sets
 defined by $\neg H_i  \wedge \phi_i$, $1 \leq i \leq m$.

 \begin{lem} \label{v+g}  \tcb{Let $X$ be a definable subset of a variety $V$ over a valued field. Let $W$ be a definable subset of $X$ which is v+g-closed in $X$.
 Then $\std{W}$ is  closed in $\std{X}$.  
  More generally, if $W$ is g-closed in $X$, then $cl(\std{W}) \cap \std{X}\nsubseteq \std{cl_v(W)}\cap \std{X}$,  with $cl$ and $cl_v$ denoting respectively the closure and the v-closure.}  \end{lem}

\prf  Let $M$ be a model, $p \in \tcb{\std{X}}(M)$, with $p \in cl(\std{W}(M))$.  We will show that $p \in \std{cl_v(W)}$.
Let $(p_i)$ be a net in  $\std{W}(M)$ approaching $p$.  Let $a_i \models p_i |M$.  Let $\tp(a/M)$
be a limit type in the logic topology (so $a$ can be represented by an ultraproduct of the $a_i$).  
For each $i$ we have $\G(M(a_i))=\G(M)$, but $\G(M(a))$ may be bigger. 

 Consider the subset $C$ of
$\G(M(a))$ consisting of those elements $\g$ such that
$-\a < \g < \a$ for all $\a > 0$ in $\G(M)$.
  Thus  $C$ is a convex subgroup of $\G(M(a))$; let $N$
be the valued field extension of $M$ with the same underlying     $M$-algebra structure as $M(a)$,
obtained by factoring out $C$.  Let $\bar{a}$ denote $a$ as an element of $N$.   We have $a_i \in W$,
so $a \in W$; since $W$ is g-closed \tcb{in $X$} it is clear that $\bar{a} \in W$.   (This is the easy direction of 
\lemref{gcriterion0}.)

Let $b \models p|M$.
  For any regular function $f$ in $M [U]$, with $U$ Zariski open in $V$,  
  we have: $(*)$ \  $\val (f(a_i)) \to \val (f(b))$ in $\G_\infty(M)$  (since $p_i \to p$).

Let $R = \{x \in N: (\exists m \in M)( \val(x) \geq \val(m)) \}$.  Then $R$ is a valuation ring of $N$ over $M$.
By  $(*)$, for large enough $i$,   $\val (f(a_i))$ is bounded \tcb{below by} some element of $\G(M)$ (namely any element below
$p_*(f)$).  So $\val (f(a))$ and $\val (f(\bar{a}))$ must lie above the same element.  Thus $\bar{a} \in R$.  
Also by $(*)$,  if $\val (f(\bar{a}))= \infty $, or just if  $\val (f(\bar{a})) > \val (M)$, then $f(b)=0$.  Thus we have a well-defined
map from the residue field of $R$ to     $M(b)$, with $\res  \, \bar{a} \mapsto b$.   
Since $\bar{a} \in W$, it follows that $b \in cl_v(W)$ \tcb{(cf. 
the last part of the proof of \lemref{lem:v-closed})}, hence $p \in \std{cl_v(W)}$.
  \eprf 
 
\begin{prop}\label{hausdorff}
\tcr{Let $V$ be an algebraic variety over a valued field $F$. Then $\std{V}$ is Hausdorff.}
\end{prop}

\prf \tcr{Let us first consider the case when $V$ is affine. In this case we may assume $V = \Aa^n$. Let $p \in \std{\Aa^n}$ and 
$p' \in \std{\Aa^n}$. Let $M$ be a model of $\ACVF$ such that $p \in \std{\Aa^n} (M)$ and 
$p' \in \std{\Aa^n} (M)$. There exists a polynomial
$F \in M [x_1, \ldots, x_n]$ such that
$\val (F (p)) \not= \val (F (p'))$. Thus, for some $\alpha \in \G (M)$, the disjoint open sets
defined by the conditions
$\val (F) < \alpha$ and $\val (F) > \alpha$ will separate $p$ and $p'$.}

\tcr{In general,  since $V$ is separated as an algebraic variety,  $V$
is
Hausdorff for the valuation topology. Indeed,  the diagonal $\Delta$ of $V$ in $V \times V$ being Zariski closed, it is closed for the valuation topology, and the valuation topology on
$V \times V $ is the product topology.}

\tcr{Let $p \in \std{V}$. Fix $U$ an affine open subset of $V$ such that $p \in \std{U}$.
Let $Z$ be the intersection of the Zariski closed subsets $W$ of $U$ such that $p \in \std{W}$. Choose a closed embedding of $U$ in some affine space
$\Aa^n$. There exists a definable v-closed subset $C$ of $Z$ such that $p \in \std{C}$. 
We may further assume $C$ to be bounded, meaning that there exists some $\gamma$ in $\G$
 such that for any $1 \leq i \leq n$, $\val (x_i) \geq \gamma$ on $C$.
We denote by  $d : U \times U \to \G_\infty$  the restriction of the function
$\min_{1 \leq i \leq n} (\val (x_i -y_i))$ on $\Aa^n \times \Aa^n$.
Now let $p' \in \std{V}$ with $p' \not= p$. Fix $U'$ an affine open subset of $V$ such that $p' \in \std{U'}$
and a closed embedding of $U'$ in some affine space
$\Aa^{n'}$. Define $Z'$ and $d'$ similarly as $Z$ and $d$ and fix a 
bounded definable v-closed subset $C'$ of $Z'$ such that $p' \in \std{C'}$.}

\tcr{If $Z \meet U' \neq \varnothing$, $p$ and $p'$ both lie in $\std{U'}$ and we are done by the affine case.
Thus, we may assume $Z \meet U'= \varnothing$ and $Z' \meet U = \varnothing$.  In particular, $Z \meet Z' = \varnothing$ and so $B \meet B' = \varnothing$.  
For $c \in U$ and $\alpha \in \G$, let $B_{\alpha} (c) = \{x \in U: d (x, c) > \alpha\}$. One defines similarly
$B'_{\alpha} (c')$ for $c' \in U'$. Since $V$ is Hausdorff for the valuation topology, for any $(c, c') \in C \times C'$ there exists some $\alpha \in \G$ such that
$B_{\alpha} (c) \meet B'_{\alpha} (c') = \varnothing$.}

\tcr{Let us now prove we can take a single $\alpha$ to work for all $(c, c') \in C \times C'$. Assume this is not the case and fix a valued field extension $F \leq M$ with
$M$
a model of $\ACVF$.
Then, for any $\alpha \in \G (M)$, there would exist $c_{\alpha} \in C (M)$,
$c'_{\alpha} \in C' (M)$ and $w_{\alpha} \in (U \meet U') (M)$ such that
$d (c_{\alpha}, w_{\alpha}) > \alpha$ and
$d' (c'_{\alpha}, w_{\alpha}) > \alpha$. Take a non-principal ultrafilter on the set $\G (M)$ and let $\widetilde M$ be the corresponding ultrapower.
Set $R = \{x \in \widetilde M: (\exists \alpha \in \G (M)) (\val (x) > \alpha)\}$
and
$I = \{x \in \widetilde M: (\forall \alpha \in \G (M))  (\val (x) > \alpha)\}$. The quotient
$M' = R / I$ is an elementary extension of $M$.
Denote by $\tilde c$, $\tilde{c}'$ and $\tilde w$ the class of 
$(c_{\alpha})$,
$(c'_{\alpha})$ and $(w_{\alpha})$ in $V(\widetilde M)$.
The boundedness assumptions imply that
$\tilde c \in C \meet V (R)$, $\tilde{c}' \in C' \meet V (R)$, and
$\tilde w \in U \meet U' \meet V (R)$, so we can consider their images 
$c$, $c''$ and $w$  in $V (M')$. Since $C$ and $C'$ are v-closed, it follows from
\lemref{lem:v-closed} that $c \in C (M')$ and $c' \in C' (M')$.
But then  $c = w= c'$ so $C$ and $C'$ are not disjoint.}

\tcr{Fix such an $\alpha$ and set
$O = \{x \in U: (\exists c \in C)  (d (x, c) > \alpha)\}$. Define similarly $O' \subset U'$.
By construction, $O$ and $O'$ are disjoint, thus $\std{O}$ and $\std{O'}$ are disjoint.
We have $p \in \std{O}$ and $p' \in \std{O'}$.
Let us check that $\std{O}$ is open. Let $q \in \std{O}$ and $M$ a valued field extension of $F$ which is 
a model of $\ACVF$.
Let  $a$ such that $\tp (a/M) = q {\vert M}$.
Since $a$ belongs to $O$,  it belongs to $B_{\alpha} (c)$ for some $c \in C$.
Thus $q$ belongs to $\std{B_{\alpha} (c)}$ which is open and contained in $\std{O}$.  It follows that $\std{O}$ is open. For similar reasons $\std{O'}$ is open, and thus
$\std{V}$ is Hausdorff.}
\eprf

 \section{Canonical extensions}\label{canon-ext}
  Let $V$ be a definable set 
 over some $A$
and let $f:  V \to \std{W}$ be an  
$A$-pro-definable morphism,
where $W$ is an 
$A$-definable 
subset of \tcb{$Z\times \G_\infty^m$, with $Z$ an algebraic variety defined over $A$}.   We can define a {\em canonical extension} \index{canonical extension} to $F: \std{V} \to \std{W}$,
 as follows.

  If $p \in \std{V} (M)$, say $p|M = \tp(c/M)$, let
 $d \models f(c) | M(c)$.
 By transitivity of stable domination, \propref{transitive}, $\tp(cd/M)$ is stably dominated, and hence so is 
$\tp(d/M)$.  Let 
$F(c) \in \std{W} (M)$ be such that $F(c)|M = \tp(d/M)$; this does not depend
on $d$.  Moreover $F(c)$ depends only on $\tp(c/M)$, so we can let $F(p)=F(c)$.  
\tcb{By \lemref{cepd}}, $F : \std{V} \to \std{W}$ is an $A$-pro-definable 
 morphism.
 Sometimes the canonical extension $F$ of $f$ will be denoted by
 $\std{f}$ or even by $f$.

\begin{lem}\label{cepd}\tcb{Let $f:  V \to \std{W}$ be an  
$A$-\tcb{pro-}definable map as above. Then the canonical extension $F : \std{V} \to \std{W}$ is an $A$-pro-definable morphism.
}
\end{lem}

\prf \tcb{Let $g : W \times Z \to \G_{\infty}$ be a definable map
and let $Y_{Z, g}$ be
the
corresponding definable set  of definable functions $Z \to \Gamma_{\infty}$ considered in
the proof of \thmref{prodef}.
The composition of $f$ with the projection $\std{W} \to Y_{Z, g}$ yields a definable map
$\bar f: V \to Y_{Z, g}$. Let $\bar g : V \times Z \to \G_{\infty}$ be the definable map
sending $(v, z) \in V \times Z$ to $\bar f (v) (z)$.
For any $p \in \std{V}$ we have $p_* (\bar g) = F (p)_* (g)$, hence
there is a definable inclusion $Y_{Z, \bar g} \hookrightarrow Y_{Z, g}$.
Since the composition of $F$ with the projection $\std{W} \to Y_{Z, g}$
factors through that inclusion,  it follows that $F$ is an $A$-pro-definable morphism.}
\eprf

 \begin{lem}  \label{basic0} Let $f:  V \to \std{W}$ be a   \tcb{pro-}definable morphism, where $V$ is \tcb{a definable subset of}
 an algebraic variety and   $W$ is a definable subset of $\Pp^n \times \G_\infty^m$.
Let  $X$ be a definable subset of $V$.   Assume $f$ is {g}-continuous \tcb{on $V$} and v-continuous at each point of $X$; i.e.
  $f \inv (G)$ is g-open whenever $G$ is open, and   $f \inv(G)$ is v-open at $x$
whenever $G$ is open, for any $x \in X \meet f \inv (G)$.   
 Then the canonical extension $F$ is continuous at each point of  $\std{X}$.   \end{lem}
 
 \prf 
 The topology on $\std{\Pp^n}$ may be described as follows, \tcb{cf.  \ref{reppn}}. It is generated by
 the preimages of open sets of $\G_\infty^N$ under continuous definable functions $\Pp^n \to \G_\infty^N$ of the form
 \[[x_0: \ldots: x_n] \longmapsto  
 [\val(x_0^d):\ldots:\val(x_n^d):
 \val (h_1): \ldots: \val (h_{N-n})]\] for some homogeneous polynomials $h_i(x_0, \ldots,  x_n)$ of degree $d$; where 
 in $\G^{N}_{\infty} \m \{\infty\}^N$ we define $[u_0: \ldots: u_m]$
 to be $(u_0 - \min u_i,\ldots,u_m - \min u_i)$.
 Composing with such functions we reduce to the case of $\G_\infty^m$, and hence to the case of $f: V \to \G_\infty$.  
 
 Let $U= f\inv(G)$ be the pullback of  a definable open subset $G$ of  $\G_\infty$.   Then 
  $F \inv (G) = \std{U}$.   Now $U$ is g-open, and v-open at any $x \in X \meet U$.     By \lemref{v+g} applied
  to the complement of $U$ in $V$, it follows that $\std{U}$ is open at any $x \in \std{X}$.        \eprf


 \begin{lem} \label{hbasic0}   Let $K$ be a valued field and $V$ be an algebraic variety over $K$.
Let $f: I \times V \to \std{V}$ be a g-continuous $K$-\tcb{pro-}definable morphism,  where $I = [a, b]$ is a closed interval.  
Let $i_I$ denote one of $a$ or $b$ and
$e_I$ denote the remaining point.
Let $X$ be a $K$-definable subset of $V$. Assume $f$ restricts
to a definable morphism
$g : I \times X \to \std{X}$ and that $f$ is $v$-continuous at every point of $I \times X$.
Then $g$ extends uniquely to a continuous $K$-pro-definable morphism $G: I \times \std{X} \to \std{X}$. 
 If moreover, for every $v \in X$,  $g(i_I,v)=v$
and $g(e_I,v) \in Z$, with $Z$ a $\G$-internal \tcb{iso-definable} 
subset, then \tcb{for every $x \in \std{X}$,}
$G(i_I,x)=x$, and $G(e_I,x) \in Z$.   
\end{lem}

\prf   
Since $\std{I \times V} = I \times \std{V}$ by \lemref{simple0}, the first statement
follows from \lemref{basic0}, by considering the pullback of
$I$ in $\Aa^1$.  The equation  $G(i_I,x)=x$ extends by continuity from the dense set of simple points to $\std{X}$.
We have by construction $G(e_I,x) \in Z$, using the fact that any stably dominated type on $Z$ is constant.   \eprf 

 \begin{lem}  \label{basic}Let $K$ be a valued field and $V$ be a definable subset of 
an  algebraic variety over $K$.
 Let $f:  V \to \std{W}$ be a   $K$-\tcb{pro-}definable morphism, with   $W$  a $K$-definable subset of $\Pp^n \times \G_\infty^m$. Assume $f $ is v+g-continuous.  
 \tcb{Then $F:   \std{V} \to \std{W}$ is continuous and it is the unique extension of 
$f$  to a continuous $K$-pro-definable morphism $  \std{V} \to \std{W}$}. 
\end{lem}
 
 \prf   \tcb{Let us prove the continuity of $F$. As in the proof of \lemref{basic0}, it is enough to consider the case $W = \G_\infty$ which follows directly from
 \lemref{v+g}.}
   There is clearly at most one continuous extension, because of  the   density in $\std{V}$ of the set of simple points $V(\Uu)$, cf. \lemref{simple}.\eprf 
 
\begin{lem} \label{hbasic}   \tcb{Let $K$ be a valued field and $V$ be a definable subset of 
an  algebraic variety over $K$.
Let $f: I \times V \to \std{V}$ be a v+g-continuous $K$-\tcb{pro-}definable morphism,  where $I = [a, b]$ is a closed interval.  
Let $i_I$ denote one of $a$ or $b$ and
$e_I$ denote the remaining point.
Then $f$ extends uniquely to a continuous $K$-pro-definable morphism $F: I \times \std{V} \to \std{V}$. 
 If moreover, for every $v \in V$,  $g(i_I,v)=v$
and $g(e_I,v) \in Z$, with $Z$ an  \tcb{iso-definable}  $\G$-internal subset, then, for every $x \in \std{V}$,
$G(i_I,x)=x$, and $G(e_I,x) \in Z$.}
\end{lem}
 
\prf \tcb{Follows from \lemref{basic} similarly as \lemref{hbasic0} follows from \lemref{basic0}.}
 \eprf
    

  \section{Paths and homotopies}\label{ph}
By an {\em interval} \index{interval} we mean a 
sub-interval of $\G_{\infty}$.
Note that intervals of different length are in
general not definably homeomorphic, and that the gluing of two intervals (e.g. $[0,1]$ coming to the right of $[0,\infty]$)  may not 
result in an interval. 
We get around the latter issue by formally introducing a more
general notion, that of a generalized interval.    


Given an interval $I$ in $\Gamma_{\infty}$, we may consider it either with the \tcb{induced} order or with the opposite order. The choice
of one of these orders will be \tcb{called}
an orientation of $I$.
\tcb{Let $I_1$, $\dots$, $I_n$ be oriented sub-intervals of $\Gamma_{\infty}$. Assume $I_1$ is right-closed (i.e. contains its largest endpoint), $I_n$ is left-closed (i.e. contains its smallest endpoint), 
and that
each $I_j$ is closed for $1<j<n$. Then one may glue end-to-end the intervals $I_i$ in a way respecting the orientations by identifying the largest endpoint of $I_i$ with the smallest endpoint of $I_{i + 1}$ for $1 \leq i <n$,
and obtain a definable space.
Any definable space $I$ that may be obtained this way will be called a {\em  generalized interval}.}\index{generalized interval}

If \tcb{the generalized interval} $I$ is closed, we denote by $i_I$ \nomenclature{$i_I$}{smallest element of $I$} the smallest element of $I$ and by 
$e_I$ its largest element. \nomenclature{$e_I$}{largest element of $I$}
Note that if $I$ is obtained by gluing intervals
$I_1$, $\dots$, $I_n$,
a function
$I \times V \to W$ is definable, resp. continuous, resp. v+g-continuous, if and only if it is obtained by gluing
definable, resp. continuous, resp. v+g-continuous, functions
$\varphi_i : I_i \times V \to W$.



Let $V$ be a definable set.
By a {\em path} on $\std{V}$ \index{path}
we mean a continuous definable map
$I \to \std{V}$ with $I$ some generalized interval.

\begin{defn} \label{homotopy-defn} Let $X$ be a pro-definable subset of $\std{V} \times \G_\infty^n$.  A {\em homotopy} is \index{homotopy}
 a continuous pro-definable map
$h : I \times X \to X$ with $I$ a closed generalized interval. 
The maps $h_{i_I}$ and $h_{e_I}$
are then said to be homotopic  \tcb{(one denotes by $h_t$ the map sending $x \in X$ to $h (t, x)$)}.
\tcb{The homotopy $h$ is 
called a  {\em deformation retraction} \index{deformation retraction} to $A \nsubseteq \std{V}$ if
$h_{i_I} = \mathrm{id}_X$, $h (t, a) = a$ for all $t$ in $I$ and $a$ in $A$ and furthermore
$h_{e_I} (x) \in A$ for each $x$.}
(In the literature, this is sometimes referred to as a {\em strong} deformation
retraction.) We say \tcb{$A = h_{e_I} (X)$}
is the {\em image} of $h$. \index{image of a deformation retraction}
If $\varrho = h_{e_I}$,
we say  that $(\varrho, \varrho(X))$ 
is a deformation retract.
Sometimes, we shall also call $\varrho$ or $\varrho(X)$  a deformation retract, the other member of the pair being understood implicitly. \index{deformation retract}

If $W$ is a definable subset of $V \times \G_\infty^n$, we will also refer to a v+g-continuous pro-definable map
$h_0: I \times W \to \std{W}$ as a homotopy; by \lemref{hbasic}, $h_0$   extends uniquely to a homotopy $h: I \times \std{W} \to \std{W}$.  
\tcb{One defines similarly a deformation retraction $h_0: I \times W \to \std{W}$ and its image.} 
\end{defn}

\tcb{By \lemref{hbasic} if $h_0$ is a deformation retraction 
with image an iso-definable $\G$-internal subset
then its canonical extension 
is a deformation retraction with the same image.}

\begin{example} Generalized intervals may in fact be needed to connect points of $\std{V}$.  For instance
let $V$ be a cycle  of $\tcb{2} n$  copies of $\Pp^1$, with consecutive pairs meeting in a point.  \tcb{By gluing 
$2n$ copies of the  homotopy $\psi_{\{0, \infty\}}$ as 
defined in  \ref{ss7.6},  one gets
a deformation retraction $[0,\infty] \times \std{V} \to \std{V}$
with image} a cycle made of $\tcb{2}n$ copies of $[0,\infty] \nsubset \G_\infty$.  However  
it is impossible to connect two points at extreme ends of this topological circle without gluing together $\tcb{n}$ intervals. \end{example}

\begin{defn}\label{star}Let $X$ be a pro-definable subset of $\std{V} \times \G_\infty^n$. A homotopy
$h :  I \times X \to X$
is said to satisfy condition $(*)$ if \index{condition $(*)$}
$h (e_I, h (t,x)) = h (e_I, x)$ for every $t$ and $x$.
One defines similarly condition $(*)$
for a homotopy
$h_0: I \times W \to \std{W}$ when $W$ is a definable subset of $V \times \G_\infty^n$. Note that $h_0$ satisfies $(*)$ if and only if its canonical extension does.
\end{defn}

Let $h_1: I_1 \times X \to X$
and $h_2 : I_2 \times X \to X$
two homotopies. Denote by 
$I_1 + I_2$ the (generalized) interval obtained by gluing
$I_1$ and $I_2$  at $e_{I_1}$ and
$i_{I_2}$. Assume
$h_2 (i_{I_2}, h_1 (e_{I_1}, x)) = h_1 (e_{I_1}, x)$ 
for every $x$ in $X$.
Then one denotes
by 
$h_2 \circ h_1$ the homotopy
$(I_1  + I_2) \times X\to X$
given by 
$h_1 (t, x)$ for $t \in I_1$ and by
$h_2 (t, h_1 (e_{I_1}, x))$ for $t$ in $I_2$, \tcb{and one calls $h_2 \circ h_1$ the composition (or concatenation) of $h_1$ and $h_2$}.
\index{composition of homotopies}


\begin{defn}\label{defn:defclosed}
\tcr{Let $X$ be a pro-definable subset of $\std{V} \times \G_\infty^n$
and let $X'$
be a pro-definable subset of $\std{V'} \times \G_\infty^{n'}$. A pro-definable map $f : X \to X'$ is said to be
{\em definably closed} if 
for any closed pro-definable subset $Z$ of $X$,
$f (Z)$ is closed in $X'$.}
\end{defn}
 \index{definably closed map}

 \begin{rem}\label{injcl}
\tcr{Note that an injective pro-definable map $f$ is definably closed if and only if it is closed, since in this case taking the image under $f$ commutes with arbitrary intersections.}
\end{rem}

 \begin{lem} \label{homotopy-descent} \tcr{Let $V$ be an algebraic variety over a valued field, 
and let $X$ and $X_1$ be pro-definable subsets of $\std{V}\times \G_{\infty}^N$.
Let $f: X_1 \to X$ be a continuous, definably closed and surjective pro-definable map.  Let
$I$ be a closed generalized interval and  $h_1: I \times X_1 \to X_1$   be a    homotopy.
Assume $h_1$  respects the fibers
of $f$, in the sense that
$ f( h_1(t,x))$ depends only on $t$ and $f (x)$.
Then $h_1$  descends to a   homotopy of $X$.}   \end{lem}

\prf     \tcr{Define $h: I \times X \to X$
 by
  $h (t, f(x)) =  f( h_1(t,x))$ for $x \in X_1$; then $h$ is well-defined and pro-definable.    We denote the map 
  $\mathrm{Id} \times f : I \times X_1 \to X$ by $f_2$.   Clearly, $f_2$ is a continuous, definably closed and surjective pro-definable map   (the topology on $I \times X_1$, $I \times X$
  being the product topology). 
   To show that   $h$ is continuous, it suffices therefore to show that  $h \circ f_2$ is continuous.
Since $h \circ f_2=f \circ h_1$  this is clear.}  \eprf 

 \begin{rem}\label{rem:homotopy-descent}In particular,  let 
$f: V_1 \to V$ be a  proper surjective morphism of algebraic varieties over a valued field. 
   Let $h_1$ be a    homotopy $h_1: I \times \std{V_1} \to \std{V_1}$, and assume $h_1$   respects the fibers of
$\std{f}$.    Then $\std{f}$ is surjective by \lemref{extend}, 
and \tcr{definably} closed by \lemref{closedmap3}; so  $h_1$  descends to a   homotopy of $X$. 
\end{rem}

\begin{defn}\tcr{Let $X$ be a pro-definable subset of $\std{V} \times \G_\infty^n$
and let $X'$
be a pro-definable subset of $\std{V'} \times \G_\infty^{n'}$.}
A continuous pro-definable map $F : X \to X'$ is said to be
a {\em homotopy equivalence} if there exists a \index{homotopy equivalence}
 continuous pro-definable map $G : X' \to X$
 such that $G \circ F$ is homotopic to $\Id_{X}$
 and
 $F \circ G$ is homotopic to $\Id_{X'}$.
\end{defn}

\section{Good metrics}
\label{Smetrics}

By a {\em definable metric} on an algebraic variety $V$ over a valued field $F$, we mean an $F$-definable function $d: V^2 \to \G_\infty$ \index{definable metric}
which is v+g-continuous and such that
\begin{enumerate}
\item  $d(x,y)=d(y,x)$; $d(x,x) = \infty$;
\item  $d(x,z) \geq \min(d(x,y),d(y,z))$;
\item   if $d(x,y) = \infty$ then $x=y$.  
 \end{enumerate}  
 
 Note that given a definable metric on $V$,
 for any $v \in V$,
$B(v;d,\g) := \{y: d(v,y) \geq \g \}$ is a family of g-closed, v-clopen sets whose intersection is $\{v\}$. Iit follows by a definable compactness argument that $d$ induces the v-topology on $V$; this is anyhow clear
for the metrics we will use.

We call $d$ a {\em good metric} \index{good metric} if  
 there exists a v+g-continuous \tcb{$F$-}definable function
 $\rho: V \to \G$ (so $\rho(v)<\infty)$, such that
 for any $v \in V$ and any $\alpha \tcb{\geq} \rho(v)$, 
$B(v;d,\alpha)$ is affine and has a unique generic type,
i.e. a definable type $p$ such that
for any  Zariski closed \tcb{subset $V' $ of  $V$ not containing $B(v;d,\alpha)$}
and any regular $f$ on $V \m V'$, 
$p$ concentrates on $B(v;d,\alpha) \m V'$,  
and $p_{*}(f)$ attains the minimum valuation of $f$ on $ B(v;d,\alpha) \m V'$.   Such a type is  
orthogonal to $\Gamma$, hence stably dominated. 
 
 \tcb{The continuity of $\rho$ can be replaced by local boundedness in this definition, using \lemref{majorize}.}
 
 \begin{lem} \label{metric} \leavevmode \begin{enumerate}
 \item  $\Pp^n$ admits a good  metric, with $\rho=0$.
 \item Let $F$ be a valued field, $V$ a   quasi-projective variety over $F$.  Then there exists a
definable metric on $V$.  
\end{enumerate}
\end{lem} 

\prf  Consider first the case of $\Pp^1 = \Aa^1 \union \{\infty\}$.  
If $x,y \in \Oo$, set
$d(x,y) = d(x \inv, y \inv) = \val(x-y)$ and let $d(x,y) = 0$ if $v(x)$ and $v(y)$ have
different signs.  This is easily checked to be consistent, and to satisfy the conditions (1-3).
It is also clearly {v-}continuous.  
\tcb{Let us now prove g-continuity. By \propref{gcriterion}, it is enough to check that}
if $F \leq K$ is a valued field extension, 
$\pi: \G(K) \to \bar{\Gamma}$ a homomorphism of ordered $\Qq$-spaces extending
$\G(F)$,  and $\bK = (K,\pi \circ v)$, then $ \pi (d_K(x,y)) = d_{\bK}(x,y)$.
If $x,y \in \Oo_K$ then $x,y \in \Oo_{\bK}$ and $\pi (d_K(x,y)) = \pi (v_K(x-y)) = d_{\bK}(x,y)$.
Similarly for $x \inv, y \inv$.  If $v(x)<0<v(y)$, then $v(x-y)<0$ so $\pi (v(x-y)) \leq 0$,
hence $d_{\bK}(x,y) =0 = d_K(x,y)$.  This proves  g-continuity.   It is clear that the metric is good, with $\rho=0$.

Now  consider $\Pp^n$ with homogeneous coordinates
$[X_0 : \ldots : X_n]$.
For $0 \leq i \leq n$ denote by $U_i$ 
the subset $\{x \in \Pp^n : X_i \not=0 \wedge \inf \val (X_j / X_i) \geq 0 \}$.
If $x $ and $y$ belong both  to
$U_i$, one sets $d (x, y) = \inf \val (X_j / X_i - Y_j /Y_i)$.
If $x \in U_i$ and $y \notin U_i$, one sets
$d (x, y) = 0$. One checks that this definition is unambiguous
and reduces to the former one when $n = 1$. The proof it is v+g-continuous is similar to the case $n = 1$ and the fact it is good with $\rho = 0$ is clear.
This metric restricts to a   metric on any subvariety of $\Pp^n$.
\eprf

\section{Zariski topology}\label{zt}
We shall occasionally use the Zariski topology on $\std{V}$.
If $V$ is an algebraic variety over a valued field, a subset of $\std{V}$ of the form
$\std{F}$ with $F$ Zariski closed, resp. open, in $V$ is said to be {\em Zariski closed}, resp. {\em open}.
Similarly, a subset $E$ of $\std{V}$ is said to be {\em Zariski dense} in $\std{V}$ if
$\std{V}$ is the only Zariski closed set containing $E$.  For $X \subset\std{V}$, the Zariski topology on $X$ is the one induced 
from the Zariski topology on $\std{V}$. \index{Zariski topology} \index{Zariski closed} \index{Zariski open} \index{Zariski dense}

\section{Schematic distance}\label{schedis}
Let $f (x_0, \ldots, x_m)$ be a homogenous polynomial
with coefficients in the valuation ring $\Oo_F$ of a valued field $F$.
One defines a function
$\val (f) : \Pp^m \to [0, \infty]$
by
$\val (f) ([x_0 : \ldots : x_m]) = \val (f ({x_0}/{x_i}, \ldots, {x_m}/{x_i}))$
for any $i$ such that $\val (x_i) = \min_j (\val (x_j))$.

Now let  $V$ be a projective variety over a valued field $F$
and let $Z$ be a closed $F$-subvariety of $V$.
Fix an embedding $\iota : V \hookrightarrow  \Pp^m$ and a family $f$ of
homogenous polynomials
$f_i$, $1 \leq i \leq r$, in $\Oo_F [x_0, \ldots, x_m]$
such that $Z = V \meet (f_1 = \cdots = f_r = 0)$.
For $x$ in $V$ set $\varphi_{\iota, f} (x) = \min_i (\val (f_i (x)))$.
The function $\varphi_{\iota, f} : V \to [0, \infty]$
is clearly $F$-definable and v+g-continuous and
$\varphi_{\iota, f} \inv (\infty) = Z$.
Any function  $V \to [0, \infty]$ of the form
$\varphi_{\iota, f} $ for some $\iota, f$ will be called
a {\em schematic distance function to} $Z$.
\index{schematic distance}


\chapter{Definable compactness}\label{definable-compactness}

{\small \noindent \textbf{Summary.}
This chapter is devoted to the study of definable compactness for subsets of
$\std{V}$. One of the main results is
\thmref{compactiffboundedclosed} which establishes the equivalence between
being definably compact and being closed and bounded.
\par\bigskip}

\section{Definition of definable compactness}
We will use definable types as a replacement for the curve selection lemma, whose purpose is often to use the definable type associated with a curve at a point.
Note that the curve selection lemma itself is not true for $\G_{\infty}$, e.g. in $\{(x,y) \in \G_\infty^2:
y>0, x<\infty\}$  there is no curve approaching $(\infty,0)$.

\tcb{As we already observed, one can consider definable types in infinitely many variables, thus the notion of a definable type on a pro-definable set makes sense.}  \index{definable type}
Let $X$ be a definable or pro-definable topological space
in the sense of \ref{deftop}.  Let $p$ be a definable type on $X$.

\begin{defn}  \label{deflimit} A point $a \in X$ is a {\em limit} of $p$ if for any definable neighborhood
$U$ of $a$ (defined with parameters), $p$ concentrates on $U$.\end{defn} \index{limit of  a definable type}

When $X$ is Hausdorff, it is clear that a limit point is unique if it exists.
\tcr{Recall that by \propref{hausdorff}, $\std{V}$ is Hausdorff for any variety $V$ over a valued field (hence also $\std{V} \times \G_{\infty}^n$).}

\begin{defn}  \label{defcomp}Let $X$ be a definable or pro-definable topological space. One says    $X$ is {\em definably compact} if   any definable type $p$ on $X$ has a limit point in $X$.
 \end{defn}
 \index{definably compact}
   
For subspaces of $\G^n$ with $\G$ o-minimal, our definition of definable compactness
 lies between the definition of \cite{peterzil-steinhorn} in terms of curves, and
the property of being closed and bounded; so all three are equivalent.   This will be treated in more detail later.

\section{Characterization of definable compactness} \label{4.2}
A  subset of   $\mathrm{VF}^n$ is said to be 
{\em bounded} if for some $\gamma$ in $\G$ it is contained 
\index{bounded subset}
in $\{(x_1, \dots, x_n) : v (x_i) \geq \gamma, 1 \leq i \leq n\}$.
This notion extends to varieties $V$ over a valued field, cf., e.g., \cite{serre} p.~81:  $X \nsubseteq V$ is
defined to be bounded if there exists an affine cover $V=\union_{i=1}^m U_i$, and bounded subsets $X_i \nsubseteq U_i$, with $X \nsubseteq \union_{i=1}^m X_i$.  Note that the projective space $\Pp^n$ is bounded within itself, and
so any subset of a projective variety $V$ is bounded in $V$.
We shall say a subset of $\G_\infty^m$ is {\em bounded} if it is contained in $[a,\infty]^m$ for some
$a$.   More generally a subset of $V \times \G_\infty^m$ is bounded if its
pullback to $V \times \VF^{m}$ is bounded. 
\tcb{We shall say a subset $Y$ of $\std{V}$, resp. $\std{V \times \G_\infty^m}$,  is bounded if there exists a bounded definable subset $X$ of $V$, resp. $V \times  \G_\infty^m$, such that $Y \nsubseteq \std{X}$.}

 Let $Y$ be a  definable subset of $\G_\infty$.  Let $q$ be a definable
type   on $Y$. 
\tcr{If $Y$ is bounded there is a  unique $\a \in \tcb{\G_\infty}$,  such that $q$ concentrates
on any neighborhood of $\a$.  Indeed, consider the $q(x)$-definition of the formula $x>y$; it must have the form $y<\a$ or $y \leq \a$.}


 Let 
$V$ be a definable set and let $q$ be a definable type on $\std{V}$. 
Assume there exists $r \in \std{V}$ such that
for any continuous pro-definable map
$f : \std{V}\rightarrow \G_{\infty}$,
$\lim f_{*}(q)$ exists and $f (r)  = \lim f_{*}(q)$.
Then $\lim q$ exists and $r = \lim q$.

\begin{lem}\label{specialdef}
Let $V$ be an affine variety over a valued field and let
$q$ be a definable type on $\std{V}$.  
We have $\lim q = r $ if and only if
for any regular function $H$ on $V$, setting $h = \val \circ H$,
\[\tcb{h_{*}(r)} = \lim h_{*}(q).\]
\end{lem}

\prf One implication is clear, let us prove the reverse one. Indeed, by hypothesis, for any pro-definable neighborhood $W$ of \tcb{$r$},
$\tcb{q}$ implies $x \in W$. In particular, if $U$ is a definable neighborhood of
$\tcb{h_{*}(r)}$, $\tcb{q}$ implies $x \in f^{-1}(U)$, hence $\tcb{h_{*} (q)} $ implies $x \in U$.
It follows that $\tcb{\lim h_{*}(q) = h_* (r)}$.
\eprf

\begin{lem}  \label{bddcompact} Let $X$ be a  bounded definable subset of an algebraic variety $V$ over a valued field and let $q$ be a definable type on $\std{X}$.   Then $\lim q $ exists in $\std{V}$. 
\end{lem}

\prf   It is possible to partition $V$ into open affine subsets $V_i$
and $X$ into bounded definable subsets $X_i \nsubseteq V_i$.  We may thus assume $V$ is affine; and indeed that  $X$ is a bounded subset of $\Aa^n$.
For any regular \tcb{function} $H$ on $V$, 
setting $h = \val \circ H$,
$h(X)$ is a bounded
subset of $\G_{\infty}$ and 
$h_{*}(q)$ 
is a definable type on $h(X)$, hence has a limit $\lim h_{*}(q)$.

Now let $K$ be an algebraically closed valued field
containing the base of definition of $V$ and $q$.
Fix $\d \models q |K$ and  $d \models p_\d | K(\d)$, where
$p_\d$ is the type coded by the element $\d \in \std{V}$.     
Let $B=\G(K)$,  
$N = K(\d,d)$ and
$B'=\G(N)$. 
So $B$ is a 
divisible ordered abelian group.
We have $\G(N) = \G(K(\d))$ by orthogonality to $\G$ of $p_\d$.
Since $q$ is definable,
for any $e \in B'$,  $\tp(e/B)$ is definable; in particular the cut of $e$ over $B$ is definable. Set
 $B'_0 = \{b' \in B': (\exists b \in B) b<b' \}$.  
It follows  that if $e \in B'_0$ there exists an element $\pi(e) \in B \union \{\infty\}$ which is nearest
$e$.  Note $\pi: B'_0 \to B_\infty$ is an order-preserving retraction and a homomorphism in the obvious sense.  The ring
$R = \{a \in K(d): \val(a) \in B'_0 \}$ is a valuation ring of $K(d)$, containing $K$. 
Also $d$ has its coordinates
in $R$, because of the boundedness assumption on $X$.   Consider
 the maximal 
ideal  $M = \{a \in K(d): \val(a) > B \}$ and set $K'=R/M$. We have a canonical homomorphism
$R  \to K'$; let $d'$ be the image of $d$.  We have a valuation on $K'$
extending the one on $K$,
namely $\val(x+ M) = \pi (\val(x))$.  So $K'$
is a valued field extension of $K$, embeddable in some elementary extension.
Let $r = \tp(d'/K)$.  Then $r$ is  
definable and stably dominated; the easiest way to see that is to assume
 $K$ is
maximally complete (as we may); in this case stable domination follows
from $\G(K(d'))=\G(K)$ by \thmref{maxcomp}.    The fact that, for any $h$ as above, $r_{*}(h) = \lim h_{*}(q)$ is a direct consequence from the definitions.
\eprf

\begin{rem}Let $V$ be a definable set. According to \defref{defcomp}
a pro-definable $X \nsubseteq \std{V}$ is  definably compact if for any
definable type $q$ on $X$ we have $\lim q \in X$.  
Under this definition, any intersection of definably compact sets is definably compact, in particular an interval such as $\meet_n [0,1/n]$ in $\G$.
However we mostly have in mind strict pro-definable sets.
\end{rem}

\begin{lem}\label{closedcompact}  Let $V$ be an algebraic variety over a valued field, $Y$   a closed pro-definable subset of $\std{V}$.  
Let $q$ be a definable type on $Y$, and suppose $\lim q$ exists.  Then $\lim q \in Y$.

Hence if $Y$ is bounded and closed in $\std{V}$, then $Y$ is definably compact.
\end{lem}

\prf The fact that $\lim q \in Y$ when $Y$ is closed pro-definable follows from the definition of the topology on $\std{V}$.
The second statement thus follows from \lemref{bddcompact}.  \eprf

\begin{defn}
Let $T$ be a theory with 
universal domain
$\Uu$.
Let $\Gamma$ be a stably embedded sort with a $\varnothing$-definable linear ordering.
Recall
$T$ is said to be {\em metastable} over $\Gamma$ if for any small \index{metastable}
$C \nsubset \Uu$, the following   condition is satisfied:
 
(MS) For some small $B$ containing $C$, for any $a$ belonging to a finite product of sorts, $\tp (a /B, \Gamma (Ba))$ is stably dominated.
\end{defn}

Such a $B$ is called a {\em metastability base}. \index{metastability base}
It follows from \thmref{maxcomp} that $\ACVF$ is
metastable \tcb{over $\G$ and that any maximally complete algebraically closed valued field is a metastability base}.

\medskip
Let $T$ be any theory, $X$ and $Y$ be pro-definable sets, and $f: X \to Y$ a  pro-definable map.
Then $f$ induces a map $f_{\tcb{*}}: S_{def, X} \to S_{def, Y}$ from the set of definable types
on $X$ to the set of definable types on $Y$.
We remark that if $f$ is injective, then so is $f_{\tcb{*}}$.    This reduces to the case of definable $f: X \to Y$, where it is clear.
We now deal with surjectivity.

\begin{lem} \label{extend}   Let  $f: X \to Y$ a surjective pro-definable map between pro-definable sets.  
\begin{enumerate}
\item
  Assume $T$ is o-minimal.  Then $f_{\tcb{*}}: S_{def, X} \to S_{def, Y}$ is surjective.
\item
Assume $T$ is metastable over some o-minimal $\G$.  Then   $f_{\tcb{*}}: S_{def, X} \to S_{def, Y}$ is surjective.
\item  Assume $T$ is metastable over some o-minimal $\G$.    Then $f_{\tcb{*}}$  restricts to   a surjective map $\std{X} \to \std{Y}$.
\end{enumerate}
\tcb{In \textup{(2)}, if $X$, $Y$  and $f$ are defined over a metastability base $M$,  
then for any $M$-definable type $r$ on $Y$ there exists an 
$M'$-definable type $p$ on $X$, with $M'$  generated over $M$ by elements of $\G$, such that $r = f_* (p)$. More precisely, 
there exists a set  $A$, possibly infinite, an $M'$-definable type  $p'$ on the pro-definable set $W = \G^A$ and an
$M$-pro-definable map $H : W \to  \std{X}$ such that
$p = \int _{p'} H$ verifies $f_* (p) = r$.}
\end{lem}

\begin{rem}  \label{extend-r}   \begin{enumerate}  
\item  In general  surjectivity over 
a given base set  does not hold in (2) and (3) 
  (e.g. take $X$ a finite set, $Y$ a point).
\item It would also be possible to prove the C-minimal case analogously to the o-minimal one, as below.
  \end{enumerate}
  \end{rem}

\prf 
\tcb{Let us prove (1).}
First note it is enough to consider the case where $X$ consists of real elements.
 Indeed if $X$, $Y$ consist of imaginaries,  
find a set $X'$ of real elements and a surjective map $X' \to X$; then it suffices
to show $S_{def, X'} \to S_{def, Y}$ is surjective.

The statement reduces to the case that $X \nsubseteq U \times Y$ is a complete type, $f: X \to Y$ is the projection, and $U$ is one of the basic sorts.   Indeed, we can first let $U=X$ and
replace $X$ by the graph of $f$.  Any given definable type $r(y)$ in $Y$ restricts to some
complete type $r_0(y)$, which we can extend to a complete type $r_0(u,y)$ \tcb{over some model} implying $X$.
Thus we  can take $X \nsubseteq U \times Y$ to be complete.  
\tcb{Recall  that when $X = \limproj_j X_j$, we have $S_{def, X}= \limproj_j S_{def, X_j}$ 
naturally. Thus,  by transfinite induction, it is enough to consider the case
of 1-variable $U$.}

We can take $X,Y$ to be complete types with 
  $X \nsubseteq \G \times Y$,   and $f: X \to Y$ the projection.  It follows from completeness
  that for any $b \in Y$, $f \inv(b)$ is convex.  
 Let $r(y)$ be a definable type in $Y$.  
  Let $M$ be a model with $r$ defined over $M$, let $b \models r |M$, and consider
  $f \inv(b)$.

If   for any $M$, $x \in X \, \wedge \, f(x) \models r |M$ is a complete type $p|M$ over $M$, then
$x \in X \union p(f(x))$  already generates a definable type by \lemref{defgen} and we are done.
So, let us assume for some $M$, and $b \models p|M$, $x \in X \, \wedge \,  f(x) =b$ does not
generate a complete type over $M(b)$. 
Then there exists an $M(b)$-definable
set $D$ that splits $f \inv(b)$ into two pieces.   We can take $D$ to be an interval.
Then since $f \inv(b)$ is convex, one of the endpoints of $D$ must fall in $f \inv(b)$.
This endpoint is $M(b)$-definable, and can be written $h(b)$ with $h$ an $M$-definable
function.  In this case   $\tp(h(b), b /M)$ is $M$-definable, and has a unique extension
to an $M$-definable type.   
In either case  we found $p \in S_{def, X}$ with $f_{*}(p)=r$.  
Note  that the proof works when only $X$ is contained in the definable closure of an o-minimal definable set, for any pro-definable $Y$.

For the proof of (2) consider $r \in S_{def, Y}$.  Let $M$ be a metastability base, with $f$, $X$, $Y$, and $r$ defined over $M$.   Let
   $b \models r |M$, and let $c \in f \inv(b)$.    Let $b_1$ enumerate $\G(M(b))$.
   Then $\tp(b/M(b_1)) = r' | M(b_1)$ with $r'$ stably dominated, and $\tp(b_1/M) = r_1 |M $ with $r_1$ definable.
   Let  $c_1$ enumerate $\G(M(c))$; 
   then  $\tp(cb / M(c_1)) = q'|M(c_1)$ with $q'$  stably dominated.   \tcb{We have $q'=\tau(c_1)$ for some $M$-definable function
   into   the stable dominated types,  and $r'=\si(b_1)$ similarly.
By (1) (and stable embeddedness of $\G$), it is possible to extend
   $\tp(c_1b_1 /M ) \union r_1$ to a definable type $q_1(x_1,y_1)$ over some $M'$, where $M'$ can be taken to be generated over
   $M$ by elements of $\G$. Let  $c_1b_1 \models q_1 | M'$,
 and $cb \models q' | M'(c_1b_1)$.  So $\tp(b/M') = r|M'$.    Now $\tp(bc / M')$ extends to a definable type $p=\int_{q_1} \tau$ 
  by transitivity, and 
 $ f_{*}(p)= \int_{r_1} \si = r$.}
 
\tcb{Note that  the proof in \cite{hhm} 10.7 and 10.8 holds verbatim in 
the metastable setting, yielding that a definable type $p$ is stably dominated if and only if it is orthogonal to $\G$, as in \propref{equivstd}.
Thus,  the proof of (3) is similar to the proof of (2);}
in this case there is no $b_1$,
 and $q_1$ can be chosen so that $c_1 \in  M'$.  Indeed $\tp(c_1/M)$
 implies $\tp(c_1/M(b))$ so it suffices to take $M'$ containing $M(c_1)$.
 \eprf

\begin{remark}\label{extend+} 
There should be no difficulty to give an abstract version of \lemref{extend}; let us just mention one more case that 
we will require.   Say $T$ has the extension property if $f_{\tcb{*}}$ is always surjective, in the situation of \lemref{extend}.
First, let $T=   \Th( A)$, where $A$ is a linearly ordered  group with a definable convex subgroup $B$, such that ($*$) $B$ and $A/B$ are o-minimal.
Then $T$ has the extension property.  This is proved exactly as in the beginning of the proof of (1) in \lemref{extend}, by reduction to 1-types; here we can reduce
to $B$ and cosets of $B$ (all o-minimal) and to $A/B$.  Secondly, assume $T$ is metastable with respect to a linearly ordered group with ($*$);
then the proof of (2) shows that $T$ has the extension property.  

In particular, the theory $\ACV2F$ obtained from $\ACVF$ by expanding $\G$ by a predicate for a convex subgroup considered in \ref{ss8.3}  has the extension property.
\end{remark}

\begin{prop} \label{compactimage}
\tcb{Let $V$ and $V'$ be algebraic varieties over a valued field.}
Let $W$ be a definably compact pro-definable subset of $\std{V \times \Gamma_{\infty}^m}$,
and
let $f: \tcb{W} \to \tcb{\std{V' \times \Gamma_{\infty}^{m'}}}$ be a continuous  pro-definable morphism.  Then $f(W)$ is definably
compact. \end{prop}

\prf  
 Let $q$ be a definable type on $f(W)$.  By \lemref{extend} there exists
a definable type $r$ on $W$, with $f_{*}(r)=q$.    Since $W$ is  definably
compact, $\lim r$ 
exists and belongs to $W$.  But then $\lim q = f ( \lim r) $   belongs to $f(W)$
(since this holds after composing with any continuous morphism to $\G_{\infty}$). 
So $f(W)$ is definably compact. \eprf

 \begin{lem}\label{bounded}Let $V$ be an algebraic variety over a valued field, and let $W$ be a definably compact
 pro-definable subset of $\std{V \times \G_{\infty}^m}$.  Then $W$ is contained in $\std{X}$ for some
 bounded definable v+g-closed subset $X$ of
 $V \times \G_{\infty}^m$.
  \end{lem}

\prf By using \propref{compactimage} for the projections
$\std{V \times \G_{\infty}^m} \rightarrow \std{V}$ and
$\std{V \times \G_{\infty}^m} \rightarrow \G_{\infty}$, one may assume
$W$ is a pro-definable subset of $\G_{\infty}$ or $\std{V}$.
The first case is clear.    For the second one, one may assume $V$ is affine
contained in $\Aa^n$ with coordinates
$(x_1, \ldots, x_n)$.
Consider the function $\min  \val (x_i)$ on $V$, extended to $\std{V}$;
it is a continuous function on $\std{V}$.  
The image of $W$ is a definably compact subset of $\G_\infty$, hence is bounded
below, say by $\alpha$.  Let $X=\{(x_1,\ldots,x_n): \val(x_i) \geq \alpha \}$. 
Then $W \nsubseteq \std{X}$. 
\eprf

By a {\em countably pro-definable} set  we mean a  pro-definable set  isomorphic to one with   a countable inverse limit system.  Note that $\std{V}$ is countably pro-definable.  \index{countably pro-definable}

\begin{lem} \label{cc0}  Let $X$ be a strict, countably pro-definable set over a model $M$, $Y$ a relatively definable subset of $X$ over $M$.
If $Y \neq \varnothing$ then $Y(M) \neq \varnothing$.  \end{lem}

\prf   Write $X = \limproj_n X_n$ with transition morphisms $ \pi_{m,n}: X_m \to X_n$,  and  $X_n$ and $\pi_{m,n}$ definable. Let $\pi_n: X \to X_n$ denote the projection. Since $X$ is strict pro-definable, the image of $X$ in $X_n$ is definable; replacing
$X_n$ with this image, we may assume $\pi_n$ is surjective.     Since $Y$ is relatively definable,  it has the form $\pi_n \inv (Y_n)$ for some nonempty $Y_n \nsubseteq X_n$.
We have $Y_n \neq \varnothing$, so there exists $a_n \in Y_n(M)$.  Define
inductively $a_m \in Y_m(M)$ for $m >n$, choosing $a_{m} \in Y_m(M)$
with $\pi_{m,m-1}(a_m)=a_{m-1}$.    For $m<n$ let $a_m = \pi_{n,m}(a_n)$.  Then $(a_m)$ is an element of $\tcb{Y}(M)$.
\eprf

Let $X$ be a pro-definable set with a definable topology (in some theory).  
Given a model $M$, and an element
$a$ of $X$ in some elementary extension of $M$, we say that $\tp(a/M)$ has a limit $b$ if  $b \in X(M)$, and for any 
$M$-definable open neighborhood $U$ of $b$, we have $a \in U$.   This  extends the notion of a limit of a definable
type; if $a \models q|M$ with $q$ an $M$-definable type, the limits have the same meaning.  In the o-minimal setting
of $\G_\infty$, we show however that in fact limits appear only for definable types.

\begin{lem} \label{cc1} Let $M$ be an elementary submodel of $\tcb{\G}$. \tcb{Let $A$ be a set and let $\tcb{a \in \G_\infty^A}$}. Let  $p_0= \tp(a/M)$ and assume
$\lim p_0 $ exists.  Then there exists a \textup{(}unique\textup{)} $M$-definable type $p$
extending $p_0$.  \end{lem}

\prf  \tcb{It is enough to consider the case when $A$ is finite, so we may assume $\G_\infty^A
= \G_\infty^n$.}
In case $n=1$, $\tp(a/M)$ is determined by a cut in $\G_\infty(M)$.  If this cut is irrational then by definition
there can be no limit in $M$.  So this case is clear.  
We have to show that for any formula $\phi(x,y)$ over $M$, with $x=(x_1,\ldots,x_n)$ and  $y=(y_1,\ldots,y_m)$, the set
$\{c \in M: \phi(a,c) \}$ is definable.   Any formula is a Boolean combination of unary formulas and of formulas of 
the form:   
$\sum \alpha_i x_i + \sum  \beta_j y_j + 
\gamma \diamond 0$,
 where $i,j$ range over some subset of $\{1,\ldots,n\}$, $\{1,\ldots, m\}$ respectively, 
 $\alpha_i, \beta_j \in \Qq, \gamma \in \G(M)$,
and $\diamond \in \{=,<\}$.  This case follows from the case $n=1$ already noted, applied to $\tp(\sum\a_ia_i/M)$.  
\eprf

 \begin{prop}  \label{cc3} Let $X$ be a pro-definable subset of $\std{V} \times \G_\infty^m$ with $V$ an algebraic variety over a valued field.  Let $a$ belong to the closure of $X$.  Then there exists a definable
 type on $\std{V} \times \G_\infty^m$ concentrating on $X$, with limit point $a$.     \end{prop}
 
 \prf  \tcb{We may assume $V$ is affine, and by lifting the $\G$-coordinates to $\std{\Aa^m}$ and absorbing in the field coordinates,   that $m=0$.  
 Let $M$ be a maximally complete model of $\ACVF$ over which the data is defined. It is a metastability base.  Let $\mathcal F=(f_i)_{i\in I}$ list all functions  on $\std{V}$ of the form $\val (F)$ for a regular function $F$ on $V$, defined over $M$.
Since $\mathcal F_*(a)$ is a limit point of $\mathcal F(X)$, there exists a  type $q_M$ on $\mathcal F(X)$ over $M$, with limit $\mathcal F_*(a)$.    
By \lemref{cc1},  $q_M$ extends to an $M$-definable type $q$ with limit $\mathcal F_*(a)$.
By \lemref{extend}, there exists a definable type $p$ on $X$ such that 
 $\mathcal F_*(p) = q$.
 Furthermore one can assume $p$ is defined over  $M'= M \union E$ with $E \nsubset \G$
 and $p = \int _{p'} H$, where $p'$ is an $M'$-definable type on a pro-definable set $W= \G^A$ and $H$ is an $M$-pro-definable
 map $W \to \std{X}$.    
Let us prove that $\lim p = a$.  
Recall the canonical map
$ \vartheta : \doublewidehat{V} \to \std{V}$ from  \remref{doublest}
sending
a stably dominated
type $q$ on 
$\std{V}$ to $\vartheta (q) = \int_{q} \text{id}_{\std{V}}$.
By composing $H$ with $\vartheta$, one obtains
an $M$-pro-definable map 
$h: W    \to \std{V}$.
To prove that $\lim p = a$, it is enough to check that
$\lim h_* (p') = a$.
Now assume $V$ is embedded in $\mathbb{A}^m$ and consider  the morphisms   $J_d : \std{\mathbb{A}^m} \to L (H_d)$
defined in  \ref{modules}. For every $d \geq 0$, set $h_d = J_d \circ h$. Note that $h_d (W)$ is a $\G$-internal subset of $L (H_d)$ defined over $M$.
By \lemref{diagonal}, 
 there exists a finite number of 
bases  of $H_d$ over $M$ such that each  semi-lattice in $h_d (W)$ is diagonal for one of these bases.  
It follows that
there exists a common basis $B_d$, defined over $M$, that diagonalizes all semi-lattices $h_d(t)$ for $t \models p'$.
Since for any basis element
 $e \in B_d$, the
 valuative norm of $e$ according to the semi-lattice $h_d(t)$ is given by  the functions in $\mathcal F$ and $\lim q= \mathcal F_*(a)$, it follows that 
 $\lim h_d$  exists for all $d$. 
 Since by \thmref{protop}
 the morphism $J : \std{\Aa^m} \longrightarrow  \limproj_d L (H_d)$ induced by the system $(J_{d})$ is injective
 and induces a homeomorphism between
$\std{\Aa^m}$ and its image, it follows that
$\lim (J \circ h)_* (p') = J (a)$ and
$\lim h_* (p') = a$.}
\eprf

\begin{cor} \label{compactclosed}\tcb{Let $X$ be a pro-definable subset of $\std{V} \times \G_\infty^m$,  with $V$ an algebraic variety over a valued field.  If $X$ is definably compact, then $X$ is closed in
$\std{V} \times \G_\infty^m$.}
\end{cor}

\prf  The fact that $X$ is closed  is immediate from \propref{cc3} and the definition of definable compactness. 
\eprf

\begin{rem}\label{omin7}  Let  $\G$ be a Skolemized o-minimal structure, $a \in \G^n$.  Let $D$
be a definable subset of $\G^n$ such that  $a$ belongs to the topological closure $cl(D)$ of $D$.  Then there exists a definable type 
$p$ on $D$ with limit $a$, in the sense of  \defref{deflimit}.
Indeed, consider the family $F$ of all rectangles (products of intervals) whose interior
contains $a$.  This is a definable family, directed downwards under containment.   
By \lemref{cofinal-def}
 there exists a definable type $q$ on $F$ concentrating, for each $b \in F$,
 on $\{b' \in F: b' \nsubseteq b \}$.   Since $a \in cl(D)$, there exists a definable (Skolem) function $g$ such that
for  $u \in F$,  $g(u) \in u \meet D$.  To conclude it is enough to set $p =g_{*}(q)$. An alternative proof is provided, in our case, by \propref{cc3}.
It follows that if the  limit of any definable type on $D$ exists and lies in $D$, then $D$ is closed.
Conversely, if $D$ is bounded,  any definable type on $D$  will have a limit, and if $D$ is closed then this
limit is necessarily in $D$. 
  \end{rem}

 Even for $\Th(\G)$, definability of a type $\tp(ab/M)$ does not imply that $\tp(a/M(b))$ is definable.
For instance $b$ can approach $\infty$, while $a \sim \alpha b$ for some irrational real $\alpha$, 
i.e. $qb < a < q'b$ if $q,q' \in \Qq$, $q < \alpha < q'$.    However we do have:

\begin{lem} \label{cc5} Let $p$ be a definable type of $\G^n$, over $M$.  Then up to a definable change of coordinates,
$p$ decomposes as the join of two orthogonal definable types $p_f, p_i$, such that $p_f$ has a limit in $\G^m$, 
and $p_i$ has limit point $\infty^{\ell}$.  \end{lem}

\prf 
If $\alpha \in  \Qq^n$ and $x \in \G^n$, we write $\alpha \cdot x$ for the scalar product
$\sum_i \alpha_i x_i \in \G$.
  Let 
$\a_1,\ldots,\a_k$ be a maximal set of linearly independent vectors in $\Qq^n$ such that
the image of $p$ under $x \mapsto \alpha_i \cdot x$ has a limit point in $\G$.
Let $\b_1,\ldots,\b_{\ell}$ be a maximal set of vectors in $\Qq^n$ such that for $ x \models p|M$,
$\a_1 \cdot x,\ldots,\a_k \cdot x,\b_1 \cdot x,\ldots,\b_{\ell} \cdot x$ are linearly independent over $M$.  
 If $a \models p|M$, let $a'=(\a_1 \cdot a, \ldots, \a_k \cdot a)$, $a'' = (\b_1 \cdot a, \ldots,\b_{\ell} \cdot a)$.    For $\a \in \Qq(\a_1,\ldots,\a_k)$
 we have  that
 $\a \cdot a$ is bounded between elements of $M$.  On the other hand each $\b\cdot  a$, with $\b \in \Qq(\b_1,\ldots,\b_{\ell})$, satisfies
 $\b \cdot a > M$ or $\b \cdot a < M$.
For if $m \leq \b \cdot a \leq m'$ for some $\tcb{m, m'} \in M$, since $\tp(\b \cdot a/M)$ is definable it must have a finite limit,
contradicting the maximality of $k$.  It follows that   $\tp(\a \cdot a/  M) \union \tp(\b \cdot a /M)$ extends uniquely  to a complete 2-type, namely $\tp((\a \cdot a,\b  \cdot a) / M)$; in particular $\tp((\a \cdot a) + (\b \cdot a) / M)$
is determined; from this, by quantifier elimination, $\tp(a'/M)  \union \tp(a''/M)$ extends to a unique type in $k+\ell$ variables.
So $\tp(a'/M)$ and $\tp(a''/M)$ are orthogonal.  After some sign changes in $a''$, so that each coordinate is $>M$,
the lemma follows. 
\eprf

\begin{rem} \label{gv-alternate}
It follows from \lemref{cc5} that to check for definable compactness of $X$,  it suffices to check definable maps from definable types on $\G^k$ that either have limit $0$, or limit $\infty$.  From this an alternative proof of the  g- and v-criteria of Chapter \ref{specializations} for closure in $\std{V}$  can be deduced.   
\end{rem}

\tcb{For the sake of completeness we shall provide the proof of  the following lemma  from \cite{mst} here.}
\begin{lem}[\cite{mst} Lemma 2.19] \label{cofinal-def} \tcb{Let $P$ be a definable directed partial ordering in an o-minimal structure $\G$.  Then there exists a definable type $p$ cofinal in $P$.}  \end{lem}

\prf  \tcb{We assume $P$ is 0-definable, and work with 0-definable sets; we will find a 0-definable type 
with this property.  
Note first that we may replace $P$ with any 0-definable cofinal subset.   Also if $Q_1,Q_2$ are non-cofinal 
subsets of $P$, there exist $a_1,a_2$ such that  no element of $Q_i$ lies above $a_i$; but
by directedness there exists $a \geq a_1,a_2$; so no element of $Q_1 \union Q_2$ lies above $a$,
i.e. $Q_1 \union Q_2$ is not cofinal.  In particular
if $P = P' \union P''$,  at least one of $P'$, $P''$ is cofinal in $P$ (hence also directed).}

\tcb{If $\dim(P) = 0$ then $P$ is finite, so according to the above remarks we may assume it is
one point; in which case the lemma is trivial.  We use here the fact that in an o-minimal theory, 
any point of a finite 0-definable set is definable.}

\tcb{If $\dim(P)=n>0$, we can divide
$P$ into finitely many 0-definable sets $P_i$, each admitting a map $f_i: P_i \to \G$
with fibers of dimension $<n$.  We may thus assume
that there exists a 0-definable map $f: P \to \G$ with fibers of dimension $<n$.  Let $P(\g) = f \inv (\g)$, and $P(a,b) = f \inv (a,b)$.}

\begin{claim1}  \tcb{One of the following holds:
\begin{enumerate}
\item For any $a \in \G$, $P(a,\infty)$ is cofinal in $P$.
\item For some 0-definable $a \in \G$, for all $b > a$, $P(a,b)$ is cofinal.
\item For some 0-definable $a \in \G$,  $P(a)$ is cofinal.
\item For some 0-definable $a \in \G$, for all $b<a$, $P(b,a)$ is cofinal. 
\item  For all $a \in \G$, $P(-\infty,a)$ is cofinal.
\end{enumerate}}
\end{claim1}

\begin{proof}[Proof of the claim]  \tcb{Suppose (1) and (5) fail.  Then $P(a,\infty)$ is not cofinal in $P$ for some $a$; 
so $P(-\infty, b)$
must be cofinal, for any $b>a$.  Since (5) fails, the set $\{b: P(-\infty,b) \hbox{ is cofinal}\}$ is a nonempty proper definable subset of $\G$, closed upwards, hence of the 
form $[A,\infty)$ or $(A,\infty)$ for some 0-definable $A \in \G$.  In the former case, $P(-\infty,A)$ is cofinal,
but $P(-\infty,b)$ is not cofinal for $b<A$, so $P(b,A)$ is cofinal for any $b<A$; thus (4) holds.
In the latter case, $(-\infty,b)$ is cofinal for any $b>A$, while $(-\infty,A)$ is not; so $P([A,b))$
is cofinal for any $b>A$.  Thus either (2) or (3) hold.}  \end{proof}

\tcb{Let $p_1$ be a 0-definable type of $\G$, concentrating on sets $X$ with $f \inv(X)$ cofinal.
(For instance in case (1) $p_1$ concentrates on intervals $(a,\infty)$.)}

\begin{claim2}  \tcb{For any $c \in P$,
if $a \models p_1 | \{c\}$ then there exists $d \in P(a)$ with $d \geq c$.}
\end{claim2}
\begin{proof}[Proof of the claim]  \tcb{Let $Y(c) = \{x: (\exists y \in P(x)) (y \geq c) \}$.  Then $P \inv (\G \m Y(c))$ is not cofinal
in $P$, so it cannot be in the definable type $p_1$.  Hence $Y(c) \in p_1 | \{c\}$.}  \end{proof}

\tcb{Now let $M \models T$.  Let $a \models p_1 | M$.
By induction, let $q_a$ be an $a$-definable type, cofinal in $P(a)$, and let $b \models q_a | Ma$.
Then $\tp(ab/M)$ is definable.  If $c \in M$ then by Claim 2, there exists $d \in P(a)$ with $d \geq c$.
So the set $\{y \in P(a): \neg (y \geq c) \}$ is not cofinal in $P(a)$.  Therefore this set is not in $q_a$.  Since
$b \models q_a | Ma$, we have $b \geq c$.  This shows that $\tp(ab/M)$ is cofinal in $P$.}  \eprf

\begin{lem}  \label{compactclosedgamma}  Let $S$ be a definably compact definable subset  of  an o-minimal structure.   If ${\mathcal D}$ is a uniformly definable family of nonempty closed definable subsets of $S$, and ${\mathcal D}$ is directed \textup{(}the intersection of any two elements of ${\mathcal D}$ contains a third one\textup{)}, then $\meet {\mathcal D} \neq \varnothing$.  \end{lem}
\prf   
By \lemref{cofinal-def}
there exists a cofinal definable type $q(y)$ on ${\mathcal D}$ concentrating, for each $U \in  {\mathcal D}$,  on $\{V \in  {\mathcal D}: V \nsubset U \}$.     

    Using the lemma on extension of definable types,  \lemref{extend},   let $r(w,y)$
 be a definable type extending $q$ and implying $w \in U_y \meet S$. 
Let $p(w)$ be the projection of $r$ to the $w$-variable.  By definable compactness $\lim p=a$
exists.  Since $a$ is a limit of points in $D$, we have $a \in D$ for any $D \in {\mathcal D}$.  So $a \in \meet {\mathcal D}$. 
\eprf

\lemref{compactclosedgamma} gives another proof that a definably compact set is closed:  let   ${\mathcal D} = \{S \m U \}$, 
where $U$ ranges over basic open neighborhoods of a given point $a$ of the closure of $S$.

\begin{thm}\label{compactiffboundedclosed}
Let $V$ be an algebraic variety over a valued field, and let $W$ be a
 pro-definable subset of $\std{V \times \G_{\infty}^m}$.
 Then $W$ is definably compact if and only if it is closed and bounded.
\end{thm}

\prf If $W$ is  definably compact
it is closed and bounded by \corref{compactclosed} and \lemref{bounded}.
If $W$ is closed and bounded,
its preimage $W'$ in
$\std{V \times \Aa^m}$ under $\mathrm{id}\times \val$ is also
closed and bounded, hence definably compact by \lemref{closedcompact}.
It follows from \propref{compactimage} that $W$ is definably compact.
\eprf

For $\G^n$, 
\thmref{compactiffboundedclosed} is a special case of \cite{peterzil-steinhorn},
 Theorem 2.1.

\begin{prop} \label{closed=closed}Let $V$ be a  variety over a valued field $F$, and let  $W$ be an $F$-definable subset of
$V \times \Gamma_{\infty}^m$.
Then $W$ is v+g-closed \textup{(}resp. v+g-open\textup{)}
if and only if $\std{W}$ is closed \textup{(}resp. open\textup{)} in $\std{V}$.
\end{prop}
\prf   A Zariski-locally v-open set is v-open, and similarly for g-open,  hence for v+g-open.  So we may assume $V = \Aa^n$ and by pulling back to
$V \times \Aa^m$ that $m = 0$.    It  is enough to prove the statement about closed subsets.
Let $V_\a = (c \Oo)^n$ be the closed polydisc of valuative radius $\a=\val(c)$.
Let $W_\a = W \meet V_\a$, so $\std{W_\a} = \std{W} \meet \std{V_\a}$.  
Then $W$ is v-closed if and only if $W_\a$ is v-closed for each $\a$; by \lemref{gcriterion1}, the same holds for g-closed;
also   $\std{W}$ is closed
if and only if $\std{W_\a}$ is closed for each $\a$.  This reduces the question to the case of bounded $W$.

By
\lemref{v+g}, if  $W$ is v+g-closed then $\std{W}$ is closed.  
In the reverse direction, if $\std{W}$ is closed it is definably 
compact.    It follows that $W$ is v-closed.  For otherwise there exists an accumulation point $w$ of $W$,
with $w=(w_1,\ldots,w_m) \notin W$.  Let $\d(v) = \min_{i=1}^m \val(v_i - w_i)$.  Then $\d(v) \in \G$
for $v \in W$, i.e. $\d(v)< \infty$.  Hence the induced function $\d: \std{W} \to \G_\infty$ also has image
contained in $\G$; and $\d(\std{W})$ is definably compact.  It follows that $\d(\std{W})$ has a maximal point 
$\g_0 < \infty$.  But then the $\g_0$-neighborhood around $w$ contains no point of $W$, a contradiction.

\tcb{To conclude it is enough} to show that if $\std{W}$ is  definably compact, then $W$ must be g-closed.  This follows from
 \corref{gcriterion1.1}. 
\eprf

\begin{cor}\label{compact=closed+bounded} Let $V$ be an algebraic variety over a valued field, and let $W$ be a definable subset of $ V  \times \G_\infty^m$.  Then $W$ is bounded and v+g-closed  {\em if and only if}
 $\std{W}$ is definably compact.  \end{cor}
  
\prf Since $W$ is v+g-closed if and only if $\std{W}$ is closed by \propref{closed=closed}, this  is a special case of
\thmref{compactiffboundedclosed}. \eprf

 \begin{lem} \label{closedmap2}  Let $V$ be an algebraic variety over a valued field and let  
 $Y$ be 
 a v+g-closed, bounded \tcb{definable} subset of $V \times \G_{\infty}^m$.   Let
 $W$ be a definable subset of $V' \times \G_\infty^m$, with $V'$ another variety.
Let $f : Y \to W$ be a definable map.
 Assume 
 $\std{f}: \std{Y} \to \std{W}$ is continuous.   Then $\std{f}$ is   a \tcr{definably closed map in the sense of \textup{Definition \ref{defn:defclosed}}}.
   \end{lem} 

\prf  \tcb{We may assume $f$ is surjective, in which case $\std{f}$ is also surjective 
by \lemref{extend}.}
By \propref{closed=closed}  and \thmref{compactiffboundedclosed}, 
$ \std{Y}$ is definably compact and any  \tcb{pro-definable} closed subset  of
$ \std{Y}$ is definably compact, so the result follows from \propref{compactimage} and \corref{compactclosed}.  
\eprf

\begin{lem}  \label{complete2}  Let $X$ and $Y$ be
v+g-closed, bounded definable subsets of a product of an algebraic variety over a valued field with
some $\G_{\infty}^m$. Then, \tcr{the image of any closed pro-definable subset of $\std{X} \times \std{Y}$ under}
the projection $\pi: \std{X} \times \std{Y} \to \std{Y}$ is \tcr{closed}.
\end{lem}  

\prf   \tcr{Let $h$ denote the canonical morphism $\std{X \times Y}  \to   \std{X} \times \std{Y}$.
The map $\pi \circ h: \std{X \times Y} \to \std{Y}$ is definably closed by \lemref{closedmap2}.
Let $Z$ be a closed pro-definable subset of $\std{X} \times \std{Y}$. Its preimage $h^{-1} (Z)$ is a closed pro-definable subset of
$\std{X \times Y}$, hence $\pi (Z) = (\pi \circ h) (h^{-1} (Z))$ is closed.}
\eprf

\begin{cor}  
Let $U$ and $V$ be v+g-closed, bounded definable subsets of a product of an algebraic variety over a valued field with
some $\G_{\infty}^m$.
If  $f: \std{U} \to \std{V}$ is a pro-definable  morphism with closed graph, then $f$ is continuous.  
 \end{cor}
  
\prf  By \lemref{complete2} \tcr{and Remark \ref{injcl}}, the projection  $\pi_1$ from the graph of $f$ to $\std{U}$ is a homeomorphism onto its image.
The   projection $\pi_2$  is continuous.  Hence $f=\pi_2 \tcr{\circ} \pi_1 \inv$ is continuous.
\eprf
 
  \begin{lem} \label{closedmap3}  Let $f: V \to W$ be a proper morphism of algebraic varieties.  Then $\std{f}$ is a \tcr{definably} closed map.
   So is $\std{f} \times \mathrm{Id} : \std{V} \times \G_\infty^m \to \std{W} \times \G_\infty^m$.  
   \end{lem} 
  
  \prf    
\tcb{Note that} $\std{V \times \G_\infty^m}$ can be identified with a subset $S$ of $\std{V} \times \tcb{\std{\Aa^m}}$ (projecting on generics of balls around zero in the second coordinate); with this identification, $\std{f} \times \mathrm{Id}$ identifies with the restriction of $\std{f \times \mathrm{Id}_{\Aa^m}}$  to $S$.  Thus  
 the second statement, for  $V \times \G_\infty^m$, reduces to first for the case of the map $f \times \mathrm{Id}: V \times \Aa^m \to W \times \Aa^m$.   

To prove the statement on $f: V \to W$, let $V', W'$ be complete varieties containing $V,W$, and let $\bar{V}$ be the closure
of the graph of $f$ in $V' \times W'$.  In the Zariski topology, the map $\mathrm{Id} \times f:  V' \times V \to V' \times W$ is closed by properness (universal closedness). 
The image of the diagonal on $V$, under this map, is the graph of $f$, a subset of $V \times W$; so it is Zariski closed as a subset of $V' \times W$.  
Let $g : V \to V \times W$ given by $g(v)=(v,f(v))$; so $g$ is the composition of   the isomorphism $v \mapsto (v,f(v))$ of $V$ onto the graph of $f$, with the inclusion of the graph of $f$
in $V \times W$.  Both of these induce \tcr{definably} closed morphisms on $\std{ \ }$-spaces, so $\std{g}$ is \tcr{definably} closed.  

Let $\pi: \bar{V} \to W'$ be the projection.  Now $\std{\pi}$ is
a \tcr{definably} closed map by \lemref{closedmap2}.  So $\std{\pi}\circ \std{g} = \std{ \pi \circ g} = \std{f} $ is \tcr{definably} closed.   
(We could also obtain the result directly from \propref{cc3}.)   \eprf

 \begin{rem}   The previous lemmas apply also to $\infty$-definable sets. \end{rem}

 \begin{cor}\label{radicial}Let $f: V \to W$ be a \tcb{finite} radicial surjective morphism of algebraic varieties over a valued field $K$.
 Then $\std{f} : \std{V}\to \std{W}$ is a homeomorphism.
 \end{cor}
 
 \prf \tcb{Since for  $f$ to be radicial means that for any field extension $K'$ the induced map $V (K') \to W (K')$ is injective,} $f$ is an isomorphism in the category of definable sets.
 Thus, $\std{f} : \std{V} \to \std{W}$ is a bijection, say by \lemref{extend}.
 On the other hand, $f$ being a \tcb{universal} homeomorphism for the Zariski topology, it is proper,
 thus $\std{f}$ is  \tcr{definably} closed by  \lemref{closedmap3}, hence a homeomorphism, \tcr{cf. Remark \ref{injcl}}.
 \eprf

  \begin{lem} \label{max}Let $X$ be a v+g-closed bounded definable subset of an algebraic variety $V$ over a valued field. 
     Let $f: X \to \G_{\infty}$ \tcr{be a definable map which is}
v+g-continuous.  Then the maximum of $f$ is attained on $X$, similarly if $X$ is a closed bounded pro-definable
subset of $\std{V}$.
\end{lem}

\prf By \lemref{basic}, $f$ extends continuously to $F: \std{X} \to \G_{\infty}$.  
By \propref{closed=closed} and \thmref{compactiffboundedclosed}
$\std{X}$ is definably
compact. 
It follows from \propref{compactimage} that
$F(\std{X})$ is a definably
compact subset of $\G_{\infty}$ and hence has a maximal point $\g$.  Take
$p$ such that  $F(p)=\g$,
let $c \models p$,   then $f(c) = \g$.  
\eprf

\begin{prop}\tcr{Let $V$ be an algebraic variety over a valued field. Then $V$ is complete if and only if $\std{V}$ is definably compact.}
\end{prop}

\prf \tcr{If $V$ is projective, $\std{V}$ is definably compact by \thmref{compactiffboundedclosed}. If $V$ is complete, there exists a surjective morphism
of algebraic varieties $f : W \to V$ with $W$ projective by Chow's lemma. Hence $\std{V}$ is definably compact, since
$\std{V} = \std{f} (\std{W})$ by \lemref{extend}  
and
$\std{V}$ is Hausdorff  by \propref{hausdorff}. 
Conversely,  assume 
$\std{V}$ is definably compact. One may embed $V$ as a Zariski dense open subset of a complete variety $W$ by Nagata's theorem. In particular,
$\std{V}$ is dense in  $\std{W}$. On the other hand, $\std{V}$ is closed in  $\std{W}$ by \corref{compactclosed}. Thus 
$\std{V}=\std{W}$ and $V=W$.}
\eprf


\chapter{A closer look at the stable completion}\label{sec5}  

{\small \noindent \textbf{Summary.}
In \ref{modules} we give a description of $\std{\Aa^n}$ in terms of spaces of  semi-lattices which will be used in \ref{ss6.2}.
This is provided by constructing  a topological embedding of $\std{\Aa^n}$ into the inverse limit of
a system of spaces of  semi-lattices $L (H_d)$ endowed with the linear topology, where $H_d$ are finite-dimensional vector spaces.
The description is extended in \ref{reppn} to the projective setting.
In \ref{ss5.3} we relate the linear topology to the one  induced by the finite level morphism 
$\std{\Aa^n} \to L (H_d)$.\par\bigskip}

\section{$\std{\Aa^n}$ and spaces of  semi-lattices}  \label{modules}
Let $K$ be a model of $ \ACVF$ and let
$V $ be a $K$-vector space of dimension $N$.
By a {\em lattice} in $V$ we mean a free $\Oo$-submodule of rank $N$. \index{lattice}
By a {\em semi-lattice} in $V$ we mean an $\Oo$-submodule $u$ of $V$, such that for some \index{semi-lattice}
$K$-subspace $U_0$ of $V$ we have $U_0 \nsubseteq u$ and $u/U_0$ is a lattice in $V/U_0$.  
Note that  every semi-lattice is uniformly definable with parameters and that
the set $L(V)$ of semi-lattices in $V$ is definable.  \nomenclature{$L(V)$}{space of lattices}
Also,
a definable $\Oo$-submodule $u$ of $V$ is a semi-lattice if and only if there
is no $0 \neq v \in V$ such that $Kv \meet u = \{0\}$ or $Kv \meet u = \Mm v$ where $\Mm$ is the maximal ideal.

We define a topology on $L(V)$ as follows.  The pre-basic open
 sets are those of the form
$\{u: h \notin u \}$ and those of the form $\{u: h \in \Mm u \}$, where $h$
is any element of $V$.   The finite intersections of these sets clearly form an ind-definable family. 
 We call this family the {\em linear topology} on $L(V)$.  \index{linear topology}

Any finitely generated $\Oo$-submodule of $K^N$ is generated by $\leq N$ elements; hence the  intersection of any finite number of open sets of the second type is the intersection of $N$ such open sets.  However this is not the case for the first kind.

  Another description can be given in terms of linear seminorms.  By a {\em linear seminorm} on a finite-dimensional $K$-vector space $V$ \index{linear seminorm}
we mean a definable map $w: V \to \G_\infty$ with $w(x_1+x_2) \geq \min (w(x_1)+w(x_2))$ and $w(cx) = \val(c)+w(x)$.  
\begin{lem}\label{snsl}
\tcb{If $w$ is a linear seminorm on $V$,
$\Lam_w = \{x \in V: w(x) \geq 0 \}$ is a semi-lattice.
Conversely, any semi-lattice
$\Lam \in L(V)$ has the form $\Lam=\Lam_w$ for a unique $w$, namely
$w(x) = - \inf \{\val (\lambda): \lambda x \in \Lam\}$.} \qed
\end{lem}

 We may thus identify
$L(V)$ with the set of linear seminorms on $V$.  On the set of linear seminorms there is a natural topology, with   basic open sets of the form $\{w: (w(f_1),\dots,w(f_k) ) \in O\}$,
with $f_1,\ldots, f_k \in V$ and $O$ an open subset of $\G_\infty^k$.   
  The linear {topology} on $L(V)$ coincides with the linear seminorm topology.

Finally, a description as a quotient by a definable group action:
Fix a  basis for $V$, and let $\Lambda_0$ 
be the $\Oo$-module generated by this basis.
Given $M  \in \mathrm{End}(V)$, let $\lambda(M) = M \inv (\Lambda_0)$.  
We identify $\Aut (\Lambda_0)$ with the group of automorphisms $T$ of $V$ with
 $T(\Lambda_0)=\Lambda_0$.   So $T \cong \Aut (\Lambda_0) \cong GL_n(\Oo)$.
 We give $ \mathrm{End}(V) = M_n(V)$ the valuation topology, viewing $M_n(V)$ as a copy of $K^{n^2}$.

\begin{lem}\label{linear-top1}   The mapping  $\lambda:  M \mapsto \lambda (M) = M \inv (\Lambda_0) $ is surjective and continuous.  It  induces
a bijection between $\Aut (\Lambda_0) \backslash \mathrm{End}(V) $  and
 $L(V)$. \end{lem} 

\prf  It is clear that
$M \mapsto \tcb{\lambda} (M)$ is a surjective map from
$\mathrm{End}(V)$ to $L(V)$, and also that $\Lambda(N)=\Lambda (TN)$ if $T \in \Aut(\Lambda_0)$.
Conversely suppose $\Lambda(M)=\Lambda(N)$.  Then  $M$ and $N$ have the same kernel
$E = \{a: Ka \nsubseteq M \inv (\Lambda_0) \}$.   So $N M \inv$ is a well-defined
homomorphism $M V \to  N V$.    Moreover, 
$MV \meet \Lambda_0$
 is a free $\Oo$-submodule of $V$, and
$(NM \inv)( MV \meet \Lambda_0) = (NV \meet \Lambda_0)$.
Let $C$ (resp. $C'$) be a free $\Oo$-submodule of $\Lambda_0$ complementary
to $MV \meet \Lambda_0$ (resp. $NV \meet \Lambda_0)$, and let
$T_2: C \to C'$ be an isomorphism.    Let $T = (NM \inv) | (MV \meet \Lambda_0) \oplus T_2$.  Then  $T \in \Aut(\Lambda_0)$, and $NM \inv \Lambda_0 = T \inv \Lambda_0$,
so (using $\ker M = \ker N$) we have $M \inv \Lambda_0 = N \inv \Lambda_0$. 
This shows the bijectivity of the induced map $\Aut (\Lambda_0) \backslash \mathrm{End}(V) \to 
L(V)$.

Continuity is  clear:  the inverse image of   $\{u: h \notin u \}$ is $\{M:  Mh \notin \Oo^n \}$, while  
the inverse image of $\{u: h \in \Mm u \}$ is $\{M: Mh \in \Mm^n \}$.  These are in fact v+g-closed.
 \eprf

The mapping $\lambda$ is far from being closed or open, with respect to the v-topology on $\mathrm{End}(V)$.  In that topology,
$\Aut (\Lambda_0)$ is open, so $\Oo^n$ is an isolated point in the  pushforward topology.

We say a subset of $L(V)$  is  {\em bounded} if  its pullback with respect to the map above \index{bounded subset of $L(V)$}
to $\mathrm{End}(V)$ is bounded. 
 Note that if $X \nsubset M_n(K)$ is bounded then so is $\GL_n(\Oo) X$ (even $M_n(\Oo) X$);
so the image of a bounded set is bounded.  Thus a bounded subset of $L(V)$ is a set of semi-lattices admitting bases in a common  bounded ball in $V$.
  In terms of linear seminorms, if $\Lambda_w$ ranges over a bounded set, then for any $h \in V$,
$w(h)$ lies in a bounded subset of $\G_\infty$, i.e. bounded on the left.

\begin{lem} \label{top0}  The space $L(V)$ with the linear {topology} is Hausdorff.  Moreover, any 
definable type on a bounded subset of $L(V)$ has a \textup{(}unique\textup{)} limit point in $L(V)$.
\end{lem}

\prf  Let $u' \neq  u'' \in L(V)$.  One, say $u'$, is not a subset of the other.  Let $a \in u', a \notin u''$.  Let $I = \{c \in K: ca \in u'' \}$.  Then $I = \Oo c_0$ for some $c_0$ with $\val(c_0)>0$.  Let $c_1$ be such that $0<\val(c_1)<\val(c_0)$ and let $a' =c_1a$.
Then $a' \in \Mm u'$ but $a' \notin u''$.  This shows that $u'$ and $u''$ are separated by the disjoint open sets $\{u: a' \notin u \}$ and $ \{u: a' \in \Mm u \}$.

For the second statement, let   $Z$ be a bounded  set of   linear seminorms.   
  Let $p$ be a definable type on $Z$.
Let
  $w(h) = \lim_p w_x(h)$, where $w_x$ is the norm corresponding to $x \models p$.  This limit is not $-\infty$ since $Z$ is bounded.
  It is easy to see that $w$ is a linear seminorm.  Moreover  any pre-basic open set containing $\Lambda_w$ must also contain a generic
point of $p$.   
\eprf

Let $H_{m;d}$ be the space of polynomials of degree $\leq d$ in $m$
variables.   
For the rest of this section $m$ will be fixed;   we will hence suppress the index and write $H_d$. 
\tcb{For $p$ in $\std{\Aa^m}$, consider the definable $\Oo$-submodule of $H_d$
\[J_d(p) = \{h \in H_d: p_{*}( \val (h)) \geq 0 \}.\] \nomenclature{$J_d(p)$}{semi-lattice attached to $p$}
Since $h \mapsto p_{*}( \val (h))$ is a linear seminorm,
$J_d(p)$ belongs to $L (H_d)$.}
Hence we have a mapping $J_d= J_{d,m}:  \std{\Aa^m} \to L(H_d)$ given by
$p \mapsto J_d(p)$.  It is clearly a continuous map, when $L (H_d)$ is given the linear {topology}:  $f \notin J_{d}(p)$ if and only 
if $p_{*}(\val (f)) <0$, 
and $f \in \Mm J_d (p)$ if and only if $p_{*}(\val (f)) >0$.

\begin{thm} \label{protop}  The system $(J_{d})_{d=1,2,\ldots}$ induces
a continuous morphism of pro-definable sets
\[
J : \std{\Aa^m} \longrightarrow  \limproj_d L (H_d),
\]
the transition maps \nomenclature{$J$}{a morphism $\std{\Aa^m} \rightarrow  \limproj_d L (H_d)$}
$L(H_{d + 1}) \to L (H_d)$ being the natural maps induced by the
 inclusions $H_d \nsubseteq H_{d + 1}$.
 The morphism $J$ is injective and induces a homeomorphism between
$\std{\Aa^m}$ and its image.
\end{thm}

\prf  Let $f : \Aa^m \times  H_d  \to \G_{\infty}$ given by
$(x, h) \mapsto \val (h (x))$.
Since $J_{d}$ factors through
$Y_{H_d, f}$, $J$  is a morphism of pro-definable sets (here
$Y_{H_d, f}$ is defined as in the proof of \thmref{prodef}).

For injectivity, recall that types on $\Aa^n$ correspond to equivalence classes of $K$-algebra
morphisms
$\varphi : K [x_1, \ldots, x_n] \to F$ with $F$ a valued field, with $\varphi$ and $\varphi'$ equivalent if they are restrictions of a same $\varphi''$.
In particular, if $\varphi_1$ and $\varphi_2$ correspond to different types,
one should have
\begin{equation*}
\begin{split}
\{f \in K [x_1, \ldots, x_m] :  \val (\varphi_1 (f))& \geq 0 \}
\not=  \\
            \{f \in  & K [x_1, \ldots, x_m] :   \val (\varphi_2 (f)) \geq 0 \},
\end{split}
\end{equation*}
whence the result.

We noted already continuity.  Let us prove that $J$ is an open map onto its image.
The topology on $\std{\Aa^n}$ is generated by sets $S$ of the form $\{p: p_{*}(\val (f)) > \g \}$ or $\{p: p_{*}(\val(f)) < \g \}$,
where $f \in H_d$ for some $d$. 
For such an $S$,
$J (S) = \pi_d^{-1} (J_d (S))$,
with $\pi_d : \limproj L (H_d') \to L (H_d)$
the natural map.
Thus, it is enough to check  that $J_{d}(S)$ is open.
   Replacing $f$ by $cf$ for appropriate $\tcb{c}$, it suffices to consider $S$ of the form
 $\{p: p_{*}(\val(f)) > 0 \}$ or $\{p: p_{*}(\val(f)) < 0 \}$.  
Now the image of these sets is precisely the intersection with the image of $J$
of the open sets    $p_{*}(\val(f))  \notin \Lambda$ or $p_{*}(\val(f)) \in \Mm \Lambda$. 
 \eprf

\begin{rem}\label{rem:5.1.5}
\tcb{Note that the image of $J$ consists of all sequences $(\Lam_d \in L(H_d))_{d=1,2,\ldots}$, with corresponding
linear seminorms $w_d$ on $L(H_d)$, 
such that for any $f_i \in  H_{d_1}, f_2 \in H_{d_2}$ we have, $w_{d_1+d_2}(f_1f_2) = w_{d_1}(f_1)+w_{d_2}(f_2)$.}
 \end{rem}

\begin{rem} 
\thmref{protop} describes the $\std{V}$-topology in terms of the  linear {topology}  
when one takes all ``jets'' into account. \tcb{It remains interesting to describe the topology induced on the $S_n$
by the  individual 
maps  $J_{d}$.  The image of $J_d$ is described in section 7 of \cite{imagdef}; it may hint at the induced topology as well.}
 \end{rem}

\section{A representation of   $\std{\Pp^n}$} \label{reppn}

Let us define the tropical projective space
$\mathrm{Trop} \, \Pp^n$, for $n \geq 0$, as the quotient \nomenclature{$\mathrm{Trop} \, \Pp^n$}{tropical projective space}
$\G_\infty^{n+1} \m \{\infty\}^{n+1} / \G$ where $\G$ acts
diagonally by translation.  This space may be topologically embedded in $\G_\infty^{n + 1}$
since it  can be identified with 
\[\{(a_0,\ldots,a_n) \in \G_\infty^{n+1}:   \min a_i = 0 \}. \]

 Over a valued field $L$, we have a canonical definable map $\tau: \Pp^n \to  \mathrm{Trop} \, \Pp^n$,
sending $[x_0: \ldots : x_n] $ to \[[v(x_0): \ldots: v(x_n)] = ( (v(x_0) - \min_i v(x_i), \ldots,  v(x_n)-\min_i v(x_i))).\]

Let us denote by $H_{n+1;d, 0}$ the set of  homogeneous polynomials in $n+1$ variables of degree $d$ with coefficients in the valued field sort.  Again we view $n$ as fixed and omit it from the notation, 
letting $H_{d,0} = H_{n+1;d,0}$.  Denote by $H_{d,m}$ the definable subset of $H_{d, 0}^{m + 1}$ consisting of 
$m+1$-tuples of homogeneous polynomials with no common zeroes other than the trivial zero.   
Hence, one can consider 
the image $PH_{d,m}$  of $H_{d,m}$ in the projectivization $P(H_{d, 0}^{m + 1})$.
We have a morphism $c : \Pp^n \times  H_{d,m} \to \Pp^{m}$, given by 
$c ([x_0: \ldots: x_n], (h_0, \ldots, h_m)) =  [h_0(x):\ldots :h_m(x)]$.   
Since  $c (x, h)$ depends only on the image of $h$ in $P H_{d,m}$, we obtain a morphism $c: \Pp^n \times  PH_{d,m} \to \Pp^{m}$.
Composing $c$  with the map $\tau:  \Pp^m \to  \mathrm{Trop} \, \Pp^m$, we obtain
 $\tau:  \Pp^n \times  PH_{d,m} \to \mathrm{Trop} \,\Pp^m$.
 For $h $ in $ PH_{d,m} $ (or in $H_{d,m}$), we denote by
 $\tau_h$ the map 
 $x \mapsto \tau (x, h)$. Thus  $\tau_h$ extends
 to a map 
 $\std{\tau_h} : \std{\Pp^n} \to \mathrm{Trop} \,\Pp^m$.

 Let $T_{d, m}$ denote the set of  functions $PH_{d, m} \to \mathrm{Trop} \, \Pp^n$  of the form
 $h \mapsto  \std{\tau_h} (x)$ for some $x \in \std{\Pp^n}$.
 Note that $T_{d, m}$ is definable.

 \begin{prop}\label{reppnprop}The  space $\std{\Pp^n} $ may be identified 
 via the 
 canonical mappings $\std{\Pp^n} \to T_{d,m}$
 with  the projective limit of the spaces $T_{d,m}$.
 If one endows $T_{d, m}$ with the topology  induced from the Tychonoff topology, this identification is a homeomorphism. \qed
 \end{prop}
  
 The proof of the proposition is a straightforward reduction to the affine case, by using standard affine charts, that we omit.
 
 \begin{rem}\label{reppnrem}By composing with the embedding
  $\mathrm{Trop} \, \Pp^m \to \G_\infty^{m+1}$, one gets a definable map
 $\std{\Pp^n} \to \G_\infty^{m+1}$.  The topology on $\std{\Pp^n}$ can be defined directly using the above maps
into $\G_\infty$, without an affine chart.   
\end{rem}

\section{Relative compactness}\label{ss5.3}
  Let $H$ be a finite-dimensional $K$-vector space. In this section we take $L(H)$ to have the linear topology.
  
 We say that \tcr{an $\infty$-definable} subset $X$ of $L(H)$
is {\em relatively compact for the linear topology} if for any definable type $q$ on $X$, if $q$ has a limit point $a$ in $L(H)$, then $a \in X$.  \index{relatively compact for the linear topology}
The closed sets of the linear {topology} are clearly relatively compact.   
  
  \begin{lem} \label{Jrc}    The image of a closed \tcr{pro-definable} set by the  morphism $J_d : \std{\Aa^m} \to L(H_d)$  is  relatively compact. 
\end{lem}

  \prf Let $Y$ be a closed subset of $\std{\Aa^m}$.  Let $q$ be a definable type
on $J_d(Y)$, and let $b$ be a limit point of $q$ for the linear {topology}.   We have to show that $b \in J_d(Y)$. The case $d=0$ is easy as $J_0$ is a constant map, so assume $d \geq 1$.  We have in $H_d$ the monomials $x_i$.  For some nonzero $c_i' \in K$ we have
$c_i'x_i \in b$, since $b$ generates $H_d$ as a vector space.  Choose a nonzero $c_i$ such that $c_ix_i \in \Mm b$.  
 Let $U = \{b':  c_ix_i \in \Mm b', i=1,\ldots,m \}$.  Then $U$ is a pre-basic open neighborhood of $b$; 
 as $b$ is a limit point of $q$, it follows that $q$ concentrates on $U$.   Note that $J_d \inv (U)$ is contained
 in $\std{B}$ where $B$ is the polydisc $\val(x_i) \geq -\val(c_i), i=1,\ldots,m$.  Thus $J_d \inv (U)$ is bounded. 
 By \lemref{extend} we may lift $q$ to a definable type $p$ on $Y \meet \std{B}$. Then as $Y \meet \std{B}$ is closed and bounded, $p$ has a limit point $a$ on $Y\cap \std{B}$.  By continuity we have $J_d(a)=b$, hence $b \in J_d(Y)$.  
 \eprf 

It follows, writing  $X=J_d (J_d \inv (X))$, that if a definable set in $L (H_d)$ is an intersection of relatively compact sets, then it is  itself relatively compact.
Thus the relatively compact sets are the closed sets of a certain topology.   

 For $b \in L (H)$, we denote by 
 $v_b$  the linear seminorm associated with $b$.

We consider definable metrics in a different sense than in  \ref{Smetrics}.   
  Namely a {\em definable g-metric} \index{definable g-metric}
on a definable set $X$
is a map $d: \tcb{X^2} \to \G_{\geq 0}$, satisfying symmetry, the triangle law $d(x,z) \leq d(x,y)+d(y,z)$, and $d(x,y)=0$ iff $x=y$.  It induces a topology
in the obvious way (from the g-topology on $\G$).
 
\begin{rem} \begin{enumerate} 
\item  Let $L^*(H)$ be the set of lattices on $H$.    This is easily seen to be a dense subset of $L(H)$ for the linear topology.
\item  On $L^*(H)$,  we have a definable g-metric defined as follows.
Each lattice corresponds to an actual linear norm on $H$,
i.e. a linear seminorm such that $w (h) = \infty$ iff $h = 0$.  We obtain a definable g-metric  between norms:
\[ \tcb{d} (w,w') = \sup \{ |w(h) - w'(h) |:  h \in H \m (0) \} \]
\item This g-metric induces a definable topology on $L^*(H)$ (in the sense of Ziegler), finer than the linear topology.
\item  The space $L(H)$ fibers over the (Grassmannian) space of linear subspaces of $H$, 
and each fiber admits a similar metric.   
\item  $L^*(H)$ is not linearly open in $L (H)$ when $H$ is of dimension $ > 1$.  
\tcb{Fix a lattice $\Lambda$ in $H$.}
Given a finite number of vectors $h_1,\ldots,h_k$ and $h_1',\ldots,h_l'$
with $h_i \notin \tcb{\Lambda}$, $ h_i' \in \Mm \tcb{\Lambda}$, let $f: H \to K$ be a   linear map so that $\ker(f)$ does not pass through any of the vectors $h_i$ or $h_i'$; renormalize it so  such that $f(\tcb{\Lambda}) =\Oo$.  Then $\val (f(h_i)) <0$ and $\val ( f(h'_j)) >0$.   So $h_i \notin f \inv (\Oo)$, $h_j' \in f \inv (\Mm) = \Mm f \inv (\Oo)$.  Hence
$f \inv (\Oo)$ belongs to a prescribed neighborhood of $\tcb{\Lambda}$ in $L(H)$, but it is not a lattice as soon as $H$ is of dimension $ > 1$.  
\item    Let $-1 \in \G$ be negative, let $m \geq 1$, and let $Y$ be the set of lattices in $L(K^m)$   of volume $-1$:
$Y= \{ M \Oo^n:    \val (\det (M)) = -1\}$. 
  Then $Y$ is relatively definably compact, $\Oo^n \notin Y$, but $\Oo^n \in cl(Y)$ in the linear topology.
To see this last point  consider the lattice     $M \Oo^n$, where $M$ is a  lower-triangular \tcb{matrix with rows} $(a,0), (c,d)$, where $\val(a)=\val(c)<0$,
$\val(d)<0$ and $\val(a)+\val(d)=-1$.   
\end{enumerate}\end{rem} 


\chapter{$\G$-internal  spaces}         \label{secgammaint}

{\small \noindent \textbf{Summary.}
This chapter is devoted to the topological structure of
$\G$-internal  spaces. 
The main results about the topological structure
of
$\G$-internal  spaces
are proved in \ref{ss6.2}. In \ref{ss6.1} several related issues are discussed.
The rather technical results in \ref{ss6.3} are used in \ref{ss6.4} which deals with the study of the topology of relatively
$\G$-internal  spaces. 
\par\bigskip}

\section{Preliminary remarks}\label{ss6.1}\tcb{Let $V$ be an algebraic variety over a valued field.
Recall an iso-definable subset $X$ of $\std{V}$ is
said to be  $\G$-internal if it is in pro-definable bijection with a 
definable set which is  $\G$-internal. \index{$\G$-internal iso-definable}
A number of delicate issues arise here.  Let us say a pro-definable subset $X$ of $\std{V}$ is {\em $\G$-parameterized} if there \index{$\G$-parameterized}
exists \tcb{a definable subset $Y$  of $ \G^n$, for some $n$, and a pro-definable map
 $g: Y \to \std{V}$ with image $X$}.     
By the following example, there exists $\G$-parameterized subsets of $\std{V}$ which are not iso-definable, whence not  $\G$-internal.}

\begin{example} \label{nonisodefinable}\tcb{Let $A$ be a base structure consisting of a trivially valued field $F$ and a value group
containing $\Zz$.  
Let $\varphi = \sum_{i = 0}^{\infty} a_i x^i$ be a formal series with coefficients
$a_i \in F$. Assume $\varphi$ is not algebraic. For any nonnegative integer $n$, set $\varphi_n(x) = \sum_{i \leq n} \a_i x^i$.
For any $\gamma \in \Gamma_\infty$,
consider the complete type $p_{\gamma}$, in the variables $x$ and $y$, generated over $A (\gamma)$ by
the generic type of the closed ball $\val (x) \geq 1$ together with
the formulas
\[
\val(y - \varphi_n (x)) \geq \min(n+1,\gamma).
\]
If  $\gamma \leq n_0$, for some $n_0 \in \Nn$,
 this is the image under $(x,z) \mapsto (x,z+\varphi_{n_0} (x))$ of the generic type of the polydisc $\val (x)\geq 1,\val (z) \geq \gamma$.  If $\gamma > \Nn$, it is the type described in \exref{exotic}.  
Consider the continuous pro-definable map $g: \G_\infty \to \std{\Aa^2}$
 sending
 $\gamma$ to $p_\gamma$.  With the notation of \ref{ss6.5}, we have $p_\gamma \in  \stda{\Aa^2}$ iff
 $\gamma \leq n_0$ for some $n_0 \in \Nn$.
 The composition of $g$ with the projection to the space of lattices on polynomials of degree $\leq n$ is constant for $\gamma \geq n$.
 The image of $g$ is $\G$-parameterized, but is not iso-definable and hence not
 $\G$-internal.}  \end{example}

 \tcb{However, the image of a complete type is    iso-$\infty$-definable, as the following lemma shows.}

 \begin{lem}   \label{embeddedgammatypes} \tcb{Let $P \nsubset \G^n$ be the solution set of a complete type over some base structure
$A$.  Let $Y$ be a pro-definable set and $f: P \to Y$  be  a pro-definable map.  Then the kernel
of $f$ is a definable equivalence relation $E$ on $P$.}

\tcb{In case $Y = \std{V}$,   with $V$ an algebraic variety over a valued field, we have $\dim(P/E) \leq \dim(V)$, where the former
 dimension is  the o-minimal dimension and the latter, the dimension of the algebraic variety $V$.}
 \end{lem} 
\prf  \tcb{Write $Y = \limproj_i Y_i$, where $(Y_i)_{i \in I}$ is a directed system of definable sets, and denote by $\pi_i: Y \to Y_i$
the natural projection.
Let $E_i$ be the kernel of $f_i = \pi_i \circ f$.  Since $\dim(P/E_i)$ is non-decreasing, there exists some element $0$ of $I$
such that, for $i \geq 0$, $E_i$ splits each $E_{0}$-class 
into finitely many classes.  Using  elimination of imaginaries
for $\G$, there exists an $A$-definable map $\phi_i: P \to \G^n$ such that $xE_i y$ iff $\phi_i(x)=\phi_i(y)$.  
So the image under $\phi_i$ of each $E_0$-class is finite.   In particular for each $E_0$-class $X$, 
some element $c \in X$ has smallest possible image $\phi_i(c)$, under 
  the lexicographic ordering on $\G^n$.    But all elements $c \in X$ have the same type:  if $c,c' \in X$,  
  then $\tp(c/A)=\tp(c'/A)$
  since $c,c' \in P$;  let $d = \phi_0(c) = \phi_0(c')$,   so $\tp(c/Ad)=\tp(c'/Ad)$.  Thus all elements $c \in X$ have smallest possible image under $\phi_i$, i.e. they have the same image under $\phi_i$,
  so $X$ is a single $E_i$-class.  This shows that $E_i=E_0$ for all $i \geq 0$.     It follows that the kernel of $f$ is $E_0$, 
  and $f(P)$ is iso-$\infty$-definable.}
  
 \tcb{Now assume $Y = \std{V}$.  By \corref{g9} and \remref{rg9} there exist finitely many polynomials
$h_1,\ldots,h_r$ such that $h =(\val (h_1),\ldots,\val (h_r))$ induces an injective map on 
$f(P)$.  The image of $h$ in $\G^r$ has dimension $\leq \dim(V)$, proving the dimension inequality.}\eprf

The above discussion referred to the pro-definable category;  we will now move to topological questions.  When
concerned with the definable category alone,   there is no point mentioning $\G_\infty$, since $\infty$ has the same role as any other element.
But from the point of view of the definable topology, the point $\infty$ does not have the same properties as  any points of $\G$, nor of the point $0$ of $[0,\infty]$; $\G_\infty$ does
not (even locally) embed into $\G^n$, and the point $\infty$ must be taken into account.

\begin{defn}Let $V$ be an algebraic variety over a valued field and let $X$ be an iso-definable $\G$-internal subset of $\std{V}$ (or of $\std{V} \times \G_{\infty}^s$, for some $s$).
We say  $X$ is {\em topologically $\G$-internal}  \index{topologically $\G$-internal}
if $X$ is pro-definably homeomorphic to a definable subset of $\G_{\infty}^r$, for some $r$.
\end{defn}

\begin{rem}\label{tgi}In   \ref{ss6.2}, we shall prove  that, when $V$ is quasi-projective, for any 
iso-definable $\G$-internal subset $X$ of $\std{V}$ there exists a pro-definable continuous injection $f : X \hookrightarrow \G_{\infty}^r$, for some $r$.
In particular, if $X$ is definably compact, $f$ is a homeomorphism and $X$ is topologically $\G$-internal.
\tcr{Indeed, in this case the image of a closed pro-definable subset of $X$ is a closed subset of $\G_\infty^r$. Thus
$f$ is definably closed and being injective it is closed, cf. Remark \ref {injcl}.}
 In general, we do not know whether every
 iso-definable $\G$-internal subset of $\std{V}$ is topologically $\G$-internal.
 The ones that will occur in our constructions will always be contained in some definably compact  iso-definable $\G$-internal set, thus will be topologically $\G$-internal.
 \end{rem}

\medskip

 We \tcb{now} discuss briefly the role of parameters.
We fix a valued field $F$.   The term ``definable'' refers to $\ACVF _F$.  Varieties are assumed defined over $F$.
At the level of definable sets and maps,  $\G$ has elimination of imaginaries.   
\tcb{Let us say that $\G$ admits {\em topological elimination of imaginaries} if whenever} \index{topological elimination of imaginaries}
$X \nsubseteq \G_\infty^n$ and $E$
is a closed definable equivalence relation on $X$, 
  there exists a definable
map $f: X \to \G_\infty^n$ inducing a homeomorphism between the topological quotient $X/E$,
and $f(X)$ with the topology induced from $\G_\infty^n$. 
 \tcb{It seems that any o-minimal expansion of $\RCF$   \nomenclature{$\RCF$}{the theory of real closed fields}
 admits elimination of imaginaries in the topological sense.}
 
In another direction, the pair $(\kk,\G)$ also eliminates imaginaries (where $\kk$ is the residue field, with induced structure), and   so does $(\RES,\G)$, where $\RES$ denotes the generalized residue structure of \cite{HK}. 
 \nomenclature{$\RES$}{generalized residue structure}
 However,  $(\kk,\G)$ or $(\RES,\G)$ do not eliminate imaginaries topologically.  One reason for this, due to Eleftheriou (cf. \remref{omin-tri2},  \cite{elef}) and valid already for $\G$, is that the theory $\DOAG$ of divisible ordered
abelian groups 
is not sufficiently flexible to identify simplices of different sizes.  
\nomenclature{$\DOAG$}{the theory of divisible ordered abelian groups}
 A more essential reason for us is the existence of spaces with nontrivial Galois action on cohomology.  For instance take $\pm \sqrt{-1} \times [0,1]$ 
with $\pm \sqrt{-1} \times \{0\}$ and $\pm \sqrt{-1} \times \{1\}$ each collapsed to a point.   However for connected spaces topologically embedded in $\G_\infty^n$, the
Galois action on cohomology is trivial.   Hence there is no embedding of the above circle  in $\G_\infty^n$ compatible with  the Galois action.
The best we can hope for is that it may be embedded in a twisted form $\G_\infty^w$, for some finite set $w$;
after base change to $w$, this becomes isomorphic to $\G_\infty^n$.  \tcb{It will follow from \thmref{G-embed-0}} that
 such an embedding in fact exists for \tcb{topologically} $\G$-internal sets.

It would be interesting to study more generally the definable spaces occurring
as closed iso-definable subsets of $\std{V}$ parameterized by a subset of 
$\VF^n \times \G^m$.   In the case of $\VF^n$ alone, a  key example should be the set of generic points of
subvarieties of $V$ lying in some constructible subset of the Hilbert scheme.  This includes
the variety $V$ embedded with the valuation topology via  the simple points functor (\lemref{simple}); 
possibly other components of the Hilbert scheme obtain the valuation topology too, but the different components
(of distinct dimensions) are not topologically disjoint.

\section{Topological structure of $\G$-internal subsets}\label{ss6.2}

\begin{lem} \label{affineimage}  
Let $V$ be a quasi-projective variety over an infinite valued field $F$, and let $f: \G^n \to \std{V}$ be 
$F$-definable.   There exists an affine open $V' \nsubseteq V$ with $f (\G^n) \nsubseteq \std{V'} $.   If $V=\Pp^n$, there exists a linear hyperplane $H$ such that $f(\G^n) \meet \std{H} = \varnothing$.   \end{lem}

\prf  Since $V$ embeds into $\Pp^n$, we can view $f$ as a map into $\std{\Pp^n}$; so we may assume   $V=\Pp^n$.  
 For $\g \in \G^n$, let $V(\g)$ be the linear Zariski closure of $f(\g)$; i.e. the intersection of all hyperplanes $H$ 
 such that $f(\g)$ concentrates on $H$.  The intersection of $V(\g)$ with any $\Aa^n$ is  the zero set of all linear
polynomials $g$  on $\Aa^n$ such that $f(\g)_{*} (g) = 0$.  So $V(\g)$ is  definable uniformly in $\g$.  
  Now $V(\g)$ is an $\ACF_F$-definable set, with canonical parameter \nomenclature{$\ACF$}{the theory of algebraically closed fields}
  $e(\g)$; by elimination of imaginaries in $\ACF_F$, we can take $e(\g)$ to be a tuple
  of field elements.  But functions
 $\G^n \to \VF$ have finitely many values (every infinite definable subset of $\VF$ contains an open subset,
 and admits a definable map onto $\kk$).  So there are finitely many sets $V(\g)$.  Let $H$ be any hyperplane
 containing none of these.  Then no $f(\g)$ can concentrate on $H$.    
\eprf

\tcb{Let $K$ be a model of $\ACVF$ and let $H$ be a $K$-vector space of dimension $n$.
 We shall make use of}
 the space $L(H)$ of semi-lattices considered in  \ref{modules}. Given a basis
$v_1,\ldots,v_n$ of $H$, we say that a semi-lattice is diagonal  if it is a direct sum $\sum_{i=1}^n I_i v_i$, 
with $I_i$ an ideal of $K$ or $I_i = K$.

\begin{lem} \label{diagonal}  Let $Y$ be a  $\G$-internal subset of $L(H)$.  Then there exists a finite number of 
bases $b^1,\ldots,b^{\ell}$ of $H$ such that each $y \in Y$ is diagonal for some $b^i$.  
If $Y$ is defined over a valued field $F$,
these bases can be found over $F^{\alg}$.  
\end{lem}  

\prf  For $y \in Y$, let $U_y = \{h \in H: K h \nsubseteq y \}$.  Then $U_y$ is a subspace of $H$, definable from $Y$.
The Grassmannian of subspaces of $H$ is an algebraic variety, and has no infinite $\G$-internal definable subsets.
Hence there are only finitely many values of $U_y$.  Partitioning $Y$ into finitely many sets we may assume
$U_y=U$ for all $y \in Y$.  Replacing $H$ by $H/U$, and $Y$ by $\{y/U: y \in Y \}$, we may assume $U=(0)$.
Thus $Y$ is a set of lattices.

Now the lemma follows from Theorem 2.4.13 (iii) of \cite{hhmcrelle}, except that in this theorem one considers
$f$ defined on $\G$ (or a finite cover of $\G$) whereas $Y$ is the image of $\G^n$ under some definable function $f$.   In fact the proof of 2.4.13 works for functions from $\G^n$; however we will indicate how to deduce
the $n$-dimensional case from the statement there, beginning with 2.4.11.   We first formulate a relative version of 2.4.11.  Let $U=G_i$ be one of the unipotent groups considered in 2.4.11 (we only need the case of $U=U_n$, the full strictly upper triangular group).  Let $X$ be a definable set, and let $g$ be a definable map on $X \times \G$,
with $g(x,\g)$ a subgroup of $U$, for any $(x,\g)$ in the domain of $g$.  Let $f$ be another definable map on
$X \times \G$, with $f(x,\g) \in U/g(x,\g)$.  Then there exist finitely many definable functions $p_j: X \to \G$,
with $p_j \leq p_{j+1}$,   definable functions $b_j$ on $X$, such
 that letting $g_j^*(x) = \meet_{p_j(x) <\g < p_{j+1}(x)} g(x,\g)$ we have $b_j(x) \in U/g^*_j(x)$, and
 \[ (*) \ \ \ \   f(x,\g) = b_j(x) g (x,\g)\] whenever  $p_j(x) < \g < p_{j+1}(x)$.
 This relative version follows immediately  from 2.4.11 using compactness, and noting that $(*)$ determines
 $b_j(x)$ uniquely as an element of $U/g_j^*(x)$. 
 
By induction, we obtain the multidimensional version of 2.4.11:

Let  $g$ be a definable map on a  definable 
subset $I$ of  $\Gamma^n$, with 
$g(\gamma)$ a subgroup of $U$ for each $\gamma\in I$. Suppose 
$f$ is also a definable map on $I$, with $f(\gamma)\in U/g(\gamma)$. 
Then there is a partition of $I$ into finitely many definable 
subsets $I'$ such that for each $I'$ there is $b\in U$ with 
$f(\gamma)=bg(\gamma)$ for all $\gamma \in I'$.

   To prove
this for $\G^{n+1} = \G^n \times \G$, apply the case $\G^n$ to the functions $b_j,g_j$ as well as
$f,g(x,p_j(x))$ (at the endpoints of the open intervals).  

Now the lemma  for the multidimensional case follows as in \cite{hhmcrelle} 2.4.13.  Namely, 
each lattice $\Lambda$ has a triangular $\Oo$-basis; viewed as a matrix, it is an element of the triangular group $B_n$.   So there exists an element $A \in U_n$ such that $\Lambda$ is diagonal for $A$, i.e.
$\Lambda$ has a basis $DA$ with $D \in T_n$ a diagonal matrix.  If $D'A'$ is another basis for $\Lambda$ of the same form,
we have $DA = ED'A'$ for some $E \in B_n(\Oo)$.  Factoring out the unipotent part, we find that $D\inv D' \in T_n(\Oo)$.   So $D$  
is well-defined modulo $T_n(\Oo)$, the group $D \inv B_n(\Oo) D$ is well-defined, we have $D \inv E D' \in D \inv B_n(\Oo) D \meet U_n$, and the matrix $A$ is well-defined up to translation by an element of $g(\Lambda) =  D \inv B_n(\Oo) D \meet U_n$.   By the multidimensional 2.4.11, since $Y$ is  $\G$-internal, it admits a finite partition into definable subsets $Y_i$, such that for each $i$, there exists a basis $A$ diagonalizing each $y \in Y_i$.  

Moreover, $A$ is uniquely defined up to $\meet_{y \in Y_i} g(y)$.    The rationality statement now follows from \lemref{diagonal+}.                   \eprf

\begin{lem}  \label{diagonal+} Let $F$ be a valued field,   let $h$ be an $F$-definable subgroup of the unipotent group $U_n$, and let $c$ be an $F$-definable coset of $h$. 
 Then   $c$ has a point in $F^{\alg}$.  If $F$ has residue characteristic 0, or if $F$ is trivially valued and perfect,
 $c$ has a point in $F$.  \end{lem}

\prf 
As in \cite{hhmcrelle}, 2.4.11, the lemma can be proved for all unipotent algebraic groups by induction on dimension, so we are reduced to the case of the one-dimensional unipotent group $G_a$. 
In the nontrivially valued case the statement is clear for $F^{\alg}$, since $F^{\alg}$ is a model.  
If $F$ is nontrivially valued and has equal characteristic $0$,   any definable ball has a definable point, obtained by averaging a definable finite set of points.

There remains the case of trivially valued, perfect  $F$.  In this case  the subgroup must be $G_a,(0), \Oo$ or $\Mm$.  The group $\Oo$ has no other $F$-definable cosets.  As for $\Mm$ the definable cosets correspond to definable elements of the residue field; as the residue field
(isomorphic to $F$) is perfect, the definable closure is just the residue field itself;  but each element of
the residue field of $F$ is the residue of a (unique) point of $F$.   \eprf 

\begin{rem} \label{parameters}  Is the rationality statement in \lemref{diagonal+} valid in positive characteristic,
for the groups encountered in \lemref{diagonal}, i.e. intersections of conjugates of $B_n(\Oo)$ with $U_n$?
This is not important for our purposes since the partition of $Y$ may require going to the algebraic closure
at all events.  \end{rem}

\begin{cor} \label{g9}  
Let $X \nsubseteq \std{\Aa^N}$ be iso-definable over 
an algebraically closed valued field $F$ and $\G$-internal. Then for some 
$d$, and finitely many polynomials $h_i$ of degree $\leq d$, the map $p \mapsto (p_{*}(\val (h_i)))_i$ is injective on $X$.
\end{cor}  

\prf  By \thmref{protop}, the maps 
\[p \mapsto J_d(p) = \{h \in H_d: p_{*}(\val (h)) \geq 0 \}\]  
 separate points on $\std{\Aa^N}$ and hence on $X$.  So for each $x \neq x' \in X$,  
 for some $d$, $J_d(x) \neq J_d(x')$.  Since  $X$ is iso-definable, for some fixed $d$, $J_d$ is injective on $X$.
 Let $F$ be a finite set of bases as in \lemref{diagonal}, and let $\{h_i\}$ be the set of elements of these bases.
 Pick $x$ and $x'$ in $X$; if $x_{*}(\val (h_i)) = x'_{*}(\val(h_i))$ for all $i$, we have to show that $x=x'$, or equivalently that $J_d(x)=J_d(x')$;  by symmetry it suffices to show that $J_d(x) \nsubseteq J_d(x')$.  Choose a basis, say $b=(b^1,\ldots,b^m)$, such that
 $J_d(x)$ is diagonal with respect to $b$; the $b^i$ are among the $h_i$, so $x_{*}(b^i) = x'_{*}(b^i)$ for each $i$.
 It follows that $J_d(x) \meet Kb^i = J_d(x') \meet Kb^i$.  But since $J_d(x)$ is diagonal for $b$, it is generated
 by $\union_i (J_d(x) \meet Kb^i)$; so $J_d(x) \nsubseteq J_d(x')$ as required.  
 \eprf

 \begin{rem}\label{rg9}\tcb{Let us observe that the proof goes through for iso-$\infty$-defina\-ble sets $X$, definably parameterized by an $\infty$-definable subset of $\G$.    
 In quoting \lemref{diagonal},
 note that an  $\infty$-definable subset of a definable set such as $L(H)$ is always contained in a definable set,
 and in the present case in a $\G$-internal one.}
 \end{rem}
  
\begin{prop}\label{G-embed-1} Let $V\subset \Pp^N$ be a quasi-projective variety over a valued field $F$.
  Let $X \nsubseteq \std{V}$ be  $F$-iso-definable and $\G$-internal.  Then   there exist $m$, $d$ 
  and $h \in H_{d, m}(F^{\alg})$ such that, with the notations
  of  \textup{\ref{reppn}}, the restriction $\std{\tau_h} : \std{\Pp^n} \to \mathrm{Trop} \,\Pp^m$ to $X$ is injective.
  If $V$ is projective and $X$ is closed,  $\std{\tau_h}$ restricts to a homeomorphism between $X$ and its image. 
 \end{prop}

 \prf  We may take $V= \Pp^N$.   Note that if $\std{\tau_h}$ is injective, and $g \in \Aut(\Pp^n)= \mathrm{PGL}(N+1)$, it is clear
 that $\std{\tau_{h \circ g}}$ is injective too. 
 By \lemref{affineimage}, there exists a linear hyperplane $H$ with $\std{H}$ disjoint from $X$.  We may 
 assume $H$ is  the hyperplane $x_0 =0$.    Let $X_1 = \{(x_1,\ldots,x_N): [1 :x_1 :\ldots :x_N] \in X \}$.
 By  \corref{g9},   there exist finitely many polynomials $h_{1},\ldots,h_r$ 
 such that $p \mapsto (p_{*}(h_i))_i$ is injective on $X_1$.  Say $h_i$ has degree $\leq d$.
 Let $H_i(x_0,\ldots,x_d) = x_0^d h_i(x_1/x_0,\ldots,x_d/x_0)$, and let $h = (x_0^d,\ldots,x_N^d,H_1,\ldots,H_r)$, $m=N +r$.  
 Then $h \in H_{d, m}$, and it is clear that $\std{\tau_h}$ is injective on $X$.  
 \eprf 

\begin{thm}   \label{G-embed-1c} Let $V$ be a quasi-projective variety over a valued field $F$.
  Let $X \nsubseteq \std{V}$ be  $F$-iso-definable and $\G$-internal.
Then there exists an $F$-definable map
\tcb{$\b: {V} \to [0,\infty]^w$, for some finite set $w$ definable over $F$,
such that 
$\std{\b}: \std{V} \to [0,\infty]^w$ is continuous and}
 restricts to an injective $F$-definable continuous map 
$\a : X \to [0,\infty]^w$.
\end{thm}

\prf By \propref{G-embed-1}, such  map $\b_a$ exists over a finite Galois extension $F(a)$ over $F$ with values in $[0,\infty]^n$.
Let $w_0$ be the set of Galois conjugates of $a$ over $F$
and set $w = w_0 \times \{1, \dots, n\}$.   Define $\b: \tcb{V}  \to[0,\infty]^w$ by
taking all  the conjugates of the function $\b_a$.  Then the statement is clear.  
\eprf 

\thmref{G-embed-1c} applies only when the base structure is a valued field; it may not have elements of $\G$ other than $\Qq$-multiples of valuations of field elements.
We now extend the result to the case when the base structure may \tcb{contain additional} elements of $\G$.  

\begin{thm}\label{G-embed-0}  Let $A$ be a base structure consisting of a field $F$, and a set $S$ of elements of $\G$.  Let $V$ be a projective variety over $F$,
and let
$X$ be an $A$-iso-definable and $\G$-internal subset of $\std{V}$.   Then there exists 
 an $A$-definable continuous injective map  $\phi: X \to [0,\infty]^w$ for some finite $A$-definable set $w$.
 If  \tcb{furthermore $X$ is closed}, then $\phi$ is a  topological 
embedding. \end{thm}

 \prf   
   We have $\acl(A) = \dcl(A \union F^{\alg})  = F^{\alg}(S)$ 
  \tcb{by \lemref{3.4.12}}.  
  It suffices to show that  a continuous, injective $\phi: X \to [0,\infty]^n$ is definable  over $\acl(A)$, for then 
  the descent to $A$ can be done as in \thmref{G-embed-1c}.  So we may assume $F=F^{\alg}$, hence $A=\acl(A)$.
  We may also assume $S$ is finite, since the data is defined over a finite subset. 
  Say $S=\{\g_1,\ldots,\g_n\}$.  Let $q$ be the generic type
  of field elements $(x_1,\ldots,x_n)$ with $\val(x_i)=\g_i$.  Then $q$ is stably dominated.   
If $c \models q$, then by   \propref{G-embed-1} 
there exists an $A(b)$-definable  \tcb{continuous injective map} $f_b: X \to \G^n$ for some $n$ and some $b \in F(c)^{\alg}$.
Since $q$ is stably dominated, and $A=\acl(A)$,  $\tp(b/A)$ extends to a stably dominated $A$-definable type $p$.
If $(a,b) \models p^2 |A$ then $f_a f_b \inv: X \to X$; but $\tp(ab/A)$ is orthogonal to $\G$ while
$X$ is $\G$-internal, so the canonical parameter of  $f_a f_b \inv $ is defined over $A \union \G$
and also over $A(a,b)$, hence over $A$.  Thus $f_a f_b \inv  = g$.  If $(a,b,c) \models p^3$
we have $f_b f_c \inv = f_a f_c \inv = g$ so $g^2=g$ and hence $g=\mathrm{Id}_X$.  So $f_a=f_b$, and $f_a$ is $A$-definable, as required.
\tcr{If $X$ is closed  in $\std{V}$, it is definably compact, hence
$\phi$ is a  topological 
embedding by Remark \ref{tgi}.}
 \eprf

 \begin{rem}\label{rem:basetopint}
 \tcb{With the notation in \thmref{G-embed-0}, 
 if $X$ is topologically $\G$-internal, the morphism $\phi: X \to [0,\infty]^w$ induces a homeomorphism
 between $X$ and its image $Y$. Indeed, $X$ is definably homeomorphic
 to a definable subset $Y'$ of $\G_{\infty}^s$ for some $s$, and any definable continuous map $Y' \to Y$ is a homeomorphism.}
 \end{rem}
  
 \section{Guessing definable maps by regular algebraic maps}\label{ss6.3}

\begin{lem} \label{g2}  Let $V$ be  a normal, irreducible, complete variety, $Y$ an  irreducible variety, $X$
a closed subvariety of $\tcb{V}$,
$g:Y \to X \nsubseteq V$ a dominant constructible map \tcb{\textup{(}i.e. $\ACF$-definable\textup{)}} with finite fibers, all defined over a   field $F$.  Then there exists 
a  pseudo-Galois covering $f: \tV \to V$  such that each component $U$  of $f^{-1} (X)$
dominates $Y$ rationally, i.e. there exists a dominant rational map $g: U \to Y$
over $X$.   
\end{lem}

\prf    First an algebraic version.  Let $K$ be a field, $R$ an integrally closed subring, $G: R \to k$ a  
ring homomorphism onto a field $k$. Let $k'$ be a finite field extension.  Then there exists a finite normal 
 field extension $K'$
and a homomorphism $G': R' \to k''$ onto a field, where $R'$ is the integral closure of $R$ in $K'$, 
such that $k''$ contains $k'$.

Indeed we may reach $k'$ as a finite tower of 1-generated field extensions, so we may assume $k'=k(a)$
is generated by a single element.   Lift the monic minimal polynomial of $a$ over $k$ to a monic polynomial $P$
over $R$.  Then since $R$ is integrally closed, $P$ is irreducible.  Let $K'$ be the splitting field of $P$.  The kernel of
$G$ extends to a maximal ideal $M'$ of the integral closure $R'$ of $R$ in $K'$, and $R'/M'$ is clearly
a field containing $k'$.  

To apply the algebraic version let $K=F(V)$ be the function field of $V$.  Let $R$ be the local ring of $X$, i.e. the ring
of regular functions on some Zariski open set not disjoint from $X$, and let $G:R \to k$ be the evaluation homomorphism to the function field  $k=F(X)$ of $X$.   Let $k'=F(Y)$ the function field of $Y$, and  $K'$, $R'$, $G'$, $M'$ and $k''$ be as above.   Let  $f: \tV \to V$ be the normalization of $V$ in $K'$.  Then $k''$ is the function field of a component 
$X'$ of $f^{-1} (X)$, mapping dominantly to $X$.  Since $k'$ is contained in $k''$ as extensions of $k$ there exists a dominant rational map $g: X' \to Y$ over $X$.   But $\Aut(K'/K)$ acts transitively on the components of 
$f \inv(X)$, proving the lemma.   \eprf

\begin{lem} \label{g3}   Let $V$ be an algebraic variety over a field $F$,   $X_i$ a finite number of \tcb{locally} closed subvarieties, $g_i: Y_i \to X_i$ a surjective constructible
map with finite fibers.    Then there exists a surjective finite morphism of varieties $f: \tV \to V$
and 
a finite number of  \tcb{locally closed} subsets $U_{ij}$ of \tcb{$f \inv (X_i)$} and 
 morphisms $g_{ij}: U_{ij} \to Y_i$
 such that,  \tcb{for every $i$,  and every $a \in X_i$, $b \in Y_i$, $c \in \tV$}
with $g_i(b)=a$ and $f(c)=a$, 
we have  $c \in U_{ij}$ and $b = g_{ij}(c)$ for some $j$.
Furthermore, if $V$ is normal, 
 we may take $f: \tV \to V$
to be a pseudo-Galois covering.
  \end{lem}

\prf  If the lemma holds for each irreducible component $V_j$ of $V$, with   $X_{j,i} = V_j \meet X_i$
and $Y_{j,i} = g_i \inv (X_{j,i})$, then it holds for $V$ with $X_i,Y_i$:  
 assuming $f_j: \tV_j \to V_j$ is as in the conclusion of the lemma, let $f$ be the disjoint
union of the $f_j$.  In this way we may assume that $V$ is irreducible.  Clearly we may assume $V$ is  complete.
Finally, we may assume $V$ is normal, by lifting the $X_i$ to the normalization  $V_{n}$ of $V$,
and replacing $Y_i$ by $Y_i \times_{g_i} V_{n}$.    We thus assume
$V$ is irreducible, normal and complete.  We may also assume the varieties $Y_i$ and $X_i$ to be irreducible.   

 Let $X_1,\ldots,X_{\ell}$ be the varieties of maximal dimension $d$ among the \tcb{locally closed} subvarieties  $X_1,\ldots,X_n$.  We use induction on $d$.  By \lemref{g2}
there exist 
 pseudo-Galois coverings $f_i: \tV_i \to V$ such that each component of $f_i \inv (X_i)$
of dimension $d$
dominates $Y_i$ rationally. Let ${V^*}$ be an irreducible subvariety of the fiber product $\Pi_V \tV_i$
with dominant (hence surjective) projection to each $\tV_i$.  (The function field of 
${V^*}$ is an amalgam of the function fields of the $\tV_i$, finite extensions of the function field of $V$.)  
Define $f: V^* \to V$, $f(x)= f_1(x_1)=\ldots=f_n(x_n)$ for $x=(x_1,\ldots,x_n) \in V^*$.   
\tcb{Take
$a \in X_i$ and $b \in Y_i$,
with $g_i(b)=a$.
Let $F'$ be a   field extension 
such that $a \in X_i(F')$ (hence $b \in Y_i((F')^{\alg})$).}
If $a$ is sufficiently generic in $X_i$,
then there exists $c \in {V^*}((F')^{\alg})$, $c=(c_1,\ldots,c_n)$ with $f(c)=a$.  Since $f_i$ is a pseudo-Galois covering, and $\tV_i$ dominates $Y_i$, \tcr{we have} $b \in F'(c_i)$.   So there exists a  dense
  open subset $W_i \nsubseteq X_i$ such that for any $a$, $b$, $c$ and $F'$ as above, 
  with $a \in W_i$, $f (c)=a$, $g_i(b)=a$, 
 we have $b \in F'(c)$.

\tcb{We may apply the above to the generic point $a$ of $X_i$, with $F' =F (X_i)$. For any point $c \in {V^*}$
with $f(c)=a$, any $b \in Y_i$
with $g_i(b)=a$ may be expressed as a rational function of $c$ with coefficients in $F'$. 
Each of these rational functions extends to a rational morphism $g_{ij} $ defined on some dense affine Zariski open subset $U_{ij}$ of $f_i^{-1} (W_i)$.
After shrinking $W_i$, we may assume that $g_{ij}$ is in fact  a regular morphism $g_{ij}: U_{ij} \to Y_i$ such that,  for  any $a \in W_i$, $b \in Y_i$ and $c \in \tV$
with $g_i(b)=a$ and $f(c)=a$, 
we have  $c \in U_{ij}$ and $b = g_{ij}(c)$ for some $j$}.

 Let $C_{i}$ be the complement of $W_i$ in $X_i$;  so $\dim(C_i)<d$.  
 We now consider  the  family  $\{X'_\nu \}$
 of subvarieties of ${V^*}$ consisting of components of the preimages of 
 the $X_i$   for $i> \ell$ and of the $C_{i}$ for $i\leq \ell$, and the
 $\{Y'_\nu \}$ consisting of the pullback of $Y_i$ to
 $X_i$   for $i> \ell$ and to  $C_{i}$ for $i\leq \ell$.
  By induction, there exists a finite morphism 
$f':\tV' \to {V^*}$   dominating the $Y'_\nu$ in the sense of the lemma.    Let $\tV$ be the normalization of $\tV'$ in 
the normal hull over $F(V)$ of the function field  $F({V^*})$.  
To insure that 
$\tV$ is pseudo-Galois, one may proceed as follows. One replaces
${V^*}$ by its normalization, and one chooses
$\tV'$ to be pseudo-Galois over 
${V^*}$, which is possible by induction.
Then $\tV \to V$ is
pseudo-Galois, and clearly satisfies the conditions of the lemma.   
\eprf

Note that since finite morphisms are projective (cf. \cite{ega} 6.1.11), if $V$ is projective then so is $\tV$.

\begin{lem}\label{g5}  Let $V$ be a
normal 
projective variety and $L$ an ample line bundle on $V$.   Let $H$ be a finite-dimensional vector space, and let $h: V \to H$ be a rational map.   
Then for any sufficiently large integer $m$ there exists sections $s_1,\ldots,s_k$ of $\fL=L^{\tensor m}$
such that there is no common zero of the $s_i$ inside the domain of definition of $h$, and such that for each $i$, 
$s_i \tensor h$ extends to a morphism $V \to  \fL   \tensor H$.     \end{lem}

\prf  Say $H=\Aa^n$.  We have $h=(h_1,\ldots,h_n)$.  Let $D_i$ be the polar divisor of $h_i$ and $D=\sum_{i=1}^n D_i$.   Let $L_D$ be the associated line bundle.  
Then  $h \tensor 1$ extends to a section of $H \tensor L_D$.  
 Since $L$ is ample, for some $m$, 
$L^{\tensor m} \tensor L_D \inv$ is generated by global sections $\si_1,\ldots,\si_k$.
Since $1$ is a global section of $L_D$, $s_i = 1 \tensor \si_i$ is a section of
 $L_D \tensor (L^{\tensor m} \tensor L_D \inv) \cong  L^{\tensor m}$.   Since away from the support of the divisor $D$, the common zeroes of the $s_i$ are also common zeroes of the $\si_i$, 
 they have  no common zeroes there.   Now 
 $h \tensor s_i = (h \tensor 1) \tensor (1 \tensor s_i)$  extends to  a section of 
 $(H \tensor L_D) \tensor (L_D \inv \tensor L^{\tensor m}) \cong H \tensor L^{\tensor m}$.
\eprf

A theory of fields is called  an {\em algebraically bounded theory}, cf. \cite{winkler}  or \cite{vddries}, if for any \index{algebraically bounded}
subfield $F$ of a model $M$, $F^{\alg} \meet M$ is model-theoretically algebraically closed in $M$.
By \propref{2.1.3} (4), $\ACVF$ is algebraically bounded.
The following lemma is valid for any algebraically bounded theory.  
 We work over a base field $F=\dcl(F)$.

\begin{lem}  \label{guess}
Let $F$ be a valued field.
Let $V$ be an irreducible normal $F$-variety and \tcb{let $H$ be a finite-dimensional $F$-vector space.} Let $\phi$ be an $\ACVF_{\tcb{F}}$-definable subset of $V \times H$ whose
projection to $V$ has finite fibers.  Then there exists a 
 pseudo-Galois covering 
$\pi: \tV \to V$, a finite family of
Zariski open subsets $U_i \nsubseteq V$, $\tU_i = \pi \inv(U_i)$, and
morphisms $\psi_i : \tU_i \to H$ such that for any 
$\tv \in \tV$, if $(\pi(\tv),h) \in \phi$
then $\tv \in \tU_i$ and $h = \psi_i(\tv)$ for some $i$.   
\end{lem}

\prf  For $a$ in $V$ write $\phi(a) = \{b: (a,b) \in \phi\}$; this is a finite subset of $H$.  
Let $p$ be 
 an $\ACVF$-type  over $F$ (located on $V$) and 
  $a \models p$.   By the algebraic boundedness of $\ACVF$, 
  $\phi(a)$ is contained in a finite normal  field extension $F(a')$ of $F(a)$.
Let $q=\tp_{\ACF} (a'/F)$, and let $h_p: q \to V$ be a rational map with $h_p (a') = a$.

 We can
also write each element $c$ of $\phi(a)$ as $c=\psi(a')$ for some rational function $\psi$ over $F$.   This gives 
a finite family $\Psi=\Psi(p)$ of rational functions $\psi$; enlarging it, we may take it to 
 be Galois invariant.  For any $c' \models q$ with $h_p (c')=a$, we have $\phi(a) \nsubseteq \Psi(c'):= \{\psi(c'): \psi \in \Psi \}$.

The type  $q$ can be viewed
as a type of elements 
of an algebraic variety ${W}$, and after shrinking ${W}$ we can take $h_p$ to be 
a   quasi-finite morphism  on ${W}$, and assume    each $\psi \in \Psi: {W} \to H$  is defined on ${W}$;  moreover we can find ${W}$ such that:

$(\ast)$  for any $c' \in {W}$ with $h_{\tcb{p}} (c')=a \models p$, we have $\phi(a) \nsubseteq \Psi(c')$.

By compactness, there exist finitely many triples  \tcb{$({W}_i, \Psi_i, h_i)$} such that for any $p$, some triple has $(\ast)$ 
  for $p$. \tcb{We may now use  \lemref{g3} to conclude.
  Indeed, let $Y \subset V \times H$ be the set of points $(x,y)$ such that for some $w \in W_i$, $x = h_i (w)$ and
  $y = \psi  (w)$ for some $\psi \in  \Psi_i$.
  Let $X_i \subset V$ be the image of $Y_i$ under the projection to $V$.
  We may assume $X_i$ and $Y_i$ are locally closed subvarieties and we denote by 
  $g_i : Y_i \to X_i$ the morphism induced by the projection to $V$.
Applying  \lemref{g3},  we obtain  a  pseudo-Galois covering
  $f: \tV \to V$, 
a finite number of  locally closed subsets $U_{ij}$ of \tcb{$f \inv (X_i)$} and 
 morphisms $g_{ij}: U_{ij} \to Y_i$ satisfying the conditions of \lemref{g3}.
  There is no harm in assuming that
  each $U_{ij}$ is closed in some affine nonempty open $\widetilde \Omega_{ij} = \pi \inv (\Omega_{ij}) \subset \tV$,
  with $\Omega_{ij}$ Zariski open in $V$. Let $\phi_{ij} : U_{ij} \to H$ be the morphism obtained by composing $g_{ij}$ with the projection to $H$.
  We may extend $\phi_{ij}$ to a morphism $\psi_{ij} : \widetilde \Omega_{ij} \to H$.
  Now the  pseudo-Galois covering $\tV \to V$ together with the family of open subsets
  $\Omega_{ij}$ and morphisms $\psi_{ij}$ does the job.}
    \eprf

If $H$ is a vector space, or a vector bundle over $V$, let $H^n$ be the $n$-th direct power of $H$, 
and let $P(H^n)$ denote the projectivization of $H^n$.
Let $h \mapsto :h: $ denote the natural map $H \m \{0\} \to PH$.
Let $r_k: P(H^n) \to PH$ be the natural rational map, $r_k(h_1: \ldots: h_n) = (:h_k:)$. 
For any vector bundle $L$ over $V$, there is a canonical isomorphism 
$L \tensor H^n \cong (L \tensor H)^n$.  When $L$ is a line bundle, we have $P(L \tensor E) \cong P(E)$
canonically for any vector bundle $E$.  Composing, we obtain an   identification  of
 $P((L \tensor H)^n)$ with $P(H^n)$.

\begin{lem}  \label{guess2} Let $F$ be a valued field. 
Let $V$ be a normal irreducible quasi-projective $F$-variety, $H$ a vector space with a basis of $F$-definable points,
and $\phi$ an $\ACVF_F$-definable subset of $V \times (H \m (0))$ whose
projection to $V$ has finite fibers.
 Then there exists  a 
  pseudo-Galois covering $\pi: \tV \to V$,
a   regular morphism  $\theta: \tV \to P(H^m)$ for some $m$, 
such that for any 
$\tv \in \tV$, if $(\pi(\tv),h) \in \phi$ then for some $k$,  $ r_k( \theta(\tv))$ is defined
and equals $:h:$.
\end{lem}

\prf   Replacing $V$ by the normalization of the closure of $V$ in some projective embedding, 
we may assume $V$ is projective and normal.
Let $\psi_i$  be as in \lemref{guess}.  Let $L, s_{ij}$ be as in \lemref{g5},
applied to $\tV$, $\psi_i$;  
choose $m$ that works for all $\psi_i$.  
Let 
$\theta_{ij}$ be the extension to $\tV$ of $s_{ij} \tensor \psi_i$.  Define
$\theta = (\ldots : \theta_{ij}  : \ldots )$,
 using the identification above the lemma.
\eprf

\section{Relatively $\G$-internal subsets}\label{ss6.4}We proceed towards a   relative version of  \propref{G-embed-1}.  
\tcb{First let us clarify some relations of $V$ with $\doublewidehat{V}$, where $V$ is any pro-definable set.  We have
an embedding $s_V: V \to \std{V}$ of $V$ in $\std{V}$ as simple points.    We can thus form
 two natural embeddings $\std{V} \to \doublewidehat{V}$, namely $\std{s_V} $ and $s_{\std{V}}$.  If $a \in \doublewidehat{V}$,
$b \models a \vert A(a)$, and $c \models b \vert A(a,b)$, then $a$ lies in the image of $\std{s_V}$ iff $c \in A(a,b)$, while it is in the image of 
$s_{\std{V}}$ iff $b \in A(a)$.  In other words, the image of $s_{\std{V}}$ consists of the types on 
$\std{V}$ that concentrate on a point of $\std{V}$, while the image of $\std{s_V}$ consists of  the types on $\std{V}$ concentrating on the set of simple points of $\std{V}$. 
Thus the intersection of the two images  is equal to the image of $V$ in $\doublewidehat{V}$, where $v$ is mapped to the type
concentrating on the type concentrating on the single point $v$.  So, away from degenerate cases, when $V=\std{V}$ already,
the two images are distinct and neither contains the other.  It is $\std{s_V}$ that will concern us below.}

Let $\pi: V \to U$ be a morphism of algebraic varieties over a valued field $F$.
We denote by  $\std{V/U}$  the subset of $\std{V}$ consisting of types
$p \in \std{V}$ such that $\std{\pi}(p)$ is a simple point of
$\std{U}$. \tcb{Note that it follows from \lemref{simple} (1) that  $\std{V/U}$ is a relatively definable subset of $\std{V}$.}

We say $X \nsubset \std{V/U}$ is  {\em relatively $\G$-internal} over $U$, if  $X$ is a relatively definable subset of $\std{V}$, \index{relatively $\G$-internal} 
 and the fibers $X_u$ of $X \to U$
are \tcb{iso-definable and} $\G$-internal, uniformly in $u \in U$.

\begin{lem}\label{simple0-2}Let $\pi : V \to U$ be a morphism of algebraic varieties over a valued field $F$, and let  $X \nsubset \std{V/U}$ be relatively $\G$-internal over $U$.  
Then there exists a natural   embedding $\theta: \std{X} \to \std{V}$ over $\std{U}$, \tcb{determined by:  $\std{s_V} \circ \theta = \std{j}$,
where $j$ is the inclusion map $X \to \std{V}$.}
Over a simple point $u \in \std{U}$, $\theta$ restricts to the identification of  $\std{X_u}$ with $X_u$. 
\end{lem}

\prf   
Let $\pi_X: X \to U$ be the natural map.  Let $p \in \std{X}$; let $A=\acl(A)$ be such that $p$ is $A$-definable; and let $c \models p|A$, $u=\pi_X(c)$.  Since $\tp(c/A(u))$ is \tcb{contained in an  $\acl(A(u))$-iso-definable $\G$-internal set}, 
by \lemref{omin0} (5) there exists an $\acl(A(u))$-definable injective map $j$ with $j(c) \in \G^m$.
But $\acl(A(c)) \meet \G = \G(A)$.  So $j(c)=\a \in \G(A)$, and $c=j \inv(\a) \in \acl(A(u))$.  
Let $v \models c | \acl(A(u))$, and let $\theta(p)$ be the unique stably dominated, $A$-definable type 
extending  $\tp(v/A)$.  So $\theta(p) \in \std{V}$, and $\tcb{\std{\pi}}_X(p) = \pi_{*} \theta(p)$.  
 \eprf

Assume now that $X  \nsubseteq \std{V/U}$ is iso-definable and relatively $\G$-internal.
By \lemref{simple0-2}   we may identify $\std{X}$  with a pro-definable subset of $\std{V}$;
namely the set $\int_U X$ of $p \in \std{V}$ such that if $p$ is  $A$-definable  and $c \models p | A$,
 then $\tp(c/A(\pi(c))) = q|A(\pi(c))$ for some $q \in X$.  
 It is really this set that we have in mind when speaking of $\std{X}$ below.    
 In particular, it inherits a topology from 
 $\std{V}$.

 \begin{thm} \label{relative} 
Let $V \to U$ be a projective morphism of \tcb{quasi-projective} varieties over a valued field $F$.  Let 
 $X \nsubseteq \std{V/U}$ be \tcb{$F$-}iso-definable and relatively $\G$-internal.
 Then  there exists a  finite pseudo-Galois
covering ${U'} \to U$, such that letting 
$X' = U' \times_U X$ 
and  $V'=U' \times_U {V}$, 
there exists \tcb{an $F$-}definable morphism $g: V' \to {U'} \times \G_\infty^N$ over ${U'}$, such that the induced map $g: \std{V'} \to \std{U'} \times \G_\infty^N$ is continuous, and  such that the restriction of $g$ to $\std{X}'$ is injective.  
  In fact Zariski locally each coordinate of $g$ is obtained as a composition of regular maps and the valuation map.  
\end{thm}

\prf  \tcb{After pulling back the data to some $\Pp^n$ we may assume $U$ is irreducible and normal.}
By \propref{G-embed-1}, for each $u \in U$, there exists 
$h \in H_{d, m}(F(u)^{\alg})$ such that $\tau_h$ is injective on the fiber $X_u$
above $u$.  By compactness, a finite number of pairs $(m,d)$ will work for all $u$; by taking a large enough 
$(m,d)$, we may take it to be fixed.  Again by compactness, 
 there exists \tcb{an $F$-}definable $\phi \nsubseteq U \times H_{d, m}$ whose projection to $U$ has finite fibers, such that
 if $(u,h) \in \phi$ then $\tau_h$ is injective on $X_u$.  By \lemref{guess2}, there is a finite pseudo-Galois covering
 $\pi: U' \to U$, and 
 a  regular morphism  $\theta: U' \to P(H'{}_{d, m}^M)$ for some $M$, with 
 $H'_{d, m}$ the vector space generated by $H_{d, m}$,
such that for any 
$u' \in U'$, if $(\pi(u'),h) \in \phi$ then, for some $k$,  $ r_k( \theta(u'))$ is defined
and equals $:h:$.   Note that since $h \in H_{d, m}$,  it follows that $\theta(u') \in PH_{Mm,d}$.  
Let
$g(u', v) = (u', \tau_{\theta(u')}) (v)$.  
Then it is clear that $g$ is continuous and that its restriction  to $X'$ is   injective. 
It follows that its restriction to $\std{X}'$ is injective.  
\eprf


Note that the proposition has content even when the fibers of $X/U$ are finite.
Under certain conditions, the continuous injection of \thmref{relative} can be seen to be a homeomorphism.  
This is clear when $X$ is definably compact, but we will need it in somewhat greater generality.

Let $X$ be a pro-definable subset of $\std{V}  \times \Gamma_{\infty}^{\ell}$, for $V$ an algebraic variety.
If $\rho: X \to \G_\infty$ is a definable continuous  function, we shall say $X$ is {\em compact at $\rho=\infty$} if any definable type $q$ on \index{compact at $\rho=\infty$}
$X$ with $\rho_*q$ unbounded has a limit point in $X$. 
  Compactness at $\rho=\infty$ implies that $\rho \inv(\infty)$ is definably compact.  If $X$ is a subspace of a definably compact space $Y$, $\rho$ extends to a continuous definable function $\rho_Y$ on $Y$, and $\rho_Y \inv(\infty) \nsubset X$, then $X$ is 
  compact at $\rho=\infty$.    In the applications,
this will be the case.  
We say  $X$ 
is {\em $\si$-compact via} a continuous definable function $\xi: X \to \Gamma$, if \index{$\si$-compact via $\xi$}
for any $\g \in \Gamma$, $\{x\in X: \xi(x) \leq \gamma \}$ is definably compact.

More generally, let $\rho, \xi: X \to \G_\infty$ be definable continuous functions.  We
say that $X$ is {\em $\si$-compact   via} $(\rho,\xi)$ if $\xi \inv(\infty) \nsubseteq \rho \inv(\infty)$,   $X$ is compact at $\rho=\infty$, and $X \m \xi \inv(\infty)$ is $\si$-compact via $\xi$.
\index{$\si$-compact via $(\rho,\xi)$}

Assume $f : V \to U$ is a morphism of algebraic varieties, $\rho : V \to \G_\infty$
and
$\xi : U \to \G_\infty$ are definable v+g-continuous functions. We say that a pro-definable subset 
$X$ of $\std{V}$ is $\si$-compact over $U$ via $(\rho,\xi)$ if 
$X$ is $\si$-compact via $(\rho, \xi \circ f)$, where we omit the \,  $\std{}$  \, on morphisms.

\begin{lem} \label{relative+}   In \textup{\thmref{relative}}, assume $\std{X}$ is $\si$-compact over $U$ via $(\rho,\xi)$, where  $\rho: V \to \G_\infty$ and $\xi: U \to \G_{\infty}$ are definable and v+g-continuous.  Then one can find $g$ as in  \textup{\thmref{relative}} inducing a homeomorphism of $\std{X'}$ with its image in $\std{U'} \times \G_\infty^N$. 
 \end{lem}

\prf  Let $f : V' \to V$ denote the projection and
$f: \std{V'} \to \std{V}$ its extension. After replacing
$g : V'  \to U' \times \G^N_{\infty}$  in the construction of \thmref{relative}
by $V'  \to U' \times \G^{N +1}_{\infty}$ sending $x$ to $(g (x), \rho \circ f)$,
one may assume that $\rho \circ f= \rho' \circ g$ with
$\rho'$ the projection on the last factor;  and $\xi \circ \pi \circ f = \xi' \circ g$, with $\xi'$ the penultimate projection,
and $\pi: V \to U$.
As in 
\thmref{relative} we still denote by $g$ its extension
$\std{V'} \to \std{U'} \times \G^N_{\infty}$.
The restriction 
$g_{ \vert \std{ X'}} $ of $g$ to $\std{ X'}$ is injective and continuous. 
We have to  show that its inverse $g_{ \vert \std{ X'}} \inv$ is continuous too, or 
equivalently that $g_{ \vert \std{ X'}} \inv \circ \phi$ is continuous for any continuous definable $\phi: \std{X'} \to \G_\infty$.  It suffices thus to show that if $W$ is a closed relatively definable subset of $\std{X'}$,  then $g(W)$ is closed.   
By \propref{cc3}, it suffices to show this:
if $p$ is a definable
type on $W$, and $g(w)$ is a limit of $g_*p$ in $\std{U'} \times \G_\infty^N$ for $w \in W$,
then $w$ is the limit of $p$ in $\std{X'}$.
As $g$ is injective and continuous on $\std{X'}$, it suffices to show that $p$ has a limit in $\std{X'}$. 

 Let us first show that if $f_* (p)$ has a limit point in $\std{X}$, then  $p$ has  a limit point in $\std{X'}$.
Since $V' \to V$ is a finite morphism, it is proper, so \tcr{the morphism} $\std{V'} \to \std{V}$ is \tcr{definably} closed by \lemref{closedmap3}.
 It follows that the morphism $f': \std{X'} \to \std{X}$ induced by $f$ is \tcr{definably} closed. Furthermore
it is   surjective since $X' \to X$ is surjective, by \lemref{extend}.
 Let $\alpha$ be the limit of $f_* (p)$. Its fiber $f'{}^{-1} (\alpha)$ is finite and nonempty, say equal
 to $\{\beta_1, \ldots, \beta_n\}$. If $p$ has a limit in $\std{X'}$, by continuity of $f'$, it should be one of the $\beta_i$.
 Hence, if $p$ does not have a limit in $\std{X'}$, there exists  open  relatively definable subsets $O_i$  of
 $\std{X'}$ containing $\beta_i$, such that $O_i \cap O_j = \varnothing$ if $i \not=j$, and such that
 $p$ is on  $Z = \std{X'} \m \cup_{1 \leq i \leq n} O_i$. Since $Z$ is closed, its image
 $f' (Z)$ is closed, hence $\Omega = \std{X}\m f'(Z)$ is open and contains $\alpha$.
 Thus $f_* (p) $ is on $\Omega$. But $f'{}^{-1}  (\Omega) \nsubseteq \cup_{1 \leq i \leq n} O_i$, which contradicts the fact that
 $p$ is on $Z$.  Hence it suffices to show that $f_* (p)$ has a limit point in $\std{X}$.

 Assume first $\rho_* (f_* (p))$ is not bounded. Then $f_* (p)$ has a limit point in $\std{X}$ by
compactness at $\rho=\infty$. 

 Otherwise,   $\rho' $ is bounded on $g_*p$, hence as $\rho'$ is continuous,
$\rho'(g(w))< \infty$.  So $\rho(f (w)) \in \Gamma$.  It follows that 
$\xi' (g(w)) = \xi (\pi (f (w)) \in \Gamma$ also.  Since $g(w)$ is a limit of $g_*p$, the type $(\xi' \circ g)_*p$ concentrates
on a bounded subset of $\G$.  
Hence the type
$f_* (p)$  includes a formula $\xi \circ \pi \leq \alpha$ for some $\alpha \in \Gamma$.
Thus, by $\si$-compactness, $f_*p$ concentrates on a definably compact relatively definable subset of $\std{X'}$,  containing $ f(w)$;
so $f_*p $ has a limit in this set, hence in $\std{X}$.  
   \eprf


\medskip 
 
The following lemma shows that o-minimal covers may be replaced by finite covers carrying the same information, at least as far as homotopy lifting goes.  

Given a morphism $g: U' \to U$  
and  homotopies $h: I \times U \to \std{U}$ and  $h':   I \times U' \to \std{U'}$, we say $h$ and $h'$ are {\em compatible} or that $h'$ {\em lifts} $h$  \index{compatible homotopies} if $\std{g} (h'(t,u'))=h(t,g(u'))$ for all $t \in I$ and $u' \in U'$.
Here, $I$ refers to any closed generalized interval, with final point $e_I$.  Let $H$ be the canonical  homotopy $I \times \std{U} \to \std{U}$ extending $h$, 
cf. \lemref{hbasic}.
 Note that if $h(e_I,U)$ is iso-definable and $\G$-internal, then $h(e_I,U)=H(e_I,\std{U})$.


\begin{thm} \label{omin-finite} Let $\phi: V \to U$ be a  projective morphism of algebraic varieties with  $U$ normal and quasi-projective,  
over a valued  field $F$.  Let  $X \nsubseteq \std{V/U}$ be iso-definable over $F$  and  relatively $\G$-internal over $U$. 
Assume $\std{X}$ is $\si$-compact over $U$ via $(\rho,\xi)$, where  $\rho: V \to \G_\infty$ and $\xi: U \to \G_{\infty}$ are definable and v+g-continuous.
Then there exists a pseudo-Galois covering $U'$ of $U$, and an \tcb{$F$-}definable function $j: U' \times_U X \to  U' \times \G_\infty^m$ over $U'$, 
 inducing a homeomorphism between  $\std{U' \times_U X}$ and its  image in $\std{U'} \times \G_\infty^m$.  Moreover:     \begin{enumerate}

\item  There exists a  finite \tcb{family} of $F$-definable functions $\xi'_i: U \to \G_{\infty}$, such that, for  any compatible pair of \tcb{$F$-}definable homotopies 
$h: I \times U \to \std{U}$ and $h':   I \times U' \to \std{U'}$, if $h$ respects the functions  $\xi'_i$,  then $h$  lifts to \tcb{an $F$-}definable  homotopy $H_X: I \times \std{X} \to \std{X}$.
Furthermore, 
if $h'$ is a deformation retraction
with iso-definable $\G$-internal image $\Sigma'$, and
$h$ is a deformation retraction
with iso-definable $\G$-internal image $\Sigma$, then 
one may impose that $H_X$ is also
a deformation retraction with iso-definable $\G$-internal image 
$\Upsilon = \std{\phi}^{-1}(\Sigma) \cap \std{X}$.

\item   Given a finite number 
of $F$-definable functions $\tcb{{\tilde \xi}_j} : X \to \G_{\infty}$ on $X$, and a finite group action on $X$ over $U$, 
one can choose the functions $\xi'_i : U \to \G_{\infty}$ 
such that  the lift 
$ I \times \std{X} \to \std{X}$
respects the given functions $\tcb{{\tilde \xi}_j}$ and the group action.

\item  If $h'$   satisfies condition $(*)$ of \textup{Definition \ref{star}}, one may also impose that 
$H_X$ satisfies $(*)$. 
\end{enumerate}
 \end{thm}

 \prf  We take $U'$ and $j$ as given by \thmref{relative} and \lemref{relative+} (that is, $j$ is the restriction of $g$).   
First consider the case when $X \nsubseteq U \times \G_\infty^N$.   
There exists a finite number of $\G_{\infty}$-valued $F$-definable functions $\xi''_i$
on $U$  such that the set of values
$\xi''_i (u)$ determine the fiber $X_u = \{x: (u,x) \in X \}$, as well as the functions $\tcb{{\tilde \xi}_j} | X_u$ (with $\tcb{{\tilde \xi}_j}$ as in (2)),
and the group action on $X_u$.  In other words if $\xi''_i(u)=\xi''_i(u')$ for simple points $u,u'$
then $X_u=X_{u'}$, $\tcb{{\tilde \xi}_j}(u,x) = \xi(u',x)$ for $x \in X_u $, and $g(u,x) = (u,x')$ 
iff $g(u',x)=(u',x')$ for $g$ a group element from the group acting in (2).  Clearly  any  homotopy $h: I \times U \to \std{U}$ respecting the functions $\xi''_i$  
lifts to a homotopy $H_X :  I \times \std{X} \to \std{X} \nsubseteq \std{U}\times \G_\infty^N$ given by 
$(t, (u, \gamma)) \mapsto (H (t, u), \gamma)$, 
where $H$ is the canonical  homotopy $I \times \std{U} \to \std{U}$ lifting $h$ provided by \lemref{hbasic}.   Moreover $H_X$ respects 
the functions of (2) and the group action.

This applies to  $X' = U' \times_U X$, via the homeomorphism induced by $j$; so for any pair $(h,h')$ as in (1), 
if $h'$ respects the functions $\xi''_i$, then $h'$ lifts to a definable homotopy $H': I \times \std{X'} \to \std{X'}$, respecting the data of (2), in particular the Galois action on
$X'$. 
As already noted in the proof of \lemref{relative+},
  \tcr{the morphism} $\std{X'}=\int_{U'} X' \to \int_U X=\std{X}$ is \tcr{definably} closed and surjective.   
  Moreover $H'$ respects the fibers of $\std{X'} \to \std{X}$ 
in the sense of \lemref{homotopy-descent}.
Hence by this lemma, $H'$ descends
to a homotopy $H_X:    I \times \std{X} \to \std{X}$.

By \corref{constant-lift}, the condition that $h'$ respects the $\xi''$ can be replaced with the condition that $h$ respects
certain other definable \tcb{$\xi' : U \to \G_{\infty}$}.  
 
Since $X$ is iso-definable uniformly over $U$, \corref{omin3} applies to the image of $H'$; so this image is 
 iso-definable  and $\G$-internal.  The image of $H_{\tcb{X}}$ is obtained by factoring out the action of the Galois group of $U'/U$; by
    \lemref{b3cc},   the image of $H_{\tcb{X}}$ is also iso-definable, and hence $\G$-internal.

The statement regarding condition $(*)$ is verified by construction, using density of simple points and continuity.
\eprf

 \begin{example}    \label{g-embed-example}
 In dimension $>1$ there exist   definable topologies on definable subsets
of $\G^n$, induced from function space topologies, for which \tcb{\thmref{G-embed-1c}} fails.   
For instance, let $X=  \{(s,t): 0 \leq s \leq t \}$.  For $(s,t) \in X$   consider the continuous  function  $f_{s,t}$ on $[0,1]$ supported on $[s,t]$,   with slope $1$ on $(s,s+\frac{s+t}{2})$, and slope $-1$
on $(s+\frac{s+t}{2},t)$.  The topology induced on $X$ 
 from the Tychonoff topology on the space of functions $[0,1] \to \G$ is a definable topology, and definably compact.  
 Any   neighborhood of the function $0$ (even if defined with nonstandard parameters)
 is a finite union of bounded subsets of $\G^2$, but 
  contains a ``line'' of functions $f_{s,s+\varepsilon}$
whose length is at least $1/n$ for some standard $n$, so this topology is not induced from any definable
embedding of $X $ in $\G_\infty^m$.  By \tcb{\thmref{G-embed-1c}}, such topologies do not   occur within $\std{V}$ for an algebraic variety $V$.
 \end{example}





\chapter{Curves}\label{sec7}

{\small \noindent \textbf{Summary.}
In \ref{ss7.1} we prove the iso-definability of $\std{C}$ when $C$ is a curve.
This is done using Riemann-Roch. 
In \ref{ss7.3}  we explain how definable types on $C$ correspond to germs of paths
on $\std{C}$. The remainder of the chapter is devoted to the construction of the retraction
on skeleta for curves. A key result is the finiteness of forward-branching points proved in \propref{br3}.
\par\bigskip}

\section{Definability of $\std{C}$    for a curve $C$}\label{ss7.1}

Recall that a pro-definable set is called iso-definable if it is isomorphic, as a pro-definable set,  to a definable set.
\begin{thm} \label{f9}  Let $C$ be an algebraic curve defined over a valued field $F$.  Then 
$\std{C}$ is an iso-definable set.   
The topology on $\std{C}$ is definably generated, that is, generated by a definable family of
\textup{(}iso\textup{)}-definable subsets. \tcb{In other words,  there is a definable family giving a pre-basis of the topology.}
  \end{thm}

\prf     
One may assume $C$ is a projective curve. 
There exists a finite purely inseparable extension $F'$ of $F$ such that the normalization of $C \otimes F'$ is
smooth over $F'$. Since this does not change the notion of definability over $F$, we may assume $F' = F$.
Hence we may assume $C$ is projective and smooth over $F$, and that it is irreducible. Let $g$ be its genus.
Let $L \tcb{= F (C)} $ be the function field of $C$ 
and let $Y$ be the set of 
elements $f \in L$ with at most $g+1$ poles (counted with multiplicities). 
\medskip

\begin{claim}
Any element of $L^{\times}$ is a product of finitely many elements of $Y$.  
\end{claim}

\begin{proof}[Proof of the claim] We use induction on the number of poles of $f \in L^{\times}$.  If this number is $\leq g + 1$, 
then $f \in Y$.  Otherwise, let $a_1,\ldots,a_H$ be poles of $f$, not necessarily distinct,
and let $b$ be a zero of $f$.  By Riemann-Roch, any divisor of degree $\geq g$ has a nontrivial global section, which provides one a function $f_1$ with poles at most at
$a_1,\ldots,a_{g + 1}$, and a zero at $b$.  Then $f_1 \in Y$, and $f/f_1$ has fewer poles than $f$
(say $f_1$ has $m$ poles; they are all among the poles of $f$; and $f_1$  has at most $m-1$ zeroes other than $b$).
The statement follows by induction.  \end{proof}

Choose an embedding of $i : C \to \Pp^m$ in some
projective space. Thus, for every positive integer $N$,
the line bundle $i^* \Oo (N) $ has degree $N d$ with $d$ the degree of the embedding.
By Riemann-Roch, if $N$ is large enough,
for every line bundle $\mathcal{L}$ on $C$ of degree $\leq  g + 1$,
$i^* \Oo (N) \otimes \mathcal{L}^{-1}$ is generated by its global sections.
Also, for $N$ large enough, the restriction mapping
$H^0 (\Pp^m, \Oo (N)) \to 
H^0 (C,i^* \Oo (N))$
is surjective.
It follows that,
for $N$ large enough, any function on $C$ with at most  $g+1$ poles 
is the quotient of two homogeneous polynomials of degree $N$.  

Fix such an $N$. 
Let $W$ be the set of pairs of homogeneous polynomials of degree $N$.
 We consider the morphism 
$f : C \times W \rightarrow \G_{\infty}$
mapping $(x, \varphi,\psi)$ to $v (\varphi (x))- v(\psi(x))$ or to $0$ if $x$ is a zero of both $\varphi$ and $\psi$.

With notations from the proof of \thmref{prodef},
$f$ induces a mapping
$\std{C} \rightarrow Y_{W, f}$ with
$ Y_{W, f}$ definable.
Now, let us remark that any type $p$ on $C$ induces
a valuation on $L$ in the following way:
let $c \models p$ send
$g$ in $L$ to $v (g (c))$ (or say to the symbol $-\infty$ if $c$ is a pole of $g$),
and that different types give rise to 
different valuations.
It follows that
 the map $\std{C} \to Y_{W, f}$ is  injective, since if two valuations agree on $Y$ they agree
on $L^{\times}$.  This shows that $\std{C}$ is an iso-$\infty$-definable set.   
Since $\std{C}$ is strict pro-definable by
\thmref{prodef} it follows it is iso-definable. 
The statement on the topology is clear.
   \eprf

Let $h: C \to V$ be a relative curve over an algebraic variety $V$, that is, $h$ is flat with fibers of dimension one. Let $\std{C/V}$ be  
the set of $p \in \std{C}$ such that $\std{h}(p)$ is a simple point of $\std{V}$.   Then we have   
the following relative version of \thmref{f9}:

\begin{thm} \label{f11}   Let $h: C \to V$ be a relative curve over an algebraic variety $V$. Then  $\std{C/V}$ is iso-definable.   
  \end{thm}
\prf  
The proof is the obvious relativization of the proof of  \thmref{f9}.   
Indeed, after 
replacing $V$ by a dense open subset we may assume that $h$ is projective, and that there exists a finite purely inseparable
morphism $V' \to V$ such that
the normalization $h' : C' \to V'$ of the pullback of  $C$ to $V'$ is a smooth morphism. 
Thus, one may assume $h : C \to V$ is projective and smooth.
 Furthermore, by Stein factorization,
 $h$ factors as the composition of a morphism $g: C \to U$ with connected fibers and a finite \tcb{surjective} morphism
$U \to V$. Since $\std{C/U}$ may be canonically identified with $\std{C/V}$, one may assume each fiber $C_a$ of $h$
to be
\tcb{connected}. We embed $C$ in $\Pp^m_V$ and note that for $N$ large enough, for any $a \in V$, any function on 
$C_a$ with $\leq g+\tcb{1}$ poles is the quotient of two homogeneous polynomials of degree $N$.  
Let $W_1$ be the 
set of pairs of homogeneous polynomials of degree $N$, $W_2$ be
the set of characteristic functions of points of $V$, and set $W = W_1 \cup W_2$.
Let 
$f : C \times W \rightarrow \G_{\infty}$
mapping $(x, \varphi,\psi)$ to $v (\varphi (x))- v(\psi(x))$ or to $0$ if $x$ is a zero of both $\varphi,\psi$,
for 
$( \varphi,\psi)$ in $W_1$ and
mapping
$(x, \varphi)$ to $v (\varphi (h (x)))$ for $\varphi$ in $W_2$.
 The map   $\std{C} \to Y_{W,f}$ is 
 injective, and we may proceed as in \thmref{f9}. 
\eprf

\begin{remark}\label{dimgeq2}The statement of \thmref{f9} is specific to dimension one.
Indeed, assume we work over a  base valued field of equicharacteristic zero.
By \exref{exotic}, $\std{\Oo^2} (\Qq(t))$ is uncountable, \tcb{when $\Qq(t)$ is endowed with the $t$-adic discrete valuation,}
thus $\std{\Oo^2}$ cannot be iso-definable. By rescaling, it follows that for any nontrivial closed ball $b$,
$\std{b^2}$ is not iso-definable and thus also 
$\std{D}$ for  $D$ a definable subset of  $\Aa^2$ of  dimension two.
By projecting  to $\Aa^2$ and using \lemref{extend}, it follows that for any definable set $X$ in the $\VF$-sort of dimension two,
$\std{X}$ is not iso-definable. Clearly the same holds in any dimension $\geq 2$, over any nontrivially valued field of any residue characteristic (by a similar argument involving, e.g., the construction in Example 13.1 in \cite{hhm} instead of the one in \exref{exotic}).
\end{remark}


\begin{question} \label{f7q}If $f: U \to V$ is a    
 finite morphism of algebraic varieties, is the inverse image of an iso-definable subset of
  $\std{V}$  iso-definable? 
\end{question}

 When the answer is positive,  
the definability of $\std{C}$ follows from that of 
$\std{\Pp^1}$ which is clear by \exref{exp1}.  

 \section{Definable types on curves}\label{ss7.3}

Let $V$ be an algebraic variety \tcb{and  $a, b \in \G_{\infty}$}.
Two pro-definable functions $f,g: [a,b) \to \std{V}$ are said to have the same germ \tcb{at $b$} if $f | [a',b) = g | [a',b)$ for some $a'$.


\begin{prop} \label{c3}  Let $C$ be a   curve, defined over $A$.  There is a canonical bijection between:
\begin{enumerate}\item  $A$-definable types on $C$. 
\item $A$-definable germs at $b$ of   paths $[a,b) \to \std{C}$, up to reparameterization.
\end{enumerate}
Under this bijection, the stably dominated types on $C$ correspond to the germs of constant paths on $\std{C}$.
\end{prop}

\prf  A  constant path, up to reparameterization, is just a point of $\std{C}$.  In this way
the stably dominated types correspond to germs of constant paths into $\std{C}$.  Let $p$
be a definable type on $C$, which is not stably dominated.  Then, by \lemref{defty3}, for some definable
$\d: C \to \G$, $\d_{*}(p)$ is a nonconstant definable type on $\G$.  
Changing sign if
necessary, either $\d_{*}(p)$ is 
the type of very large elements of $\G$, or else for some $b$,  $\d_{*}(p)$ concentrates on elements
in some interval $[a,b]$; in the latter case there is a smallest $b$ such that $p$ concentrates
on $[a,b)$,
 so that it is the type of elements just $<b$, or else dually.  Thus we may 
assume $\d_{*}(p)$ is the generic at $b$ of   an interval $[a,b)$ (where possibly $b = \infty$).

By \thmref{defty2}
there exists a $\delta_{*} (p)$-germ $f$ of  definable function to
$\std{C}$ 
whose integral is $p$.   It is the germ of a definable function $f=f_{p,\d}: [a_0,b) \to \std{C}$;
since $\std{C}$ is definable and the topology is definably generated by \thmref{f9},    for some   
(not necessarily definable) $a$,  the restriction $f=f_{p,\d}: [a,b) \to \std{C}$ is continuous.    The germ of this function $f$ is well-defined.  
A change in the choice of $\d$ corresponds to reparameterization. 
Conversely, given $f: [a,b) \to \std{C}$, we obtain a definable type $p_f$ on $C$; namely
$p_f | E =\tp( e/E)$ if $t$ is generic over $E$ in $[a,b)$, and $e \models f(t) | E(t)$.  
It is
clear that $p_f$ depends only on the germ of $f$. \tcb{Furthermore, with the above notation, $p=p_{f_{p,\d}}$. On the other hand, 
for any $\d$ as above, $f$ and $f_{p_f,\d}$ have the same germ, up to reparameterization. 
Finally,  if the germ of $f$ is $A$-definable, then each $\phi$-definition $d_{p_f} \phi$ is $A$-definable, and so $p_f$ is $A$-definable.}
\eprf

\begin{rem}  
\begin{enumerate}
\item Over a general base set $A$, the germ may not have an $M$-definable representative.    For instance assume $A$ is the canonical code for an open ball of \tcb{valuative radius} $\gamma$ (e.g. $A=\dcl(\beta)$ with $\b$ a transcendental element of the residue field, and $b = \res \inv (\beta)$; in this case $\gamma=0$).  The path
in question takes $t \in (\gamma,\infty)$ to the generic type of a 
 closed sub-ball of $M$, of \tcb{valuative radius} $t$, containing a given point $p_0$.  The germ at $b$ does not depend on $p_0$,
 but there is no definable representative over $A$. 
 \item  Assume $C$ is $M$-definable, and $p$ an $M$-definable type on $C$.  If   $M=\dcl(F)$ for a field $F$,  the germ   in \propref{c3} (2) is represented
by an $M$-definable path. 
\item  The same proof gives a correspondence between invariant types on $C$,
and germs at $b$ of paths to $\std{C}$, up  to reparameterization, where
now $b$ is a Dedekind cut in $\G$.   The analogue of  (2) remains true if $M$ is a maximally complete model.  
\end{enumerate}
\end{rem}

\section{Lifting paths}\label{ss7.4}

Let us start by an easy consequence of Hensel's lemma, valid in all dimensions, but applicable only near simple points.

\begin{lem} \label{p4} Let $f: X \to Y$ be a finite morphism between  smooth varieties, and let 
$x \in X$ be a  closed  point. 
  Assume $f$ is \'etale at $x \in X$.  Then
there exists  neighborhoods $N_x$ of $x$ in $\std{X}$ and  $N_y$ of $y$ in $\std{Y}$ such that $\std{f}: \std{X} \to \std{Y}$
induces a homeomorphism $N_x \to N_y$.
 \end{lem}

\prf  By Hensel's lemma, there exist valuative neighborhoods $V_x$ of $x$ and $V_y$ of $y$ such
that $f$ restricts to a bijection $V_x \to V_y$.   We take $V_x$ and $V_y$ to be defined by weak inequalities;
let $U_x$ and $U_y$ be defined by the corresponding strict inequalities.  
Then $f$ induces a continuous bijection $\std{V_x} \to \std{V_y}$ which is a homeomorphism by definable compactness. 
In particular, $f$ induces a homeomorphism $N_x \to N_y$, where
  $N_x=\std{U_x}$ and $N_y=\std{U_y}$.   
  \eprf
   
 In fact this gives a   notion of a small closed ball on a curve, in the following sense:   
 
\begin{lem} \label{smallballs}  Let $F$ be a valued field, $C$ be a smooth curve over $F$, and let $a \in C(F)$ be a point.  Then there exists an $\ACVF_F$-definable decreasing family $b(\gamma)$ of g-closed, v-clopen definable subsets of $C$, with intersection $\{a\}$.  Any two such families agree eventually up to
reparameterization, in the sense that if $b'$ is another such family then for some $\gamma_0, \gamma_1 \in \G$ and $\a \in \Qq_{>0}$, for all $\g \geq \g_1$ we have $b(\gamma) = b'(\a \gamma + \gamma_0)$.
\end{lem}

\prf Choose $f: C \to \Pp^1$, \'etale at $a$.  Then $f$ is injective on some v-neighborhood $U$ of $a$. 
We may assume $f (a) = 0$. Let $b_\g$ be the closed ball of radius $\g$ on $\Aa^1$ centered at $0$.   For  some $\g_1$, for $\g \geq \g_1$ we have $b_\g \nsubseteq f(U)$ since $f(U)$ is v-open.
 Let $b(\gamma) = f \inv (b_\g) \meet U$.  
 Note that $A = \{(x, y) \in C \times b_\g : f (x) = y\}$ is a 
 v+g-closed and bounded subset of $C \times \Pp^1$.
 It follows from \propref{closed=closed}, \thmref{compactiffboundedclosed} and \lemref{closedmap2} that $b(\gamma)$
 is g-closed. 
 Since $f$ is a local v-homeomorphism it is v-clopen.

   Now suppose $b'(\gamma)$
is another such family.  Let $b'_\g = f (b'(\gamma))$.  Then by the same reasoning $b'_\g$
is a v-clopen, g-closed definable subset of $\Aa^1$, with  $\meet_{\g \geq \g_2} b'_\g = \{0\}$. 
Each $b'_\g$ (for large $\g$) is a finite union $\union_{i=1} ^m c_i(\g) \m d_i(\g)$, where
$c_i(\g)$ is a closed ball and 
$d_i(\g)$ is a finite union of open sub-balls of $c_i(\g)$, whose number is uniformly bounded,
cf. Holly Theorem, Theorem 2.1.2 of  \cite{hhmcrelle}.  From \cite{hhmcrelle}
it is known
that there exists an $F$-definable finite set $S$, meeting each $c_i(\g)$ (for large $\g$) in one point $a_i$.  The valuative radius of $c_i(\g)$ must approach $\infty$, otherwise it has
some fixed radius $\g_i$ for large $\g$, forcing the balls in $d_i(\g)$ to have
eventually fixed radius and contradicting $\meet _\g b'_\g = \{0\}$.  So, for every $i$ and  large $\g$, $c_i(\g)$
are disjoint closed balls centered at $a_i$.  
   It follows that  $c_i(\g') \m d_i(\g') \nsubseteq c_i(\g) \m d_i(\g)$ for  $\g \ll \g'$.  We have $a_i \notin d_i(\g)$, or else for large $\g'$  we would have $c_i(\g') \nsubseteq d_i(\g)$. Thus
   $a_i \in \meet_\g c_i(\g) \m d_i(\g)$ and  $a_i=0$, hence $m = 1$.

    Now the balls of $d_1(\g)$ must also
   be centered in a point of $S'$ for some finite set $S'$, and for large $\g$ we have $c_1(\g)$
   disjoint from these balls; so $\tcb{b'_{\g}} = c_1(\g)$ is a closed ball around $0$.  For large $\g$ it must have valuative
   radius $\a \g + \g_0$, for some $\a \in \Qq_{>0}, \g_0 \in \G$.  \eprf

\begin{defn}  \label{p3} A continuous  map $f: X \to Y$  between topological spaces with  finite fibers is {\em topologically \'etale} if \index{topologically \'etale}
the diagonal $\Delta_X$ is open in
$X \times_Y X$.  
\end{defn}


\begin{lem} \label{unramified-unique}Let $f: X \to Y$ be a finite morphism between  varieties over a valued field. 
 Let $c: I \to \std{Y}$ be a path, and $x_0 \in \std{X}$.  
If $\std{f}: \std{X} \to \std{Y}$ is topologically \'etale above $ c(I)$, then   $c$ has at most one lift
to a path  $c': I \to \std{X}$, with $c'(i_I)=x_0$.     
\end{lem}

\prf   Let $c'$ and $c''$ be two such lifts.     
  The set $\{t: c'(t)=c''(t) \}$ is definable, 
it contains the initial point, and is closed by continuity.  So it suffices to show that if $c'(a)=c''(a)$
then $c'(a+t)=c''(a+t)$, for sufficiently small $t$, which is
clear by openness of the diagonal.  \eprf
 
\begin{example} \label{p5}  
  In characteristic $p>0$, let  $f: \Aa^1 \to \Aa^1$, $f(x)=x^p - x$.
Let $a \in \Aa^1$ be a closed point, and consider the standard path $c_a: (-\infty,\infty] \to \std{\Aa^1}$,
with $c_a(t)$ the generic of the closed ball of valuative radius $t$ around $a$.  Then
$\std{f}\inv (c_a(t))$ consists of $p$ distinct points for $t > 0$, but of a single point for $t \leq0$.
In this sense $c_a(t)$ \tcb{may be said to be} backwards-branching.
The set of backwards-branching points is the set of balls of valuative radius $0$, which is not a $\G$-internal set.    The   complement of the diagonal within
 $\std{\Aa^1} \times_f \std{\Aa^1}$   
is the union  over $0 \neq \alpha \in \Ff_p$  of 
the sets $U_\a = \{(c_a(t),c_b(t)): a-b=\a, t>0 \}$.   The closure (at $t=0$) intersects
the diagonal in the backwards-branching points.   
\end{example}  
 %


Because of \exref{p5}, we will rely on the classical notion of \'etale only near initial simple points. 
 
\begin{lem}  \label{c0} Let $C$ be an algebraic  curve defined over a valued field $F$ and let  $a$ be a closed point of $C$.  
\begin{enumerate}
\item There exists a   path
$c: [0,\infty] \to \std{C}$ with $c(\infty) = a$, but $c(t) \neq a $ for $t < \infty$.
\item If $a$ is a smooth point,  and
 $c$ and $c'$ are two such paths, then they eventually agree, up to definable  reparameterization.
\item If $a$ is in the valuative closure  of an $F$-definable subset $W$ and $a \notin W$, then for large $t \not= \infty$ one has $c(t) \in \std{W}$.  
\end{enumerate}
\end{lem}

\prf One first reduces to the case where $C$ is smooth. As in the proof of \thmref{f9},
there exists a finite purely inseparable extension $F'$ of $F$ such that the normalization of $C \otimes F'$ is
smooth over $F'$. Since this does not change the notion of definability over $F$, we may assume $F' = F$.
Let $n: \tilde{C} \to C$ be the normalization, and let $\tilde{a} \in \tilde{C}$ be a point such that,
if a $W$ is given as above, then $\tilde{a}$ is a limit point of $n \inv (W)$.  Then the lemma
for $\tilde{C}$ and $\tilde{a}$ implies the same for $C$ and $a$. So, we may assume $C$ is normal. 
For $\Pp^1$ the lemma is clear by inspection.  In general, 
find a morphism $p: C \to \Pp^1$,
with $p(c)=0$ which is  unramified above $0$.  
By \lemref{p4} and its proof, there exists
 a definable homeomorphism for the valuation topology between a 
definable neighborhood $Y$ of $c$ and a definable neighborhood $W'$ of $0$ in $\Pp^1$
which extends to a homeomorphism between $\std{Y}$ and $\std{W'}$.
If $c$ and $c'$ are two paths to $a$ then eventually they fall into $\std{W'}$.   
This reduces to the
case of $\Pp^1$.   For (3) it is enough to notice that one can assume $p (W) \cup \{0\} = W'$. 
(2) comes from \lemref{smallballs}.
\eprf

\begin{rem} \label{c1}   More generally let $p \in \std{C}$, where $C$ is a curve.  If $c \models p$, let $\res(F)(\bar{c})$ be
the set of points of $\St_F$ definable over $F(c)$.  This is the function field of a curve $\bar{C}$
in $\St_F$.   One has a definable  family of paths in $\std{C}$  with initial point $p$, parameterized
by $\bar{C}$.  And any such path eventually agrees with some member of the family, 
 up to definable  reparameterization.  \end{rem}

\section{Branching points}\label{ss7.5}
Let $C$ be a  (noncomplete) curve  over $F$ together with a 
finite morphism of algebraic varieties $f:  C \to \Aa^1$ defined over $F$. 
Given a closed ball $b \nsubseteq \Aa^1$,  let  $p_b \in \std{\Aa^1}$ be the generic type of $b$. 


By an {\em outward path} on $\Aa^1$ we mean a path $c: I \to \std{\Aa^1}$ with $I$ an interval in $\G_{\infty}$ \index{outward path}
such that $c(t) = p_{b(t)}$, with $b(t)$ a ball around some point $c_0$ of valuative radius
$t$.  

Let $X$ be a definable subset of $C$.  
By an {\em outward path} on $(X,f)$ \tcb{with initial point $p$} we mean a germ of path  $c: \tcb{(a,d]} \to \std{X}$, \tcb{with $a < d$,}
such that  $f_{*} \circ c $ is an outward path on $\Aa^1$ and \tcb{$c (d) =p$}.   We first consider the case $X=C$.

In the next lemma, we do not worry about the field of definition of the path; this will be considered later.  
 
\begin{lem} \label{c4}  Let $p \in \std{C}$.  Then $p$ is the initial point of at least one  outward path on $(C,f)$.    
\end{lem}

\prf  The case of simple $p$ was covered in \lemref{c0}, so assume $p$ is not simple.
The point $\std{f}(p)$ is a non-simple element of $\std{\Aa^1}$, i.e. the generic of a closed ball $b_p$, 
of \tcb{valuative radius} $\alpha \neq \infty$.    Fix a model  $F$ of $\ACVF$ over which  $C$, $p$ and $f $ are defined, 
$b_p(F) \neq \varnothing$, and $\alpha = \val(a_0)$ for some $a_0 \in F$.  We will show
the existence of an $F$-definable outward path with initial point $p$.   For this purpose we may
renormalize, and assume $b$ is the unit ball $\Oo$.  
 
 Let $c \models p|F$.  Then $f(c)$ is generic in $\Oo$.  Since $C$ is a curve, 
$\kk (F(c))$ is a function field over $\kk (F)$ of transcendence degree $1$.  Let
$z: \kk (F(c)) \to \kk (F)$ be a place, mapping the image of
$f (c)$ in $\kk (F(c))$ to $\infty$.
We also have a place $Z: F(c) \to \kk (F(c))$ corresponding
to the structural valuation on $F(c)$.  The composition $z \circ Z$ gives a place
$F(c) \to \kk (F)$, yielding a valuation $v'$ on $F(c)$.  Since  $z \circ Z$ agrees with $Z$
on $F$, we can take $v'$ to agree with $\val$ on $F$.  We have 
 $0 < -v'(f (c)) < \val(y)$ for any $y \in F$ with $\val(y)>0$.

Let $q = \tp(c/F; (F(c),v'))$ be the quantifier-free type of $c$ over $F$ in the valued
field $(F(c),v')$.   
 In other words, find an embedding of valued fields $\iota: (F(c),v') \to \Uu$
over $F$, and  let $q=\tp(\iota(c)/F)$.  Similarly, set $r = \tp(f(c)/F; (F(c),v')) := \tp(\iota(f (c))/F)$. 
Clearly $r$ is definable, thus, by  \lemref{findef} it follows that
$q$ is a definable type over $F$, so we can extend it to a global $F$-definable type.  
Note that $q$ comes equipped with a definable map $\delta \to \Gamma$
with $\delta_{*} (q)$ nonconstant, namely $\val(f(c))$.
  According to \propref{c3}, $q$ corresponds to a germ at $0$ of a
path   $c: (-\infty,0) \to C$.  Since for any rational function $g \in F(C)$ regular on $p$, we have $v'(g(c)) = \val (g(c)) \mod \Zz v'(f(c))$,
one may extend
 $c$ by continuity to $(-\infty,0]$ by $c(0)=p$.  
It is easy to check that \tcb{(the germ of)} $c$ is an outward path, since $f_{*} \circ c$ is a standard outward path on $\Aa^1$.  \eprf

We note immediately that the number of germs at $a$ of   paths as given in the lemma 
is finite.  \tcb{Let $p \in \std{C}$.}
\tcb{Fix an outward path $c_0: \tcb{(}-\infty,d] \to  \std{\Aa^1}$, with $c_0(d) = f_{*}(p)$.  
Let $\mathrm{OP}(p)$ be the set of germs of paths 
$c:  (a,d] \to \std{C}$  with $c(d)=p$ and $f_{*} \circ c = c_0$ on $(a,d] $ for some $a<d$.}
If  $c_1,\ldots,c_N \in \mathrm{OP}(p)$ have distinct germs at $d$, then for $d'<d$ sufficiently
close to $d$ the points $c_i(d')$ are distinct; in particular $N \leq \deg(f)$.  

\begin{defn}\label{b-def} A point $p \in \std{C}$ is  called {\em forward-branching} for $f$  if \index{forward-branching point}
 there exists more than one germ of outward paths 
  $c: (a, d] \to \std{C}$
 with $a < d$ and $c (d) = p$, above a given outward path on $\Aa^1$.  We will also say in this
 case that $f_{*}(p)$ is forward-branching for $f$, and even that $b$ is forward-branching
 for $f$ where $  f_{*}(p)$ is the generic type of $b$.
   \end{defn}
  Let $b$ be a closed ball in $\Aa^1$, $p_b$ the generic type of $b$.   Let $  M \models \ACVF$, with $F \leq M$
and $b$ defined over $M$, and 
let $a \models p_b | M$.  Define $n(f,b)$ to be the number
of types
$\{\tp(c / M(a)): \,  f(c)=a \}$.
This is also the number of types
 $\{\tp(c / \acl(F(b))(a)): \,  f(c)=a \}$ (where $M$ is not mentioned), using the stationarity  
 lemma Proposition 3.4.13 of \cite{hhmcrelle}.  Equivalently it is the number of types $q(y,x) $ over $M$
 extending $p_b(x) | M \tcb{\cup \{f (y) = x\}}$; or again:  
\[ n(f,b) = | \{\tp(c/M):  c \in C, f(c)=a \} \vert. \] 
  In other words $n(f,b)$ is the cardinal of the fiber of ${\std{f}}^{-1} (p_b)$,
 with $\std{f}: \std{C} \to \std{\Aa^1}$. In particular, the function $b \mapsto n (f, b)$
 is definable.
 
If $b$ is a closed ball of  valuative radius $\alpha$, and $\lambda > \alpha$, both defined
over $F$,
we define a {\em generic closed sub-ball of $b$ of \tcb{valuative radius} $\lambda$} (over $F$)
to be a ball of \tcb{valuative radius} $\lambda$ around $c$, where $c$ is generic in $b$ over $F$.
Equivalently, $c$ is contained in no proper $\acl(F)$-definable sub-ball of $b$.
 
\begin{lem}   \label{b1}  Assume $b$ and $\lambda$ are in $ \dcl(F)$, and let 
 $b'$ be a generic closed sub-ball of $b$ of \tcb{valuative radius} $\lambda$, over $F$.
  Then $n(f,b') \geq n(f,b)$.  \end{lem}

\prf  Let $F(b) \leq M \models \ACVF$, and $M(b') \leq M' \models \ACVF$.  
Take $a$ generic in $b'$ over $M'$.  Then $a$ is also a generic point of $b$ over $F$.  
Now $n(f,b)$ is the number of types
$\{\tp(c / M): \,  f(c)=a \}$, while $n(f,b')$  is the number of types
$\{\tp(c / M'): \,  f(c)=a \}$.  
As the restriction map sending
 types over $M'$ to  types over $M$ is well-defined and    surjective,
we get  $n(f,b) \leq n(f,b')$.    \eprf  


\begin{lem}\label{fb}The set $\mathrm{FB}'$ of closed balls $b$ such that, for some closed $b' \supsetneqq b$, for all closed
$b''$ with $b \subsetneqq b'' \subsetneqq b'$, we have $n(f,b)  < n(f,b'')$  is a finite definable set,
uniformly with respect to the parameters.
\end{lem}

\prf  The statements about definability of   $\mathrm{FB}'$ are clear since $b \mapsto n (f, b)$
 is definable.
Let us prove that for $\a \in \G$, 
the set $\mathrm{FB}'_\a $ of balls in $\mathrm{FB}'$ of \tcb{valuative radius} $\a$
is finite.
Otherwise, by the Swiss cheese
description of $1$-torsors in Lemma 2.3.3 of \cite{hhmcrelle},
$\mathrm{FB}'$ would
 contain a closed ball $b^*$ of \tcb{valuative radius} $\a' < \a$ such that every sub-ball of $b^*$ of \tcb{valuative radius} $\a$ is in $\mathrm{FB}'$. 
 For each such sub-ball $b'$, for some $\lambda$
with $\a' \leq \lambda < \a$, we have $n(f,b') < n(f,b'')$ 
for any ball $b''$ of \tcb{valuative radius} $\g$ 
with $\lambda < \g < \a$ containing $b'$. Let
$\lambda(b')$ be the infimum of such $\lambda$'s.
Now $\lambda$ is a definable function into $\G$,  so it is constant generically
on $b^*$.  
 Replacing $b^*$ with a slightly smaller ball, we may assume $\lambda$ is actually constant;
so we find $b$ of \tcb{valuative radius} $\lambda$ such that for any sub-ball $b'$ of $b$ of \tcb{valuative radius} $\a$, 
we have $n(f,b') < n(f,b)$.  
But this contradicts \lemref{b1}.

Hence $\mathrm{FB}'$ has only finitely many balls of each \tcb{valuative radius}, so it can be viewed as a       function
from a finite cover of $\G$ into the set of closed balls.  Suppose $\mathrm{FB}'$ is infinite.  Then   it must contain
all closed balls   of \tcb{valuative radius} $\g$  containing a certain point $c_0 \in \tcb{K}$, for $\g$
in some proper interval $\a < \g < \a'$ (again by Lemma 2.3.3 of \cite{hhmcrelle}).  But then by definition of $\mathrm{FB}'$ we find 
$b_1 \nsubset b_2 \nsubset \ldots$ with $n(f,b_1) < n(f,b_2) < \ldots$, a contradiction.  
\eprf

\begin{prop}\label{br3}  The set of forward-branching points for $f$ is finite.  \end{prop}

\prf   By \lemref{fb} it is enough to prove that if $p_b$ is forward-branching, then $b \in \mathrm{FB}'$.
 Let $n=n(f,b) = |\std{f}^{-1} (p_b)|$.  Let $c$ be an outward path on $\std{\Aa^1}$
beginning at $p_b$.  For each $q \in  \std{f}^{-1} (p_b)$ there exists at least
one path starting at $q$ and lifting $c$ by \lemref{c4}, and for some such $q$, there exist more
than one germ of such path.  
So in all   there are   $>n$ distinct germs of paths $c_i$ lifting $c$.  
For $b''$ along $c$ sufficiently close to $b$, the $c_i(b'')$ are distinct; so $n(f,b'') > n$.  \eprf

\begin{prop} \label{c8} Let $f: C \to \Aa^1$ be a finite morphism of curves over a valued field $F$.  Let $x_0 \in C$ be a closed point
where $f$ is unramified, $y_0 = f(x_0)$, and let $c$ be an outward path on $\std{\Aa^1}$,
with $c(\infty)=y_0$.  Let $t_0$ be maximal such that $c(t_0)$ is a forward-\tcb{branching} point
of $f$, or $t_0=-\infty$ if there is no such point. 
Then there exists a unique $F$-definable path $c': [t_0,\infty] \to \std{C}$
with $\std {f} \circ c' = c$, and $c'(\infty)=x_0$.  \end{prop}

\prf  Let us first prove uniqueness. Suppose $c'$ and $c''$ are two such paths.  By \lemref{p4} and
 \lemref{unramified-unique}, 
 $c'(t)=c''(t)$
for sufficiently large $t$.  By continuity, $\{t: c'(t)=c''(t) \}$ is closed.  Let $t_1$
be the smallest $t$ such that $c'(t)=c''(t)$. Then we have two germs of paths lifting
$c$ beginning with $c'(t)$, namely the continuations of $c',c''$.  So $c'(t)$ is a forward-branching point, and hence $t \leq t_0$.  This proves uniqueness on $[t_0,\infty)$. 

Now let us prove existence.  Since we are aiming to show existence of a unique and definable object, we may increase the base field; so we may assume the base field $F$
\tcb{is a maximally complete model of $\ACVF$}.

\begin{claim1}  Let $P \nsubseteq (t_0,\infty]$ be a complete type over $F$,
with $n(f,a)=n$ for $a \in c(P)$.  Then there exist  continuous definable  $c_1,\ldots,c_n:P \to \std{C}$ with $\std{f} \circ c_i = c$, such that $c_i(\a) \neq c_j(\a)$ for $\a \in P$ and $i \neq j \leq n$. 
\end{claim1}
\begin{proof}[Proof of the claim] The proof is similar to that of \propref{c3}, but we repeat it.   Let $\a \in P$.
\tcb{We consider the distinct preimages $\beta_1, \ldots, \beta_n$ of $c(\alpha)$ on $\std{C}$, and for each $\beta_i$ we chose a realization $b_i$ of the corresponding type.}
\tcb{The morphism $f$ is finite, so $\Gamma (F(f(b_i)))$  has finite index in $\Gamma (F(b_i))$. 
Since
$\Gamma (F(f(b_i)))$
 is generated by $\Gamma(F)$ and $\alpha$,
it follows from}  \thmref{maxcomp} that
$\tp(b_i / \acl(F(\a)))$ is stably dominated. 
By \cite{hhmcrelle}, Corollary 3.4.3  
and Theorem 3.4.4,  
 $\acl(F(\a))=\dcl(F(\a))$.  
Thus $\tp(b_i/F(\a)) \in \std{C}$ is $\a$-definable over $F$, and we can write 
$\tp(b_i/F(\a))= c_i(\a)$.
\end{proof}


\begin{claim2}  For each complete type $P \nsubseteq (t_0,\infty]$ over $F$, there exists a \tcb{half-open interval  $ (\a_P,\b_P]$, $\alpha_P,\b_P \in \G_\infty(F)$,   with $P \subset (\a_P,\b_P]$,} and                     
for each $y \in \std{f} \inv( c(\b_P))$, a \textup{(}unique\textup{)} $F(y)$-definable  path 
$c': (\a_P,\b_P\tcb{]} \to \std{C}$ with $\std{f} \circ c' = c$ and $c'(\b_P)=y$.  
\end{claim2}
\begin{proof}[Proof of the claim] For $P=\{\infty\}$ this again follows from  \lemref{p4}. 
\tcb{When $P$ is a realized type different from $\infty$, the statement for $P$ follows from the one for the $F$-type $P^{-}$ of elements infinitely close to $P$ and smaller than $P$.
Thus it  remains to consider the case when $P$ is not realized.}    Then $P$ is an intersection of open intervals defined over $F$.
Say $n(f,a)=n$ for $a \in c(P)$.   By Claim 1 there exist 
  disjoint $c_1,\ldots,c_n$ on $P$ with $\std{f} \circ c_i = c$.
 \tcb{By definability of the space $\std{C}$, and compactness, they may be extended 
  to an open interval   $(\a,\b)$ around $P$   defined over $F$, 
   such that moreover $n(f, c(a))=n$ for $a \in I$, and 
  the  $c_i(a)$ are distinct.  So $\{c_i(a): i=1,\ldots,n\} = 
\std{f} \inv (c(a))$.  Since $\b>t_0$ it is not forward-branching, so we have $n(f,c(\b))=n$ also, and the paths $c_i$ remain distinct at $c(\b)$.}
 The claim follows.   
\end{proof}
Now by compactness of the space of types \tcb{over $F$}, $(t_0,\infty]$ is covered by a finite union of open intervals  where the conclusion 
of Claim 2 holds.  It is now easy to produce $c'$, beginning at $\infty$ and gluing along these intervals.
\eprf 

\begin{rem} 
Here we continue the path till the first time $t$ such that some point of $C$
above $c(t)$ is forward-\tcb{branching}.  
It is possible to continue the path $c'$ a little further, to the first point such that
$c'(t)$ itself is forward-\tcb{branching}.  However in practice, with the continuity with respect
to nearby starting points in mind, we will stop short even of $t_0$, reaching only
the first $t$ such that $c(t)$ {\em contains} a forward-\tcb{branching} ball.  \end{rem}

\section{Construction of a  deformation retraction}\label{ss7.6}

Let $\Pp^1$ endowed with the standard metric of \lemref{metric}, dependent on a choice
of open embedding $\Aa^1 \to \Pp^1$.  Define $\psi: [0,\infty] \times \Pp^1 \to  \std{\Pp^1}$
by letting $\psi(t,a)$ be the generic of the closed ball around $a$ of valuative radius $t$, for this metric.  
By definition of the metric, the homotopy preserves $\std{\Oo}$ (in either of the standard copies of $\Aa^1$). 
We will refer to $\psi$ as the {\em standard homotopy} of $\Pp^1$.  \index{standard homotopy}

\tcb{Note that $\std{\Pp^1}$ has a natural tree structure. Given two points $x$ and $y$ in 
$\std{\Pp^1}$ there exists a unique iso-definable subset $[x, y]$ definably isomorphic to a closed generalized interval with endpoints $x$ and $y$.
If $D$ is a subset of $\std{\Pp^1}$, one defines the convex hull of $D$  as the union of all 
the sets $[x, y]$, for $x,y \in D$.}

Given a Zariski closed subset $D \nsubset \Pp^1$,        
let $\rho(a,D) = \max \{ \rho(a,d): d \in D  \}$.    
Define $\psi_D: [0, \infty]\times \Pp^1 \to \std{\Pp^1}$ by $\psi_D(t,a) = \psi{(\max(t,\rho(a,D)), a)}$.  
\tcb{We call $\psi_D$ the {\em standard homotopy with stopping time defined by} $D$}. In case $D=\Pp^1$ this is  the identity homotopy, $\psi_D(t,a)=a$;
but we will be mostly interested in the case of finite $D$.  \index{standard homotopy with stopping time} \nomenclature{$\psi_D$}{standard homotopy with stopping time defined by $D$}
 In this case $\psi_D$ has \tcb{a $\G$-internal  image, namely the convex hull of $D$.}
(Note:  it is important to use the metric minimum distance, and not schematic distance.  For instance if one uses the latter
for the subscheme on $\Aa^1$ having a double point at $0$, the image would not be $\G$-internal.)

Let $C$ be a projective curve over $F$ together with 
a finite morphism  $f:  C \to \Pp^1$ defined over $F$. 
Working in  the two standard affine charts $A_1$ and $A_2$ of $\Pp^1$,
one may extend the definition of forward-branching points of $f$ to the present setting. 
The set of  forward-branching points of $f$ is  contained in a finite  definable set, uniformly with respect to the parameters.
Factor $f$ as $C \overset{h}{\longrightarrow} C' \overset{f'}{\longrightarrow} \Pp^1$
with $h$ \tcb{finite} radicial and $f'$ generically \'etale. By \corref{radicial},
$\std{h} : \std{C} \to  \std{C'}$ is a homeomorphism. Note that $h$ induces a bijection
between the set of forward-branching points of $f$ and of $f'$.

\begin{thm}\label{ret3}Fix 
a finite $F$-definable subset $G_0$ of $\std{C'}$, including all  forward-branching points of $f'$, 
all singular points of $C'$ and all ramification points of $f'$.
Set $G = \std{f'} (G_0)$ and fix a \tcb{nonempty} divisor $D$ in $\Pp^1$ having a nonempty intersection
with all balls in $G$ \textup{(}i.e. all balls of either affine line in $\Pp^1$, whose generic point lies in $G$\textup{)}. 
\tcb{In other words,  the convex hull of $D$ contains all the aforementioned points.}
Then
$\psi_D: [0, \infty]\times \Pp^1 \to \std{\Pp^1}$
lifts 
uniquely
to
a v+g-continuous $F$-definable function
$ [0, \infty]\times C \to \std{C}$ 
extending to a
 deformation retraction
$H: [0, \infty]\times \std{C} \to \std{C}$ onto an iso-definable $\G$-internal subset of
$\std{C}$. 
\end{thm}

\prf  Since $\std{h} : \std{C} \to  \std{C'}$ is a homeomorphism we may assume $C = C'$ and $f = f'$.
Fix $y \in \Pp^1$. The function $c'_y : [0, \infty] \to \std{\Pp^1}$ sending $t$ to $\psi_D  (t, y)$ is v+g-continuous.
By \propref{c8}, for every $x$ in $C$ there exists a unique  path $c_x : [0, \infty] \to \std{C}$
lifting $c'_{f (x)}$. 
This path remains within the preimage of either copy of $\Aa^1$. 
By \lemref{homotopy-lift} with $X=\Pp^1$,   it follows that
the function
$h: [0, \infty]\times C \to \std{C}$
defined by $(t, x) \mapsto c_x (t)$ is v+g-continuous.
By \lemref{hbasic}, $h$ 
extends to a  deformation retraction
$H: [0, \infty]\times \std{C} \to \std{C}$.
To show that 
$H(0, C)$  is $\G$-internal,
it is enough to check that $\std{f} (H(0, C))$ is $\G$-internal, which is clear.
Uniqueness is clear by \propref{c8}.
\eprf

\begin{example}\label{elliptic}
Assume the residual characteristic of the valued field $F$ is not $2$. 
Fix $\lambda \in F$, $\lambda \not=0$, with $\val (\lambda) > 0$.
Let $C_{\lambda}$ be the projective model of the Legendre curve  $y^2 = x (x-1) (x - \lambda)$
and let $f : C_{\lambda} \to \Pp^1$ be the projection to the $x$ coordinate.
With the notation of  \thmref{ret3}, we may take $D$ to be the divisor
consisting of the four  points $0$, $1$, $\lambda$ and $\infty$. For $x \in F$ with $\val (x) \geq 0$, denote by
 $\eta_{x}$ the generic point of the smallest closed ball containing
$0$ and $x$. Thus, the final image of
$\Pp^1$ under $\psi_D$ is the finite graph $K$ \tcb{that} consists of the union of five segments connecting respectively
$0$ to $\eta_{\lambda}$, $\lambda$ to  $\eta_{\lambda}$,
$1$ to  $\eta_{1}$, $\eta_{\lambda}$ to $\eta_{1}$ and $\infty$ to $\eta_{1}$.
The final image of $H$ is the preimage $K'$ of $K$ under $\std{f}$ which may be described as follows:
over each point of the interior of the segment connecting
$\eta_{\lambda}$ to $\eta_{1}$ there are exactly two points in $K'$ and over all other
points of $K$ there is exactly one (note that $\std{f}^{-1} (\eta_{\lambda})$ is a forward-branching point).
Thus $K'$ retracts on the preimage of the segment connecting
$\eta_{\lambda}$ to $\eta_{1}$ which is combinatorially a circle (see \exref{ellipticber} for the translation of this example in the Berkovich setting).
\end{example}

\begin{example}\label{3lines}
Let $C$ be the union of the three lines $x = 0$, $y = 0$ and $x + y = 1$ in
$\mathbb{A}^2_F$ or its closure in $\mathbb{P}^2_F$.
On each line $L$ consider $\psi_D$ with $D$ the divisor consisting of the intersection points with the two other 
lines. They paste together to produce a retraction of $\std{C}$ to
an iso-definable $\G$-internal subset definably homeomorphic to the subset $\Sigma$
of $\G_\infty^3$ defined as follows.
 Let $Y = \{(\infty, t, 0); 0 \leq t \leq \infty\}$
be the segment connecting
$(\infty, \infty, 0)$ to $(\infty, 0, 0)$ and let the symmetric group $S_3$ act on $\G_\infty^3$ by permuting the coordinates.
Then $\Sigma$ is the hexagon $\cup_{\sigma \in S_3} \sigma (Y)$. One may check, 
similarly as in the example of \remref{necinf}, that $\Sigma$ is not homotopically equivalent to
a definable subset
of some $\G^n$ (or $\G^w$ with finite definable $w$).
In particular, there is no way to retract definably
$\std{C}$ onto an iso-definable $\G$-internal subset definably homeomorphic to a subset of some
$\G^n$ or $\G^w$.
\tcb{Note that this phenomenon detects the singularities of $C$; for instance,
a similar  statement would hold when $C$ is a nodal cubic ($\std{C}$ would retract to a ``circle'' containing the singular point and such a circle is not definably homotopy equivalent to a definable subgraph of 
some $\G^n$).}
\end{example}

  \chapter{\tcb{Strongly stably dominated points}} \label{section:ssd}

{\small \noindent \textbf{Summary.}
\tcb{In
\ref{ss6.5} we study further the properties of strongly stably dominated types over valued fields bases. 
In this setting, strong stability corresponds to a strong form of the Abhyankar property for valuations: the transcendence degrees of the extension and of the residue field extension coincide.
In \ref{ssbertini}} we prove a Bertini type result and also that the 
strongly stable points form a strict ind-definable subset $\stda{V}$  of $\std{V}$. In \ref{ss6.6} we prove a rigidity statement for iso-definable
$\G$-internal subsets of maximal o-minimal dimension of $\std{V}$, namely that 
 they cannot be deformed by any homotopy \tcb{leaving appropriate functions invariant}. This result will be used in \ref{ss10.6}.
\tcb{In \ref{ssclosure}, we study the closure of iso-definable $\G$-internal sets 
 in $\stda{V}$ and we prove that 
 $\stda{V}$ is exactly the union of all skeleta (using \thmref{1}).}
\par\bigskip}

\section{Strongly stably dominated points}\label{ss6.5}\tcb{Recall the notion of being  strongly stably dominated   from \defref{ssd}. This definition makes sense for types of arbitrary imaginaries,
but we will be interested here in the case of types on an algebraic variety.}

\tcb{Let $q$ be a  definable type  on a variety $V$ over  a valued field.  Write $\dim(q)$ for the dimension of the Zariski closure of 
$q$, i.e. of the smallest subvariety of $V$ on which $q$ concentrates.} \nomenclature{$\dim(q)$}{dimension of  the Zariski closure of $q$}

\tcb{We call a definable type {\em sequentially stably dominated} if for all $A=\acl(A)$ with $q$ based on $A$ and $q|A=\tp(c/A)$, there exist $c_1,\ldots,c_n \in A(c)$ \index{sequentially stably dominated type}
with $\tp(c_i/A(c_1,\ldots,c_{i-1}))$  stably dominated, and $ c \in \acl(A(c_1,\ldots,c_n))$.  Here each $c_i$
is a singleton from the field sort.    We will see in \propref{st-st-dom} that this is the same notion, on a variety, as being strongly stably dominated; and that it suffices to check the    property
for {\em some} $A=\acl(A)$ with $q$ based on $A$.}

\tcb{We call a type $\tp(c/A)$ over $A$  strongly stably dominated, respectively sequentially stably dominated,  if it extends to a definable type over $\acl(A)$,
with the corresponding property.    In this case, the definable type is uniquely determined by $\tp(c /{\acl}(A))$.}

\begin{lem} \label{stst-density}
\tcb{Assume $A=\acl(A)$ is generated by  $\VF (A) \union \G (A)$. Let $V$ be an algebraic variety 
defined over $\VF (A)$.
Then the set of sequentially stably dominated types on $V$ over $A$ is dense in the space of types
on $V$ over $A$.
If in addition   $\G(A) \neq (0)$, this 
  remains true if one restricts to   Zariski dense types on  $V$.}
\end{lem}

\prf    \tcb{Let $P$ be the class of sequentially stably dominated types (respectively sequentially stably dominated Zariski dense 
types) over $A$.  
To show $\tp(cd/A)$ is approximated by types of the given class $P$, we may use transitivity.
Consider a formula $\phi(x,y) \in \tp(cd/A)$.  If we know the density for $1$-types, we can find
$d'$ with $\tp(d'/A)$ in $P$, and such that $(\exists x)\phi(x,d')$.  Then we can
find $c'$ with $\phi(c',c')$ and $\tp(c'/\acl(A(d'))) \in P$, and by transitivity (\propref{ssd1} (3)) we have $\tp(c'd'/A) \in P$.}

\tcb{Let $(c_1,\ldots,c_n)$ 
be affine coordinates of $c$ in an appropriate affine embedding.   It suffices to approximate
$\tp(c_i/ \acl (A(c_1,\ldots,c_{i-1})))$ for each $i$; 
so we may assume $c \in \Aa^1$.}

\tcb{Let $D$ be a nonempty  $A$-definable subset  of $\Aa^1$.
By C-minimality, $D$ contains either
a subset $B$ which is an $A$-definable closed ball of finite radius in $\Gamma$
possibly with finitely many proper $A$-definable sub-balls removed 
or an $A$-definable point.
Moreover, if $\G (A) \neq (0)$ and $D$ is Zariski dense, $D$ always contains  such a $B$.
Note that such a definable set $B$ has a canonical definable type, namely the type of elements in this diminished ball avoiding
any proper sub-ball and that this type  yields a (sequentially) stably dominated type over $A$ within $D$. 
}    
\eprf

\begin{prop}\label{st-st-dom}\tcb{Let $q$ be an $A$-definable type  on a variety $V$ over  a valued field.
Let $F$ be a valued field with $A \leq \dcl(F)$.
The following   conditions are equivalent:  
\begin{enumerate} 
\item   $q$ is strongly stably dominated;
\item   over $F$ there exists a locally closed subvariety $W$ of  $V$ with $q \in \std{W}$ and $q$ Zariski dense in $W$, and a quasi-finite morphism   $f: W \to \Aa^n$ of varieties,
such that $f_* q = p_\Oo^{n}$ where $p_\Oo$ is the generic type of $\Oo$;
\item     $\dim(q)=\dim(g_*q )$ for some $F$-definable map $g$  into a variety over the residue field;
\item   $\dim(q)=\dim(h_*q )$ for some $A$-definable map $h$  into a stable sort; here $\dim(h_* q)$ refers
to Morley dimension;
\item there exist  singletons $c_1,\ldots,c_n \in A(c)$
with $\tp(c_i/A(c_1,\ldots,c_{i-1}))$  stably dominated, and $ c \in \acl(A(c_1,\ldots,c_n))$.  
\item  $q|A$ is sequentially stably dominated over $A$.   
  \end{enumerate}}
 \end{prop}

 \prf
 \tcb{(1) implies (2):  Assume first that $F$ is not trivially valued.  
  Let $c \models q|F$.    Then $\tp(c/\St_F(c))$ is isolated.
Now $\St_F = \acl(F \union \kk)$, where $\kk$ is the residue field (\tcb{over the model $\acl(F)$}, a $\kk$-internal set is
contained in $\dcl(\kk)$).  So $\St_F(c) = \dcl(F(c)) \meet \acl(F \union \kk)$.  
Thus $\tp(c/F(d_1,\ldots,d_{n}))$ is isolated for some $d_1,\ldots,d_{n} \in \kk(\acl(F(c))$; by taking conjugates
over $F(c)$ we may assume $d_1,\ldots,d_n \in \kk(F(c))$.
Let  $n$ be minimal, thus $d_1,\ldots,d_{n}$ are algebraically independent over $F$.  We may write $d_i = \res f_i(c)$
where $f_i$ is an $F$-definable function. In fact, upon replacing $d_i$ with $d_i^{p^m}$ for high enough
$m$, if the residue characteristic is $p>0$,  we can take $f_i$ to be a rational function.
So $\tp(c/F(f_1(c),\ldots,f_{n}(c)))$ is isolated.  But $F(f_1(c),\ldots,f_{n}(c))^{\alg} \models \ACVF$;
so $\tp(c/F(f_1(c),\ldots,f_{n}(c)))$ is realized in $F(f_1(c),\ldots,f_{n}(c))^{\alg}$, i.e.
$c \in F(f_1(c),\ldots,f_{n}(c))^{\alg}$.  It follows that $n = \dim(q)$.  Now one may easily find $W$
such that $f|W$ is quasi-finite.}

\tcb{If $F$ is trivially valued then so is $F(c)$, since $\tp(c/F)$ is orthogonal to $\G$;
this case is proved similarly to the above but more easily and is left to the reader.}  

%
%

\tcb{(2) implies (3) is clear; we may take $\Aa^n$ over the residue field, and $g = \res \circ f$.}

\tcb{(3) implies (4) and (1):  It follows from (3) that $q$ is stably dominated  via a function defined over $F$.
Indeed,
the image under a map into the residue field of a definable type $q$ on an $n$-dimensional variety 
is never more than $n$, and if it equals $n$ then the image of any definable map into $\G$ must
be constant.  As definable types orthogonal to $\G$ are stably dominated, $q$ must be stably dominated, and any dominating
function would be algebraic over the given one over $F$, so $q$ is already dominated by that function.
It follows from the Descent Theorem 4.9 in \cite{hhm}
that $q$ is stably dominated  via some $A$-definable function  $h$  into a stable sort. Thus $q|F$ is stably dominated via $h$, and hence $g_*q$ is dominated by $h_* q$.
It follows that $\dim h_* q \geq \dim g_*q = \dim q \geq \dim h_* q$, so equality holds.   This yields (4).  To prove (1),
we may assume $g_*q$ is a Zariski dense type of $\Aa^n$ over the residue field; then $g = \res \circ f$
for some $f$ as in (2).   As $\dim(q) = n$,
if $c \models q | F$, then $c \in \acl_M(f(c))$.  In particular, as $\tp(f(c)/g(c))$ is isolated and implies a type
over $M(g(c))$, $\tp(c / M,g(c))$ is isolated; so $\tp(c,h(c) / M,g(c))$ is isolated,
hence also $\tp(c/ M,g(c),h(c))$.   But $g(c) \in \acl(M,h(c))$.  So $\tp(c/M,h(c))$ is isolated, proving (1).}

\tcb{(4)   implies (5):   We may assume $V$ is affine and use affine coordinates.  Take
$c =(c_1,\ldots,c_m)$ such that $q|A=\tp(c/A)$.
Reordering the coordinates we may assume $c_1,\ldots, c_n$ are algebraically independent over $A$, while
$c \in \acl(A (c_1,\ldots,c_n))$.  So $\dim(q)=n$.  Let $C_i = \St_C(A (c_1,\ldots,c_i))$  and let $d_i$ be the 
Morley transcendence degree of $C_i$ over $C_{i-1}$, i.e. the
supremum of
the Morley rank of $\tp(e/C_{i-1})$, with $e \in C_i$.  Then $\sum_{i=1}^n d_i  = n$.  It follows that $d_i=1$ for each $i$.
Hence (this was seen in the proof of (3) implies (1), as a special case)  $\tp(c_i/A (c_1,\ldots,c_{i- 1}))$ is stably dominated.}
 
 \tcb{(5)   implies (6) is clear, since (5) holds for every base $A$.}

\tcb{(6) implies (1):   By transitivity of strong stable domination, \propref{ssd1} (3), this reduces to the case $\dim(q)=1$.
In this case, taking a maximally complete model $M$ containing $A$, it is clear that (3) holds over $M$.
The implication (3) to (1) was seen above.}  \eprf
%

\begin{example}  If some closed ball $b$ is $A$-definable and $q$ is the generic type of $b$, then $q$ is strongly stably dominated.   Extending $A$ by a realization of $q$ may not add any residue field points, but
it does add a point of a torsor of the residue field, corresponding to $b$.   \end{example}

If $V$ is a definable set, we denote the set of strongly stably dominated types on $V$ by $\stda{V}$.  
 \nomenclature{$\stda{V}$}{set of strongly stably dominated types on $V$}

\begin{lem} \label{abh-char}  Let $U$, $V$ and $W$ be varieties over a valued field, $f: V \to U$ be a definable map.  
\begin{enumerate} 
\item  If $\dim(V)=1$, then $\std{V} = \stda{V}$.
\item  Let $q \in \stda{V}$.  Then $f_* q \in \stda{U}$.
\item  If $f$ \tcb{has finite fibers}, $(\tcb{f}_*) \inv (\stda{U}) = \stda{V}$.
\item If $f$ is surjective,  then $f_* (\stda{V} ) = \stda{U}$.
\item Let $g: V \to  \stda{W}$ be a pro-definable morphism, $p \in \stda{V}$.  Then $\int_p g \in \stda{W}$.
\end{enumerate}
\end{lem} 

\prf (1)  \tcb{By (3) we may assume $V = \Pp^1$ in which case it is clear.}

(2)  follows from \propref{ssd1} (2).

(3) Clear from the characterization of being strongly stably dominated in \propref{st-st-dom}   in terms of dimensions.

(4)   \tcb{Let $p \in \stda{U}$ based on a model $M$ and write
$p|M= \tp (c/M)$. By the density statement in \lemref{stst-density}, there exists $d \in f^{-1} (c)$ such that $\tp (d/ \acl (M(c)))$
is sequentially stably dominated, hence
strongly stably generated by \propref{st-st-dom}.
Thus, by the transitivity property (\propref{ssd1} (3)),
$\tp (d/ M)$ is also strongly stably generated. This yields a definable type $q  \in \stda{V}$ such that $f_* q = p$.}

\tcb{(5) follows from \propref{ssd1} (3).}
%
 \eprf

\begin{remark}\label{nonabhy2}It follows from Example 13.1 in \cite{hhm}, already mentioned in \exref{nonabhy}, that $\stda{(\Aa^2)} \not= \std{\Aa^2}$.
Thus, for any $n \geq 2$,
$\stda{(\Aa^n)} \not= \std{\Aa^n}$.
By rescaling, it follows that for any nontrivial closed ball $b$,
 $\stda{(b^n)} \not= \std{b^n}$. Thus, the same holds for any definable subset of  $\Aa^n$ of  dimension $n$,
 hence, by projecting to $\Aa^n$ and using \lemref{abh-char} (2),
  for any definable set $X$ in the $\VF$-sort of dimension $n$.\end{remark}

\section{A Bertini theorem}\label{ssbertini}
Let $F_0$  be a valued field with infinite residue field and set $F=\acl(F_0)$

Let $p_\Oo$ denote the generic type of $\Oo$.  We will view the tensor power $p_\Oo^{mk}$ as 
 the   generic type of the matrices $M_{m,k}(\Oo)$; thus a {\em generic} matrix over $F$ is
 one realizing $p_\Oo^{mk}$.  

Since $F_0$ is a field with infinite residue field, $p_\Oo$ and hence also the generic type of $M_{m,k}(\Oo)$ are finitely satisfiable in $F_0$.   Thus a definable property that holds for a generic matrix also holds for many
matrices with entries in $F_0$.

Recall that $F(e)$ denotes $\dcl(F \union \{e\})$; this is generally bigger than the field generated by $F$ and $e$.

\begin{prop}\label{bertini}  Let $V$ be an algebraic variety over $F_0$.
Let $c \in V$ such that $ \tp(c/F_0)$ is stationary and  strongly stably dominated. Assume that  $\mathrm{trdeg}_F F(c)=m$.  
Then,  for some  locally closed subvariety  $W$ of $V$ defined over $F_0$ and containing $c$, and some $F_0$-morphism $g :W \to \Aa^{m-1}$, with $\bc = g(c)$, $\bc \models p_\Oo^{m-1}$, and 
$\tp(c/F(\bc))$ is stationary and  strongly stably dominated.  
\end{prop}

\prf  Let $f: W \to \Aa^m$ be as in \propref{st-st-dom} (\tcb{2}).    We will take $g$ of the form $L \circ f$, with $L : \Oo^m \to \Oo^{m-1}$ an $\Oo$-linear function.  In fact, we will show that a   generic  such $L$ will work.
By {\em generic}, we mean a realization of the generic type of $M_{m,k}(\Oo)$.  
  Since $g \inv(\bc)$ is a curve,   $\std{g \inv (\bc)}$ is uniformly iso-definable.  The stationarity 
  statement is equivalent to the existence of a {\em unique} element of  $\std{g \inv (\bc)}$ extending
  $\tp(c/F(\bc))$;  
  it follows easily that the required property holds not only for realizations of $p_\Oo$ but for all sufficiently
close approximations.

Thus it suffices to prove the claim below for $k=m-1$.  For the simple existence statement of $L$,
the claim for any $k$ follows   inductively from the case $k=1$; but we prefer to exhibit the genericity.
 \eprf

\begin{claim}  Let $k<m$.  For a generic  $L : \Oo^m \to \Oo^{k}$, with $\bc = L(f(c))$,  $\bc \models p_\Oo^{k}| F$ 
  and $\tp(c/ F(\bc))$ is stationary.  
  \end{claim}

 \begin{proof}[Proof of the claim] 
For $k=0$ there is nothing to prove.  Assume the claim holds for $k-1$.   
Consider a generic realization $L$ of the generic type of $M_{m,k}(\Oo)$, over $F$.  Let  $\bc = L(f(c))$.  
It is clear that  $\bc \models p_\Oo^{k}| F(L)$ and in particular $\bc \models p_\Oo^{k}| F$.
\tcb{Let us prove that, moreover,} $F(f(c)) \meet \acl(F(\bc)) = F(\bc)$. Indeed, in appropriate coordinates, over $F(L)$, $\bc$ is 
the first $k$ coordinates of a tuple $f(c)$ realizing $p_\Oo^m$; so $F(f(c)) \meet \acl(F(\bc))\subset  F(L,\bc)$;
but $L,\bc$ are independent over $F$, and $L/F$ is stationary, so $F(L,\bc) \meet \acl(F(\bc)) \subset F(\bc)$.
Now suppose $\tp(c/F(\bc))$ is not stationary,
so  \[ F(c) \meet \acl(F(\bc)) \neq  F(\bc). \]
 Let 
$G=\Aut (F(c)^{\alg} / F(f(c)))$ be the (profinite) Galois group.   We have a canonical isomorphism
$\phi: G \to \Aut(F(c)^{\alg} / (F(f(c),\bc)))$ which is the inverse of the restriction map.  
The displayed inequality above implies that  $\Aut(F(c)^{\alg} / \acl(F(\bc)))$ is a proper subgroup of 
$\Aut(F(c)^{\alg} / F(f(c),\bc))$.  Let $H$ be the pullback under $\phi$ of this subgroup.   Let $\hat{J}= \mathrm{Fix}(H)$,
so $\hat{J}$ is a proper algebraic extension of $F(f(c))$, and $\hat{J} \subset \dcl(F(\bc)^{\alg}(f(c)))$.
In fact by Galois theory, there exists an algebraic extension $E$ of $F(\bc)$ such that 
\[ \hat{J}(\bc) = H(f(c)). \]
 
Now let $L$ and $L'$ be mutually generic realizations of the generic type of $M_{m,k}(\Oo)$.  
Let $ \bc'=L'(f(c))$.  If $k \leq m/2$ then $F(\bc)$ and $F(\bc')$ are linearly disjoint over $F(L,L')$ and hence over $F$.   If $m>k>m/2$, they are linearly independent over their intersection, which is generated over $F$
 by a realization of $p_\Oo^{2k-m}$.   (To see this, it is convenient to express $L=L_1 \oplus L_2, L'=L_1' \oplus L_2'$ where $L_1=L_1' \circ L_3$ for some invertible $L_3$ so that they have the same image, and   $L_1,L_2,L_2',L_3$ are generic.)   At any rate, 
 \[ \hat{J}(\bc,\bc') = H(f(c),\bc') = H'(f(c),\bc). \]
   Now $\tp(f(c) / F(\bc,\bc'))$ is strongly stably dominated and stationary.
 It follows that there exist finite extensions $J$ of $F(\bc)$ and $J'$ of $F(\bc')$, with $F(\bc,\bc',J)=F(\bc,\bc'J')$.
 This contradicts the inductive hypothesis.   \end{proof}

%

 \begin{rem}\label{bertini-r}  \leavevmode
\begin{enumerate} \item From the fact that $\tp(c/A(\bc))$ extends to an $A(\bc)$-definable type, it follows that 
\[\acl(A(\bc)) \meet \dcl(A(c)) = \dcl(A(\bc)).\]
 \item The same argument within $\ACF$ shows that for almost all $L$ (outside of a proper Zariski closed subset of $M_{n,k}$),     we have $\acl(A(\bc)) \meet \dcl(A(c)) = \dcl(A(\bc))$ {\em in the sense of $\ACF$}.   Hence this can be required at the same time,
     i.e. we can require $W_{\bc}$ is an irreducible curve.\end{enumerate}\end{rem}
     
     We briefly digress to mention a geometric picture for \propref{bertini}, that should be developed elsewhere.   Let $F$ be a valued field, algebraically closed for simplicity.
  Consider  a subset of affine space
of the form $A=\{x: \val (f_i(x)) \geq 0, i \in I \}$, where $(f_i)_{ i \in I}$ is a set of polynomials over  $F$.  These are $\infty$-definable
sets in $\ACVF_{F}$ that we will call {\em polynomially convex}. \index{polynomially convex set}
   If $W$ is the Zariski closure
of $A$, we   prefer to write $A= \{x \in W: \val (f_i(x)) \geq 0, i \in I \}$.  
Any $p \in \std{\Aa^n}$ has an associated
polynomially convex set $A(p)$, where $f_i$ is the set of polynomials over $F$ such that $p_* (\tcb{\val} (f_i)) \geq 0$; call polynomially convex sets arising in this way {\em irreducible.} \index{irreducible polynomially convex set}

The generically stable type can be recovered from $A(p)$, via   $p_*(\tcb{\val}(f)) = \inf_{a \in A(p)} \{\val (f(a))\}$.
If $p$ is strongly stably dominated, call $A(p)$ a {\em strictly algebraic irreducible affinoid}. \index{strictly algebraic irreducible affinoid}
Note that $(f_i)_{i \in I}$ may be taken to have finitely many polynomials of any given degree (generators of the appropriate lattice).

It probably follows from results in \cite{mst} that if  one can take $I$ to be finite, then 
$A$ is a  strictly algebraic irreducible affinoid.
  The (close) relation between these two notions should be clarified.  
\medskip

In this language, the proof of \propref{bertini} can be adapted to show:

\begin{prop}\label{bertinip}  A strictly algebraic irreducible affinoid of
dimension $>2$  admits strictly algebraic irreducible hyperplane sections. \end{prop}

\begin{remark}  It may be possible to approximate any affinoid (possibly including   analytic affinoids in the Berkovich setting) by a strictly algebraic one, leading to a more general Bertini theorem.   Strict irreducibility is roughly the same as having a Shilov boundary consisting of
a single element.
\end{remark}
     \medskip

 \propref{bertini} will  allow us to think of a strongly stably dominated type of dimension $n$ as the integral over $p_\Oo^{n-1}$ of a definable function into $\std{V/\Aa^{n-1}}$, where $\dim(V)=n$.

\begin{prop}\label{stda-char-2} Let $V$ be an algebraic  variety over a valued field and let $q \in  \stda{V}$ such that  $\dim(q)=m$.  
Then there exists a Zariski open subvariety $W$ of the  Zariski closure of $q$,   a morphism $ W \to \Aa^{m-1}$ making $W$ a relative curve over 
an open subset of $\Aa^{m-1}$, and a definable map $j: \Oo^{m-1} \to \std{W/ \Aa^{m-1}}$,
such that $q = \int_{{p_\Oo}^{m-1}} j$.  Conversely for any such $W$ and $j$, 
$ \int_{{p_\Oo}^{m-1}} j$ lies in  $\stda{V}$. 
\end{prop}

\prf  Let $A$ a base for $q$, $c \models q|A$, and let notation  ($W,m,g,\bc,q_\bc$) be  as in \propref{bertini}.  By \remref{bertini-r} the generic fiber of $g$ can be taken to be an irreducible curve.  Restricting to a Zariski open subset of $W$, we can arrange that $g: W \to U \nsubset \Aa^{m-1}$ is a relative curve.
  We view $q_\bc$ as an element of the iso-definable set $ \std{W_\bc}$ 
\tcb{(cf. \thmref{f9})}.  As $q_\bc \in \dcl(A,\bc)$, and $\bc \models p_\Oo^{m-1}$, there exists an $A$-definable
$j: \Oo^{m-1} \to \std{W/ \Aa^{m-1}}$ such that $j(\bc)=q_\bc$.  Now $c \models q_\bc | A(\bc)$; by definition, $ \int_{{p_\Oo}^{m-1}} j  $ is the unique stably dominated type based on $A$ and extending $\tp(c/A)$;   but $q$ has these properties, so $ \int_{{p_\Oo}^{m-1}} j  =q$.  

The converse statement is a special case of \lemref{abh-char} (5). \tcb{It holds for any definable $j: \Oo^{m-1} \to \std{W/ \Aa^{m-1}}$, though the natural case is when $j$ is a section of $\std{W/ \Aa^{m-1}} \to \Aa^{m-1}$.}
\eprf


\tcb{For a binary map $R(x,v)$, we write $R_x$ for the unary map defined by $R_x(v) = R(x,v)$.}

\begin{defn}  \label{uniformly-stably-dominated} 
 \tcb{A {\em uniform parameterization} \index{uniform parameterization}
  is a definable set $X$ with a pro-definable map $p: X \to \std{V}$, along with a definable map
$R$ on $ X \times V$ such that for any $x \in X$, $R_x$  is a definable map $V \to \St_x$,
and $p(x)$ is stably dominated via $R_x$.}

\tcb{If in addition there exist formulas $\phi_{\nu}$, $1 \leq \nu \leq n$, and a definable partition $X=\union_{\nu=1}^n X_i$,
such that $\dim(p(x))$ is constant on $X_\nu$, and for any $1 \leq \nu \leq n$ and $x \in X_\nu$, $p(x)$ is strongly stably dominated via $\phi_{\nu}$ and  $R_{x}$,  
we say that $p$ is  a {\em strong uniform parameteri\-za\-tion}. \index{strong uniform parameterization}
}

\tcb{A {\em uniform ind-parameterization}, \index{uniform ind-parameterization} resp. a {\em strong uniform ind-pa\-ram\-e\-ter\-i\-za\-tion}, \index{strong uniform ind-pa\-ram\-e\-ter\-i\-za\-tion}
is a morphism $p: X \to \std{V}$ with $X$ an ind-definable set, 
along with an ind-definable $R$ on $X \times V$, such that the restriction to any definable $X' \nsubset X$ is a 
uniform parameterization, resp. a strong uniform parameterization.}

\tcb{We say a subset $W$ of $\std{V}$ is {\em uniformly stably dominated}  \index{uniformly stably dominated}
 (resp. {\em strongly uniformly stably dominated}, {\em ind-uniformly stably dominated}, {\em strongly ind-uniformly stably dominated})
 if there exists a uniform parameterization
 (resp. a strong uniform parameterization, a uni\-form ind-parameterization, a strong uniform ind-parameterization)
$p: X \to \std{V}$
with $p(X) =W$.}    \end{defn}\index{strongly uniformly stably dominated}\index{ind-uniformly stably dominated}\index{strongly ind-uniformly stably dominated}


 \tcb{For $p \in \stda{V}$, note that $p$ is stably dominated via $r$
iff for any base model $M$ for $p$ and any $c \models p |M$, $r(c)$ algebraically generates the residue field of $M(c)$ over $\res(M)$;
while if $p$ is strongly stably dominated via $\phi$ and $r$, then $r(c)$  generates the residue field of $M(c)$ as a field over $\res(M)$.}

\begin{lem} \label{usd0} \tcb{Let $\pi: W \to V$ be a relative curve, $X$ a definable set, and let $j: X   \to \std{W/V}$ be a definable map.  Then $j$
 is a strong uniform parameterization.}
 \end{lem}
 
 \prf    \tcb{First suppose $j: X \to \std{\Pp^1}$; then it is easy to see explicitly that $j$ is a strong uniform parameterization.  In the general case,
 for $x \in X$, let $W_x = \pi \inv (\pi (j(x)))$.  After partitioning $X$ into definable pieces, we may assume that for some morphism $h: W \to \Pp^1$, and some fixed $k$, for any $x \in X$, $\std{h} (j(x))$ has exactly $k$ preimages in  $\std{W_x}$.   The lemma follows by a standard
 compactness argument.} 
 \eprf

\tcb{We denote by  $\stda{V}_m$  the set of elements $p \in \stda{V}$ of dimension $\dim(p)=m$.}\nomenclature{$\stda{V}_m$}{points of dimension $m$ in $\stda{V}$}

\begin{lem} \label{usd2} \tcb{Let $p: X \to \std{V}$  be  a uniform parameterization with image contained in $\stda{V}_m$.  Then $p$ is a strong uniform parameterization.}\end{lem}

\prf
\tcb{By compactness it suffices to show that for each $x \in X$, $p(x)$ has a definable neighborhood where the parameterization is strong.
Fix $x \in X$.   By \propref{stda-char-2} there exists a Zariski open subvariety $W_x$ of the  Zariski closure of $p(x)$,   a morphism $ f_x:W_x \to \Aa^{m-1}$ making $W_x$ a relative curve over 
an open subset of $\Aa^{m-1}$, and a definable map $j_x: \Oo^{m-1} \to \std{W_x/ \Aa^{m-1}}$,
such that $p(x) = \int_{{p_\Oo}^{m-1}} j_x$. The fact that $(f_x)_* (p(x)) = {{p_\Oo}^{m-1}}$ is equivalent
to $(f_x)_* (p (x)) \nsubset \Oo^{m-1}$ along with  $(\res \circ f)_*p(x)$ having transcendence degree $m-1$; the latter is equivalent
to $R_x$ having transcendence degree $\leq 1$ over $(\res \circ f)_*p(x)$; so it can be witnessed in a definable neighborhood of $x$.
On the other hand, by \lemref{usd0}, $j_x$ is a strong uniform parameterization over $\Oo^{m-1}$.    Now isolation is transitive,
in a uniform way:  if $\tp(c/Eb)$ is isolated via $\phi(y,b,e)$, and $\tp(b/E)$ is isolated via $\psi(x,e')$, then $\tp(bc/E)$ 
is isolated via $\psi(x,e') \wedge \phi(y,x,e)$, and $\tp(c/E)$ is isolated via $(\exists x)(\psi(x,e') \wedge \phi(y,x,e))$, so that the form of the isolating
formula is fixed.   Putting this together, using
transitivity of isolation, we see that $p$ is a  strong uniform parameterization as well.}
\eprf

\begin{lem} \label{usd1}  \tcb{Let $p: X \to \std{V}$ and $q: Y \to \std{W}$ be   strong uniform parameterizations.  Let  $H: V \to W$ be a definable
map.  Then the set \[\{(x,y) \in X \times Y:  H_* (p(x)) =q(y) \}\] is definable.}
 \end{lem}
 
 \prf  \tcb{Say the data is defined over $C_0$.  As $\{(x,y) \in X \times Y:  H_* (p(x)) =q(y)  \}$ is clearly $\infty$-definable, it suffices to show that it is also ind-definable.    
 We may again work in a definable neighborhood of a given type over the base set;
 in particular we may assume  $p(X) \nsubset \stda{V}_m$ and $q(Y) \nsubset \stda{W}_{m'}$.   
  As $p$ is a strong uniform parameterization, there exists a definable map $R(x,v)$ such that $R_x(v)$ generates
 $\St_{C_0(x)} (v)$ over $C_0(x)$, whenever  $v \models p(x) | C_0(x)$.  Let $R'$ and $\phi$ witness that $q$ is a strong uniform parameterization
 (partition again so that one $\phi$ works).
 Find a formula $\theta(y,z)$ such that for any $y \in Y$,   $\theta(y,w)$ is a formula of    
  Morley dimension $m'$ and multiplicity $1$ in the stable definable type $(R'_y)_* q(y)$.  (Note that Morley dimension and multiplicity
  vary definably in definable families of formulas of $\St$; this reduces to the case of ACF.)
    Then  $H_*(p(x)) =q(y)$ iff for some $C_0$-definable $h$, $h(x,R_x(v)) = R'(y,H(v))$, and
 $H(v) \models y | C_0(y)$ whenever  $v \models x | C_0(x)$.   The latter condition reduces to the following three conditions:
 \begin{enumerate}
 \item   $\phi(H(v), h(x,R_x(v)))$;
 \item  $\theta(y, h(x,R_x(v)))$;
 \item  $R_x(v)$ has Morley rank $\geq m'$ over $y$.
 \end{enumerate}
 The first two conditions are clearly definable, and the third can be ascertained ind-definably using a formula that shows $R_x(v)$ to
 have Morley rank $\leq m-m'$ over $h(x,R_x(v))$.}\eprf

\begin{rem}\label{8.2.9} \tcb{Applying \lemref{usd1} in the case $V=W, H= \mathrm{Id}$, we see that an  ind-uniformly strongly stably dominated set $X \nsubset \stda{V}$
 admits a strict ind-definable structure.  Moreover by the same lemma, the strict ind-definable structure induced
 from any other ind-uniformly strong parameterization is the same.}\end{rem}

\begin{prop}\label{strongisind} \tcb{Let $V$ be an algebraic variety over a valued field.  
Then  $\stda{V}_m$ admits a unique strict ind-definable
structure, so that it becomes ind-uniform\-ly stably dominated.   With this structure, it is in fact
ind-uniform\-ly strongly stably dominated.}     \end{prop}

\prf
\tcb{The set   $\Ss_1$ of subvarieties of $V$ is a strict ind-definable set, already in the theory $\ACF$.
The same is true of the set $\Ss_2$ of pairs $(W,f)$ where $W$ is a locally closed subvariety of $V$ of dimension $m$,
and $f: W \to U \nsubseteq \Aa^{m-1}$ is a morphism to an open subset of $\Aa^{m-1}$, whose fibers are absolutely irreducible curves.  
Let $\Ss_3$ be the set of triples $(W,f,g)$, where $(W,f) \in \Ss_2$, $U=f(W)$,  and $g: U \to \std{W/U}$ is a definable section of $f$ (in $\ACVF$ now).
It is clear that $\Ss_3$ is an ind-definable set (recall that $\std{W/U}$ is iso-definable  by  \thmref{f11}; this is uniform in $(W,f) \in \Ss_2$).
Define a map $h: \Ss_3 \to \std{V}$ by $h(W,f,g) = \int_{{p_\Oo}^{m-1}} g$.  By \propref{stda-char-2}, the image of $h$ is $\stda{V}_m$.  This is clearly an ind-uniform parameterization.   By \lemref{usd2}, it is strong.  By \lemref{usd1}
the kernel of $h$ is definable on definable pieces, and so
a strict ind-definable structure is induced.   Uniqueness similarly follows by comparing to another parameterization, which will
also be strong by \lemref{usd2}, and so isomorphic to the given one by \lemref{usd1} and \remref{8.2.9}.}
\eprf

\tcb{As $\stda{V}$ is the disjoint union of  $\stda{V}_m$ over $m \leq \dim(V)$,  \lemref{unionstrict} endows $\stda{V}$ with a  strict ind-definable structure;  it is the unique such structure such that  the dimension
$\dim(p)$  is an ind-definable   function and $\stda{V}$ is
ind-uniformly strongly stably dominated. From now on this will be the way we shall view $\stda{V}$ as ind-definable.}

 
\section{$\G$-internal sets and strongly stably dominated points} \label{ss6.6}

\tcb{Let $V$ be a variety over a valued field and let $W$ be an iso-definable $\G$-internal subset of $\std{V}$.
By the o-minimal dimension $\dim (W)$ of $W$\index{o-minimal dimension} we mean the dimension of any definable subset of $\G^m$, for some $m \geq 0$,
pro-definably isomorphic to $W$. Note that by \lemref{embeddedgammatypes}, $\dim (W) \leq \dim V$.
If $W'$ is an iso-$\infty$-definable subset of $W$, we set $\dim(W') = \inf \dim(Z)$, where
$Z$ ranges over  all iso-definable $\G$-internal subsets containing  $W'$.  Note that if $\dim(W')=n$ then $W'$ extends
to a complete type of dimension $n$ over any model over which $V$ and $W'$ are defined.}

\tcb{For a point ${x}$ of $\std{V}$, we define $\dim_{x}(W)$ to be the infimum of $\dim (W \meet O)$,
where $O$ runs over  all relatively definable neighborhoods   of ${x}$.} \nomenclature{$\dim_{x}(W)$}{o-minimal dimension at $x$}
 \tcb{Assume that $\dim_{{x}}({W})=n \geq 0$ and that $V$ and $W$ are defined over some base structure $A$.   Then  there  exists a complete type  $q$ over $A$, whose solution set is a subset $W' \nsubseteq {W}$,  such that ${x}$ lies in the closure $cl (W')$ of $W'$
(i.e. every definable neighborhood of ${x}$ intersects   $W'$) and $W'$ has o-minimal dimension $n$.  Indeed, the collection $C_A(x)$ of $A$-definable subsets $W''$ of $W$ such that
 ${x} \notin cl(W'')$ is closed under finite unions.  By assumption,  for $W'' \in C_A({x})$, 
 $\dim(W \m W'') \geq n$.   Hence  $C'_A({x}) = \{W'' \union W''': W'' \in C_A(x), W''' \nsubset W, \dim(W''') < n \}$
 is also  closed under finite unions and 
does not include $W$.  So there exists a type over $A$, on $W$, avoiding each element of  $C'_A({x})$.}

\medskip

We shall say $W$ is of {\em pure dimension $n$} if it has o-minimal dimension $n$ at every point. \index{pure dimension}

\medskip

\tcb{In  \thmref{1} (7) we will prove the existence of skeleta of pure dimension $n$ for varieties of pure dimension $n$.    
 By \thmref{1} (5) \tcb{(or by  \thmref{sgc} (3))}
the skeleton points will be strongly stably dominated.   The following proposition will   permit us to find homotopies fixing such a given skeleton; the idea is roughly that when the skeleton already has
dimension $n$,  there is no room for the homotopy to move things around.}

\begin{prop}  \label{maxdimabh}
\tcb{Let $V$ be a  variety over a valued field 
 and let ${W} \nsubset \std{V}$ be iso-definable and $\G$-internal. Assume $V$ is of dimension $n$.
\begin{enumerate}  \item 
Away from a countable union of iso-definable subsets of dimension $<n$,  all points of ${W}$ are strongly stably dominated \textup{(}see \textup{\thmref{sgc}  (3)} for a stronger statement\textup{)}.
\item    Let  $\phi:  V  \to \G_\infty^r$  be a definable function inducing a finite-to-one map $\tcb{W} \to \G_{\infty}^r$.   Let  $p \in \std{V}$ with $\dim_p W =n $, and let 
 $h: I \times \std{V} \to \std{V}$ be a homotopy respecting   $\phi$.  Then $h$ fixes  $p$.  In particular
 if $W$ has pure dimension $n$, then $h$ fixes pointwise $W$.
\end{enumerate}}
 \end{prop}
 
\prf   
\tcb{(1) For $\a \in {W}$, let $p_\a$ denote the associated stably dominated definable type.    Let $A$ be a countable base model such that 
$V$ and $W$ are defined over $A$, and there exist  $A$-definable functions $\phi_i: V \to \G_\infty$, $1\leq i\leq r$,   
such that  the restriction of $(\phi_1,\ldots,\phi_r) : V \to \G_{\infty}^r$ to $W$ is finite-to-one, cf. \propref{G-embed-1}.}

\begin{claim} \tcb{Let  $W'  \nsubset W$ be the solution set of a type over $A$ with $\dim(W')=n$.  Then  for any $\a \in W'$, $p_\a$ is strongly stably dominated.}
\end{claim}

\begin{proof}[Proof of the claim]
\tcb{Pick $\a \in W'$.    Let  $M$ be a maximally complete model  containing $A$.
 There exists $\a' \models \tp(\a/A)$ with $\tp(\a'/M)$ of o-minimal dimension $n$.
 Without loss of generality (applying an automorphism of the
  universal domain, say) we may assume $\a=\a'$; so $\tp(\a/M)$ has o-minimal dimension $n$.
 Let $c \models p_\a | M(\a)$.  Let $\b$ be a basis for $\G(M(c))$ over $\G(M)$.  So $\b \in M(\a)$.  Also, as $M$ is maximally complete, $\tp(c/M(\b))$ extends to a stably dominated type $r$;
so $r|M(\b)$  generates a complete type over
 $M(\b) \union \G$, and in particular over $M(\a)$.  It follows that $r | M(\a) = p_\a | M(\a)$, so $r=p_\a$.  Thus $\a \in M(\b) \nsubseteq M(c)$.  
 After multiplying each $\a_i$
 by some positive integer,  we can write $\phi_i(\a) = \val (f_i(c))$,  where $f_i$ is a rational function over $M$.
 Reordering if necessary,  we may assume $\phi_1(\a),\ldots,\phi_n(\a)$ are $\Qq$-linearly independent.   Let $\phi=(\phi_1,\ldots,\phi_n)$, $\g_i: = \phi_i(\a), \g=\phi(\a)$.  Let $\psi_i(x) = \rv (f_i(x))$, $\psi=(\psi_1,\ldots,\psi_n)$;  note $\rv(u)$ lies in a $\val(u)$-definable stable sort.  
As $(\g_1,\ldots,\g_n)$ are linearly independent modulo $\G(M)$, the type of $f_1(c),\ldots, f_n(c)$ over $M$ is determined; 
in particular,   
$\rv (f_1(\a)), \ldots, \rv (f_n(\a))$ are algebraically independent over $M(\g)$.  By \propref{st-st-dom} (4), $p_\a$ is  strongly stably dominated;
in fact 
   $p_{\a}$ is dominated by $\psi_* p_\a$, over $A(\g)$.}
 \end{proof}

\tcb{Thus all points of $W$ are strongly stably dominated, apart from ones lying in an $A$-definable
$n-1$-dimensional set. 
As there are only countably many such $A$-definable sets,  this proves (1). 
}
 
 \medskip

\tcb{(2) Let 
 $h: I \times \std{V} \to \std{V}$ be a homotopy respecting the $\phi_i$.  Let $W' \nsubset W$ be the solution
 set of a complete type over $A$, with $\dim(W')=n$, such that $p$ lies in the closure of $W'$.  
 It suffices to prove that the elements of $W'$ are fixed by $h$. 
 Pick $\a \in W'$, let  $M$ be a maximally complete model  containing $A$ and set
 $\g_i = \phi_i (\a)$. As above, we may assume  
$(\g_1,\ldots,\g_n)$ are linearly independent modulo $\G(M)$.
 Let $t \in I$ be non-algebraic over $M(\g_1,\ldots,\g_n)$ and  set $\a' = h_t (\a)$.
Since $h$ respects the levels of the $\phi_i$, 
we have $\phi_{i}(\a') = \g_i$ for each $i$.     Again by the linear independence of  $(\g_1,\ldots,\g_n)$ over $\G(M)$, 
$\rv (f_1(\a')), \ldots, \rv (f_n(\a'))$ are algebraically independent over $M(\g)$.   So $\psi_* p_{\a'} = \psi_* p_\a$  is the generic type of $\RV(\g)= \Pi_i \RV({\g_i})$ (which is the unique
type over $M(\g)$ in $\RV(\g)$).  As above it follows that $p_{\a'}$ is defined over $M(\g)$, and so 
does not depend on $t$.  Thus for non-algebraic $t$, $ h_t(\a) $ takes a constant value; since non-algebraic values of $t$ are
dense, and $h_t$ is continuous, this constant value must be $\a$, and we must have $h_t(\a)=\a$ for all $t \in I$.}  \eprf

\begin{rem}The proof of \propref{maxdimabh} (1) shows also the following.   Let $V$ be a quasi-projective variety \tcb{of dimension $n$} over a valued field and let $\rho: W \to \std{V}$
be pro-definable, continuous, and injective (or finite-to-one) where $W \nsubset \G_\infty^m$ is a definable set of pure dimension $n$.   Then almost   all points of ${\rho(W)}$ are strongly stably dominated. 

On the other hand,
a non-strongly stably dominated point $p$ can always be deformed in at least one direction, at least in the weaker
sense of the existence of a path from $p$ to a strongly stably dominated one.  If $p$
cannot be moved by a homotopy, it must belong to (every) skeleton in the sense of \defref{defskel} and is thus strongly stably dominated after all by \thmref{1}.
 \end{rem}
 
 \section{Topological properties of $\stda{V}$}\label{ssclosure}
\tcb{Let us call a pro-definable subset $X$ of $\std{V}$ {\em pro-$\G$-internal} if the image of $X$ in each definable quotient of $\std{V}$ is $\G$-internal. \index{pro-$\G$-internal set}
For $X \subset \std{\Aa^n}$, this is equivalent to the statement that for each $m$,  letting $H_m$ be the space of polynomials
on $\Aa^n$ of degree $\leq m$ modulo the polynomials that vanish on $X$,  there exists a finite set $F_m$ of bases of $H_m$
such that  for any $p \in X$, the lattice $\Lambda_m(p)$ corresponding to $p$ in $H_m$ is diagonal with respect to one of the bases
in $F_m$.  A $\G$-parameterized pro-definable set, in particular an iso-definable $\G$-internal one,  is pro-$\G$-internal.}


\begin{prop}\label{closed-gamma-strongly2}\tcb{Let $P$ be a pro-$\G$-internal subset of $\std{V}$.  Then the closure of $P \meet \stda{V}$ is contained in $\stda{V}$. 
In particular the closure of a $\G$-parameterized subset of $\stda{V}$ is contained in $\stda{V}$.}  \end{prop}

\prf  \tcb{Let $q \in \std{V}$ and assume every neighborhood of $q$ contains a  point of $P \meet \stda{V}$.  
Say $P$ and $q$ are $N$-definable with $N \models \ACVF$ somewhat saturated.  
Find a net $p_i \in (P \meet \stda{V})(N)$ with $p_i \to q$.  We have to show that   $q \in \stda{V}$.} 

\tcb{We may assume the dimension  of the Zariski closure of $p_i$ is a fixed integer $n_0 \leq \dim(V)$, and that the
stable dimension of all $p_i$, i.e. the maximal dimension of an image in a stable sort, 
 is a fixed number $d_0 \leq n_0$.  By  \propref{st-st-dom} (4) we have in fact $d_0=n_0$, and
 it suffices to show the analogous fact for $q$.   It is enough to prove  that $d(q) \geq d_0$ and $n(q) \leq n_0$.}
  
\tcb{We may assume $V$ is affine, and even $V=\Aa^n$.  Let $H_d$ be the vector space of polynomials in $n$-variables
of degree $\leq d$; so $\union_d H_d$ is a $K$-algebra.  For $p \in \std{V}$, let $J_d(p) = \{h \in H_d:   p_* (\val (h)) \geq 0 \}$.  Then $J_d(p)$ is a lattice
in $H_d/K_d(p)$, where $K_d(p)$ is the maximal $K$-space contained in $J_d(p)$.   Since all definable maps
from $\G$ to varieties are piecewise constant, the set of possibilities for $K_d(p)$ is an $\infty$-definable
set of bounded cardinality, so it is finite; we may assume that $K_d(p_i)$ is constant for large $i$, say equal to $K_d$.
Since $p_i \to q$, we   have $K_d \nsubseteq K_d(q)$.
It follows that $n(q) \leq n_0$.}

\tcb{For some $d$,
we may find  $x_1,\ldots,x_{d(q)}$ in $H_d$  whose $q$-residues are algebraically independent elements of the residue field over $\kk (N)$,
and thus form a transcendence basis.  For $p \in \std{V}$, let $\g_{\nu}(p)= p_* (\val (x_{\nu}))$.  
Thus,
$\g_{\nu} (p_i) \to \g_{\nu} (q) =: \g_\nu$.
Let $r_{\nu}(p)$ the image of $p$ under $\rv (x_{\nu})$
 in $\RV_{\g_{\nu}(p)} = \{u: \val(u) = \g_{\nu}(p) \} / (1+\Mm)$.  Then for any $y \in H_{d'}$, $d' \geq d$,  $y + \Mm J_d(q)$ depends algebraically on
$r_1(q),\ldots, r_{d (q)}(q)$.   For a fixed such  $y$, 
 the dependence is witnessed by a strict valuation inequality  of the form
\[  \val (f(y,x_1\ldots,x_{d(q)})) _* (q)> \min ( \val (c_\mu) + \mu \cdot \g_\nu)\]  
for some polynomial
$f=\sum   c_\mu x^\mu \in N[x_0,x_1,\ldots,x_{d(q)}]$ with coefficients $c_\mu$ in $\Oo$.
Thus, for large enough $i$, one has
\[\val (f(y,x_1\ldots,x_{d(q)}))_* (p_i) > \min ( \val (c_\mu) + \mu \cdot \g_\nu (p_i)).\]  
This  shows that the $r_i(p)$ algebraically span the
image of $p$ in the stable sorts too.  Thus $d_0=d(p_i) \leq d(q)$.} 
\eprf

\tcb{We borrow from  \defref{defskel} the notion of a 
skeleton of $\std{V}$. It is an 
 iso-definable $\G$-internal subset 
$\Upsilon$ of $\std{V}$,
definably 
homeomorphic to a definable subset of $ \G_\infty^w$, for some finite definable set $w$, such that there exists
 a definable deformation retraction 
 $h : I \times \std{V} \to \std{V}$ 
 with image  
  $\Upsilon$, and such that for each irreducible component $W$ of $V$,
  $\Upsilon \cap W$ is of pure dimension $\dim (W)$. In particular skeleta are topologically $\G$-internal.
  Let us conclude our study of $\stda{V}$ with the following theorem.
  Note that   the proof of    (4) relies on \thmref{1}, that we permit ourselves to quote here.}

\begin{thm}  \label{sgc} \tcb{Let $V$ be a quasi-projective variety over a valued field.}
  \begin{enumerate}
\item  \tcb{Let $X \subset \stda{V}$ be  iso-definable and $\G$-internal.    Then the closure of $X$ in $\std{V}$ is contained in $\stda{V}$.}
\item \tcb{In any   iso-definable $\G$-internal subset of $\std{V}$, the strongly stably dominated points form a {\em closed} iso-definable subset.}
\item  \tcb{Let $X$ be an  iso-definable $\G$-internal subset of $\std{V}$ of pure dimension $n=\dim(V)$. Then $X \subset \stda{V}$.}
\item \tcb{The set $\stda{V}$ is exactly the union of all skeleta of $\std{V}$.}
  \end{enumerate}
\end{thm}

\prf  \tcb{Clause (1)  follows directly from \propref{closed-gamma-strongly2}.
Iso-defina\-bili\-ty in (2) follows from  \ref{strongisind} and closedness follows from (1).}

 \tcb{For (3), observe that by (2), if $X$ is an  iso-definable $\G$-internal subset of $\std{V}$, $X \cap \stda{V}$
 is closed. But it follows from \propref{maxdimabh} (1) that if $X$ is of pure dimension $n=\dim(V)$,
then  $X \cap \stda{V}$ is dense in $X$.}

 \tcb{Let us prove (4).  The fact that any skeleton of $\std{V}$ is contained in $\stda{V}$ follows from (3) (note that it is enough to consider the case when $V$ is irreducible).
For the converse, we shall use that $\stda{V}$ has a canonical strict ind-definable structure.
Let $a$ be a point of $\stda{V}$, 
 and fix $M \models \ACVF$, over which $a$ is realized.  
 Since we may require that our retractions are  Zariski generalizing in the sense of
 \thmref{1} (3), there is no harm in assuming that  $a \in \std{W}$ for any irreducible component $W$ of $V$.
Set $n = \dim(V)$.
  We shall prove by descending induction on $k \leq n$ that if $b \in \stda{V}$ is such that  $P = \tp(b/M)$ has o-minimal dimension $k$, then $b$ belongs to some
 skeleton.  For $k=0$ this includes the case of $b=a$.  
 Let $\a : V  \to\G^\ell$ be an $M$-definable
 function  which is injective on $P$ as provided by
 \thmref{G-embed-1c}.
  Let $h :  I \times \std{V} \to \std{V}$ be an $M$-definable deformation retraction 
 as in \thmref{1}, preserving the levels of $\a$.  
When $k = n$, it follows from \propref{maxdimabh} (2) that all points of $P$ are fixed by $h$ hence belong to the image of $h$ which is a  skeleton, and thus $b$ too.
Suppose now $k < n$.  For $c \in P$, let $\tau(c)$ be the maximal point $\tau \in I$ such that $h(c,t)=c$ for $t< \tau$.  If $\tau(c)$ is the final 
 point of $I$ for one (hence for all) $c \in P$, then all $c \in P$ are contained in the final image of $h$
 and so (4) holds.  Otherwise, 
  let $q_c$ be the type over $M(c)$ of elements of $I$ just greater than $\tau(c)$.  
 Consider the set $S=\{h(c,t): c \in P, t \in q_c\}$.  This is a type-definable
subset of some $\G$-internal definable subset of $\stda{V}$.  If $\dim(S) \leq k = \dim(P)$, find an $M$-definable set $S'$ of dimension $k$ containing $S \union P$; for $c \in P$ and $t-\tau(c)$ sufficiently small, $h(c,t)$ lies in $S'$.
Note that
    $P$ is open in $S'$ since $P$ is a complete type, hence
 for $t -\tau(c)>0$   small enough,
$h(c,t)$ must still lie in $P$.  However, as the levels of $\a$ are preserved, $\a (h(c,t))=\a (c)$
so $h(c,t)=c$, contradicting the definition of $\tau$.
Thus $\dim(S)=k+1$.  Clearly
each  point of $P$ lies in the 
closure of $Q$ (consider a path reversing the homotopy).  By induction, any point of $Q$ lies in the closure of an iso-definable $\G$-internal set of dimension
$n$; hence so does each point of $P$. Thus,  it follows from \thmref{1} together with
\propref{maxdimabh} (2), similarly as when $k = n$, that $b$ lies on some skeleton.}
 \eprf

 \begin{rem} \label{ssd3}    \tcb{Let $W$ be an o-minimal subset of $\std{V}$ of pure dimension $n$.  The fact that
 every point of $W$ is strongly stably dominated also follows from     \thmref{1}.  Indeed by \thmref{1} (1) and (5) and by \propref{maxdimabh} (2)
we can find a homotopy fixing $W$ and with strongly stably dominated final image.}     \end{rem}

\begin{rem}\label{localstate}\tcb{Modulo \thmref{1}, \thmref{sgc} (4) is equivalent to a converse to (3) that does not mention 
  retractions, namely that the local o-minimal dimension of $\stda{V}$ is everywhere equal to the 
  local dimension of $V$:   e.g. if $V$ has pure dimension $n$, then every point $p$ of $\stda{V}$ is contained in a   $\G$-internal set of local dimension $n$ at $p$. 
However
 we do not know how to prove this local statement without using \thmref{1}.}
 \end{rem}

  \begin{rem}\tcb{It would be natural to consider $\stda{V}$ with the {\em direct limit} topology, rather than the topology
  induced from $\std{V}$.  We saw that $\stda{V}$  has a canonical ind-definable structure; 
  we topologize each definable subset according to the embedding in $\std{V}$, but then topologize
  $\stda{V}$ as a direct limit.  This is another, and probably better, {\em canonical} topology on $\stda{V}$.
  \thmref{sgc} (1)  
  implies that any ind-o-minimal subset of $\stda{V}$ becomes an ind-o-minimal space, i.e.  a direct limit of o-minimal spaces under a system of closed embeddings.}
\end{rem}

  \chapter{Specializations and $\ACV2F$}\label{specializations}

{\small \noindent \textbf{Summary.} 
We  introduce the theory $\ACV2F$ of   iterated places in \ref{ss8.3}.
It provides us with algebraic  criteria for v- and g-continuity. 
Some applications of the continuity criteria are given in  \ref{ss8.7} and  \ref{ss8.8}.
The result on definability of v- and g-criteria in \ref{ss8.9} will be used in \ref{ss10.7} to handle uniformity with respect to parameters.
Compare to \cite{huber-knebusch}.
\par\bigskip}

\section{g-topology and specialization}   \label{ss8.1}
 Let $F$ be a valued field, and consider pairs $(K,\Delta)$,
with $(K,v_K)$ a valued field extension of $F$, and $\Delta$ a proper convex subgroup of $\G(K)$,
with $\Delta \meet \G(F) = (0)$.  Let   $\pi: \G(K) \to \G(K) / \Delta$ be the quotient 
homomorphism.   We extend $\pi$ to $\G_\infty (K)$ by $\pi(\infty)=\infty$.  
Let 
$\bK$ be the field $K$ with valuation $\pi \circ v_K$.  We will refer to pairs $(K, \bK)$ as  {\em g-pairs} over $F$.  \index{g-pair}

\tcb{The convention of  \ref{ss2.1} shall be in use:   any $\ACVF_F$-definable set  or function   will be assumed to be defined by a quantifier-free formula. This will allow us to evaluate them on g-pairs.
Note that if $F$ has characteristic $(0,p)$, i.e. $0< v(p)< \infty$, then as $p \in F$, $v(p) \notin \Delta$, so $\bK$ has characteristic $(0,p)$ as well.
The  residue field of $\bK$ is thus a valued field of characteristic $(p,p)$, with the same residue field as the one of $K$.}

\begin{lem} \label{gcriterion0}  Let $F$ be a valued field, $V$ an $F$-variety,  and let $U$ and $X$ be $\ACVF_F$-definable subsets with $U \nsubseteq X \nsubseteq V$.  Then the following conditions are equivalent:
\begin{enumerate}
\item $U$ is g-open in $X$;
\item   \tcb{$U$ is  the intersection of  $X$ with a positive Boolean combination 
of Zariski closed and open sets defined over $F$ and sets of the form $\{w \in  W : \val(f(w)) > \val(g (w))\}$, with $f$ ang $g$ regular functions on a Zariski open $W$ in $U$, all defined over $F$;}
\item  for any g-pair
$(K,\bK)$ over $F$, we have $U(\bK) \meet X(K) \nsubseteq U(K)$;
\item   same as \textup{(3)}, with    $K$ \textup{(}as a field\textup{)} of the   form $F(a)$, with $a \in U$.
\end{enumerate}
 \end{lem}

\prf \tcb{One verifies immediately that each of the conditions is satisfied if and only if it holds on every $F$-definable Zariski open subset of $V$.  So we may assume $V$ is affine.}

\tcb{Let us prove that (1) implies (3). Assume $U$ is g-open in $X$, and let $(K,\bK)$ be a g-pair over $F$.  If $a \in X(K)$ and $a \in U(\bK)$, we have to show that $a \in U(K)$.   
 We may pass for this to algebraic closures of $K$ and $\bK$; thus we may assume $K=K^{\alg}$.}   Let $U'$ be g-open, with $a \in U'$ and such that $\ACVF_F \models U' \meet X \nsubseteq U$.
\tcb{As $U'$ is g-open, it is defined by a positive Boolean combination of strict inequalities 
$\val(f)<\val(g)$, with $f$ and $g$ regular functions on $V$ and algebraic equalities and inequalities.   Since $K$ is a model,   all these data can be chosen to be defined over $K$.
 Since $\pi$ is order-preserving on $\G_\infty$, if $\pi \circ v_K(f) < \pi \circ v_K(g)$ then $v_K(f) < v_K(g)$.  The algebraic
equalities and inequalities are preserved since the fields are the same.  Hence
$U'(\bK) \nsubseteq U'(K)$, so $a \in U(K)$.}

\tcb{Since trivially (2) implies (1) and (3) implies (4), it remains to prove  that (4) implies (2). Let $W=X \m U$.  Then $W \nsubseteq \VF^n$ is $\ACVF_F$-definable,  and  for any g-pair
$(K,\bK) $ over $F$,   $X(\bK) \meet  W(K)  \nsubseteq W(\bK)$.  We have to show that  $W$ is cut out of $X$ by a finite disjunction of finite conjunctions of   weak valuation inequalities $\val (f) \leq \val (g)$,   
equalities $f=g$ and inequalities $f \neq g$ involving regular functions defined over $F$.}
It suffices to show
that any complete quantifier-free type $q$ over $F$ extending $W$ implies a finite conjunction  of this form,
which in turn implies $W$.   Let $q'$ be the set of all 
equalities, inequalities and  weak valuation inequalities in $q$, along with the formula defining $X$.  By 
compactness, it suffices to show that $q'$ implies $W$.  Let $a \models q'$, and let $\bK$
be the valued field $F(a)$.   We have $a \in X$, and we are done if $a \in W$; so suppose $a \in U$.
 Let $b \models q$, and let $K=F(b)$.  Since $q'$ is complete
inasfar as $\ACF$ formulas go, $F(a),F(b)$ are $F$-isomorphic, and we may assume $a=b$ and $K$ \tcb{and} $\bK$ coincide as fields.  Any element $c$ of $K$ can be written as $f(a)/g(a)$ for some 
polynomials $f,g$.  Let $c, c' \in K$; say $c=f(a)/g(a), c'=f'(a)/g'(a)$.  
If $v_K(c) \geq v_K(c')$ then $v_K(f(a)g'(a)) \geq v_K(f'(a)g(a))$;   the weak valuation 
inequality $v_K(f(x)g'(x))\geq v_K(f'(x)g(x)) $ is thus in $q$, hence in $q'$, so  
$v_{\bK}(f(a)g'(a)) \geq v_{\bK}(f'(a)g(a))$, and hence $v_{\bK}(c) \geq v_{\bK}(c')$.  It follows
that the map  $v_K(c) \mapsto v_{\bK}(c)$ is well-defined, and weak order-preserving; it
is clearly a group homomorphism $\G(K) \to \G(\bK)$, and is the identity on $\G(F)$.  
By the hypothesis, $W(K) \meet X(\bK) \nsubseteq W(\bK)$.  Since $b \in W(K)$, we have $a \in W(\bK)$.
But $a$ was an arbitrary realization of $q'$, so $q' $ implies $W$.  \eprf 

\begin{rem} 
We can now see that  the family of g-open sets is {\em definable in definable families}.  In other words, if $\{U_a: a \in P\}$
is an $F$-definable family of definable subsets of $V$, and $C$ is the set of elements $a \in P$ with $U_a$ g-open,
then $C$ is a definable subset of $P$.    (We may take $P$ to be affine $k$-space.)   Indeed it is clear from the definition that $C$ is a union of definable sets; so it suffices to show
that if $a \notin C$, then for some formula $\phi \in \tp(a/F)$, any realization of $\phi$ is not in $C$. Recall the theory $\ACV2F$
(cf.   \ref{ss8.3}).  Here we take the sorts to be the valued field sort, 
 and the value group; the latter is enriched with a predicate for a convex subgroup $\Delta \leq \G$.   If $(K,\Delta)$ is the data for a g-pair,
 with $\val: K \to \G$ surjective and $K$ algebraically \tcb{closed}, then $(K,\Delta) \models \ACV2F$.  
 Let $T= \Th(K,\Delta,c)_{c \in F}$.  The complete $T$-type of $a/F$ is then generated by the  $ACVF_F$-diagram $D$ of $a$, along with the  set $S$ of
  sentences:
  $\val(f(a))>0 \to f(a) \notin \Delta$ (for every rational function $f$ over $F$, defined \tcr{at} $a$).  
  Equivalently, we can take $S$ to be the set of sentences:  $\val (f(a)) > \val (g(a)) \to \val (f(a) )- \val (g(a)) \notin \Delta$, with $f,g$ polynomials in $k$ variables over   $F$.  (This makes it clear that $S$ is independent of the type of $a$.)
   By    \lemref{gcriterion0}, as $a \notin C$, 
 $T +D+S \vdash   U_a(\bK) \not \nsubseteq U_a(K)$.  So for some $\ACVF$-formula $\psi \in D$,   already $T + \psi(a) + S  \vdash U_a(\bK) \not \nsubseteq U_a(K)$.  Hence again by the criterion, as soon as $\psi(a',b)$ holds,
 $a' \notin C$.  
 \end{rem}

\begin{lem} \label{gcriterion1}  Let $F_0$ be a valued field, $V$ an $F_0$-variety, and let $W \nsubseteq V$ be $\ACVF_{F_0}$-definable.   Then $W$ is g-closed if and only if   for any $F  \geq F_0$ with $F$ maximally complete and algebraically closed, and any g-pair
$(K,\bK )$ over $F$ such that $\G(K)=\G(F)+\Delta$ with $\Delta$ convex and $\Delta \meet \G(F)=(0)$, 
 we have $W(K) \nsubseteq W(\bK)$.  
 
 When  $V$ is an affine variety,    $W$ is g-closed iff $W \meet E$ is g-closed for every bounded, g-closed, definable subset $E$ of $V$.   
   \end{lem}

\prf  The ``only if'' direction follows from \lemref{gcriterion0}.  For the ``if'' direction, suppose $W$ is not g-closed. 
By \lemref{gcriterion0} there exists a g-pair $(K,\bK)$ over $F_0$ with $W(K) \not \nsubseteq W(\bK)$; 
\tcb{furthermore, one may assume $K$ is finitely generated over $F_0$}, so that $\G(K) \tensor \Qq$ is finitely generated over $\G(F_0) \tensor \Qq$ 
as a $\Qq$-space.  Let $c_1,\ldots,c_k \in K$ be such that $\val(c_1),\ldots,\val(c_k)$ form a $\Qq$-basis 
for $\G(K) \tensor \Qq / (\Delta + \G(F_0)) \tensor \Qq$.  Let $F=F_0(c_1,\ldots,c_k)$.  Then $(K,\bK)$ is  
a g-pair over $F$, $\G(K) = \G(F) + \Delta$, and $W(K) \not \nsubseteq W(\bK)$.    We continue to modify $F$, $K$, and $\bK$.
As above we may replace $F$ by $F^{\alg}$.  Next, let $K'$ be a maximally
complete immediate extension of $K$, $F'$ a maximally complete immediate extension of $F$, and embed
$F'$ in $K'$ over $F$.  Let $\bK '$ be the same field as $K'$, with valuation obtained by composing $\val: K' \to \val (K')= 
\val (K)$
with the quotient map $\val (K) \to \val (K) / \Delta$.  Then $\bK$ embeds in $\bK '$ as a valued field.   We have now the same
situation but with $F$ maximally complete.  
This proves the criterion.

For the statement regarding bounded sets, suppose again that $W$ is not g-closed; let $(K,\bK)$ be a g-pair as above,
$a \in W(K)$,  $a \notin W(\bK)$.   Then $a \in V \nsubset \Aa^n$; say $a=(a_1,\ldots,a_n)$ and let $\g = \max_{i \leq n} -\val(a_i)$.  
Then $\g \in \Delta+ \G(F)$ so   $\g \leq \g'$ for some $\g' \in \G(F)$.    Let $E=\{(x_1,\ldots,x_n) \in V:   
\val(x_i) \geq - \g' \}$.  Then $E$ is $F$-definable, bounded, g-closed, and $W \meet E$ is not g-closed, by the criterion.
 \eprf

As  pointed out by an anonymous referee,    if $W$ is not g-closed, there may still be no bounded subset $E$ defined over $F_0$ with $W \meet E$ non-g-closed; for instance this happens when $F_0$ is trivially valued and $W=\{x: \val(x)<0\}$.  On the other hand since the family of g-closed
sets is definable in definable families, if $F_0$ is nontrivially valued, then such a set $W$ will be definable over $F_0^a$ (a model of $\ACVF$);
and it follows that one will also be definable over $F_0$.

\begin{cor}  \label{gcriterion1.1}  Let $W$ be a definable subset of a variety $V$.  Assume whenever a definable type $p$  
on $W$, viewed as a set of \textup{(}simple\textup{)} points on $\std{W}$, has a limit point $p' \in \std{V}$, then $p' \in \std{W}$.  
Then $W$ is g-closed.
\end{cor}   

\prf    
We will verify the criterion of \lemref{gcriterion1}.     
 Let $(K,\Delta)$ give rise to a g-pair $(K, \bK)$ over $F$ with $K$ finitely generated over $F$, and $\G(K)=\Delta+ \G(F)$,
 $F$ maximally complete.  
 Let $a \in W(K)$.   
 Let $a'$ be the same point $a$, but viewed as a point of $V(\bK)$. 
 We have to show that $a' \in W(\bK)$.    Let $d=(d_1,\ldots,d_n)$ be a basis for $\Delta$.
Note $\tp(d/F)$  has $0=(0,\ldots,0)$ as a limit point, in the sense of \lemref{cc1}.  Hence $\tp(d/F)$ 
extends to an $F$-definable type $q$.  Now $\tp(a/F(d))$ is 
stably dominated by \thmref{maxcomp} (2), so in particular definable;
hence $p=\tp(a/F)$ is definable.   Since $F$ is maximally complete and $\G(\bK)=\G(F)$, 
  $p' = \tp(a'/F)$ is stably dominated by \thmref{maxcomp}.  
  Furthermore, $p'$ is  a limit of $p$.   To check this, since $F$ is an elementary submodel 
  and $p,p'$ are $F$-definable, it suffices to consider $F$-definable open subsets of $\std{V}$, of the form
  $\val (g ) <\infty$, $\val (g) < 0$ or $\val (g) > 0$ with $g$ a regular function on a Zariski open subset of $V$.  
  If $p'$   belongs to such an open set, the strict inequality holds of $g(a')$, and hence clearly of $g(a)$; so $p$ belongs to it too.
    By assumption, $p' \in \std{W}$, so $a' \in W$.   \eprf

\begin{lem} \label{gcriterion02}  Let $F$ be a valued field, $V$ an $F$-variety, and let $Z \nsubseteq V \times \G^{\ell}$ be $\ACVF_F$-definable.  Then $Z$ is g-closed if and only if  for any g-pair
$(K,\bK )$ over $F$,   
 $\pi(Z(K)) \nsubseteq Z(\bK)$.
   \end{lem}
   
   \prf If $Z$ is g-closed then  the condition on g-pairs
   is also clear, since $\pi$ is order-preserving.  In the other direction, 
   let $\widetilde Z$ be the pullback of $Z$ to $V \times \VF^{\ell}$.  Then $Z$ is g-closed if and only if $\widetilde Z$
   is g-closed.   The condition 
   $\pi(Z(K)) \nsubseteq Z(\bK)$ implies $\widetilde Z(K) \nsubseteq \widetilde Z(\bK)$.  By \lemref{gcriterion0},
   since this holds for any g-pair $(K,\bK)$, $\widetilde Z$ is indeed g-closed.
   \eprf 

\section{v-topology and specialization}\label{ss8.2} 
Let $F$ be a valued field, and consider pairs $(K,\Delta)$,
with $(K,v_K)$ a valued field extension of $F$, and $\Delta$ a proper convex subgroup of $\G(K)$,
with $\G(F) \nsubseteq \Delta$.  Let   $R = \{a \in K: v_K(a) > 0 \ \text{or} \ v_K(a) \in \Delta \}$.
Then $M = \{a \in R: v_K(a) \notin \Delta \}$ is a maximal ideal of $R$ and we
may consider the field 
$\tK = R / M$,  with valuation  $v_{\tK}(r) = v_K(a)$  for nonzero $r = a+M \in \tK$.
We will refer to $(K,\tK)$ and the related data    as a {\em v-pair} over $F$.   For an affine $F$-variety $V \nsubseteq \Aa^n$, \index{v-pair}
let $V(R) = V(K) \meet R^n$.  If $h: V \to V'$ is an isomorphism between $F$-varieties,
defined over $F$, then since $F \nsubseteq R$ we have $h(V(R)) = V'(R)$.  Hence
$V(R)$ can be defined independently of the embedding in $\Aa^n$, and the notion can be extended to an arbitrary $F$-variety.  We have a residue map $\pi: V(R) \to V(\tK)$.  We will write $\pi(x')=x$ to mean:  $x' \in V(R)$ and $\pi(x')=x$, and say:  $x'$ specializes to $x$.
Note that $\G(\tK) = \Delta$.  If $\g = v_K(x)$ with $x \in R$, we  also write $\pi(\g)= \g$ if $v_K(x) \in \Delta$, and $\pi(\g) = \infty$ if $\g > \Delta$.
\tcb{Note also that if $F$ has characteristic $(0,p)$, i.e. $p \neq 0$ but $v(p)>0$ in $F$, then $v(p) \in \Delta$, so $p \notin M$, 
and hence $\tK$ also has characteristic $(0,p)$.}

\begin{lem} \label{lem:v-closed}     Let $V$ be an $F$-variety, $W$ an $\ACVF_F$-definable subset of $V$.
Let $(K,\tK)$ be any v-pair over $F$, with $\tK \models \ACVF$.   
Then $W$ is v-closed if and only if $\pi(W(R) ) \nsubseteq W(\tK)$. \end{lem}

\prf  Since $\ACVF_F$ is complete and eliminates quantifiers, we may assume $W$ is defined without quantifiers.  By the discussion above, we may take $V$ to be affine; hence we may assume $V = \Aa^n$. 

 Assume the criterion  holds.  Let $b \in V(\tK) \m W(\tK)$.  
 If $a \in V(R), b = \pi(a)$,
then $a \notin W$.  
Thus   there exists a $K^{\alg}$-definable open ball containing $a$ and disjoint from
$W$. Since $F \tcb{\subset} \tK$, we may view $\tK$ as embedded in $R$, hence take $a = b$.
It follows that the complement of $W$ is v-open, so $W$ is v-closed.    

 Conversely, assume $W$ is v-closed, and let $a \in W(R)$, $b = \pi(a)$.  Then $b \in V(\tK)$.  If $b \notin W$,
there exists $\g \in \G(F)$ such that, in $\ACVF_F$, the $\g$-polydisc $D_\g(b)$ is disjoint from $W$.  However we
have $a \in D_\g(b)$, and $a \in W$, a contradiction.  
\eprf

\begin{lem} \label{vcc}   Let $U$ be a variety over a valued field $F$, let 
$f: U \  \to \G_\infty$ be an $F$-definable function, and let $e \in U(F)$.
Then $f$ is v-continuous at $e$ if and only if for any v-pair $(K,\tK)$ over $F$ 
and any $e' \in U(R)$,  with $\pi(e') = e$,  we have $f(e) = \pi (f(e'))$.  
Furthermore,
if $F$ is nontrivially valued, one can take $\tK=F$, and if
 $f(e) \in \Gamma$ then in fact $f$ is v-continuous  at $e$ if and only if it is constant on some v-neighborhood of $e$. 
\end{lem}

\prf    
Embed $U$ in affine space; then we have a basis of v-neighborhoods $N(e,\delta)$
of $e$ in $U$ parameterized by elements of $\Gamma$, with $\delta \to \infty$.  

First suppose $\g = f(e) \in \G$. 
Assume for some nontrivial
v-pair $(K, F)$ 
and for every  $e' \in U(R)$ with $\pi(e') = e$,  we have $f(e) = \pi (f(e'))$.  
To show that $f \inv(\g)$
contains an open neighborhood of $e$, it suffices,
since $f \inv(\g)$ is a definable
set, to show that it contains an open neighborhood defined
over some set of parameters.  Now if we take $\delta > \G(F)$, $\delta \in \G(K)$,
then any element $e'$ of $N(e,\delta)$ specializes to $e$,
i.e. $\pi(e')=e$, hence $f(e) = f(e')$ and $f \inv(\g)$ contains an open neighborhood.

Conversely if $f \inv(\g)$ contains an open neighborhood of $e$, this
neighborhood can be taken to be $N(e,\delta)$ for some $\delta \in \Qq \tensor \G(F)$.
It follows that the criterion holds, i.e. $\pi(e') =e$ implies $e' \in N(e,\delta)$ so
$f(e')=f(e)$, for any v-pair $(K,\tK)$.

Now suppose $\g = \infty$.  Assume for some nontrivial
v-pair $(K, F)$ 
and for every  $e' \in U(R)$ with $\pi(e') = e$,  we have $f(e) = \pi (f(e'))$.  
We have to show that for any $\g'$,
$f \inv((\g',\infty])$  contains an open neighborhood of $e$.    
In case $F$ is nontrivially valued, it suffices to take $\g' \in \G(F)$.  Indeed as above, any
element $e'$ of $N(e,\delta)$ must satisfy $f(e') > \g'$, since $\pi (f(e')) = \infty$.
Conversely, if continuity holds, then for some definable function $h: \Gamma_{> 0} \to \Gamma_{> 0}$,
if $e' \in N(e,h(\g'))$ then $f(e') > \g'$; so if $\pi(e')=e$, i.e.
$e' \in N(e,\delta)$ for all $\delta > \G(F)$, then $f(e') > \G(F)$ so $\pi (f(e')) = \infty$. \eprf

\begin{rem} Let  $f: U \to \G$ be as in \lemref{vcc}, but suppose it is  merely
(v-to-g-)-continuous at $e$,  i.e. the inverse image of any interval around 
$\g=f(e) \in \G$ contains a v-open neighborhood of $e$.  Then $f$ is v-continuous
at $e$. \end{rem}

\prf It is easy to verify that, under the conditions of the remark, the criterion 
holds:  $\pi(f(e'))$ will be arbitrarily close to $f(e)$, hence they must be equal.  (Let us also sketch   a
more geometric proof.  
We have to show that $f \inv(\g)$ contains an open neighborhood of $e$.
If not then there are points $u_i$ approaching $e$ with $f(u_i) \neq \g$.  By
curve selection we may take the $u_i$ along a curve; so we may replace $U$
by a curve.  By pulling back to the resolution, it is easy to see that we may take
$U$ to be smooth.  By taking an \'etale map to $\Aa^1$ we find an isomorphism
of  a v-neighborhood of $e$ with a neighborhood of $0$  in $\Aa^1$; so we may
assume $e=0 \in U \nsubseteq \Aa^1$.  For some neighborhood $U_0$ of $0$
in $U$, and some rational function $F$, we have $f(0) = \val(F)$ for $u \in U_0 \m 0$.
By (v-to-g-)-continuity we have $f(0) = \infty $ or $f(0) = \val(F) \neq \infty$ 
also.  But by assumption $\g \neq \infty$.  Now $f  = \val(F)$ is v-continuous, a contradiction.)
\eprf 

\begin{lem} \label{v-density}   
 Let $V$ be an $F$-variety with $F$ algebraically closed, $W' \nsubseteq W$ two $\ACVF_F$-definable subsets of $V$.
Then $W'$ is v-dense in $W$ if and only if for any $a \in W(F)$, for some v-pair $(K,F)$ and
$a' \in W'(K)$, $\pi(a')=a$.  
\end{lem}

\prf  Straightforward, but this and \lemref{vcc-f} will not be used and are left as remarks.\eprf

\begin{lem} \label{vcc-f}   Let $U$ be an algebraic variety over a valued field $F$, and let 
$Z$ be an $F$-definable family of definable functions $U  \to \G$.   
Then the following are equivalent:
\begin{enumerate}
\item There exists an $\ACVF_F$-definable, v-dense subset $U'$ of $U$ such that each $f \in Z$ is v-continuous at each point;
\item  for any $K,\tK $ such that $(\tK,F)$ and $(K,\tK)$ are 
both v-pairs over $F$,  for any $e \in U(F)$, for some $e' \in U(\tK)$ specializing to $e$,
for any $f \in Z(\tK)$ and any  $e'' \in U(K)$ specializing to $e'$, we have  
 $f(e'') = f(e')$.
 \end{enumerate}
\end{lem}

\prf   Let $U'$ be the set of points where each $f \in Z$ is v-continuous.  Then $U'$
is  $\ACVF_F$-definable, and by \lemref{vcc}, for $\tK \models \ACVF_F$   we have
that  $e' \in U'(\tK)$ if and only if for any $f \in Z(\tK)$, any v-pair $(K,\tK)$ 
and any  $e'' \in U(K)$ specializing to $e'$,    $f(e'') = f(e')$.
Thus (2) says that for any v-pair $(\tK,F)$, and any $e \in U(F)$,   some $e' \in U'(\tK)$ specializes to $e$.  By \lemref{v-density} this is equivalent to $U'$ being dense.  \eprf

\medskip

Let $U$ be an $F$-definable v-open subset of a smooth quasi-projective variety $V$ over a valued field $F$, let $ W$ be an $F$-definable 
open subset of $\G^m$,  let $Z$ be an algebraic variety over $F$, 
and let $f: U \times W \to  \std{Z}$ or $f: U \times W \to \G_\infty^k$
be an $F$-definable function.  We consider $\G^m$ and $\G_\infty^k$
with the order topology. We say
$f$ is {\em \textup{(}v,o\textup{)}-continuous} at $(a,b) \in U \times W$ if the preimage of every \index{(v,o)-continuous}
open set containing $f (a, b)$ contains
the product of a v-open containing $a$ and an open containing $b$.

\begin{lem}  \label{cont3}
Let $U$ be an $F$-definable v-open subset of a smooth quasi-projective variety $V$ over a valued field $F$, let $ W$ be an $F$-definable 
open subset of $\G^m$,  let $Z$ be an algebraic variety over $F$, 
and let $f: U \times W \to  \std{Z}$ or $f: U \times W \to \G_\infty^k$
be an $F$-definable function.  
Then $f$ is \textup{(}v,o\textup{)}-continuous 
  if and only if it is continuous
separately in each variable.   More precisely $f$ is \textup{(}v,o\textup{)}-continuous at $(a,b) \in U \times W$
provided that $f(x,b)$ is v-continuous at $a$, and $f(a',y)$ is continuous at $b$ for any $a' \in U$,
or dually that  $f(a,y)$ is continuous at $b$, and $f(x,b')$ is v-continuous at $a$ for any $b' \in W$.  
\end{lem}

\prf  Since a base change will not affect continuity, we may assume $F \models \ACVF$. 
The case of maps into $\std{Z}$ reduces to the case of 
maps into $\G_\infty$,
by composing with continuous definable maps into $\G_\infty$, which determine the  topology on $\std{Z}$.  For maps into  ${\G_\infty^k}$,
since the  topology on ${\G_\infty^k}$ is the product topology,  it suffices also to check
for maps into $\G_\infty$.
 So assume $f: U \times W \to \G_\infty$ and $f(a,b)=\g_0$.  Suppose $f$ is not continuous at $(a,b)$.   So for some neighborhood $N_0$ of $\g_0$ (defined over $F$)
there exist $ (a',b')$ arbitrarily close to $(a,b)$ with $f(a',b') \notin N_0$.  
Fix a  metric on $V$ 
near $a$, and write $\nu(u)$ for the valuative distance of
$u$ from $a$.  Also write $\nu'(v)$ for $\min | v_i-b_i|$, where $v=(v_1,\ldots,v_m), b=(b_1,\ldots,b_m)$.  For any $F' \supset F$, let $r_0^+|F' $ be the type of elements $u$ with
  $\val(a) < \val(u)$ for every nonzero $a$ in $F'$, and 
  let $r_1^- |F'$ be the type of elements $v$ with 
 $0 < \val(v) < \val(b)$ for every  $b$ in $F'$ with $\val (b) >0$. 
  Then $r_0^+,r_1^-$ are definable types, and they are orthogonal
 to each other, that is,  $r_0^+ (x) \union r_1^-(y)$ is a complete definable type.  
 Consider $u, v \in \Aa^1$ with $u \models r_0^+|F, v \models r_1^- |F$.  
Since $F(u,v)^{\alg} \models \ACVF$, there exist $a' \in U(F(u,v)^{\alg})$ and $b' \in W(F(u,v)^{\alg})$ such that $\nu(a')  \geq \val(u)$, $\nu'(b') \leq \val(v)$, and $f(a',b') \notin N_0$.  Note
that any nonzero coordinate of $a'-a$ realizes $r_0^+$; since $r_0^+$ is orthogonal to
$r_1^-$ and $v \models r_1^- | F(u)$, we have $a'-a \in F(u)^{\alg}$, so $a' \in F(u)^{\alg}$.
Similarly $b' \in   \G(F(v)^{\alg})$.   Say two points of $\G_\infty$ are very close over $F$
if the interval between them contains no point of $\G (F)$.  
   By the continuity assumption (say the first version),
$f(a',b')$ is very close to $f(a',b)$ (even over $F(u)$) and $f(a',b)$ is very close to $f(a,b)$
over $F$.  So $f(a',b')$ is very close to $f(a,b)$ over $F$.  But then $f(a',b') \in N_0$,
a contradiction.  
\eprf  

\begin{cor}  \label{cont3c}
More generally, let $f: U \times \G_\infty^{\ell} \times W \to \std{Z}$ be $F$-definable, and let 
$a \in U \times \G_\infty^{\ell}$, $b \in W$.  Then  $f$ is \textup{(}v,o\textup{)}-continuous at $(a,b)$ if  
$f(a,y) $ is continuous at $b$, and $f(x,b')$  is \textup{(}v,o\textup{)}-continuous at $a$ for any $b' \in W$. 
 \end{cor}

\prf  Pre-compose with $\mathrm{Id}_U \times \val \times \mathrm{Id}_W$.  \eprf

\begin{rem}   It can be shown that a definable function $f: \G^n \to \G$, continuous in each variable, is continuous.  But this is not the case for $\G_\infty$.  For instance, $|x-y|$ is continuous in each variable, 
if it is given the value $\infty$ whenever $x=\infty$ or   $y=\infty$.  But it is not continuous
at $(\infty,\infty)$, since on the line $y=x+\b$ it takes the value $\b$.   By pre-composing 
with $\val \times \mathrm{Id}$ we see that \lemref{cont3} cannot be extended to $W \nsubseteq \G_\infty^m$. \end{rem}
 
 \section{$\ACV2F$}\label{ss8.3}

 We consider the theory
$\ACV2F$  of triples $(K_2,K_1,K_0)$ of fields with surjective, non-injective places $r_{ij}: K_i \to K_j$ for $i>j$, $r_{20}=r_{10} \circ  r_{21}$,  such that $K_2$ is algebraically closed.  
\tcb{We shall denote by $\ACV2F_{p_2, p_1, p_0}$ the theory of such triples with $K_i$ of characteristic $p_i$.}
We will work in $\ACV2F_{F_2}$, i.e. over constants for some 
subfield of $K_2$, but will suppress $F_2$ from the notation. \nomenclature{$\ACV2F$}{the theory of  valued algebraically
closed fields with valued residue field}
 The lemmas below should be valid over imaginary    constants too, at least from $\G$.
\nomenclature{$\ACV2F_{p_2, p_1, p_0}$}{completion of $\ACV2F$}

We let $\G_{ij}$ denote the value group corresponding to $r_{ij}$ \tcr{and we write $\G_{ij, \infty}$ for $\G_{ij}$ with an element $\infty$ added with usual conventions}.
Then we have a natural exact sequence
\[ 0 \to \G_{10} \to \G_{20} \to \G_{21} \to 0.\]
The 
inclusion $\G_{10} \to \G_{20} $ is given as follows:  for $a \in \Oo_{21}$, 
$\val_{10}(r_{21}(a)) \mapsto \val_{20}(a)$.  Note that if $\val_{10}(r_{21}(a))=0$
then $a \in \Oo_{20}^*$ so $\val_{20}(a)=0$.  The surjection on the right is
$ \val_{20}(a) \mapsto \val_{21}(a)$.  

Note that $(K_2,K_1,K_0)$ is obtained from $(K_2,K_0)$ by expanding the value group
$\G_{20}$ by a predicate for $\G_{10}$.  On the other hand it is obtained from
$(K_2,K_1)$ by expanding the residue field $K_1$.  

\tcb{We will use characteristics $(0,0,0)$, resp. $(p,p,p)$, when starting with a value field of characteristic $(0,0)$, resp. $(p,p)$; when the field we start with
has characteristic $(0,p)$, we will use $\ACV2F_{0,p,p}$ for the g-criterion, and $\ACV2F_{0,0,p}$ for the v-criterion.}

\begin{lem} \label{2v1} \leavevmode \begin{enumerate} \item \tcb{The theory $\ACV2F_{p_2,p_1,p_0}$ is complete.}
\item 
\tcb{The induced structure on $(K_1,K_0)$ is just the valued field structure; moreover $(K_1,K_0)$ is stably embedded.}
\item \tcb{The set of stably dominated types 
$\std{V}$ is unambiguous for $V$ over $K_1$, whether interpreted in $(K_1,K_0)$ or in $(K_2,K_1,K_0)$.}
\item  
 \tcb{The sorts $(K_2,\G_{20})$ admit quantifier elimination in the language with the ring operations on $K_2$,
 the valuation map into $\G_{20}$, the group operations on $\G_{20}$ and  the predicate $\G_{10}$ on $\G_{20}$.}
\end{enumerate}  \end{lem}

\prf \tcb{(1) and (2) are special cases of \lemref{stablyemb}. Indeed, take $T$ to be $\ACVF_{p_2,p_1}$, with sort $D$ denoting the residue field,
and take $T^*_D$ to be $\ACVF_{p_1,p_0}$.  
 (3) is a consequence of (2).
 For (4), we use the quantifier elimination statement of \lemref{stablyemb} applied to $T=\ACVF_{p_2,p_0}$,
 $D=\G$ (which we will refer to as $\G_{20}$), $T_D^*$ the expansion of $\Th(\G_{2,0})$ by a predicate for a 
 nontrivial convex subgroup $\G_{10}$.     Quantifier-elimination for $\ACVF$ as well as the 
  stable embeddedness of $\G$ via term functions is well-known, cf.   \ref{ss2.6}; quantifier elimination for $(\G_{20},\G_{21},+,-,0,<)$ is easy and left to the reader.}
 \eprf

\begin{lem} \label{2v2}Let $W$ be a definable set in $(K_2,K_1)$ \textup{(}possibly in an imaginary sort\textup{)}.
\begin{enumerate}
\item  Let  $f: W \to \G_{2, 1}$ be a definable function in $(K_2,K_1,K_0)$.  Then 
there exist $(K_2,K_1)$-definable functions $f_1,\ldots,f_k$ such that on any
$a \in \dom(f)$ we have $f_i(a) = f(a)$ for some $i$.
\item Let $f: \G_{21} \to W$ be a $(K_2,K_1,K_0)$-definable function.    Then $f$ is  $(K_2,K_1)$-definable
\textup{(}with parameters; see remark below on parameters\textup{)}.
\end{enumerate}
 In fact this is true for any
expansion of $(K_2,K_1)$ by relations $R \nsubseteq K_1^m$.  \end{lem}

\prf   We may assume $(K_2,K_1,K_0)$ is saturated.  
We shall use some basic properties of stably embedded sets for which we refer to the appendix
of \cite{ch}.

(1)  It suffices to show that for any
$a \in W$ we have 
$f(a) \in \dcl_{21} (a)$, where $\dcl_{21}$ refers to the structure $M_{21} = (K_2,K_1)$.
We have at all events that $f(a)$ is fixed by $\Aut(M_{21}/K_1,a)$.  By stable embeddedness of $K_1$
in $M_{21}$, we have $f(a) \in \dcl_{21}(e,a)$ for some $e \in K_1$.  But 
by orthogonality of $\G_{21}$ and $K_1$ in $M_{21}$ we have $f(a) \in \dcl_{21}(a)$.

(2)  Let $A$ be a base structure, and consider a type $p$ over $A$ of elements of $\G_{21}$.  Note that the induced structure
on $\G_{21}$ is the same in $(K_2,K_1,K_0)$ as in  $(K_2,K_1)$, and that $\G_{21}$ is orthogonal to $K_1$ in both senses.
For $a \models p$, $b=f(a)$, let $g(b)$ be an enumeration of the $(K_2,K_1)$-definable
closure  of $b$ within $K_1$ (over $A$).  By
orthogonality, $g \circ f$ must be constant on $p$; say it takes value $e$ on $p$.  Now $\tp_{21}(ab / A, e) \models \tp_{21}(ab/ A, K_1)$
by stable embeddedness of $K_1$ within $(K_2,K_1)$.  By considering automorphisms it follows that 
 $\tp_{21}(ab / A, e) \models \tp_{210}(ab/A, e) $, so  $\tp_{21}(ab / A, e)$ is the graph of a function on $p$; this function must
 be $f | p$.  By compactness, $f$ is $(K_2,K_1)$-definable.  
\eprf

 \begin{rem}  \label{2v2.5}  Let $D$ be definable in $(K_2,K_1,K_0)$ over an algebraically closed substructure $(F_2,F_1,F_0)$ of constants. In particular $\res_{21}(F_2) \nsubseteq F_1$ and $\res_{10}(F_1) \nsubseteq F_0$, but we do not assume equality.
If $D$ is $(K_2,K_1)$-definable with additional parameters, then $D$ is $(K_2,K_1)$-definable over $(F_2,F_1)$.    \end{rem}

\prf  We may take $(K_2,K_1,K_0)$ saturated.  Let $e$ be a canonical parameter 
for $D$ as a $(K_2,K_1)$-definable set.    Note that $e$ is fixed by the group $\Aut(K_2,K_1 / F_2,K_1)$.  Hence by 
stable embeddedness of $(K_1,K_0)$, we have $e \in \dcl_{K_2,K_1}(F_2,F_1')$ for some (small) $F_1' \nsubset K_1$.
As $K_1$ is \tcb{a pure algebraically closed field} stably embedded in \tcb{$(K_2,K_1)$} and has elimination of imaginaries, $\dcl_{K_2,K_1}(F_2,e)=
\dcl_{K_2,K_1}(F_2,c)$ for some tuple $c=(c_1,\ldots,c_m)$ of elements of $K_1$.  Now each $c_i$ is fixed
by  $\Aut(K_2,K_1,K_0 / F_2,F_1,F_0)$, hence by $\Aut(K_1,K_0 / F_1,F_0)$; it follows easily that $c_i \in F_1$
(since non-algebraic elements of a valued field cannot be definable over residue field elements).
\eprf

\begin{lem}  \label{2v3}  Let $W$ be a $(K_2, K_1)$-definable set and   let $I $ be a definable subset of $\G_{21}$ and let  $f: I   \times W  \to \G_{21, \infty} $ be a $(K_2,K_1, K_0)$-definable 
function such that for fixed $t \in I$,
$f_t(w)=f(t,w)$ is $(K_2,K_1)$-definable.  Then $f$ is  $(K_2,K_1)$-definable.
\end{lem}

\prf     
Applying compactness to the hypothesis, we see that 
there exist finitely many functions $g_k,h_k$ such that $g_k$ is $(K_2,K_1)$-definable,
$h_k$ is definable, and that for any  $t \in I$ for some $k$ we have $f(t,w) = g_k(h_k(t),w)$.  
Now by  \lemref{2v2} (2), $h_k$ is actually $(K_2,K_1)$-definable too.  So we may simplify 
to $f(t,w) = G_k(t,w)$ with $G_k$ a $(K_2,K_1)$-definable function.  But   every definable subset of $I$ is $(K_2,K_1)$-definable, in particular $\{t: (\forall w) (f(t,w) = G_k(t,w))\}$.  
From this it follows that $f(t,w)$ is $(K_2,K_1)$-definable. \eprf

\begin{lem}  \label{2v4}  Let $T$ be any theory, $T_0$ the restriction to a sublanguage $L_0$,
and let $\Uu \models T$ be a saturated model, $\Uu_0 = \Uu | L_0$.   Let $V$ be a definable set of $T_0$.
Let $\std{V}$, $\std{V}_0$ denote the spaces of generically stable types in $V$ of $T,T_0$
respectively.  Then there exists a map $r_0: \std{V} \to \std{V}_0$ such that 
$r_0(p) | \Uu_0 = (p | \Uu) | L_0$.  If $A = \dcl(A)$ \textup{(}in the sense of $T$\textup{)} 
and $p$ is $A$-definable, then $r_0(p)$ is $A$-definable. \end{lem}

\prf  In general, a definable type $p$ of $T$ over $\Uu$ need not restrict to a definable type of $T_0$.  However, when $p$ is generically stable, for any formula $\phi(x,y)$
of $L_0$ the $p$-definition $(d_p x) \phi(x,y)$ is equivalent to a Boolean combination
of formulas $\phi(x,b)$.  Hence $(d_p x) \phi(x,y)$ is $\Uu_0$-definable.  The statement
on the base of definition is clear by Galois theory.  \eprf

\begin{remark}The same holds of course when $T_0$ is interpreted in $T$ (not necessarily as a reduct).
\end{remark}

Returning to $\ACV2F$, we have:

\begin{lem} \label{2v5} Let $V$ be an algebraic variety over $K_1$.  Then the restriction map
of \textup{\lemref{2v4}} from the stably dominated types of $V$ in the sense of $(K_2,K_1,K_0)$ to those in the sense of $(K_1,K_0)$ 
is a bijection.  \end{lem}

\prf This is clear since $(K_1,K_0)$ is embedded and stably embedded in $(K_2,K_1,K_0)$.  (``Embedded'' means that the induced
structure on $(K_1,K_0)$ is just the $\ACV2F$-structure.) \eprf

We can thus write  unambiguously $\std{V}_{10}$ for $V$ an algebraic variety over $K_1$.

Now let $V$ be an algebraic variety over $K_2$.  Note that $K_1$ may be interpreted
in $(K_2,K_0,\G_{20}, \G_{10})$ (the enrichment of $(K_2,K_0,\G_{20})$ by a predicate
for $\G_{10}$).

\begin{lem}  \label{2v6} Any stably dominated type of $(K_2,K_0)$ in $V$ over $\Uu$
generates a complete type of $(K_2,K_1,K_0)$.    More generally,
assume $T$ is obtained from $T_0$ by expanding a linearly ordered sort $\G$ of $L_0$,
and that $p_0$ is a stably dominated type of $T_0$.  Then $p_0$ generates a complete definable type of $T$; over any base set $A=\dcl(A) \leq M \models T$, $p_0 |A$ generates a complete $T$-type over $A$.
\end{lem}

\prf   We may assume $T$ has quantifier elimination.  Then $\tp(c/A)$
is determined by the isomorphism type of $A(c)$ over $A$.  Now 
since $\G(A(c))=\G(A)$, any $L_0$-isomorphism $A(c) \to A(c')$ is
automatically an $L$-isomorphism.  
\eprf

\begin{lem}\label{2v7}   Assume $T$ is obtained from $T_0$ by expanding a linearly ordered sort $\G$ of $L_0$, and that in $T_0$, a type is stably dominated if and only if it is orthogonal to $\G$.
\tcb{Let $V$ be an $L_0$-definable set.} Then  the following properties of a type on 
$V$ over $\Uu$ 
 are equivalent:  \begin{enumerate}
\item  $p$ is stably dominated; 
\item  $p$ is generically stable;
\item  $p$ is orthogonal to $\G$;
\item  the restriction $p_{0}$ of $p$ to $L_0$ is
stably dominated.
\end{enumerate}
\end{lem}

\prf   The implication (1) to (2) is true in any theory, and so is (2) to (3) given that
$\G$ is linearly ordered.  Also in any theory (3) implies that $p_0$ is
orthogonal to $\G$, so by the assumption on $T_0$, $p$ is stably dominated, hence (4).  Finally, let  $p_0$ be stably dominated and generating a  
type $p$ of  \tcb{$T$} (\lemref{2v6}), let us prove this type is also stably dominated. Using the terminology from \cite{hhm} p.~37,
say $p$ is dominated via some $*$-definable functions $f: V \to D$,
with $D$ a stable ind-definable set of $T_0$.

Since $T$ is obtained by expanding $\G$, which is orthogonal to $D$, the set 
$D$ remains stable in $T$.  Now for any base $A$ of $T$ we have that
$p |A$ is  implied by $p_0 |A$, hence by $(f_{*}(p_0)|A)(f(x))$, hence by 
$(f_{*}(p)|A)(f(x))$.  So (4) implies (1).
\eprf

It follows from Lemmas \ref{2v6} and  \ref{2v7} that for any definable set $V$ in \tcb{$M_{20}^{eq}$}, the restriction map $\std{V}_{210} \to \std{V}_{20}$
is a bijection.  

\begin{lem} \label{2v8}
 For $\ACV2F$, the following properties of a type on 
$V$ over $\Uu$ 
 are equivalent:  \begin{enumerate}
\item  $p$ is stably dominated;
\item  $p$ is generically stable;
\item  $p$ is orthogonal to $\G_{20}$;
\item  the restriction $p_{20}$ of $p$ to the language of $(K_2,K_0)$ is
stably dominated.
\end{enumerate}
\end{lem}

\prf  Follows directly from \lemref{2v7} upon letting $T_0$ be the theory of $(K_2,K_0)$. \eprf

\section{The map $\rto: \std{V}_{20} \to \std{V}_{21}$}    \label{ss8.4}

Let $V$ be an algebraic variety over $K_2$.   We write $V_{210}, V_{20},V_{21}, V_2$, etc., when we wish to view $V$
as a definable set for $(K_2,K_1,K_0)$, $(K_2,K_0)$, $(K_2,K_1)$, or just the field $K_2$, respectively.

We have on the face of it three spaces:  $\std{V}_{2j}$ the space of stably dominated types for $(K_2,K_j)$ for $j=0$ and $1$,
and   $\std{V}_{210}$ the space of stably dominated types with respect to the theory $(K_2,K_1,K_0)$.   But in fact
$\std{V}_{20}$ can be identified with $\std{V}_{210}$, as  \lemref{2v6} and \lemref{2v7} show.
We thus identify $\std{V}_{210}$ with $\std{V}_{20}$.  In particular we use this identification to define a topology on $\std{V}_{210}$.

By \lemref{2v4}, we have a restriction map $\rto: \std{V}_{20} = \std{V}_{210} \to \std{V}_{21}$.
If a stably dominated type over a model $M$ is viewed as a sequence of functions into $\G$ (sending an $M$-definable function into $\G$
to its generic value), then $\rto$ is just composition with the natural homomorphism $\G_{20} \to \G_{21}$. 
Note that $\rto$ is the identity on simple points and  that $\rto$ is continuous.  

\medskip

The following \lemref{2v8.5} will not be used in the rest of the paper.
Note that in  (2) and (3) of \lemref{2v8.5}, it is important that $V$ be allowed to be made of  imaginaries of $(K_2,K_1)$. (In (1) this is 
irrelevant, since $\ACF$ eliminates imaginaries.)  This allows applying them to stable completions in (4).   

\begin{lem} \label{2v8.5} \leavevmode \begin{enumerate} 
\item Let $U$ be a variety  \textup{(}or constructible set\textup{)} over $K_1$.  Let $\std{U}_1$ be the space of stably dominated types of $U$ within $\ACF$.
Then the restriction map $\std{U}_{10} \to \std{U}_1$ is
surjective.
\item  Let $V$ be a pro-definable set over $(K_2,K_1)$.  Then the restriction $\std{V}_{210} \to \std{V}_{21}$ is surjective.
\textup{(}The same is true rationally over any algebraically closed substructure of $(K_2,K_1,K_0)$.\textup{)}
\item  Let $V$ be a pro-definable set over $(K_2,K_1)$.  Then any definable type $q$ on $V_{21}$ extends to a definable type $q'$ of $V_{210}$
 \textup{(}moreover, with $q'$ orthogonal to $\G_{10}$\textup{)}.
\item   Let $V$ be a quasi-projective variety over $K_2$.  Then $\rto$ is surjective and \tcr{definably} closed. 
\item  The topology on $ \std{V}_{21}$ is the quotient topology from $ \std{V}_{210}$.
\end{enumerate}
\end{lem}

\prf
(1)    $\std{U}_1$ is  also the space of definable types of $U_1$,
or again the space of generics of irreducible subvarieties of $U$.
Let $W$ be an absolutely irreducible variety over $K_1$.   We have to show that the  generic type of $W$
 expands to a stably dominated type of $(K_1,K_0)$.  Let $\mathcal{W}$ be a scheme over $\Oo_1$ with generic fiber $W$,
 and with special fiber of dimension equal to $\dim(W)$.  Then there are finitely many types $q$ over $K_1$ of elements of $\mathcal{W}(\Oo)$
 whose residues have transcendence degree equal to  $\dim(W)$, and all of them are stably dominated and have Zariski closure equal to $W$.
 
(2)  Let $p$ be a stably dominated type of $(K_2,K_1)$; it is dominated via some definable map $f$ to a finite-dimensional vector space over $K_1$.
 So $f_*p$ is a definable type of $K_1$.  By the previous paragraph, $f_*p$ expands to a stably dominated type $q$ of $(K_1,K_0)$.
 It is now easy to see (as in \remref{2v2.5}) that $q$ dominates a unique definable type $r$ of $(K_2,K_1,K_0)$ via $f$; and clearly $\rto(r)=\tcb{p}$.

(3) For types on $\G_{21}$ this is easy and left to the reader; in this case, note that  the type of $n$ $\Qq$-linearly  independent  elements
over $\G_{21}$ actually generates a complete type over $\G_{210}$.   Now any definable type $r$ on $V_{21}$ is 
the integral over some definable $q$ on $\G_{21}$ of a definable map into $\std{V}_{21}$; i.e. for any $M$ (over which $r$ is defined),
$r=\tp(c/M)$ where $a \models q |M$ and $s=\tp(c/M(a))_{21}$ is stably dominated.  Let $q'$ be an expansion of $q$ to $\G_{210}$;
we may assume $a \models q'|M$.  By (2), there exists a $(K_2,K_1,K_0)$  expansion $s'$
of $s$ to a stably dominated type  $s'$ over $M(a)=\acl(M(a))$.   Integrating $s'$ over $q'$ we obtain a  definable type of $V_{210}$ restricting to $r$.

(4)  Since $V$ itself is open in some projective variety, we may assume $V$ is projective.  
  Let $X$ be a closed pro-definable subset of $\std{V}_{210}$ and let 
  $q$ be a definable type on $\bX = \rto(X) \nsubset \std{V}_{21}$.  By (3), \tcb{$q$ extends to a definable type  $q'$ on $\bX_{210}$}
  (the same pro-definable set $\bX$, now viewed within the structure $(K_2,K_1,K_0)$).     Using \remref{extend+},  $q'$ lifts to 
  a definable type $\tilde q$ on $X$.   Let $\tilde c \in X$ be a limit point of $\tilde q$; it exists by definable compactness of $\std{V}_{210}=\std{V}_{20}$.
  Let $c=\rto(\tilde c)$; by continuity it is a limit point of $\bX$.  
  
  (5)  Follows from (2) and (4).  
  \eprf

 %
 %
 %
 %
 %
 %
 %

 %
 
We move towards    the  $(K_2,K_1)$-definability of the image of $(K_2,K_1, K_0)$-definable paths in $\std{V}$.

\begin{lem} \label{2v10}  Let $f: \G_{20, \infty} \to \std{V}_{20}$ be $(K_2,K_1,K_0)$-\textup{(}pro\textup{)}-definable.
Assume $\rto \circ f = \bar{f} \circ \pi$ for some $\bar{f}:  \G_{21, \infty} \to \std{V}_{21}$ with $\pi: \G_{20, \infty} \to \G_{21, \infty}$ be the natural projection. 
Then $\bar{f}$ is $(K_2,K_1)$-\textup{(}pro\textup{)}-definable.
\end{lem}  

\prf Let $U$ be a $(K_2,K_1)$-definable set, and let $g: V \times U \to \G_{21, \infty}$ be definable.
  We have to prove the $(K_2,K_1)$-definability of the map:  $(\gamma,u) \mapsto  g(\bar{f}(\alpha), u)$,
  where $g(q,u)$ denotes here the $q$-generic value of $g(v,u)$.  For fixed $\gamma$, this is
  just $u \mapsto g(q,u)$ for a specific $q = \rto (p)$, which is certainly $(K_2,K_1)$-definable. 
By \lemref{2v3}, the map :  $(\gamma,u) \mapsto  g(\bar{f}(\alpha), u)$ is $(K_2,K_1)$-definable. \eprf

 \begin{lem} \label{2v11} 
 Let $f : \Gamma_{20, \infty} \to \std{V}_{20}$ be a $(K_2,K_1,K_0)$-\textup{(}pro\textup{)}-definable path. 
 Then there exists a path $\bar{f}:  \G_{21, \infty} \to \std{V}_{21}$
such that  
   $\rto \circ f = \bar{f} \circ \pi$.
   \end{lem}

 \prf Let us first prove the existence of $\bar{f}$ as in \lemref{2v10}. 
  Fixing a point of $\G_{21, \infty}$, with a preimage $a$ in $\G_{20, \infty}$,  it suffices to show
that $\rto \circ f$ is constant on $\{\g +a: \g \in \G_{10, \infty} \}$.  Hence, for any definable family of test function 
$\phi(x,y) : V \to \G_{20, \infty}$ we need to show that $\g \mapsto \pi (f(\g +a)_* \phi)$ is constant in $\g$;
or again that for any $b$, the map $\g \mapsto \pi (f(\g +a)_* \phi (b))$ 
is constant in $\g$.  This is clear since any definable map $\G_{10} \to \G_{21}$
has finite image (due to orthogonality of $\G_{21}$  and $K_1$ inside $M_{210}$, and
since $\G_{10} \nsubseteq K_1^{eq}$),  and by continuity. 
By  \lemref{2v10}
$\bar{f}$ is definable, it remains to show it is continuous.
This amounts, as the topology on $\std{V}$
is determined by continuous functions into $\G_{20, \infty}$,  
to checking that if $g: \G_{20, \infty} \to \G_{20, \infty}$ is continuous and $(K_2,K_1,K_0)$-definable, then the induced map
$\G_{21, \infty} \to \G_{21, \infty}$ is continuous, which is easy.
\eprf



\begin{example}  \label{2v12} 
Let $a \in \Aa^1$ and let $f_a: [0,\infty] \to \std{\Aa^1}$ be the 
map with $f_a(t) =$ the generic of the closed ball around $a$ of valuative radius $t$.
Then $\rto \circ f_a ( t) = f_a( \pi(t))$, where on the right $f_a$ is interpreted
in $(K_2,K_1)$ and on the left in $(K_2,K_0)$.  
Also, if $f_{a}^\g(t)$ is defined by $f_a^\g (t)=f_a(\max(t,\g))$ for  
then $\rto \circ f_a^\g (t) = f_a^{\pi(\g)} ( \bar{t})$.  \end{example}

Let $\Pp^1$ be endowed with  the standard metric of \lemref{metric}.
Given a Zariski closed set $D \nsubset \Pp^1$ of points,    
recall the standard homotopy $\psi_D: [0, \infty]\times \Pp^1 \to \std{\Pp^1}$
defined in 
 \ref{ss7.6}.

\begin{lem}  \label{2v12p}  For every $(t, a)$ we have
$\rto \circ \psi_D (t,a) = \psi_D( \pi(t),a)$,  where on the right $\psi$ is interpreted
in $(K_2,K_1)$ and on the left in $(K_2,K_0)$.   \end{lem}

\prf  Clear, since $\pi (\rho(a,D)) = \rho_{21} (a,D)$. \eprf 

\begin{lem} \label{2v13}  Let $f: V \to V'$ be an $\ACF$-definable map of varieties over $K_2$.
Then $f$ induces $f_{20}: \std{V}_{20}\to \std{V'}_{20}$ and also
$f_{21}:   \std{V}_{21}\to \std{V'}_{21}$.  We have $\rto \circ f_{20} = f_{21} \circ \rto$. \end{lem} 

\prf Clear from the definition of $\rto$. \eprf

\section{Relative versions}\label{ss8.5}
Let $V$ be an algebraic variety over $U$, with $U$ an algebraic variety over $K_2$,
that is, 
a morphism of algebraic varieties $f : V \to U$  over $K_2$.
We have already defined the relative space $\std{V/U}$.
It is the subset of $\std{V}$ consisting of types
$p \in \std{V}$ such that $\std{f}(p)$ is a simple point of
$\std{U}$.
A map $h : W \rightarrow \std{V/U}$ will be called 
pro-definable (or definable) if the composite \index{pro-definable}
$W \to \std{V}$ is.
We endow $\std{V/U}$ with the topology induced by
the topology of $\std{V}$.
In particular one can speak of
continuous, v-, g-, or v+g-continuous maps with values in $\std{V/U}$.
Exactly as above we obtain 
$\rto :  \std{V/U}_{20} \to \std{V/U}_{21}$. Thus, for any $u_0 \in U$, the map $\rto$ 
restricts to the previous map $\rto : \std{V_{u_0}}_{20} \to \std{V_{u_0}}_{21}$
between the 
respective
fibers over $u_0$.

 The relative version of all the above lemmas holds without difficulty:

\begin{lem} \label{2v10rel}  Let $f:  \tcb{U \times}\G_{20, \infty} \to \std{V /U}_{20}$ be \tcb{a} $(K_2,K_1,K_0)$-\textup{(}pro\textup{)}-definable \tcb{map commuting with the structural maps to $U$}.
Assume $\rto \circ f = \bar{f} \circ \pi$ for some $\bar{f}:  \tcb{U \times} \G_{21, \infty} \to \std{V /U}_{21}$.  
Then $\bar{f}$ is $(K_2,K_1)$-\textup{(}pro\textup{)}-definable.
\end{lem}  
\prf
Same proof as 
\lemref{2v10}, or by  restriction.
\eprf

\begin{lem} \label{2v11rel} 
\tcb{Let $f : \tcb{U \times} \Gamma_{20, \infty} \to \std{V /U}_{20}$ be a $(K_2,K_1,K_0)$-\textup{(}pro\textup{)}-definable map
commuting with the structural maps to $U$}. Then
the assumption that  $\rto \circ f$   factors through    $\tcb{U \times} \G_{21, \infty}$ is   automatically verified.\end{lem}

\prf
 This follows from \lemref{2v11}  since 
 a function on $U \times \G_{20, \infty}$ factors through $U \times \G_{21, \infty}$ if and only if
this is true for the section at a fixed $u$, for each $u$.  
\eprf

\exref{2v12} goes through
for the relative version $\std{\Aa^1 \times U /U}$, where now $a$ may be taken
to be a section $a: U \to \Aa^1$. 

The standard  homotopy on $\Pp^1$  defined in \ref{ss7.6}
may be extended fiberwise to a
homotopy $\psi :  [0, \infty] \times \Pp^1 \times U
\to \std{ \Pp^1 \times U / U}$, which we still call standard.  
Consider now an
$\ACF$-definable (constructible) set  $D \nsubset \Pp^1 \times U$  whose projection to $U$ has finite fibers.  
One may consider as above the standard homotopy with stopping time defined by $D$
at each fiber 
$\psi_D: [0, \infty]\times \Pp^1 \times U \to \std{\Pp^1 \times U / U}$.

In this framework \lemref{2v12p} still holds, namely:
\begin{lem}  \label{2v12prel}  For every $(t, a)$ we have
$\rto \circ \psi_D (t,a) = \psi_D( \pi(t),a)$,  where on the right $\psi$ is interpreted
in $(K_2,K_1)$ and on the left in $(K_2,K_0)$.   \end{lem}

Finally  \lemref{2v13} also goes through in the relative setting:

\begin{lem} \label{2v13rel}  Let $f: V \to V'$ be an $\ACF$-definable map of varieties over $U$ \textup{(}and over $K_2$\textup{)}.
Then $f$ induces $f_{20}: \std{V / U}_{20}\to \std{V' /U}_{20}$ and also
$f_{21}:   \std{V/ U}_{21}\to \std{V' / U}_{21}$.  We have $\rto \circ f_{20} = f_{21} \circ \rto$. \qed \end{lem}

\section{g-continuity criterion}\label{ss8.6}

Let $F \leq K_2$.  Assume $v_{20}(F) \meet \G_{10} =(0)$; so $(F,v_{20} |F) \cong (F,v_{21} |F)$ and
\tcb{$((K_2, v_{20}), (K_2, v_{21}))$} is a g-pair over $F$.   In this case any $\ACVF_{F}$-definable object $\phi$ can be interpreted with respect to $(K_2,K_1)_F$ or to $(K_2,K_0)_F$.  We refer to $\phi_{20},\phi_{21}$.  
In particular if $V$ is an algebraic variety over $F$, then $V_{20}=V_{21}=V$;   $\std{V}$ is $\ACVF_F$-pro-definable, and 
$\std{V}_{20}, \std{V}_{21}$ have the meaning considered above.  If $f: W \to \std{V}$
is a definable function with
$W$ a g-open $\ACVF_F$-definable subset of $V$, we obtain $f_{2j}: W \to \std{V}_{2j}$, $j=0, 1$.
Let $W_{21}, W_{20}$ be the interpretations of $W$ in $(K_2,K_1)$, $(K_2,K_0)$.  By \lemref{gcriterion0}
we have $W_{21} \nsubseteq W_{20}$.

\begin{prop}\label{gcriterion} Let $V$ be an algebraic variety over $F$ and $W$ be a g-open $\ACVF_F$-definable subset of $V$. 
 Assume $v_{20}(F) \meet \G_{10} =(0)$. 
\begin{enumerate}
\item An $\ACVF_{F}$-definable map $g:W \to \Gamma_\infty$ is g-continuous  if and only if 
$g_{21}   = \pi \circ g_{20}$ on $W_{21}$.  
\item  An $\ACVF_{F}$-definable map $g:W \times \G_\infty^{\tcb{n}} \to \Gamma_\infty$ is g-continuous if and only if 
$g_{21} \circ \pi_2  = \pi \circ g_{20}$ on $W_{21} \times \G_{20, \infty}$,
where $\pi_2(u,t) = (u,\pi(t))$, $\pi$ being the projection $\G_{20} \to \G_{21}$. 
\item  An $\ACVF_{F}$-definable map $f:W \to \std{V}$ is g-continuous if and only if 
$f_{21}   = \rto \circ f_{20}$ on $W_{21}$.
\item  An $\ACVF_{F}$-definable map $f: W \times \G_\infty^{\tcb{n}} \to \std{V}$ is g-continuous if and only if 
$f_{21} \circ \pi_2  = \rto \circ f_{20}$ \tcb{on} $W_{21} \times \G_{20, \infty}$.
\end{enumerate}
\end{prop}

\prf (1) Recall that g-continuity of maps to $\G_\infty$ was defined with respect to the g-topology on $\G_\infty$ (as well as on $W$).
The function   $g$ is g-continuous with respect to $\ACVF_F$ if and only
if $g^{-1} (\infty)$ is g-open and
for    any open interval $I$
of $\G_{21}$, $g \inv (I)$ is g-open. 

Let us start with an interval of the form
$I_a = \{ x:  x > \val_{21}(a) \}$,  with $a \in K_2$.

By increasing $F$ we may assume $a  \in F$.   (We may assume $F=F^{\alg}$.  There is no problem replacing
$F$ by $F(a)$ unless $v_{20}(F(a)) \meet \G_{10} \neq (0)$.  In this case it is easy to see that $v_{21}(a)=v_{21}(a')$
for some $a' \in F$, so we may replace $a$ by $a'$.)

We view $U_a  =  g \inv (I_a)$ 
as defined by $\infty > g(u) > \val(a)$ in $\ACVF_F$.
By  \lemref{gcriterion0}, $U_a$ is g-open if and only if  $(U_{a})_{21} \nsubseteq (U_a)_{20}$, that
is,  $ \infty > g_{21} (u)    > \val_{21}(a)$ implies $ \infty > g_{20} (u) > \val_{20}(a)$.
Thus, $g \inv (I_a)$ is g-open for every $a$ if and only if
$g_{21}(u) \leq \pi (g_{20}(u))$ and $g_{20}(u) < \infty$   whenever $g_{21}(u) < \infty$. 
Let $I'_a =    \{x: x < \val_{21}(a) \}$. One gets similarly that
$g \inv (I'_a)$ is g-open for every $a$ if and only if
$g_{21}(u) \geq \pi (g_{20}(u))$  whenever $g_{21}(u) < \infty$.
Again by \lemref{gcriterion0},
$g^{-1} (\infty)$ is g-open if and only if  $g_{20}(u) = \infty$ whenever $g_{21}(u) = \infty$.
The statement follows.

(2)  Let $G(u,a) = g(u,\val(a))$.  Then $g$ is g-continuous if and only if $G$ is g-continuous. 
The statement follows from (1) applied to $G$.

For (3) and (4), we pass to affine $V$, and consider a regular function $H$ on $V$.  Let 
$g(u) = f(u)_{*} (\val (H))$.  Then $f_{21}   = \rto \circ f_{20}$ if and only if for each such $H$ we have
$g_{21} =    \pi \circ g_{20}$; and $f$ is g-continuous if and only if, for each such $H$, $g$
is g-continuous.   Thus (3) follows from (1), and similarly (4) from (2).
\eprf

\begin{remark}  A similar criterion is available when $W$ is g-closed rather than g-open; in this case we have $W_{20} \nsubseteq W_{21}$, and the equalities
must be valid on $W_{20}$.  In practice we will apply the criterion only with g-clopen $W$.  \end{remark}

\section{Some applications of the continuity criteria}\label{ss8.7}
As an example of using the continuity criteria,  assume $h: V \to W$ is a   finite surjective 
morphism of separable degree $n$ between algebraic varieties of pure dimension $d$, with $W$ normal.  
For $w \in W$, one may endow $h \inv(w)$ with the structure of a multi-set (i.e. a finite set with
multiplicities assigned to points) of constant cardinality $n$ as follows. 
One consider\tcb{s} a pseudo-Galois covering $h' : V' \to W$ of separable degree $n'$ with Galois group $G$ factoring as $h' = h \circ p$ with $p : V' \to V$ finite of separable degree $m$. If $y' \in V'$, one 
sets $m (y') = \vert G \vert / \vert \mathrm{Stab} (y') \vert$ and for $y \in V$, one sets 
$m (y) = 1/m \sum_{p (y') = y} m (y')$. The  function $m$ on $V$ is independent from the choice of the pseudo-Galois covering $h' $
(if $h''$ is another  pseudo-Galois covering,  consider a pseudo-Galois covering dominating both $h'$ and $h''$).
Also, the function $m$ on $V$ is $\ACF$-definable.
Let $R$ be a regular function on $V$
and set $r = \val \circ R$.  More generally, $R$ may be a tuple of regular functions $(R^1,\ldots,R^m)$, and $r=(\val \circ  R^1,\ldots,\val \circ R^m)$.  
The pushforward $r( h \inv(w))$ is also a multi-set
of size $n$, and a subset of $\G_\infty^m$.  Given a multi-set $Y$ of size $n$ in a linear ordering, 
we can uniquely write $Y = \{y_1,\ldots,y_n\}$ with $y_1 \leq \ldots \leq y_n$ and with repetitions equal
to the multiplicities in $Y$.  Thus, using the lexicographic ordering on $\G_\infty^{\tcb{m}}$, we can write $r(h \inv(w)) = \{r_1(w),\ldots,r_n(w) \}$; in this way we obtain definable functions $r_i: W \to \G_\infty$, $i=1,\ldots, n$.    In this setting we have:

\begin{lem}\label{finite} The functions  $r_i$ are  v+g-continuous. \end{lem}

\prf Note that if $g: A \to B$ is a weakly order-preserving map of linearly ordered set, $X$ is a multi-subset of $A$ of size $n$ and $Y=g(X)$, then $g(x_i)=y_i$ for $i \leq n$.
It follows that both the v-criterion \lemref{vcc} and the g-criterion \propref{gcriterion} (1) hold in this situation.
\eprf

\begin{cor}\label{finite-open}  Let   $h: V \to W$ be a finite surjective  morphism  between algebraic varieties of pure dimension $d$  over a valued field, with $W$ normal.  Then $\std{h}: \std{V} \to \std{W}$ is an open map. 
\end{cor}

\prf  We may assume that $W$ and hence $V$ are affine.  A  basic open subset of $\std{V}$ may be written
as $G=\{p:  (r(p)) \in U \}$ for some $r=(\val \circ R^1,\ldots, \val \circ R^m)$, $R^i$  regular functions on $V$, and some v+g-open definable subset
$U$ of $\G_\infty^n$.   Consider the functions $r_i$ 
  as in  \lemref{finite}.   By \lemref{finite} they are v+g-continuous. By \lemref{basic}, they extend to continuous functions
  $\std{r_i} : \std{W} \to \G_\infty$.
 Since $w \in \std{h}(G)$ if and only if  for some $i$ we have $\std{r_i}(w) \in U$, it follows that
 $\std{h}(G)$ is open.
 \eprf

Note the necessity of the assumption of normality.  If $h$ is a  a pinching of $\Pp^1$, identifying two points $a \neq b$,   the image of a small valuative neighborhood of $a$ is not open.
 
 We also have: 
 
 \begin{lem} \label{proj-open} Let $U$ and $V$ be algebraic varieties over a valued field and let $p: U \times V \to U$ be the projection. Then $\std{p}$ is open. \end{lem}
 
 \prf  By taking open covers, we may assume $U$, and then $V$, are affine.  Embedding $V$ in $\Aa^n$, so that an open
 subset of $\std{U \times V}$ is the restriction of an open subset of $\std{U \times \Aa^n}$, we may assume $V=\Aa^n$.
 By induction on $n$, we reduce to the case $V = \Aa^1$.   It suffices to consider open subsets $\std{H}$ of 
 $\std{U \times V}$ cut out by inequalities $\val (F_i) > 0, \val (G_j) < 0$ where $F_i,G_j$ are regular functions on $U \times V$.  
 By \lemref{extend}, $\tcb{\std{p}}(\std{H}) = \std{p(H)}$.  Since $F_i,G_j$ are continuous in the valuation topology, it is clear that
 $p(H)$ is v-open.  To see that it is g-open,  it suffices by \lemref{gcriterion0} to show that  for any g-pair
$(K,\bK)$ over the base field,  $p(H)(\bK) \nsubseteq p(H)(K)$.  This is clear since $H(\bK) \nsubset H(K)$ (strict inequalities being stronger
for $\bK$), and since $K,\bK$ have the same underlying set.   \eprf

 \begin{cor}   \label{finite-open+} Let   $h: V \to W$ be a   morphism  between algebraic varieties  over a valued field, with $W$ normal.  \tcb{Assume
 $W$ and $V$ are of pure dimension $m$ and $m + n$} and that $h = f \circ g$ where
 $f: V \to W \times \Pp^n$ is a finite surjective morphism,  $g$ is the projection map $W \times  \Pp^n \to W$.
 Then $\std{h}: \std{V} \to \std{W}$ is an open map. 
\end{cor}

\prf  Clear from \corref{finite-open} and \lemref{proj-open}.   \eprf


\begin{cor} \label{constant-lift} Let   $h: V \to W$ be a finite morphism of algebraic varieties  of pure dimension $d$ over a valued field, with $W$ normal and $V$ quasi-projective.   
Let 
$\xi: V \to \G_\infty^n$ be a definable function.  Then there exists a definable function $\xi': W \to \G_\infty^m$ such that
for any path $p: I \to \std{V}$, still denoting by $\xi$ and $\xi'$ their canonical extensions to
$\std{V}$ and $\std{W}$, if $\xi' \circ h \circ p $ is constant on $I$, then so is $\xi \circ p$.  \end{cor}

\prf  By \lemref{contfact}
we may assume $\xi$ is continuous.  Also, we can treat the coordinate functions separately, so we may as well take $\xi: V \to \G_\infty$.
Let $d=\deg(h)$, and 
define $\xi_1,\ldots,\xi_d$ on $W$ as above, so that the canonical extension of $\xi_i$ (still denoted by $\xi_i$) is continuous on $\std{W}$ and $\xi(v) \in \{\xi_1(h(v)),\ldots,\xi_d(h(v))\}$. 
Let $\xi'=(\xi_1,\ldots,\xi_d)$. 
Now if $\xi' \circ h \circ p $ is constant on $I$, then $\xi \circ p$ takes only finitely many values, so by definable connectedness of $I$, cf.  \ref{ss9.4}, it must
be constant too.  
\eprf

\begin{lem}\label{contfact}
Let $V$ be a quasi-projective variety over a valued field and let $\xi : V \to \G_\infty^n$ be a definable function.
Then there exists a v+g-continuous definable function $\xi^*: V\to \G_\infty^N$ and a definable function $d : \G_\infty^N \to \G_\infty^n$
such that $\xi = d \circ \xi^*$.
\end{lem}

\prf 
We may assume $V = \Pp^m$. The statement follows from the following remark:
if $f/g$ is a rational function on $\Pp^m$ with $f$ and $g$ homogeneous of the same degree,
the map $x \mapsto \val ((f/g) (x))$ factors through the maps
$x \mapsto \max (0, \val (f (x)) - \val (g (x)))$ and
$x \mapsto \max (0, \val (g (x)) - \val (f (x)))$.
\eprf

\section{The v-criterion on $\std{V}$}\label{ss8.8}

 %


 
 %

Let $V$ be an algebraic variety defined over a field  $F_2 \nsubseteq \Oo_{21}$.  This means
that $v_{21}(a) \geq 0$ for $a \in F_2$, so $v_{21}(a)=0$ for $a \in F_2$, equivalently
$v_{20}(F_2^{\times}) \nsubseteq \G_{10}$.   This is the \tcb{condition considered in relation with the v-criterion in  \ref{ss8.2}.}
The place $r_{21}$ induces a field isomorphism
$\res_{21}: F_2 \to F_1$.  Let $V_1$ be the conjugate of $V$ under this field isomorphism, so $(F_2,V ) \cong (F_1,V_1)$.   We can also
view $V_1$ as the special fiber of the $\Oo_{21}$-scheme $V_2 \tensor _{F_2} \Oo_{21}$.    As noted earlier, $\std{V}_1$ is unambiguous for varieties over $F_1$.  

  Recall $\std{V}_{20} = \std{V}_{210}$.  Now
$\std{V}_{210}$ has a subset $\std{V}_\Oo = \std{V(\Oo_{21})}$  consisting of types concentrating on $V(\Oo_{21})$.
We   have a definable map $\res: V(\Oo_{21}) \to V(K_1)$.  This induces a map
\[  \res_{21*}: \std{V}_\Oo \to \std{V}_1.\]  

Let $\G_{20}^+ = \{x \in \G_{20, \infty}: x \geq 0 \vee x \in \G_{10} \}$. 
Define a retraction $\pi: \G^+_{20} \to \G_{10,\infty}$ by letting
$\pi (x)=\infty$ for $x \in  \G^+_{20}  \m \G_{10}$.    Note that this is the same as the map $\pi$ in  \ref{ss8.2}.

\begin{lem}  \label{vcc2}    Let $V$ be an algebraic variety over $F_2$, let $W$ be an $\ACVF_{F_2}$-definable  subset
of  $\Pp^n \times \G_\infty^m$ and consider an $\ACVF_{F_2}$-definable map $f:V \to \std{W}$.  
\tcb{\begin{enumerate}
\item Let $x$ be a point in $V(\Oo_{21})$. Then $f$ is v-continuous at $x$ if and only if
$(\res_{21*} \circ f_{20}) (x) = (f_{10} \circ \res_{21}) (x)$.
\item Let $X$ be an $\ACVF_{F_2}$-definable subset of $V$ and
assume $\res_{21*} \circ f_{20} = f_{10} \circ \res_{21}$
at $x$ whenever $x \in V(\Oo_{21})$ and  $\res_{21}(x) \in X$.   Then $f$ is v-continuous at each point of $X$.  In particular, 
if  $f$ is also g-continuous, then  the canonical extension   $F : \std{V} \to \std{W}$   is continuous at each point of $\std{X}$.  
\end{enumerate}}
\end{lem}

\prf  
\tcb{Let $x$ be a point in $V(\Oo_{21})$.
As in  the proof of \lemref{basic0},
 $f$ is v-continuous at  $x$ if and only if
 for every
continuous definable function
$c: \std{W} \to \G_\infty^n$,
$c \circ f $ is v-continuous at $x$.
On the other hand, by    the 
 ``only if''  direction in \lemref{vcc}, 
the other condition holds for $f$ at $x$ if and only if it  holds for $c \circ f$, for any
continuous definable function
$c: \std{W} \to \G_\infty^n$. Thus, in the proof of (1), 
we may assume $f: V \to \G_\infty$, in which case the  statement  follows from
\lemref{vcc}. (2)  follows directly from (1) and \lemref{basic0}.}
  \eprf

\begin{rem} \label{roots}  Let $F(X) \in \Oo_{21}[X]$ be a polynomial in one variable, 
and let $f(X)$ be the specialization to $K_1[X]$.  Assume $f \neq 0$.
Then the map $r_{21}$ takes the roots of $F$ onto the roots of $f$.  
Indeed, consider a root of $f$; we may take it to be $0$.  
Then the Newton polygon
of $f$ has a vertical edge.  So the Newton polygon of $F$ has a very steep edge
compared to $\G_{10}$.  Hence it has a root of that slope, specializing to $0$.
\end{rem}

The following lemma states that a continuous map on $X$ remains continuous relative to a set  $U$
that it does not depend on;  i.e. viewed as a map on $X \times U$ with dummy variable $U$, 
it is still continuous.  This sounds trivial, and  the proof is indeed  straightforward  if one uses the continuity criteria; 
it seems curiously  nontrivial to prove directly. 

 For $U$ a variety and $b \in U$, let $s_b$ denote
the corresponding simple point of $\std{U}$, i.e. the definable type $x=b$.    For $V$ a variety and $q \in \std{V}$,
let $q \tensor s_b$ denote the unique definable type $q(v,u)$ extending $q(v)$ and $s_b(u)$.  

\def\baf{\bar{f}}

\begin{lem}  \label{relcont} Let $U$, $V$ and $V'$ be varieties over a valued field. \tcb{Assume $U$ and $V$ are quasi-projective}
and $X$ be a v+g-open definable subset of  $V'$, or of $V' \times \G_\infty^N$.  
Let $f: X \to \std{V}$ be 
v+g-continuous,  
and let $\baf(x,u) =   f(x) \tensor s_u$.  Then $\baf: X \times U \to \std{V \times U}$ is v+g-continuous.  \end{lem}  

\prf     For g-continuity, we use  \propref{gcriterion} (3) and (4).  We have $f_{21}=\rto \circ f_{20}$ on $X_{21}$. 
 Also for $x \in X_{21}, u \in U_{21}$, we have $\baf_{21} (x,u) = f_{21}(x) \tensor s_u$, and  
  $\baf_{20} (x,u)= f_{20}(x) \tensor s_u$.  Moreover we noted that $\rto$ is the identity on simple points, so 
  $\rto(p \tensor s_b) = \rto(p) \tensor s_b$ in the natural sense.  The criterion follows. 

For v-continuity, \lemref{vcc2} applies.  Assume $\res_{21}(x) \in X$, so $x \in X$.   
 Let $u \in U(\Oo_{21})$.    We have $\res_{21*} \circ f_{20} (x) = f_{10} \circ \res_{21} (x)$.  Now  
 $\res_{21*}(q \tensor s_u) =\res_{21*}(q) \tensor s_{\bar{u}}$, 
where $\bar{u}=\res_{21}(u)$, and $\res_{21}(x,u) = (\res_{21}(x),\bar{u})$, so the criterion follows.
\eprf  

\begin{rem}\tcb{In the context of the previous lemma}, recall that the map $\tensor: \std{U} \times \std{V} \to \std{U \times V}$ is well-defined 
but not continuous \tcb{in general}.  If $f: I \times \std{V} \to \std{V}$ is a homotopy, let $\phi: I \times V \to \std{V}$
be the restriction to simple points, and let $(\phi \tensor \Id)(t,v,u) = \phi(t,v) \tensor u$.  By \lemref{relcont}, 
$(\phi \tensor \Id)$ is v+g-continuous.  By \lemref{hbasic}, it extends to a homotopy  
$I \times \std{V \times U} \to \std{V \times U}$, which we denote $\std{f \times \Id}$.   We easily compute: 
$\std{f \times \Id}(t,p \tensor q) = f(t,p) \tensor q$. 
\end{rem}

\begin{cor}  \label{product-homotopy} Let $U$ and   $V$ be \tcb{quasi-projective} varieties over a valued field and let $X$ and $Y$ be definable subsets of $U$ and $V$. 
Let $f: I \times \std{X} \to \std{X}$ and $g: I' \times \std{Y} \to \std{Y}$ two definable deformation retractions
onto iso-definable
$\G$-internal subsets $S$ and $T$, respectively.
Assume $f$ and $g$ are restrictions of homotopies
$F: I \times \std{U} \to \std{U}$ and $G: I' \times \std{V} \to \std{V}$, respectively.
Then there exists a definable deformation retraction  $h:  (I+I')  \times  \std{X\times Y} \to \std{X \times Y}$  whose image is equal to  $S \tensor T$.  
\end{cor}

\prf    Recall $I+I'$ is obtained from the disjoint union of $I$ and $I'$ by identifying the endpoint $e_I$ of $I$ with the initial point of $I'$.  
The homotopy $\std{F \times \Id}$
restricts to a homotopy $\std{f \times \Id} : I \times \std{X \times Y} \to \std{X \times Y}$
and similarly
$\std{ \Id \times G}$
restricts to a homotopy $\std{\Id \times g} : I' \times \std{X \times Y} \to \std{X \times Y}$.
Let $h$ be the concatenation of $\std{f \times \Id}$ with  $\std{\Id \times g}$, that is, defined by 
\[h(t,z)=\std{f \times \Id}   \hbox{ for } t \in I, \ \ \  h(t,z)= \std{\Id \times g} (t,\std{f \times \Id}(e_I,z)) \hbox{ for } t \in I'.\] 

So $h(t,p \tensor q) = f(t,p) \tensor q$  for $t \in I$, and $= f(e_I,p) \tensor g(t,q)$ for $t \in I'$.  In particular,
$h(e_{I'}, p \tensor q) = f(e_I,p) \tensor g(e_{I'},q)$.  

Since any simple point of $\std{X \times Y}$  has the form $a \tensor b$, we see that   $h(e_{I'},X \times Y) \nsubseteq S \tensor T$.
Hence for any $r \in \std{X \times Y}$,  $h(e_{I'},r)$ is   an integral over $r$  of a function into $S \tensor T$.  But as $S \tensor T$ is
\tcb{iso-definable and} $\G$-internal,  and $r$ is stably dominated, this function is generically constant on $r$, and the integral is an   element of $S \tensor T$.
Thus the final image of $h$ is contained in $S \tensor T$.

Using again the expression for $h(t, p \tensor q)$ we see that if 
 $f(t,s) = s$ for $s \in S$ and $g(t,y)=y$ for $y \in T$,   then      $h(t,z)=z$ for all $t$ and all  $z \in S \tensor T$.   So the final image is
 exactly equal to $S \tensor T$.  
 \eprf

The following statement is a consequence of \corref{product-homotopy} and  \thmref{1}.

\begin{cor}  \label{product-homotopy-eq} Let $U$ and $V$ be quasi-projective varieties over a valued field and
let  $X$ and $Y$ be definable subsets of $U$ and $V$.  
The canonical map $\pi: \std{X \times Y} \to \std{X} \times \std{Y}$ is a homotopy equivalence.  
\end{cor}

\prf  We may assume $U$ and $V$ are projective. 
By \thmref{1}, there exists definable deformation retractions
$F : I \times \std{U} \to  \std{U}$ and  $G:  I' \times \std{V} \to \std{V}$,  leaving  $X$ and $Y$ invariant,
 whose images $\Sigma$ and $\Theta$ are \tcb{iso-definable and} $\G$-internal.  Since $\Sigma$ and $\Theta$ are continuous definable images of
 $\std{U}$ and $\std{V}$, they are definably compact.  
The map $\pi_{\Sigma} \times \pi_{\Theta} : \Sigma \tensor \Theta \to \Sigma \times \Theta$ is continuous
 and injective,  hence a homeomorphism.  Thus the inverse map $\tensor: \Sigma \times \Theta \to \Sigma \tensor \Theta$ is continuous.
 
 Let $f: I \times \std{X} \to \std{X}$, $g: I' \times \std{Y} \to \std{Y}$ be the restrictions of $F$ and $G$, respectively, with images
\tcb{iso-definable and}  $\G$-internal subsets $S$ and $T$.
 Being the restriction of a continuous map, $\tensor: S \times T \to S \tensor T$ is 
 continuous, thus $\pi_S \times \pi_T :
 S \otimes T \to S \times T$ is
 a homeomorphism.  
 Denote by $e$ and $e'$ the endpoints of $I$ and $I'$.

By \corref{product-homotopy}, we have a 
  homotopy equivalence $h_{e'}: \std{X \times Y} \to S \tensor T$ such that the following diagram is commutative:
   \begin{equation*}\xymatrix{
     \std{X \times Y}
	\ar[d]_{\pi_X \times \pi_Y} \ar[r]^{h_{e'}} & 
	S \tensor T  \ar[d]^{{\pi_S \times \pi_T}}\\
	\std{X} \times \std{Y} \ar[r]_{f_{e} \times g_{e'}} &S \times T.
	}
	\end{equation*}
Since $f_{e} \times g_{e'}$ is a homotopy equivalence and
 $\pi_S \times \pi_T$  is a homeomorphism, $\pi_X \times \pi_Y$ is  a homotopy equivalence.  
  \eprf

 %
 


\section{Definability of v- and g-criteria}\label{ss8.9}

\tcb{We shall consider $\stda{V}$  with its canonical strict ind-definable structure defined in \ref{ssbertini}.}

\begin{prop}\label{def-vg-crit}  Let $V$ and $W$ be varieties over a valued field and let $C$ be the set of definable functions $V \to \stda{W}$ that extend to continuous functions $\std{V} \to \std{W}$.
\tcb{Assume $W$ is quasi-projective.} Then $C$ is \textup{(}strict\textup{)} ind-definable.  If $V$ and $W$ depend on a parameter $t$, then this is uniform in the parameter.
\end{prop}

\prf  \tcb{We will use the v- and g-criteria to  show that for  each definable set of definable functions $V \to \stda{W}$ the subset of those that are v-continuous, resp. g-continuous, is definable.}

We begin with v-continuity.  Let  $V$ and $W$ be defined over a field $F_2 \nsubseteq \Oo_{21}$.  Let $f=f_b: V \to \stda{W} \nsubset \std{W}$
be a definable map, with parameter $b \in F_2$.  
By  \lemref{vcc2} (1), $f$ is v-continuous iff the equation:  
\begin{equation}\tag{$\ast$}\label{eqlab}\res_{21*} \circ f_{20} = f_{10} \circ \res_{21}\end{equation}
holds on $V(\Oo_{21})$. 

There is no harm in assuming that $W$ is projective, so as to simplify notation:   $W(\Oo_{21}) = W(K_2)$.  
Now the map \[\res_{21}:   W(K_2) = W(\Oo_{21})   \to W(K_1)\] is $\ACV2F$-definable.   It induces a map
\[  \res_{21*}: \std{W} _{20} \to \std{W_{10}}.\]  

It is easy to see that $\res_{21*} (\stda{W} ) \nsubset \stda{W_1}$.  For instance, the argument for $\std{W_{20}} = \std{W_{210}}$ shows
also that the strongly stably dominated types of these structures coincide; and the image of a strongly stably dominated type under a definable map is strongly stably  dominated in $K_{210}$, and hence in $K_{10}$ which is stably embedded;
see \propref{ssd1} (2).   

The restriction $r$ of $\res_{21*}$ to $\stda{W}$ is $\ACV2F$-piecewise definable (i.e. definable on definable pieces), since $\res_{21*}$ itself is pro-definable.
Now   the set  $f_{20} (V(\Oo_{21}))$ is contained in a definable subset of $\stda{W}$ which does not depend on $b$.
 Hence the displayed equation (\ref{eqlab}) is an $\ACV2F$-definable property of $b$.  
\tcb{Now any  $\ACV2F$-definable subset  $X$ of   $K_2^n$ is defined by quantifier-free formulas in the field
language along with valuation maps $v_{20},v_{21}$, and the group operations on $\G_{20}$; 
this is \lemref{2v1} (4).  But on $F_2$, the valuation $v_{21}$ is trivial; hence $X \meet F_2^n$
is cut out by quantifier-free formulas in the field
language along with   $v_{20}$ alone, so it is already cut out by $(K_2,K_0)$-formula.    } 
 Thus the set of $b$ from $F_2$ for which $f_b$ is v-continuous   is $\ACVF$-definable.  

Similarly, we use the g-criterion \propref{gcriterion} (3)  for  proving that for
each definable set of definable functions $V \to \stda{W}$ the subset of those that are  g-continuous is definable.
Here the defining equation is
\[f_{21}   = \rto \circ f_{20} \hbox{ on } V, \]
$\rto$ is the composition of the equality $\stda{W}_{210} = \stda{W_{20}}$ with the restriction map
$\stda{W}_{210} \to \stda{W}_{21}$, and is clearly piecewise definable in $\ACV2F$.  Once again the quantifier-free induced
structure on $F_2$ is the same as the $v_{20}$-$\ACVF$-structure, which implies the statement.
 \eprf

 
 \chapter{Continuity of homotopies}\label{sec9}

{\small \noindent \textbf{Summary.} This chapter consists mostly of  preliminary material useful for the proof of the main theorem in  Chapter \ref{sec10}.
In \ref{ss9.1} and \ref{ss9.2} we use the continuity criteria of Chapter \ref{specializations} to prove the continuity of functions and  homotopies used in Chapter \ref{sec10}.
In \ref{ss9.3}, we construct inflation homotopies, which are a key tool in our approach. Finally, in \ref{ss9.4} we prove GAGA type results for connectedness
and prove additional results regarding  the  Zariski topology.
\par\bigskip}

\section{Preliminaries}\label{ss9.1}The following lemma will be used both for the relative curve homotopy, and for the inflation homotopy.  

\begin{lem}  \label{homotopy-lift} Let $f: W \to U$ be a  morphism of \tcb{quasi-projective} varieties over some valued field $F$. Let 
 $h: [0,\infty]   \times U \to \std{U}$ be $F$-definable.  Let 
 $H: [0,\infty]   \times W \to \std{W}$ be an $F$-definable lifting of $h$.  Let $H_w(t) = H(t,w)$ and $h_u(t)=h(t,u)$.  
  Assume for all $w \in W$, 
  $H_w$ 
  and $h_{f(w)}$ are   paths and that
  $H_w$ is the
  unique path  lifting $h_{f(w)}$ 
  with $H_w(\infty)=w$.   
Let $X$ be a g-open definable subset of $U$.
  Assume $h$ is g-continuous, and v-continuous on \textup{(}respectively, at each point of\textup{)}  $[0,\infty] \times X$.  \index{v-continuous}
  Then $H$ is  
  g-continuous, and is v-continuous on \textup{(}respectively, at each point of\textup{)} $[0,\infty] \times f \inv (X)$
    \textup{(}we say a function is v-continuous on a subset, if its  restriction to that subset is v-continuous\textup{)}.
\end{lem}  

\prf  We first use the criterion of \propref{gcriterion} (4) to prove g-continuity.  We may assume the
data are defined over a subfield
 $F$ of $K_2$, such that $v_{20}(F) \meet \G_{10} =(0)$; so $(F,v_{20}) \cong (F,v_{21})$.  

To show that $H_{21} \circ \pi_2  = \rto \circ H_{20}$, we fix $w \in W$.  
By \lemref{2v11},    $ \rto \circ H_{20} (w,t) = H'_w \circ \pi$ for some path $H'_w$. 
To show that $H_{21}(w,t) = H'_w(t)$, it is enough to show that  \tcb{$f_{21} \circ H'_{w, 21} = h_{f(w), 21}$}.  It is clear that $H'_w(\infty) =H_{20}(\infty) = w$
since $\rto$ preserves simple points.  To see that \tcb{$f_{21} \circ H'_{w, 21} = h_{f(w), 21}$} it suffices to check that 
$f \circ H'_w \circ \pi= h_{f(w)} \circ \pi$, i.e. $ f \circ \rto \circ H_{20} (w,t) = h_{f(w)}(\pi(t))$.   
Now $ f \circ \rto \circ H_{20}= \rto \circ h_{20} = h_{21} \circ \pi_2$.
It follows that  the g-continuity criterion for $H$ is satisfied.

Let  \tcb{us} now use the v-continuity criterion in \lemref{vcc2} above $X$, $(\res_{21*} \circ H_{20})(t,v) = (H_{10} \circ \res_{21}) (t,v)$ whenever $(f \circ \res_{21})(v) \in X$.  Fixing 
 $w=\res_{21}(v)$, $H_{10}(t,w)$, for $t \in \G_{10}$,   is the unique path lifting $h_{f(w)}$ and starting at 
 $w$, hence to conclude it is enough to prove that 
 $\res_{21*} \circ H_{20}(t,v)$ also has these properties. But continuity follows from
  \lemref{2v20} and the lifting property from \lemref{2v21}.  \eprf


In the next two lemmas we shall use the notations and assumptions in  \ref{ss8.8}.
In particular we will  assume  that
$v_{20}(F_2^{\times})\subset
\Gamma_{10,\infty}$.

\begin{lem} \label{2v20}  Let $V$ be a quasi-projective variety over $F_{2}$.
  Let $f: [0,\infty] \nsubset \G_{20,\infty}  \to \std{V}_{20}$ be a $(K_2, K_0)$-definable  path
  defined over $F_2$, 
  with 
  $f(\infty)$ a simple point $p_0$ of $\std{V}_\Oo$.  
  Then:
  \begin{enumerate}
\item For all $t$, $f(t) \in \std{V}_\Oo$.
\item 
  We have $\res_{21*} ( f(t)) = \res_{21*}(p_0)$ for positive $t \in \G_{20} \m \G_{10}$.    
\item The restriction
of 
$\res_{21*} \circ f $ to
$[0,\infty] \nsubset \G_{10,\infty}$ is a continuous
 $(K_1,K_0)$-definable path
$[0,\infty] \nsubset \G_{10,\infty}  \to \std{V}_1$.
\end{enumerate}
\end{lem}

\prf  Using base change if necessary and  \lemref{affineimage} we may assume $V\nsubseteq \Aa^n$ is affine.
So $f:  [0,\infty] \nsubset \G_{20,\infty}  \to \std{\Aa^n}_{20}$ and we may assume $V=\Aa^n$.

To prove (1) and (2),  by using
the projections to the coordinates, one reduces to the case $V=\Aa^1$.
Let $\rho(t)=v(f(t)-p_0)$. Then  $\rho$ is a continuous function $[0,\infty]
\to \Gamma_{\infty}$, which is $\tcb{F_2}$-definable (in $(K_2,K_0)$), and sends
$\infty$ to $\infty$. If $\rho$ is constant, there is nothing to prove, since $f $ is
constant, so suppose not. As $\Gamma$ is stably embedded, it follows that
there is $\alpha\in \Gamma_{20}(\tcb{F_2})\nsubset \Gamma_{10}$ such
that for all $ t\in [0,+\infty], \alpha \leq \rho(t)$.
Hence, if $t\in [0,+\infty]_{20}$, then $v_{20}(f(t)-p_0)\geq \alpha$, which
implies that $f(t)\in \std{O_{21}} $ as desired, and gives (1).
Again, by $\tcb{F_2}$-definability and since $f$ is not constant, for some $\mu>0 $ and
$\beta\in \Gamma_{20}(\tcb{F_2})$, if  $t>\beta$, then $ \rho(t)>\mu t$. Thus, when
$t>\Gamma_{10}$, then $\pi(\rho(t))=0$, i.e., $\res_{21*}(f(t))=\res_{21*}(p_0)$. 
   
(3)  Definability of
the 
restriction
of 
$\res_{21*} \circ f $ to
$[0,\infty] \nsubset \G_{10,\infty}$  follows directly  from \lemref{2v1}.   
   For continuity, note that if  $h$ is  a polynomial on $V=\Aa^n$, 
   over $K_1$ and if $H$ is a polynomial over $\Oo_{21}$ lifting 
   $h$, then $v_{20} (H(a)) = v_{10} (h( \res(a)))$.
   It follows that for $t \not= \infty$ in $[0,\infty] \nsubset \G_{10,\infty}$
   continuity of $f$ at $t$ implies  continuity of $\res_{21*} \circ f$.
   
 In fact since $(\res_{21*} \circ f(t))_{*} h $ factors through $\pi_{10}(t)$ as we have
 shown in (2), the argument in (3) shows continuity   at $\infty$ too.  
 To see this directly, one may again consider a polynomial
 $h$ on  $V=\Aa^n$
   over $K_1$ and a lift $H$ over $\Oo_{21}$, 
 and also lift an open set containing $\res_{21*}(p_0)$ to one defined over 
 a subfield $F_{2}'$ contained in $\Oo_{21}$.  The inverse image 
    contains an interval $(\g,\infty)$, and since $\g$ is definable over $F_{2}'$ we necessarily have  $\g \in \G_{10}$.  The pushforward by $\pi_{10}$ of $(\g,\infty)$
    contains an open neighborhood of $\infty$.    \eprf

  \begin{lem} \label{2v21}  Let $f: V \to V'$ be a morphism of varieties defined over $F_2$.  
 Then $f$ induces $f_{20}: \std{V}_{20}\to \std{V}'_{20}$ and also
 $f_{10}:   \tcb{\std{V}_{10}\to \std{V}'_{10}}$.   We have $\res_{21*} \circ f_{20} = f_{10} \circ \res_{21*}$
 on $\std{V}_{\Oo}$. \end{lem} 

\prf  In fact $f_{20}, f_{10}$ are just induced from
restriction of the morphism
$f \tensor_{F_2} \Oo_{21}: V \times _{F_2} \spec \, \Oo_{21} \to V' \times _{F_2} \spec \, \Oo_{21}$, to the general and special fiber respectively, and the statement is clear. \eprf

\begin{lem} \label{cutoff-divisor}    Let $U$ be a projective variety over a valued field, $D$ a divisor.  Let $m$ 
be a metric on $U$, cf. \textup{\lemref{metric}}.  Then the function $u \mapsto \rho(u,D) = \sup \{ m (u,d): d \in D \}$
is v+g-continuous on $U$. \end{lem}

\prf  By \lemref{max}, the supremum is attained.  Let $\rho(u) = \rho(u,D)$. It is clearly  v-continuous. Indeed,   
 if $\rho(u,D) = \alpha  \in \Gamma$, then
$\rho(u',D) = \alpha$ for any $u'$ with $m(u,u') > \alpha$.  If $\rho(u,D) = \infty$
then $\rho(u',D) > \alpha$ for any $u'$ with $m(u,u') > \alpha$.  Let us show g-continuity  by using the criterion in
\propref{gcriterion}.  Let $(K_2,K_1,K_0)$ and $F$ be as in that criterion.  Let $u \in U(K_2)$.
We have to show that $\rho_{21}(u) = (\pi \circ  \rho_{20}) (u)$.   Say $\rho_{20}(u) = m(u,d)$ with $d \in D(K_2)$.
Then $m_{21}(u,d) = \pi ( m(u,d))$ by g-continuity of $m$.  Let $\alpha = \pi(m(u,d))$ and suppose for contradiction
that $\rho_{21}(u) \neq \alpha$.  Then $m_{21}(u,d') > \alpha$ for some $d'$.  We have 
again $m_{21}(u,d') = \pi(m_{20}(u,d'))$ so $m_{20}(u,d') > m_{20}(u,d)$, a contradiction.
\eprf

 \begin{rem}
In the proof of \lemref{cutoff-divisor}, semi-continuity can be seen directly as follows. 
Indeed,   $\rho \inv(\infty) = D$ which is g-clopen.  It remains to show  
 $\{u:  \rho(u,D) \geq \alpha \}$ and 
 $\{u:  \rho(u,D) \leq \alpha \}$ are g-closed.  Now 
  $\ \rho(u,D) \geq \alpha $ if and only if $(\exists  y \in D)(\rho(u,y) \geq \alpha )$; this is the projection 
  of a v+g-closed subset of $U$, hence v+g-closed.  The remaining inequality seems less
 obvious without the criterion, which serves in effect as a topological refinement of quantifier elimination. 
 \end{rem}

  \begin{lem} \label{homotopy-cutoff}\tcb{Let $U$ be an algebraic variety over a valued field or a definable
  subset of such an algebraic variey.} Let $h: I \times \std{U} \to \std{U}$ \tcb{\textup{(}resp.  $h: I \times U \to \std{U}$\textup{)}} be a homotopy.  
Let $\g: \std{U} \to I$  be a   definable continuous function \tcb{\textup{(}resp.  $\g: U \to I$  be a   definable v+g-continuous function\textup{)}}.  Let $h[\g]$ be the cut-off, defined by 
$h[\g](t, u) = h(\max(t,\g(u)), u)$.  Then $h[\g]$ is a homotopy.  Also, if $h$ satisfies $(*)$ of \textup{\ref{ph}}, then so does $h[\g]$. \nomenclature{$h[\g]$}{cut-off of $h$}

 %
 %
\end{lem}

\prf Clear.   
  \eprf

 %

\tcb{Let $U$ be a quasi-projective variety, $Z$ a definable subset of $U$, $f: Z \to \G$ a definable function.  We say $f$  is {\em locally bounded} \index{locally bounded}
if every point $p \in Z$  has a neighborhood, in the valuation topology, on which $f$ is bounded.}
Say $f$ is {\em $U$-locally bounded  } if every point $p \in U$  has a neighborhood $O$ in the valuation topology, with  $f|O$ bounded.  Note that when $Z$ is v-closed, these two notions coincide.

\begin{lem}\label{locbounded}
\tcb{Let $U$ be a quasi-projective variety over a valued field, $Z$ a  $U$-definable subset of $U$, $f: Z \to \G$ a definable function. Then  $f$ is  $U$-locally bounded
on $Z$ if and only if for any bounded definable subset $W$ of $U$,  $f|(W \meet Z)$ is bounded.}
\end{lem}

\prf \tcb{Assume $U$, $Z$ and $f$ are defined over $K$ with $K \models \ACVF$. We may also assume $U$ is affine.
It is enough to prove that if  $f$ is locally bounded on $U$, then for every v-closed bounded $K$-definable subset $W$ of $Z$,
$f (W)$ is bounded. Suppose this does not hold. Then there would exist such a $W$ such that, for some elementary extension $K^* \geq K$, there exists $a \in W (K^*)$ with
$f (a) > \G (K)$.
Consider the valuation ring  \[R = \{x \in K^*: \exists b \in K \val (x) \geq \val (b)\}.\]
Since $W$ is bounded, $a \in W (R)$. The residue field $K'$ is an elementary extension of $K$.
Denote by $\pi : R \to K'$ the canonical projection and set $b = \pi (a)$.
Since $W$ is v-closed, $b \in W (K')$   by \lemref{lem:v-closed}.
We claim that $f$ is not locally bounded at $b$.
Otherwise, there would exist $\gamma, \delta \in \G (K')$ such that, denoting by $B_{\gamma} (b)$ the open polydisc of polyradius $(\gamma, \ldots, \gamma)$ around
$b$, for every $y$ in $Z \cap B_{\gamma} (b)$, $f (y) \leq \delta$. After increasing $\gamma$ and $\delta$ we may assume they belong to $\G (K)$.
Now consider a Hahn field extension
$L = K' \llp t ^{\mathbb{Q}}\rrp$ with $\val (t) > \G (K)$.
Let $k'$ denote the residue field of $K'$.
Since $(K^*, K', k')$ and $(L, K', k')$ are models of $\ACV2F$ with the same characteristics, they are elementary equivalent with parameters in $K'$ by \lemref{2v1}. It follows 
there exists $a' \in W (L)$ such that $a' \in B_{\gamma} (b)$ and
$f (a') > \delta$, leading to a contradiction with the definition of $\delta$ and $\gamma$, since $L$ is an elementary extension of $K'$.}
\eprf

\begin{lem} \label{majorize} \tcb{Let $V$ be a  projective variety over some valued field $F$, $V'$ a Zariski locally closed subset, $U$ a v-closed 
definable subset of $V'$, $f: U \to \G$ be an $F$-definable function.  Assume $f$ is locally bounded on $U$.   Then there exists a v+g-continuous $F$-definable
function $G: V \to \G_\infty$  such that $f(x) \leq G(x) \in \G$ for $x \in U$.}
\end{lem}  
 
\prf  \tcb{
By embedding $V$ as a  closed subset of projective space, we can find a v+g-continuous function
$g: V \to [0,\infty]$ (distance to the boundary),  such that $g$ is finite on $V'$ and for $\a \in \Gamma$, 
 \[V_\alpha = \{x \in V: g(x) \leq \alpha \}\]
 is a v+g-closed and bounded subset of $V'$.   Let $U_\alpha = V_\alpha \meet U$.  Then
  $f$ is bounded on $U_\alpha$ by \lemref{locbounded}; let $f_1(\alpha) $ be the least upper bound.   Since $f_1$ is 
a piecewise affine function, one can find $m \in \Nn $  and $c_0 \in \G$ such that $f_1(\alpha) \leq  m \alpha + c_0$
for all $\alpha \geq 0$ and the function $G(x) = m g(x) + c_0$ does the job.}
\eprf

\section{Continuity on relative $\Pp^1$}   \label{ss9.2}
We fix three points $0$, $1$, $\infty$ in $\Pp^1$. In particular,
the notion of a ball and the standard homotopy are well-defined, cf. \lemref{metric}, \thmref{ret3}.
Let $U_0$ be a normal variety and set $E_0 = U_0 \times \Pp^1$. In practice, $U_0$ will be a dense open subset of $U=\Pp^{n-1}$.
Let $D$ be a divisor on $E_0$  containing the divisor at $\infty$ at each fiber.

Write $D = D' \cup D''$, with $D'$ finite over $U_0$ and
$D''$ the preimage of a closed divisor $Z$ in $U_0$.
Set $U'_0 = U_0 \m Z$ and $E'_0 = E_0 \m D''$.
Let  $\psi_{D'} : [0,\infty] \times E'_0 \to \std{E'_0/U'_0}$ be the   standard homotopy  with stopping time defined by \tcb{$D'$} at each fiber, as defined above \lemref{2v12prel}.  
 We extend $\psi_{D'}$ to a map $\psi_D : [0,\infty] \times E_0 \to \std{E_0/U_0} \nsubset \std{E_0}$
 by $\psi_D(t,x)=x$ for $x \in D''$.  


\begin{lem}    \label{P1-homotopy} Assume $D$ is finite over $U_0$ \tcb{\textup{(}thus $D''$ is empty\textup{)}}. Then the pro-definable map
 $\psi_{D} : [0,\infty] \times E_0 \to \std{E_0/U_0}$ is  v+g-continuous.   \end{lem}

\prf    
Thanks to the g-criterion in \propref{gcriterion}, one deduces from \lemref{2v12prel} and \lemref{2v13rel},  
that $\psi_{D}$ is g-continuous.

We clearly have
v-continuity for the basic homotopy on $\Pp^1$, applied fiberwise on $\Pp^1 \times U_0$.   
\tcb{Let 
$\varrho_D^F : \Pp^1 \times U_0 \to \G_{\infty}$
be the fiberwise  distance to $D$: 
 $\varrho_D^F (y, u)$ is the maximum of all $d (y, z)$ with $(z, u) \in D$, with $d$ the metric on $\Pp^1$.
 Let us check 
$\varrho_D^F$ is v-continuous.
There is no harm in assuming $U_0$ is projective. 
Fix  $(y, u) \in \Pp^1 \times U_0$ and let $\alpha = \varrho_D^F (y, u)$.
Fix $\varepsilon > 0$ in $\G$ and set $W_{\varepsilon} = \{x \in \Pp^1 : d (y, x) \geq \alpha + \varepsilon\}$.
Fix a metric on $\Pp^1 \times U_0$ and
consider as in \lemref{cutoff-divisor} the distance function to $D$. By 
 \lemref{cutoff-divisor} it is v+g-continuous on
$\Pp^1 \times U_0$. Thus, by \lemref{max}, on the bounded v+g-closed definable set
$W_{\varepsilon} \times \{u\}$  its  maximum is attained. Since $D \cap (W_{\varepsilon} \times \{u\} )= \varnothing$,
this maximum is finite. It follows that on some definable v-open set containing
$(y, u)$, $\varrho_D^F  \leq \alpha + \varepsilon$,
proving upper semi-continuity. Lower semi-continuity 
follows from the fact that  the morphism $\std{D}  \to \std{U_0}$ is an open map
by \corref{finite-open}, since 
$D$ has pure codimension $1$ in $E_0$.}
Thus, by
\lemref{homotopy-cutoff}, $\psi_{D}$ is  v-continuous  on $[0,\infty] \times E_0$.  
 \eprf

\begin{lem}  \label{preserve-xi}
 Let $\xi: \Pp^1 \times U \to \G_{\infty}$ be a definable map, with $U$ an algebraic variety over a valued field.  Then there exists a divisor $D_\xi$ on $\tcb{\Pp^1 \times} U$ 
such that, \tcb{for any divisor $D$ containing  $D_\xi$},  the standard homotopy  with stopping time defined by $D$ 
preserves $\xi$.  
\end{lem}

\prf \tcb{If $U$ is not affine, there exists a divisor $D_0$ in $U$ whose complement is affine. By making $\Pp^1 \times D_0$ a component of $D_\xi$, we reduce to the 
case that $U$ is affine. Write $\Pp^1 = \Aa^1 \union \{\infty\}$;
by adding $\{\infty\} \times U$ to $D_\xi$ we can ensure that $\xi$ is preserved there, 
and so it suffices to preserve $\xi | \Aa^1 \times U$. 
Since $\xi | \Aa^1 \times U$ factorizes through a finite number of
functions of the form $\val (g)$, with $g$ regular on $\Aa^1 \times U$, we may assume 
 $\xi | \Aa^1 \times U$ is actually of the form 
$\xi(u) = \val (g)$ with  $g$ regular on $\Aa^1 \times U$.}
Write $g=g(x,u)$, so for fixed
$u \in U$ we have a polynomial $g(x,u)$; let $D_\xi$ include the divisor of zeroes of $g$.
Now it suffices to see for each fiber $\Pp^1 \times \{u\}$ separately that the standard homotopy
with stopping time defined by  a divisor containing the roots of $g$ must preserve $\val (g)$.  This is clear
since this standard homotopy fixes any ball containing a root of $g$; while on a ball containing
no root of $g$, $\val (g)$ is constant. 
\eprf

\begin{lem}\label{nocite}
Let $f: W \to U$ be a generically finite morphism of varieties over a valued field $F$, with $U$ a normal variety, and $\xi: W \to \G_{\infty}$ an $F$-definable map.  
Then there exists a divisor $D$ on $U$ and $F$-definable maps $\xi_1,\ldots,\xi_n: U \to \G_{\infty}$ such that 
any homotopy $I \times W \to \std{W}$ lifting a homotopy of $I \times U  \to \std{U}$ fixing pointwise
$D$ and the levels of the functions $\xi_i$ also preserves $\xi$.  
\end{lem}

\prf  There exists a divisor $D_0$ of $U$ such that $f$ is finite above the complement of $D_0$, and such that $U \m D_0$ is affine.  By making $D_0$ a component of $D$, we reduce to the 
case that $U$ is affine, and $f$ is finite. 
So $W$ is also affine, and $\xi$ factorizes through
functions of the form $\val (g)$, with g regular. We may thus assume  $\xi$ is of this form and in particular that it is  v+g-continuous,
so that it induces a continuous function on $\std{W}$.  
Let $\xi_i(u)$, $i=1,\ldots,n$, list the values
of $\xi$ on $f \inv(u)$ and let also  $\xi_{n + 1}$ be the function \tcb{given by the valuation of} the characteristic function of $D$.
Let $h$ be a  homotopy of $W$  lifting a homotopy of $U$ fixing $D$ and the levels of the $\xi_i$.  Then for fixed $w \in W$, $\xi ( h(t,w))$ can only take  finitely many
values as $t$ varies.  On the other hand $t \mapsto \xi (h(t,w))$ is continuous, so it must be constant.
\eprf

\section{The inflation homotopy}  \label{ss9.3}

\begin{lem} \label{separate} Let $V$ be a quasi-projective variety over a valued field $F$ and let  $W$  be a closed and bounded $F$-pro-definable subset of $\std{V}$. 
Let $D$ and $D'$ be closed $F$-subvarieties of $V$, and suppose $W \meet \std{D'} \nsubseteq \std{D}$.
Then there exists a v+g-closed, bounded $F$-definable subset $Z$ of $V$ with $Z \meet D' \nsubseteq D$, and 
$W \nsubseteq \std{Z}$.
\end{lem}

\prf We may assume $V$ is affine. Indeed,  
we may assume $V=\Pp^n$;  then find finitely many affine open $V_i \nsubset V$ and closed bounded definable subsets $B_i \nsubset V_i$ such that $\tcb{W} = \union_i B_i$; given $Z_i$ solving the problem for $V_i$, set $Z=\union_i (B_i \meet Z_i)$.

Choose a finite generating family $(f_i)$ of
the ideal of regular functions vanishing on 
$D$ and set
$d (x, D) =  \inf \val (f_i (x))$ for $x$ in $V$.
Similarly, fixing 
a finite generating family of
the ideal of regular functions vanishing on 
$D'$, one defines a distance function
$d (x, D')$ to $D'$.
Note that the functions
$d (x, D)$ and $d (x, D')$ may be extended 
to $x \in \std{V}$.
 
For $\alpha \in \G$, 
let  $V_\alpha$ be the set of points $x$ of $V$ with $d (x, D) \leq \alpha$. Let $W_\alpha = W \meet \std{V_\alpha}$.  Then $W_\alpha \meet \std{D'} = \varnothing$.  So $d(x,D') \in \Gamma$ for $x \in W_\alpha$.  By \lemref{max} there exists $\delta(\alpha) \in \Gamma$ such that $d(x,D') \leq \delta(\alpha)$ for 
$x \in W_\alpha$.  We may take $\delta : \G \to \G$ to be a continuous
non-decreasing definable function.     Since any such function $\G \to \G$ extends to a 
 continuous function $\G_\infty \to \G_\infty$, we may extend
 $\delta$ to a 
  continuous function $\delta : \G_\infty \to \G_\infty$.
  Also, since any such function is bounded by
  a continuous function with value $\infty$ at $\infty$ we may assume
  $\delta (\infty) = \infty$.
 Let \[Z_1 = \{x \in V:  d(x,D') \leq \delta( d(x,D)) \}.\]
 This is a v+g-closed set.
 Let $c$ be a realization of $p \in W$.
 We have $c \in Z_1$ and
 $Z_1  \meet D' \nsubseteq D$.
 Since, by \lemref{bounded}, $W$ is contained in
 $\std{Z_2}$ with $Z_2$ a bounded v+g-closed definable subset of
 $V$, we may set $Z = Z_1 \meet Z_2$.
\eprf


\begin{lem}  \label{inflation}Let $D$ be a closed subvariety of a projective variety $V$ over a valued field $F$, and assume there exists an \'etale \index{inflation homotopy}
map $e:  V \m D \to U$,  with $U$ a \tcb{Zariski} open subset of $\Aa^n$.  
Then there exists an $F$-definable 
homotopy
 $H: [0,  \infty] \times \std{V}\to \std{V}$  fixing $\std{D}$ \textup{(}that is, such that $H (t, d) = d$ for $t \in [0,  \infty]$ and $d \in \std{D}$\textup{)}, 
with image $\bZ = H (0, \std{V})$, such that for any   subvariety $D'$ of $V$ of dimension $< \dim(V)$ we have $\bZ \meet \std{D'} \nsubseteq \std{D}$.  
Moreover given a finite family of $F$-definable  v-continuous 
 functions $\xi_i: V \m D \to \G$, $i \in I$,  one can choose the homotopy such that  
 the levels of the $\xi_i$ are preserved. The same statement remains true if instead of being $F$-definable, the functions $\xi_i$ are only assumed
 to be $F'$-definable, with $F'$ a finite Galois extension of $F$, and the functions $\xi_i$ are permuted by the action of
 $\mathrm{Gal} (F' / F)$.
If a finite group $G$ acts on $V$ over $U$, inducing a continuous action on $\std{V}$ and  leaving $D$ and the fibers of $e$ invariant, then $H$ may be chosen to be $G$-equivariant.
  \end{lem} 

\prf    Let $I = [0, \infty]$ and
let  $h_0: I \times  \Aa^n \to \std{\Aa^n}$  be the standard homotopy sending
$(t, x)$ to the generic type of the closed polydisc of polyradius
$(t, \ldots, t)$ around $x$. Denote by $H_0 : I \times \std{\Aa^n} \to \std{\Aa^n}$
its canonical extension (cf. \lemref{hbasic}). 
Note the following fundamental inflation property of $H_0$: if $W$ is closed subvariety of $\Aa^n$ of dimension 
$< n$,
then, for any $(t, x)$ in $I \times \std{\Aa^n} $, if $t \not= \infty$, then $H_0 (t, x) \notin \std{W}$.

By \lemref{p4} and \lemref{unramified-unique}, for each 
$u \in U$ there exists $\g_0(u) \in \Gamma$ such that $h_0(t, u)$ lifts uniquely
to $V \m D$ beginning with any $v \in e \inv(u)$, up to $\g_0(u)$.    By \lemref{majorize}
we can take $\g_0$ to be v+g-continuous.  For $t \geq \g_0(u)$, let $h_1(t, v)$ be the unique
continuous lift.

Since $\xi_i$ is v-continuous outside $D$, $\xi_i \inv (\xi_i(v))$ contains a v-neighborhood of $v$.
So for some $\g_1(u) \geq \g_0(u)$, for all $t \geq \g_1(u)$ we have $\xi_i( h_1(t, v)) = \xi_i(v)$.  
We may take $\g_1 (u) = \min \{\alpha \in \Gamma_{\geq 0} : 
\xi_i (h_1 (t, v)) = \xi_i (v), \forall t \in [\alpha, \infty), \forall v \in e^{-1} (u), \forall i\}$, which is locally bounded
and $F$-definable, not only when the functions $\xi_i$ are assumed to be $F$-definable, but also 
when they are only assumed
 to be $F'$-definable, with $F'$ a finite Galois extension of $F$, and permuted by the action of
 $\mathrm{Gal} (F' / F)$.
Thus,  we may use \lemref{majorize} again to replace $\g_1$ by a v+g-continuous $F$-definable function.
By \lemref{homotopy-cutoff}, the cut-off
$h_0[\g_1]$  is \tcb{v+g-}continuous,
and  by \lemref{homotopy-lift}, $h_1[\g_1 \circ e]$  is \tcb{v+g-}continuous on $V \m D$.  However we would like
to fix $D$ pointwise and have continuity on $D$.

Let $m$ be a metric on $V$, as provided by \lemref{metric}. Given $v \in V$ let 
$\rho(v) = \sup \{m(d,v): d \in D \}$.  By \lemref{max} we have $\rho: V \m D \to \G$.   
Let $\g_2: \tcb{U} \to \G$,
$\g_2 \tcb{\geq} \g_1$, such that for $t \geq \g_2 (u)$ we have $\tcb{m}( h_1(t, v),v) \tcb{\geq} \rho(v)$
for each $v$ with $e(v) = u$.  
By \lemref{majorize} we can take $\g_2$ to be v+g-continuous.

  Let  $H$ the canonical extension of $ h_1[\g_2 \circ e]$
  to $\std{V \m D}\times I$ provided by \lemref{hbasic}.  We
  extend $H$ to
  $\std{V}\times I$ by setting $H (t, x) = x$ for $x \in \std{D}$.
  We want to show that $H$ is continuous on $\std{V}$. Since we already know it is
  continuous at each point of the open set $\std{(V \m D)} \times I$, it is enough to prove $H$ is   continuous at each point of $\std{D} \times I$.    

 Let $d \in \std{D}$, $t \in I$.
Then $H(t, d)=d$.  Let $\bG$ be an open neighborhood of $d$.     One may assume $\bG$  to have the form  $\{x \in G_0: \val (r(x)) \in J \}$, with $J$  open in $\G_\infty$, and $r$ a regular function 
 on a Zariski open neighborhood $G_0$ of $d$ (which is just a Zariski open subset of $V$ supporting $d$).  So $\bG = \std{G}$ where $G$ is a v+g-open subset of $V$.  
 
     We have to find an open neighborhood
$W$ of $(t, d)$ such that $H(W) \nsubseteq \bG$.   We may take $W \nsubseteq \bG \times \G_{\infty}$,
so we have $H(W \meet \std{D})   \nsubseteq \bG$.   Since the simple points of $W \m \std{D}$ are
 dense in $W \m \std{D}$ and by construction of the canonical extensions in  \ref{canon-ext},
 it suffices to show that for some neighborhood $W$, the simple points are mapped
 to $\std{G}$.    
 
  View $d$ as a type (defined over \tcb{some model} $M_0$); if $z \models d | M_0$, then for some $\varepsilon \in \G$, $H( B(z; m,\varepsilon)) \nsubseteq \bG$. 
  Fix $\varepsilon$, independently of $z$.  The set
  \[W_0 = \{v \in V:  B(v;m,\varepsilon) \nsubseteq G\}\]  
   is v+g-open since its 
complement is 
\[\{v \in V: (\exists y) m(x,y) \leq \varepsilon \wedge y \in (V \m G) \}.\]  Now the projection of a (bounded) 
  v+g-closed set  is also v+g-closed.  

 If there is no neighborhood $W$ as desired,    there exist \tcb{a net $(t_i, v_i)$ with $t_i \to t$ and $v_i \in V \m D$ simple  points with   $v_i \to d$}, such that
$H(t_i, v_i) \notin G$.

\tcb{
Since \[H(t_i, v_i) = h_1(\max (\g_2(e(v_i)), t_i), v_i),\] we have $m(H(t_i, v_i), v_i) \geq \rho(v_i)$.
As $\rho(v_i) \to \rho(d) =m(d,D) = \infty$, it follows that $m(H(t_i, v_i), v_i) \to \infty$.}  
So, for large $i$,  we have $H(t_i, v_i) \in \std{B(v_i; m,\varepsilon)}$, and  
also $v_i  \in W_0$.  So $B(v_i,m,\varepsilon) \nsubseteq G$, hence $H(t_i, v_i) \in \std{G} = \bG$, a contradiction.  This shows that $H$ is continuous. 

It remains to prove that if $\bZ = H (0, \std{V})$, then,  for any   subvariety $D'$ of $V$ of dimension $< \dim(V)$, we have $\bZ \meet \std{D'} \nsubseteq \std{D}$.  
This follows from the inflation property of $H_0$ stated at the beginning, applied to 
$e (D' \meet (V \m D))$.

The  statement on the group action follows from the uniqueness of the continuous lift.  
\eprf

 \begin{rem}\label{reminflation}  \lemref{inflation} remains true if one supposes only that $D$ contains the singular points of $V$. Indeed,
one can find divisors $D_i$ with $D = \meet_i D_i$, and \'etale morphisms $h_i: V \m D_i \to \Aa^n$,
and iterate the lemma to obtain successively $\bZ \meet \std{D'} \nsubseteq \std{D_1} \meet \ldots \meet \std{D_i}$.
In particular, when $V$ is smooth, \lemref{inflation} is valid for $D = \varnothing$.
\end{rem}

\section{Connectedness and the  Zariski topology}\label{ss9.4}

Let $V$ be an algebraic variety over some valued field. 
We say a
strict pro-definable subset $Z$ of $\std{V}$ is {\em definably connected} if \index{definably connected}
it contains no clopen strict pro-definable subsets other than $\varnothing$ and $Z$.  
We say that 
$Z$ is {\em definably path connected} if for any two points $a$ and $b$ of  $Z$ there exists \index{definably path connected}
a definable path in $Z$ connecting $a$ and $b$.
Clearly definable
path connectedness implies definable connectedness. When $V$ is quasi-projective and
$Z = \std{X}$ with $X$ a definable subset of $V$, the reverse implication will eventually follow from \thmref{1}.

If $X$ is a definable subset of $V$, $\std{X}$ is definably connected if and only if $X$ contains no 
v+g-clopen definable subsets, other than $X$ and $\varnothing$. Indeed, 
if $U$ is a clopen strict pro-definable subset of $\std{X}$,
the set $U \cap X$ of simple points of $U$ is a  v+g-clopen definable subset of $X$, and $U$
 is the closure of 
 $U \cap X$. 
When $X$ is a definable subset of $V$, we shall say $\std{X}$ has a finite number of connected components 
if $X$ may be written as a finite disjoint union of 
v+g-clopen definable subsets $U_i$ with each $\std{U_i}$ definably connected.
The $\std{U_i}$ are called {\em connected components} of $\std{X}$. \index{connected component}

\begin{lem}\label{codimconnected}Let $V$  be a smooth quasi-projective variety over a valued field and let $Z$ be a nowhere dense Zariski closed subset of $V$. 
Then $\std{V}$ has a finite number of connected components if and only if
$\std{V \m Z}$ has a finite number of connected components. Furthermore,
if $\std{V}$ is a finite disjoint union of connected components
$\std{U_i}$ then the $\std{U_i} \m \std{Z}$  are the connected components of
$\std{V \m Z}$.
\end{lem}

\prf By \remref{reminflation}, there exists a homotopy
$H : I \times \std{V}
\to \std{V}$ such that its final image $\Sigma$ is contained in
$\std{V \m Z}$. Also, by construction of $H$, the simple points of $V \m Z$ move within $\std{V \m Z}$, 
and so  $H$ leaves $\std{V \m Z}$ invariant.
Thus, we have a continuous morphism of strict pro-definable spaces 
$\varrho : \std{V} \to \Sigma$.
If $V$ is a finite disjoint union of 
v+g-clopen definable subsets $U_i$ with each $\std{U_i}$ definably connected,
note that each $\std{U_i}$ is invariant by the homotopy $H$. Thus,
$\varrho  (\std{U_i}) =  \Sigma \cap \std{U_i}$ is definably connected.
Since $\Sigma \cap \std{U_i} = \Sigma \cap (\std{U_i} \m  \std{Z})$
and any simple point of
$U_i \m Z$ is connected via $H$ within $\std{U_i \m Z}$ to
$\Sigma \cap \std{U_i}$ it follows that $\std{U_i} \m \std{Z}$ is definably connected.
For the reverse implication, assume
 $V \m Z$ is a finite disjoint union of 
v+g-clopen definable subsets $V_i$ with each $\std{V_i}$ definably connected.
Then $\varrho  (\std{V_i}) =  \Sigma \cap \std{V_i}$ is definably connected.
Let $U_i$ denote the set of simple points in $\varrho^{-1}(\Sigma \cap \std{V_i})$.
Then $\std{U_i}$ is definably connected.
\eprf

\begin{thm}\label{connectedness}  Let $V$ be an \tcb{algebraic}  variety over a valued field $F$. Assume $V$ is geometrically connected for the Zariski topology.  Then $\std{V}$ is definably connected.     \end{thm} 

\prf  We may assume $F$ is algebraically closed and $V$ is  irreducible.  It follows from \tcb{the version of Bertini's theorem given in} \cite{mumford} p.~56,
that any two points of $V$ are contained in an  irreducible curve $C$ on $\tcb{V}$.
So, since simple points are dense, the lemma reduces to the case of irreducible curves, and by normalization, to the case of normal irreducible curves $C$.  
As in the beginning of \thmref{f9}, one may thus assume $C$ is smooth and irreducible.
By \lemref{codimconnected} one may assume that  $C$ is also projective. 
The case of genus $0$ is clear using the standard homotopies of $\Pp^1$.  So assume
$C$ has genus $g>0$.  
By \thmref{ret3} there is a retraction
 $\varrho : \std{C} \to \Upsilon$ with  $\Upsilon$ \tcb{an iso-definable} $\G$-internal subset.
 It follows from
\thmref{G-embed-1c}  that $\Upsilon$ is a finite disjoint union of 
 connected iso-definable $\G$-internal subsets $\Upsilon_i$. Denote by
 $C_i$  the set of simple points in $C$ mapping to $\Upsilon_i$.
 Each $C_i$ is a 
 v+g-clopen definable subset of $C$ and
  $\std{C_i}$ is definably connected, thus $\std{C}$ has a 
finite number of connected components.
Assume this number is $> 1$. Then $\std{C^g} /  \Sym(g)$ has also a
finite  number $> 1$ of connected components, since $\std{C^g}$ may be written as a disjoint union
of the definably connected sets $\std{C_{i_1} \times \cdots \times C_{i_g}}$.

 Let $J$ be the Jacobian variety of $C$.
There exist proper subvarieties $W$ of $C^g$ and $V$ of $J$, with $W$ invariant under $ \Sym(g)$,
and a biregular isomorphism of varieties $(C^g \m W)/ \Sym(g) \to J \m V$.  
By \lemref{codimconnected},
$\std{(C^g \m W)/ \Sym(g)}$ has a
finite  number $> 1$ of connected components, hence also $\std{J \m V}$.
By \lemref{codimconnected} again, $\std{J}$ would have a finite  number $> 1$ of connected components.
The group of simple points of $J$ acts by translation on $\std{J}$, homeomorphically, and so acts also on 
the set of connected components. Since it is a divisible group,  the action must be trivial.
On the other hand, it is transitive on simple points, which are dense, hence on
connected components.  This leads to a contradiction, hence $\std{C}$ is connected, which finishes the proof.
\eprf

\begin{lem}  \label{infty}  Let $V$ be an algebraic variety over a valued field $F$ and let  $f: V \to \G_\infty$ be a v+g-continuous $F$-definable function.  Then $f \inv (\infty)$ is a  subvariety of $V$.     \end{lem}

\prf   \tcb{Note that, for constructible sets, the
Zariski closure and the v-closure coincide.
Hence,} since $f \inv (\infty)$ is v-closed, it suffices to show that it is constructible.  We may assume $F$ is algebraically closed.
By noetherian induction we may
assume $f \inv(\infty) \meet W$ is a subvariety of $W$, for any proper subvariety $W$ of $V$.
So it suffices to show that $f \inv(\infty) \meet V'$ is an algebraic variety, for some Zariski open $V' \nsubseteq V$.
In particular we may assume $V$ is affine, smooth and irreducible.  Since any definable set is v-open away from some proper subvariety,
we may also assume that   $f \inv(\infty)  $ is v-open.
 On the other hand $f \inv(\infty) $ is v-closed.
The point $\infty$ is an isolated
point in the g-topology, so  $f \inv(\infty)$ is  g-closed and g-open.  By \lemref{v+g} it follows that
$\std{f \inv(\infty)}$ is a clopen subset of $\std{V}$.  Since $\std{V}$   is definably connected by \thmref{connectedness}, one deduces  that $f \inv(\infty) = V$ or $f \inv(\infty)= \varnothing$, proving the lemma.  \eprf 

Let $w$ be a finite definable set.
It will be convenient to use the following terminology. 
By a {\em z-closed} subset of  $\Gamma_\infty^w$ we mean one of the \index{z-closed}
form $[x_i=\infty]$,   an intersection of such sets, or a finite union of such intersections.     Note that such sets are not automatically
defined over the given base (but some of them are).  Let $Y \nsubseteq \G_\infty^w$ be a definable set.  A z-closed subset of $Y$ is the intersection with $Y$ of a z-closed subset.  
(If $Y$ is $A$-definable, an $A$-definable z-closed subset of $Y$ can be written as $Y \meet Z$, where $Z$ is z-closed and $A$-definable; this can be done
by taking unions of Galois conjugates.) By a {\em z-irreducible} subset we mean a z-closed subset \index{z-irreducible}
which cannot be written as the union of two proper z-closed subsets.  Any z-closed set can be written
as a finite union of z-closed z-irreducible sets; these will be called {\em z-components}. \index{z-component}
A z-open set   is the complement of a z-closed set $Z$.  A z-open set is {\em dense} if  \index{z-dense}
its complement   does not contain any z-component of $Y$.

Let $Y$ be a definable subset of $\G_\infty^w$.  
Define a {\em Zariski closed} subset of $Y$ (resp. {\em Zariski open}) to be  a clopen definable subset of a z-closed subset of $Y$ (resp. z-open).  \index{Zariski closed}  \index{Zariski open}
By o-minimality,  there are finitely many such clopen subsets, the unions of the definably connected components.   A definable set $X$ thus has only finitely many  
Zariski closed subsets; if $X$ is  connected and z-irreducible, there is a maximal proper one.

\medskip
This has nothing to do with the topology on $\G^n$ generated by translates of  subspaces defined by $\Qq$-linear equations, for which the name 
Zariski would also be natural.  We will use this latter topology little, and will refer to it as the linear Zariski topology on $\G^n$, when required.

\medskip

\lemref{infty} can be strengthened as follows:

\begin{lem} \label{geninfty}Let $V$ be an algebraic variety over a valued field $F$, let $w$ be a finite $F$-definable set
and let  $f: V \to Y \nsubseteq \G_\infty^w$ be a v+g-continuous $F$-definable function.
 Then $f \inv(U)$ is Zariski open \textup{(}resp. closed\textup{)} in $V$, whenever $U$ is Zariski open \textup{(}resp. closed\textup{)} in $Y$.  
\end{lem}

\prf  It suffices to prove \tcr{the statement about closed sets}.  We may assume $F$ is algebraically closed. So $U$ is a clopen subset of $U'$, with $U'$ z-closed.  By \lemref{infty}, $f \inv(U')$ is Zariski closed;
write $f \inv(U') = V_1 \union \ldots \union V_m$ with $V_i$ Zariski irreducible.
It suffices to prove the lemma for $f|V_i$, for each $i$; so we may assume $V_i=V$ is Zariski irreducible.  By \thmref{connectedness},
$f \inv(U) = V$.  \eprf

 Here is a converse:
 
 \begin{lem} \label{z-iso1}  
Let $X \nsubseteq \G_\infty^w$ and let $\beta: X \to \std{V}$ 	be a continuous, pro-definable map.  Let $W$ be a Zariski closed subset of $\std{V}$. 
Then $\beta \inv(W)$ is Zariski closed in $X$.    \end{lem}

\prf   
Let $F_1,\ldots,F_{\ell}$ be the nonempty, proper Zariski closed subsets of $X$.  
Removing from $X$ any $F_i$ with $F_i \nsubseteq \beta \inv(W)$, we may assume no such $F_i$ exist. 
By working separately in each component, we may assume $X$ is connected, and in fact z-irreducible.  
 Moreover by induction on z-dimension, 
we can assume the lemma holds for proper z-closed subsets   of $X$.

\begin{claim} $\beta \inv(W) \meet F_i = \varnothing$ for each $i$.
\end{claim}

\begin{proof}[Proof of the claim]Otherwise, let $P$ be a minimal $F_i$ with nonempty intersection with $\beta \inv(W)$.  
Let $Q$ be the z-closure of $P$; then
$Q \neq X$.  As   
Zariski closed in $Q$ implies Zariski closed in $X$,   $Q\cap \beta\inv(W) = \varnothing$.  (Thanks to Z. Chatzidakis for this argument.)
\end{proof}

Say $ \beta\inv(W)\nsubseteq \Gamma_\infty^{w_1}\times \{\infty\}^{w_2}$
 with $(w_1, w_2)$  a partition of $w$ and $\vert w_1 \vert$ minimal.
Then $\beta \inv(W) \meet (x_i=\infty) = \varnothing$ for $i \in w_1$, i.e. 
 $\beta \inv(W) \nsubseteq \Gamma^{w_1}\times \{\infty\}^{w_2}$.
Projecting homeomorphically to $\G^{w_1}$, we may assume $w_1 = w$ and $X \nsubseteq \G^w$.  However, $W$ is of the form
$\std{F}$ with $F$ g-clopen, so 
$\beta \inv(W)$ is g-clopen.
Since any g-clopen subset of $\G^w_{\infty}$ which is also closed and contained
in $\G^w$ is clopen, it follows that 
$\beta \inv(W)$ is clopen, which  implies
that it is after all Zariski closed in $X$.   \eprf

\begin{cor}  \label{z-iso} Let $\Ups$ be an iso-definable subset of $\std{V}$, $X$ a definable subset of $\G_\infty^w$, and let $\a: \Ups \to X$
be a pro-definable homeomorphism.  Then $\a$ takes the Zariski topology on $\Ups$ to the Zariski topology on $X$.  \end{cor}
\prf Follows from  \lemref{geninfty}  and \lemref{z-iso1}.   \eprf 


\chapter{The main theorem}\label{sec10}

{\small \noindent \textbf{Summary.}
The main theorem is stated in \ref{ss10.1} and
several 
preliminary  reductions are performed  in \ref{ss10.2} that allow us to 
essentially reduce to a curve fibration.
We construct a  relative curve homotopy  in \ref{ss10.3} and 
a liftable base homotopy in \ref{ss10.4}.
In  \ref{ss10.5} a purely combinatorial homotopy is constructed in the $\G$-world.
Finally in \ref{ss10.6}  we end the proof of the main theorem; the homotopy retraction is constructed by concatenating the previous
three homotopies together with an inflation homotopy. The chapter ends with
\ref{ss10.7}  which is devoted to the relative version of the main theorem.
\par\bigskip}

\section{Statement}\label{ss10.1}
 
 \begin{thm}  \label{1} Let $V$ be a quasi-projective variety over a valued field $\tcb{F}$ and let $X$ be a
  definable subset of $V \times \G_\infty^{\ell}$ over some base set $A \nsubset \VF \union \G$, \tcb{with $F = \VF (A)$}.  
Then there exists an $A$-definable deformation retraction 
 $h : I \times \std{X} \to \std{X}$ \tcb{with image} an iso-definable subset $\Upsilon$ definably 
homeomorphic to a definable subset of $ \G_\infty^w$, for some finite $A$-definable set $w$.  

One can furthermore require the following additional properties for $h$ \tcb{to hold simultaneously}:
\begin{enumerate}
\item Given finitely many  $A$-definable functions $\xi_i: X \to \G_{\infty}$, \tcb{with canonical extension $\xi_i: \std{X} \to \G_{\infty}$ as in \textup{\ref{canon-ext}}}, one can choose $h$
to respect the $\xi_i$, i.e. to satisfy $\xi_i (h(t, x) )=\xi_i(x)$ for all $(t, x) \in I \times \std{X}$.   In particular, finitely many  definable subsets $U$ of $X$ can be preserved, in the sense that the homotopy restricts to one of $\std{U}$.
\item Assume given, in addition, a finite algebraic group $G$ acting on $V$ and
\tcb{leaving $X$ globally invariant}.   
 Then the retraction $h$ can be chosen to be equivariant with respect to the $G$-action.
\item Assume $\ell=0$.  The homotopy $h$ is Zariski generalizing, i.e. for any Zariski open subset $U$ of $V$, $\std{U} \meet \tcb{\std{X}}$  is   invariant under $h$.   \index{Zariski generalizing}
 \item The homotopy $h$ satisfies condition $(*)$ of  \textup{\ref{ph}}, i.e.:
$h (e_I, h (t,x)) = h (e_I, x)$ for every $t$ and $x$.
\item The homotopy $h$ restricts to $h^{\#}: I \times \stda{X} \to \stda{X}$, cf.  \textup{Definition \ref{ssd}} and  \textup{\ref{ss6.5}}.
\item One has $h(e_I, X)=\Upsilon$, i.e. $\Upsilon$ is the image of the simple points.
\item Assume $\ell=0$ and $X= V$. Given a finite number of closed irreducible subvarieties $W_i$ of $V$, one can demand that
$\Upsilon  \meet \std{W_i}$ has pure dimension $\dim(W_i)$.
 \end{enumerate}
\end{thm}

 \begin{defn}  \label{defskel}
 \tcb{Let $V$ be a quasi-projective variety, $X$ be a
  definable subset of $V$ over some base set $A \nsubset \VF \union \G$.  
  Let $\Upsilon$ be an $A$-iso-definable subset of $\std{X}$.  
  We call  $\Upsilon$ a {\em skeleton}  of $\std{X}$ \index{skeleton} 
  if it is definably 
homeomorphic to a definable subset of $ \G_\infty^w$, for some finite $A$-definable set $w$, 
there exists an $A$-definable deformation retraction 
 $h : I \times \std{X} \to \std{X}$ 
 with image  
  $\Upsilon$, and in addition (7) holds for each irreducible component $W$ of the Zariski closure of $X$.}
  \end{defn}
  
The last condition may look inelegant, but will allow us to prove that any two skeleta are contained in a third, and more generally that the homotopy  in \thmref{1} can
be required to fix \tcb{pointwise} any given skeleton.  
\tcb{A possible alternate definition could be to replace the last condition by the condition that $\Upsilon$ is contained in
$\stda{V}$.
By \thmref{sgc} any such skeleton is contained in a  skeleton in the sense of  \defref{defskel},
and any skeleton in the sense of   \defref{defskel} lies in $\stda{V}$.}

 \medskip
 
 \begin{rem} \label{1r} 
\leavevmode
 \begin{enumerate}
\item  Without parameters, one cannot expect in general $\Upsilon$ to be definably homeomorphic  to a subset of $\G_\infty^n$, because
of the existence of  \tcb{Berkovich analytifications for which 
the Galois group acts nontrivially on the cohomology,  cf. the earlier observation in \ref{ss6.1}}.
\item   Let 
$\pi: V' \to V$ be a finite \tcb{surjective} morphism \tcb{of $F$-varieties with $V$ normal},  and  $\xi': V' \to \G_{\infty}^m$ \tcb{be an $A$-definable morphism}.  Then, when $X = V$
one can find $h$ as in the theorem lifting to  $h': I \times \std{V'} \to \std{V'}$ respecting $\xi'$.  To see this, let 
$V'' \to V'$ be such that $V'' \to V$ admits  a finite group action $G$, and $V'$ is the quotient variety
of some subgroup. An equivariant homotopy of $\std{V''}$ will induce homotopies on $\std{V'}$ and on $\std{V}$.   \tcb{The continuity of the induced homotopies
follows from \lemref{homotopy-descent}  \tcr{and Remark \ref{rem:homotopy-descent}} and  the iso-definability of their image from \lemref{b3cc}}.
\item  \tcb{In \thmref{1} (1), one 
 may demand that the homotopy preserve  a given $A$-definable map $\xi: X \to \G_\infty^w$ with $w$ a finite $A$-definable set.}
  \tcb{Indeed,} 
let  $\xi': X \to \G_\infty^m$ (where $m=|w|$) be a map 
such that for $v \in X$, $\xi'(v)$ is an $m$-tuple in non-decreasing order, enumerating the underlying set of the $w$-tuple $\xi(v)$.
  There exist definable sets $U_i$
such that $\xi | U_i$ is continuous.   
We can ask that the homotopy $h$ preserves the $U_i$ and $\xi'$.  Then along  the homotopy $h$, $\xi$ is preserved up to a permutation of $w$, hence by continuity it is preserved. 
\item Property (3) \tcb{in \thmref{1}} implies that, for any irreducible component $W$ of $V$,
$\Upsilon \meet \std{W}$ is Zariski dense in $X \meet \std{W}$ in the sense of \ref{zt} and that
$X \meet \std{W}$ is invariant under $h$.
For  the first assertion note that one cannot have
$\Upsilon \meet \std{W} \nsubseteq \std{Z}$, for some proper Zariski closed subset $Z$ of $W$, since
then a point in $W \m Z$ would have its final image in $\std{\tcb{Z}}$.
For the second one, let $W_0$ be the complement in $W$ of the other components.
By (3) $\std{W_0} \meet X$  is   invariant under $h$ and the invariance of 
$\std{\tcb{W}} \meet X$ follows by continuity.
\item  Assume $\ell = 0$. The retraction $\std{X} \to \Upsilon$ can be taken to be {\em definably proper}, i.e. so that the pullback of a definably compact set is definably compact. \index{definably proper}
 Indeed $V$ embeds in some projective
variety $V'$, in an $G$-equivariant way \tcb{as in the beginning of \ref{ss10.2}}.   We can use the theorem to find a homotopy $\std{V'} \to \Upsilon'$, preserving
the data, and also preserving $V' \m X$ and $X$.  The retraction $\std{X} \to \Upsilon$ is just the restriction of $\std{V'} \to \Upsilon'$,
and hence also definably proper.  
 %
  %
\end{enumerate}  \end{rem}
  
It is worth pointing out
  that  the fibers of $\std{X} \to \std{Y}$, over an element $\by \in \std{Y}$, 
for a definable map $X \to Y$, are not in general spaces of the form $\std{U}$.  The fiber $\std{X}_\by$ over an element $\by \in \std{Y}$
does contain a subset $\mathrm{I}X_\by$ accessible in our language, namely $\{ \int_\by g \}$ for $g: Y  \to \std{X/Y}$ a definable \tcb{section}.
But this does not exhaust the fiber.   
Nonetheless, the proof of
\thmref{1} is inductive, using appropriate fibrations.   
 What permits this   is that our homotopy is determined by its restriction to the simple
points, cf. \lemref{hbasic}.   Given relative homotopies of the fibers,   on the simple points of $X$ one obtains a map into $\std{X}$
whose image, over a fiber $\by$, does fall into the ``inductive'' subset $\mathrm{I}X_\by$ mentioned above. 
In addition, under appropriate circumstances, a homotopy of $\std{Y}$ can be extended to a homotopy of  $\std{X}$.
Though the methods can be applied more generally, it is worth pointing out that the homotopy restricts
to a homotopy of   $\stda{X}$; and that the fibers of  \tcb{$\stda{X} \to \stda{Y}$} {\em can} all be obtained as integrals, as above.

\section{Proof of  Theorem \ref{1}:  Preparation}\label{ss10.2}

The theorem reduces easily to the case $\ell = 0$ (for instance, take the projection of $X$ to $V$, and
add $\tcb{\xi'_i}$ describing the fibers, as in the first paragraph of the proof of \thmref{omin-finite}).  We assume $\ell=0$ from now on.

\tcb{It is enough to prove the theorem when $X = V$. Indeed consider the functions $\xi'_i$
on $V$ obtained by extending the functions $\xi_i$ by $0$ on $V \m X$ together with the function given by the valuation of the  characteristic function of $X$. The theorem for $X = V$
equipped with these functions implies the statement  for $X$ and the original functions $\xi_i$. We now assume $X = V$.}

\tcb{Let $G$ be a finite algebraic group acting on $V$. We may embed equivariantly $V$ in a projective, equidimensional variety $W$ with $G$-action of the same dimension.
Indeed, let $\bar V$ be a projective completion of $V$. Embed $V$ diagonally in $V^G$, via
$v \mapsto (gv)_{g \in G}$; this is equivariant with respect to the action of $G$ on $V^G$ via the regular action of $G$ on $G$.  Taking the Zariski closure of the image in $\tcb{\bar{V}^G}$ and  an equidimensional $\bar{V}'$ containing $\bar{V}$ with the same irreducible
components of dimension $\dim V$, and then considering $\union_{h \in G} \bar{V}'$, we get $W$ as required. 
On $W$ we can consider the extensions by 0 of 
 the functions $\xi_i$ together with the functions given by the valuation of the  characteristic functions
 of the lower dimensional components of $V$. It is enough to prove the theorem for  $W$ equipped with these functions
 to have it for $V$ with the functions $\xi_i$.
 Thus, we may assume from now on that $X= V$ is projective and equidimensional.}

At this point we  note that we can take the base $A$ to be a field.  Let $F=\VF(A)$ be the field part.   
Then
$V$ and $G$ are defined over $F$.   
Write $\xi = \xi_\g$ with $\g$ from $\G$.  Let $\xi'(x)$ be the function:
$\g \mapsto \xi_\g(x)$. Clearly if the fibers of $\xi'$ are preserved then so is each $\xi_\g$ (cf. \remref{1r} (4)).  By stable embeddedness
of $\G$, $\xi'$ can be coded by a function into $\G^k$ for some $k$.  And this function is $F$-definable.
Thus all the data can be taken to be defined over $F$, and the theorem over $F$ will imply the general case.

We may assume $F$ is perfect, 
since this does not change the notion of definability over $F$.  

We use induction on $n=\dim(V)$.  
For $n = 0$, take the identity deformation $h(t,x)=x$, $w=V$, and map $a \in w$ to $(0,\ldots,0,\infty,0,\ldots,0)$
with $\infty$ in the $a$-th place.

We start with a hypersurface (that is, a closed subset everywhere of dimension $n - 1$)   $D_0$ of $V$ containing the singular locus $V_{\mathrm{sing}}$.
We assume there exists an \'etale morphism $V \m D_0 \to \Aa^n$, factoring through $V/G$.
Such a $D_0$ exists using generic smoothness, after choosing a separating transcendence basis at the generic point of $V/G$.
We also assume $D_0$ is nonempty of dimension $n - 1$ in each irreducible component of $V$.
Note that the functions $\xi_i$ factor through  v+g-continuous functions into $\G_\infty^m$.   \tcb{Indeed, if} $f$ and $g$ are homogeneous polynomials of the same degree, 
then away from the common zero set of $f$ and $g$, $\val(f/g)$ is a function of $\max(0,\val(f)-\val(g))$ and $\max(0,\val(g)-\val(f))$.  The characteristic function of a set defined by
$\val (f_i) \geq \val (f_j)$ is the composition of the characteristic function of $x_i \geq x_j$ on $\G_\infty^m$, with the function $(\val (f_1),\ldots,\val (f_m))$.
Hence taking a large enough degree, and collecting together all the polynomials mentioned,  and adding more so that
$f_1,\ldots,f_m$ never vanish simultaneously, all $\xi_i$ factor through the function
$[\val (f_1): \ldots: \val (f_m)]$ of \remref{reppnrem}.  Thus 
we may take the $\xi_i$ to  be v+g-continuous.  
We denote by $x_h$ a schematic distance function to $D_0$, cf. \ref{schedis} and
we shall assume $x_h$ is one of the $\xi_i$.

 By enlarging $D_0$, we  may assume  $D_0$ contains $\xi_i \inv(\infty) \meet U$ for any irreducible component $U$
 such that $\xi_i$ is not identically \tcr{equal to} $\infty$ on $U$, cf. \lemref{infty}.  Moreover, we can demand that   $D_0$ is $G$-invariant, and that the {\em set} $\{\xi_i: i \in I \}$ is $G$-invariant, by increasing
 both if necessary.  Note that there exists a continuous function  $m=(m_1,\ldots,m_n): \G_\infty^I \to \G_\infty^n$
 whose fibers are the orbits of the symmetric group acting on $I$, namely $m((x_i)_{i \in I}) = (y_1,\ldots,y_n)$
 if $(y_1,\ldots,y_n)$ is a non-decreasing enumeration of $\{x_i\}_{i \in I}$, with appropriate multiplicities.  
   Then $\{m \circ \xi_i\}_{i \in I}$ is $G$-invariant. 
  It is clear that a homotopy preserving $m\circ \xi$ also preserves each $\xi_i$.  Thus we may assume
  that each $\xi_i$ is $G$-invariant.

\medskip 

Let $E$ be the blowing up of $\Pp^n$ at one point. Then $E$ admits a morphism $E \to \Pp^{n-1}$,
whose fibers are $\Pp^1$.
We now show one may assume $V$ admits a finite morphism to $E$, with composed morphism to $\Pp^{n-1}$
 finite on $D_0$, at least when $F$ is infinite.

\begin{lem} \label{Ps} Let $V$ be a projective variety of dimension $n$ over a field $F$.
\tcb{Assume $F$ is infinite.
  Then there exists a finite morphism $\pi : V \to \Pp^n$ and 
  a  zero dimensional subscheme $Z$ of $V$ such that if}
  $v: V_1 \to V$ denotes the blowing up at $Z$, 
     there exists a finite morphism 
  $m: V_1 \to E$
  making the diagram
	  \begin{equation*}\xymatrix{
	V_1
	\ar[d]^m \ar[r]^{v} & 
	V \ar[d]^{\pi}\\
	E \ar[r] &\Pp^n
	}
	\end{equation*}
commutative.
  Moreover, if a divisor $D_0$ on $V$ is given in advance, we may arrange that $Z$ is disjoint from $D_0$, and that 
the composition of $m$ with the projection $E \to\Pp^{n -1}$  is finite on $v \inv (D_0)$.      If a finite group $G$ acts on $V$, we may take all these to be $G$-invariant.
      \end{lem}

\prf   Let $m$ be minimal such that $V$ admits a finite morphism to $\Pp^m$.  If $m>n$, choose an $F$-rational hyperplane $H$
inside $\Pp^m$, and an $F$-rational point neither on $H$ nor on the image of $V$; and project the image of $V$
to  $H$ through this point.  Hence $m=n$, i.e. there exists a finite morphism $V \to \Pp^n$.

Given a divisor $D_0$ on $V$, choose an $F$-rational point  $z$ of $\Pp^n$ not on the image of this divisor.
The projection through
this point to a $\Pp^{n-1}$ contained in $\Pp^n$ and not containing $z$ determines a morphism
$E \to \Pp^{n-1}$.  If $V_1$ is the blowing up of $V$ at the inverse image $Z$ of $z$, we find a morphism
$V_1 \to E$; composing with $E \to \Pp^{n-1}$ we obtain the required morphism.

To arrange for $G$-invariance, we shall apply the lemma to $V':=V/G$.  Let $\phi: V \to V' $ be the natural projection.
Let $R \nsubset V'$ be the ramification \tcb{locus} of $V \to V'$. Assuming as we may that $G$ acts faithfully,
$R$ is the union over $h \in G$ of the set of fixed points of $h$; so away from $R$, $V \to V'$ is Galois and \'etale.
Let $D' $ \tcb{be a divisor containing $\phi (D_{\tcb{0}})$ and $R$}.  Applying the lemma to $(V',D')$, one
obtains $v': V_1' \to V'$, $m': V_1' \to E$, $\pi': V' \to \Pp^n$, and $Z'$ (so $v'$ is an isomorphism away from $Z'$, and $Z' \meet D' = \varnothing$).
  Let $V_1 = V_1' \times_{V'} V$.  Then $V_1 \to V$ is a blowing up of the pullback $Z$ of $Z'$ under the morphism \tcb{$\phi$ which is \'etale over $Z$},
  and all statements are clear.
  \eprf

The next lemmas provide a variant of \lemref{Ps} that works over finite fields too.  They provide a less detailed description of $V_1$, but still sufficient for our purposes; the reader who wants to assume an infinite base field may skip them. \tcb{Note that non-archimedean geometry over trivially valued fields, including finite ones, may have some relevant applications, cf. \cite{thuillier}.}
We are grateful  to Antoine Ducros for pointing out the need for a special argument in the case of a finite base field.

\begin{lem} \label{Ps1} Let $V$ be a subvariety of dimension $n$ of $\Pp^m$ over a finite field $F$.  Then there exist homogenous polynomials
 $f_1,\ldots,f_n $ in 
 $F [x_0,\ldots, x_m]$, of equal degree,
 such that $Z=V \meet (f_1=\ldots=f_n=0)$ is finite.  Given a subvariety $D$ of $V$ of dimension $<n$, we may choose $f_1,\ldots,f_n$ so that $Z$ is disjoint from $D$ and such
  that $[f_1:\ldots:f_n] : D \to \Pp^{n - 1}$ is a finite morphism. \end{lem}

\prf  Given any finite number $k$ of $F$-irreducible projective subvarieties $U_i$ of
$\Pp^m$
of positive dimension, one can always find a homogeneous polynomial $f$ in  $F [x_0,\ldots, x_m]$
\tcb{which} does not vanish on any of them.  
Indeed, by Hilbert polynomial considerations, the codimension of the space of homogeneous polynomials of degree $d$ vanishing on $U_i$ grows at least linearly with $d$.
Thus, for large enough $d$, this codimension is $> \log_q (k)$; in particular if the field $F$ has cardinality $q$,
a fraction strictly less than $1/k$ of all homogeneous polynomials of degree $d$ in  $F [x_0,\ldots, x_m]$ will vanish on $U_i$,
implying that some will vanish on no $U_i$. 
 
 On the other hand, let $w_0$ be a finite, Galois invariant,  set of points of $V(F')$,  with $F'$ a finite Galois extension of $F$.  
 We lift  $w_0 \nsubset \Pp^m$ to a finite, Galois invariant,  subset $w$ of  $\Aa^{m+1}$ in such a way that each element of $w$ has some coordinate equal to $1$.
 Let $H_d$ denote the space of homogeneous polynomials $h (x_0,\ldots, x_m)$
 of degree $d$ (with zero added),
 let $H_d(w) $, resp.  $H_d(w,1)$, denote the subspace of  $H_d$
 consisting of polynomials vanishing at each element of $w$, resp. 
 taking value $1$ on each element of $w$.
  Thus $H_d(w,1)$ is (empty or) a coset 
 of $H_d(w)$.    We now claim that $H_d(w,1) $ has a point over $F(w)$, \tcb{for $d$ large enough}.  If this is true for $w$ and $w'$ such that $w$ and $w'$ are disjoint,
 then it holds also for $w \union w'$, since $H_d / H_d(w \union w') \cong H_d/H_d(w) \oplus H_d/H_d(w')$.  So it suffices to consider a Galois
 orbit $w$.  Then for the elements $c=(c_0, \ldots, c_m)$ of $w$, we have say $c_0 =1$.  Since over the finite field $F(w)$, every function \tcb{$F(w)^m \to F (w)$}
  is represented by a polynomial, there exists over $F(w)$ a  polynomial $h$, \tcb{say of degree $d_0$,} with $h(c_1/c_0,\ldots,c_m/c_0)=1$ for $(c_0,\ldots,c_m) \in w$.  So
 $x_0^{d} h (x_1/x_0,\ldots,x_m/x_0)$ is a homogeneous polynomial \tcb{of degree $d$, for $d \geq d_0$}, as required.  Finally, by Hilbert 90, $H_d(w)$   has an $F$-basis; and since
 also  $H^1(\mathrm{Gal} (F (w)/ F), G_a^{N}) = (0)$, 
   as the affine space $H_d(w,1)$ is defined over $F$, it has a point in $F$.   
 
 We now prove the lemma.  The condition on equal degree is easily arranged afterwards, by taking appropriate powers of each $f_i$; so we ignore it.
 Inductively, we need to find $f=f_k$   that vanishes on no positive-dimensional component of $V_k=V \meet (f_1=\ldots=f_{k-1}=0)$.  
 Then it is clear that any component of $V_k$ has dimension at most $n-k$.  In particular for $k=n$ this proves the finiteness of $Z$.
 
 Further we can insist    that $f_k$  vanishes on no positive-dimensional component of $D \meet (f_1=\ldots=f_{k-1}=0)$.  As a result,
 $D \meet (f_1=\ldots=f_{n-1}=0)$ is finite, with points $c_1,\ldots,c_m \in V(F^{\alg})$.  When choosing $f=f_n$, we need also
 to insist that $f(c_i) \neq 0$.   We saw above that there exists a nonempty affine subspace of $H_d$ of codimension bounded independently 
 of $d$, whose elements satisfy $f(c_i) \neq 0$.  For large $d$, this subspace cannot be contained in the union of the linear spaces that need to be avoided in order
 to avoid vanishing on the components of $V_n$
 whose codimensions grow linearly with $d$.  
The finiteness of $[f_1: \ldots : f_n] : D \to \Pp^{n - 1}$ follows from \lemref{ital}.  
    \eprf
 
 \begin{lem}\label{ital}Let $Y$ be an irreducible quasi-projective variety of dimension $>0$ over a field $k$.
 Let $f : U \to Y$ be a dominant $k$-morphism with $U $  a Zariski open subvariety of $\Pp^m$.
 Let $X$ be a closed subvariety of $\Pp^m$ which is contained in $U$.
 Then $f \vert X$ is finite.
 \end{lem}
 
 \prf We may assume $k$ is algebraically closed and it is enough to prove $f \vert X$ is quasi-finite.
 Thus we may assume $f (X)$ is a point. Let $D$ be a divisor in $Y$ such that $f (X) \notin D$.
 Let $E$ be the Zariski closure of $f^{-1} (D)$. We have $E \nsubseteq f^{-1} (D) \cup F$, with $F = \Pp^m \m U$,
 thus $E \cap X = \varnothing$. By B\'ezout's theorem, if follows that $X$ is of dimension $0$.
  \eprf  
 
 \begin{lem}  \label{Ps2} Let $V$ be a projective variety of dimension $n$ over a finite field $F$, and let $D$ be a closed subvariety,  of dimension $< n$, containing any component of $V$ of dimension $< n$.  Then, there exists a projective variety $V_1$, a finite closed subset $Z$ of $V$, disjoint from $D$,
 a morphism   $v: V_1 \to V$ which is the blowing up of an ideal supported on $Z$
 \textup{(}in particular it is  an isomorphism above $V \m Z$\textup{)},  and a morphism $u:V_1 \to \Pp^{n-1}$ which is finite on  $v \inv (D) \cup v \inv (Z)$
 such that $v \inv (D)$ is a Cartier divisor and
 there exists a Zariski dense open subset $U_0$ of $U=\Pp^{n-1}$ such that with $V_0={u} \inv (U_0)$,  ${u} |V_0$ factors as
 $V_0 \to E_0 = U_0 \times \Pp^1 \to U_0$, with  $V_0 \to E_0$ a finite morphism, and $E_0 \to U_0$ the projection.
    If a finite group $G$ acts on $V$, we may take all these to be $G$-equivariant.
      \end{lem}

 \prf  \tcb{Fix an embedding of $V$ in $\Pp^m$.} By \lemref{Ps1} there exist homogenous polynomials
 $f_1,\ldots,f_n $ in 
 $F [x_0,\ldots, x_m]$, of equal degree,
 such that $Z=V \meet (f_1=\ldots=f_n=0)$ is finite and disjoint from $D$ and such
  that $[f_1:\ldots:f_n] : D \to \Pp^{n - 1}$ is a finite morphism.
   Let
 $V_1 \nsubset V \times \Pp^{n-1}$ be the Zariski closure of  the graph $\{ (v,(f_1(v): \ldots : f_n(v))) :  v \in V \m Z \}$.  
 Let $v$ be the first projection and
 $u$ the second projection. 
 Thus $v : V_1 \to V$ is the blowing up of $V$ along the ideal $(f_1,\ldots,f_n)$.
 \tcb{By \lemref{ital},} the restriction of $u$ to $v \inv (Z)$ is finite.
  The generic fiber of the morphism $V_1 \to U=\Pp^{n-1}$
 is a curve (possibly reducible, and possibly containing some isolated points, in $D$).  Thus it admits a finite morphism to $\Pp^1$ over $F(U)$. 
 This morphism is the generic fiber of a morphism $u: V_0 \to U_0 \times \Pp^1$, over $U_0$, for some Zariski dense open
 $U_0$ of $U$.  Equivariance is arranged by applying this construction to $V/G$ in the first place.
 \eprf


 \medskip

Let us return to the main discussion and recall our setting. We have a projective equidimensional variety $V$ together with a
hypersurface $D_0 \nsubset V$ containing the singular locus  of $V$ 
and such that there exists an \'etale morphism $V \m D_0 \to \Aa^n$, factoring through $V/G$.
Consider $v : V_1 \to V$ as provided by \lemref{Ps} and \lemref{Ps2}, respectively in the infinite and finite field case.
It is  a $G$-equivariant birational morphism  whose exceptional locus lies above a finite subset $Z$ of $V$.
By  \lemref{homotopy-descent} \tcr{and Remark \ref{rem:homotopy-descent}}, any  deformation retraction $h_1 : I \times \std{V_1} \to \std{V_1}$ leaving  the exceptional locus  invariant descends to a  deformation retraction $h : I \times \std{V}\to \std{V}$. Furthermore, if $h_1$ satisfies the theorem for $X = V_1$, so does $h$ for $X = V$.
Thus, 
pulling  back the data of \thmref{1}  to $V_1$,
and adding the above invariance requirement,
 we see that it suffices to prove the theorem for $V_1$ (which is equidimensional of dimension $n$). 
 Furthermore, 
 setting $D'_0 = v \inv (D_0) \cup v \inv (Z)$,
we have  $V_1 \m D'_0 = V \m D_0$.
In particular, $V_1 \m D'_0$ is smooth  and admits an  \'etale equivariant morphism to $\Aa^n$.
Hence, we may assume $V = V_1$ and $D_0 = D'_0$.

By construction, there is a morphism $u : V \to U=\Pp^{n-1}$, whose restriction to $D_0$ is finite, and
a Zariski dense open subset $U_0$ of $U$ such that,  setting $V_0={u} \inv (U_0)$,  ${u} |V_0 = q \circ f$ with
 $f : V_0 \to E_0 = U_0 \times \Pp^1$ a finite morphism and $q: E_0 \to U_0$ the projection.
 If a finite group $G$ acts on $V$, we may take everything to be $G$-equivariant.
Note that the hypotheses imply that $f$ is surjective.

 
 Furthermore, we may assume, after possibly shrinking $U_0$,
 that the morphism
 $f : V_0 \to E_0$ factors
through
$V_0 \overset{h}{\longrightarrow} V'_0 \overset{f'}{\longrightarrow} E_0$
with $h$ \tcb{finite} radicial and $f'$ satisfying the following condition: for every
$u$ in $U_0$, the restriction $f'_u : V'_u \to \Pp^1_u$  of $V_0 \to E_0$ over
$u$ is a generically \'etale morphism of curves.
Indeed, such a factorization exists over the generic point $\xi$ of $U_0$
and can be spread out on some dense Zariski  open set $U_0$.

 \section{Construction of a relative curve homotopy}\label{ss10.3}\nomenclature{$H_{curves}$}{relative curve homotopy}
We fix three points $0,1,\infty$ in $\Pp^1$.  We are now in the setting of  \ref{ss9.2}
 with 
$U_0 \nsubset U=\Pp^{n-1}$. 
For any divisor $D$ on $E_0$  we consider
$\psi_D : [0,  \infty] \times E_0 \to
\std{E_0/U_0}$ as in  \ref{ss9.2}.

\begin{lem}\label{pass}Let $W$ be an $A$-iso-definable subset of  
$\std{E_0/ U_0}$ such that
$W \to U_0$ has finite fibers. There exists a divisor $D'$ on
$E_0$, generically finite over $U_0$,  such that for
every $u$ in $U_0$, for every $x$ in $W$ over $u$, 
the intersection of $D'$ with the ball in $\Pp^1_u$ corresponding to
$x$ is nonempty.
\end{lem}

\prf   Recall we are working over a field base $A$.   By splitting $W$ into two parts
(then taking the union of the divisors $D'$ corresponding to each part), we may 
assume $W \nsubseteq \std{\Oo} \times U_0$ where $\Oo$ is the   unit ball.
Let $a$ be  a point in $U_0$; so  $W_a \nsubseteq \std{\Oo}$.

We claim that 
there exists a finite $A (a)$-definable subset
$D'_a$ of $\Oo$ such that 
for every $x$ in $W_a$, 
the intersection of $D'_a$ with the ball in $\Oo$ corresponding to
$x$ is nonempty.   Let $W^1$ be the set of simple points in $W$.
Thus, $W$ splits into two disjoint iso-definable sets
$W^1$ and $W^2 = W \m W^1$. 
Let $D'^{1}_a$ be the union 
of the simple points in $W^1_a$.  If  $A(a)$ is   trivially valued, any $A$-definable closed sub-ball of
$\Oo$ must have valuative radius $0$, i.e. must equal $\Oo$.  In this case 
we set $D'^{2}_a = \{0\}$.
Otherwise, $A (a)$ is a nontrivially valued field, and so $\acl (A (a))$ is a model of $\ACVF$.
Hence, if we denote by
$\widetilde W_a$ the finite set of closed  balls corresponding to the points in
$W_a$, for every $b$ in $\widetilde W_a$,
$b \meet \acl (A (a)) \not= \varnothing$. \tcb{Thus} there exists a finite $A (a)$-definable \tcb{$D'^{2}_a$}
set such that $D'^{2}_a \meet b
\not= \varnothing$ for every $b$ in
$\widetilde W_a$. Set $D'_a = D'^{1}_a \cup D'^{2}_a$.

By compactness we get a constructible set $D''$ finite over $U_0$
with the required property. Taking the Zariski closure
of $D''$ we get a Zariski closed set $D'$ generically finite over
$U_0$ with the required property. 
\eprf

\begin{lem}\label{pass2}There exists a divisor $D'$ on $E_0$
 such that,
for any divisor $D$ containing $D'$,  
$\psi_D$ lifts uniquely to an $A$-definable map  
$h : [0,  \infty] \times V_0 \to
\std{V_0/U_0}$, which is fiberwise a homotopy.
\end{lem}

\prf
We proceed as in the proof  of  \thmref{ret3}.
 By assumption the morphism
 $f : V_0 \to E_0$ factors
through
$V_0 \overset{h}{\longrightarrow} V'_0 \overset{f'}{\longrightarrow} E_0$
with $h$ \tcb{finite} radicial
and  for every
$u$ in $U_0$, the restriction $f'_u : V'_u \to \Pp^1_u$  of $V_0 \to E_0$ over
$u$ is a generically \'etale morphism of curves.
Thus, for every $u$ in $U_0$, the restriction $f_u : V_u \to \Pp^1_u$ of $V_0 \to E_0$ over
$u$ factors as
$V_u \overset{h_u}{\longrightarrow} V'_u \overset{f'_u}{\longrightarrow}  \Pp^1_u$, with $h_u$ the restriction of $h$.
Note that $V'_0 \to U_0$ is a relative curve so that
$\std{V'_0 /U_0}$ is iso-definable over $A$  by \thmref{f11}. 
There is a subset $W_0$ of $\std{V'_0 /U_0}$,  iso-definable over $A$,  
containing, for every point $u$ in $U_0$, all singular points of $C'_u$, all ramification points
of
$f'_u$ and all forward-branching points
of $f'_u$, and such that the fibers $W_0 \to U_0$ are all finite.
Such an $W_0$ exists by \lemref{fb} (uniform finiteness of the set
of
 forward-branching points).
Let $W$ be the image of $W_0$ in $E_0$. Then $D'$ provided
by \lemref{pass} does the job. 
\eprf

Let  $D$ be a divisor on $E_0$ as in \lemref{pass2},  and such that
$D$ contains the image of $D_0$ in $E_0$. Assume also $D$ contains the infinity divisor in $E_0$.
Then $\psi_D$ lifts to an $A$-definable map  
$h_{curves}^0 : [0,  \infty] \times V_0 \to \std{V_0/U_0}$.
By \lemref{preserve-xi},   after enlarging $D$, one can arrange that 
$h_{curves}^0$ preserves the functions $\xi_i$. Note that $G$-invariance follows from uniqueness of the lift.  
\tcb{After shrinking $U_0$ we may assume that the restriction of $u : E_0 \to U_0$ to $D$ is finite, that is, that $D$ has no vertical component over $U_0$.}

By \lemref{homotopy-lift} and \lemref{P1-homotopy}, $h_{curves}^0$ is
 v+g-continuous at each point of $[0,  \infty] \times \tcb{V_0}$. 
We  extend $h_{curves}^0$ to $h_{curves} : [0,  \infty] \times V \to \std{V/U}$
by setting $h_{curves} (t, x)= x$ for every $t$ in $[0,  \infty]$ and every $x$ in $V \m V_0$.

\begin{lem}\label{vcont}The mapping $h_{curves}$ is 
 g-continuous on $[0, \infty] \times V$ and v-continuous at each point of $[0, \infty] \times X$ for $X = \tcb{V_0}  \union D_0$. 
\end{lem}

\prf Since $V \m \tcb{V_0}$ is g-clopen,
g-continuity  may be shown separately on $V \m \tcb{V_0}$ and away from $V \m \tcb{V_0}$.
On $V \m \tcb{V_0}$ it is trivial since $h_{curves}$ is constant there.
Away from $V \m \tcb{V_0}$ it was already proved.

It remains to show {v}-continuity at points on $D_0$.   Let $F_2$,   $\res$ be as in  \ref{ss8.8} and in the {v}-continuity
criterion  \lemref{vcc2}.  
Let $p \in V(F_2)$ with $\res(p) \in D_0$.  If $p \notin \tcb{V_0}$ then $h_{curves}$ fixes
$p$, so assume $p \in \tcb{V_0} (F_2)$.  Set $q = \res (p)$.
Fix $t$ in $[0, \infty]$ and let $q_t = \res_{21*} (h_{curves} (t, p))$. Since $h_{curves} (t, q) = q$, it is enough to prove that
$q_t = q$.
Recall we assume one of the $\xi_i$
is  a schematic distance function $x_h$ to $D_0$, cf. \ref{schedis}.
Since
$x_h (h_{curves} (t, p)) = x_h (p)$,
it follows that
$\res_{21*} (x_h (h_{curves} (t, p))) =  \res_{21*} (x_h (p)) = \infty$.
Thus $q_t$ lies in $\std{D_0}$. Since it lies 
on the fiber of $u$ at $q$, and the intersection of this fiber with $D_0$ is a finite set $D_{0q}$, it follows that
$q_t$ is a simple point lying on $D_{0q}$.
Let $q' \not= q$ be another point of 
$D_{0q}$ and let $\vartheta$ be a regular function on some Zariski open set containing $q$ and $q'$
which vanishes at $q'$ and not at $q$. Thus $\val (\vartheta (q))$ is equal to some finite $\gamma \in \G (F_1)$
and  $\val (\vartheta (p)) = \gamma$ also. On the other hand the set of $\val (\vartheta (q_t))$ is finite.
By continuity of 
$h_{curves}$ in the $t$-variable one gets that 
$\val (\vartheta (h_{curves} (t, p)))$ cannot jump and is equal to $\gamma$ for all $t$.
Hence, for every $t$,  $q_t \not= q'$, and $q_t = q$ follows.
\eprf

By \lemref{hbasic0}
 the restriction of
$h_{curves}$ to $[0,  \infty] \times \tcb{V_0} \cup D_0$ extends to a  
 deformation retraction
$H_{curves} : [0,  \infty] \times \std{\tcb{V_0} \cup D_0}\to  \std{\tcb{V_0} \cup D_0}$.
Since $D_0$ is finite over $U$, the image
$\Upsilon_{curves}  = h_{curves} (0, \tcb{V_0} \cup D_0)$
is iso-definable over $A$ in $\std{V/ U}$ and relatively $\G$-internal.  
Thus, as above \thmref{relative}, we can  identify $\std{\Upsilon_{curves}}$ with its image in $\std{V}$.
It follows that
 the image $H_{curves} (0, \tcb{\std{V_0 \cup D_0}})$
is equal to 
 $\std{\Upsilon_{curves}}$. 
 By construction $H_{curves} (\infty, x) = x$ for every $x$ 
and $H_{curves}$ satisfies $(*)$.

Let $x_v: U \to [0, \infty]$ be a   schematic distance  
to  the image  of $\tcb{V \m V_0}$ in $U$, cf. \ref{schedis}.
We still denote by $x_v$ its pullback to $V$ (which is a schematic distance to $\tcb{V \m V_0}$)
and the corresponding extension to $\std{V}$.
Let us check  that
$\std{\Upsilon_{curves}}$ is $\si$-compact
via $(x_h,x_v)$.
Indeed, on $\std{\Upsilon_{curves}}$ the infinite locus of $x_v$ is contained in that
of $x_h$, and 
$\std{\Upsilon_{curves}}$ is compact at $x_h^{-1} (\infty)$ since
$\{x \in \std{V}: x_h (x) = \infty\}$ is contained in $\std{\Upsilon_{curves}}$.
Furthermore, 
\tcb{since for any $\gamma \in \Gamma$,
the set $\{x \in \std{V}: x_v (x) \leq \gamma\}$
is definably compact and preserved  by $H_{curves}$,
$\{x \in \std{\Upsilon_{curves}}: x_v (x) \leq \gamma\}$ is
definably compact, 
being the image by a continuous definable map of a definably compact set.}

\section{The base homotopy}\label{ss10.4}

By \thmref{omin-finite} there exists a finite pseudo-Galois covering $U'$   of $U$
and a finite number  of  $A$-definable functions $\xi'_i: U' \to \G_{\infty}$  
such that, for $I$ a generalized interval,  any  $A$-definable deformation retraction $h: I \times U \to \std{U}$ 
lifting to  a  deformation retraction
$h': I \times U' \to \std{U'}$ 
respecting the functions 
$\xi'_i$, also lifts to an  $A$-definable deformation retraction $ I \times \std{\Upsilon_{curves}} \to \std{\Upsilon_{curves}}$  
respecting the restrictions of the functions $\xi_i$
on $\Upsilon_{curves}$ and the $G$-action.

Now by the  induction hypothesis applied to $U'$ and $\mathrm{Gal}(U'/U)$,  such a pair $(h,h')$ does exist;
we can also take it to preserve $x_v$, the schematic distance to $\tcb{V \m V_0}$. Set $h_{base} = h$.
Hence, $h_{base}$ lifts to 
 a  deformation retraction
  \[H_{\widetilde{base}} :  I \times \std{\Upsilon_{curves}} \to \std{\Upsilon_{curves}},\] \tcr{which by (2) in \thmref{omin-finite} we may assume to respect the restrictions of the functions $\xi_i$ and the $G$-action}.
 \nomenclature{$H_{\widetilde{base}}$}{base homotopy}

 Recall the notion of Zariski density in $\std{U}$, \ref{zt}.
By induction $h_{base}$ has an  $A$-iso-definable $\G$-internal final image 
$\Upsilon_{base}$ and we may assume $\Upsilon_{base}$ is Zariski dense in $\std{U}$.
By \thmref{omin-finite} we may assume
$H_{\widetilde{base}}$ has an $A$-iso-definable  $\G$-internal  final image equal to 
$     \std{\Upsilon_{curves}} \cap \std{\tcb{u}}^{-1} (\Upsilon_{base})$ 
and by induction we may assume
$H_{\widetilde{base}}$ satisfies 
$(*)$.

By composing
the homotopies $H_{curves}$ and $H_{\widetilde{base}}$
one gets an $A$-definable deformation retraction  
\[H_{bc} = H_{\widetilde{base}} \circ H_{curves} : I' \times \std{\tcb{V_0} \cup D_0} \longrightarrow \std{V},\] 
where $I'$ denotes the generalized interval obtained by gluing $I$ and $[0, \infty]$. 
The image is contained in the image of $H_ {\widetilde{base}}$, but contains 
$H_{\widetilde{base}} (e_I \times \std{\Upsilon_{curves}/U)}$, the image over the simple points of $U$.    As
these sets are equal, the image   is equal to both, and is iso-definable and $\G$-internal; we denote it by
$\Upsilon_{bc}$.  Thus, $ \Upsilon_{bc} =     \std{\Upsilon_{curves}} \cap \std{\tcb{u}}^{-1} (\Upsilon_{base})$. 
In general $\Upsilon_{bc}$ is not definably compact, but it is $\si$-compact
via $(x_h,x_v)$, since $H_{\widetilde{base}}$ fixes $x_v$ and $\std{\Upsilon_{curves}}$ is $\si$-compact
via the same functions.  (Note that $ \Upsilon_{bc} \cap x_h^{-1} (\infty) = \std{D_0} \cap \std{\tcb{u}}^{-1} (\Upsilon_{base})$.)

\begin{lem}\label{Upzdense}\leavevmode\begin{enumerate}
\item The subset
$\Upsilon_{bc}$ is  a Zariski dense subset of $\std{V}$.
\item One may choose $h_{base}$ so that, for every irreducible component
$V_i$ of $V$, $\Upsilon_{bc} \cap \std{V_i}$ is of  pure dimension $n=\dim(V)$.  
\end{enumerate}
\end{lem}

\prf
Let $V_i$ denote the irreducible components of $V$, $\tcb{u} : \std{V/U} \to U$ and
$\std{\tcb{u}} : \std{V} \to \std{U}$ denote the projections.
\tcb{Since $H_{curves}$  preserves $\std{D_0}$, its complement (check it fiberwise) and the connected components of its complement by continuity, 
it preserves each of the $\std{V_i}$.
Furthermore,} there exists
an open dense subset $U_1 \nsubseteq U$ such that, for every $x \in U_1$,
$\tcb{u} \inv (x) \cap \Upsilon_{curves} \cap \std{V_i}$
is Zariski dense \tcb{in  $\std{u}^{-1} (x) \cap \std{V_i}$} for every $i$.
It follows that,
for every $x \in \std{U_1}$,
$\std{\tcb{u}} \inv (x) \cap \std{\Upsilon_{curves}} \cap \std{V_i}$
is Zariski dense \tcb{in  $\std{u}^{-1} (x) \cap \std{V_i}$} for every $i$ (recall 
$\std{\Upsilon_{curves}}$
is identified with
$\int_U \Upsilon_{curves}$).
Pick $x \in \Upsilon_{base}$ which is Zariski dense in $\std{U}$,
then $\std{\tcb{u}} \inv (x) \cap  \Upsilon_{bc}$ is 
Zariski dense in $\std{V}$.

Next, we deal with local dimension.  
Consider a component $V_i$ of $V$.    
Let $C$ be an irreducible component of a fiber of $V_i$ above $\tcb{U_0}$.
Since $D_0$ was chosen so that $D_0 \meet C \neq \varnothing$, it follows directly from the definition
that the homotopy on $C$   has   image containing more than one point. 
It follows \tcb{by construction} that the image of each irreducible component $C_{\ell}$ of $C$ over the algebraic closure of $F$
by the homotopy  also contains more than one point. 
By \thmref{connectedness},  the image of each $C_{\ell}$ under that homotopy is necessarily connected.
Since it is of dimension $\leq 1$, it follows that this image has no isolated points, so is purely one-dimensional. 
Thus the image of $C$ under the homotopy is also purely one-dimensional.

Now  $ \Upsilon_{bc} =     \std{\Upsilon_{curves}} \cap \std{\tcb{u}}^{-1} (\Upsilon_{base})$; and by the inductive assumption (7) of \thmref{1}, one may assume that
$\Upsilon_{base}$ has pure dimension $n-1$. 
\tcb{Since the morphism $V \to U$ restricts to a composition $V_0 \to U_0 \times \Pp^1 \to U_0$, where $V_0 \to U_0 \times \Pp^1$
is finite surjective, it follows from  \corref{finite-open+}  that the map $\std{V_0} \to \std{U_0}$ is open.
In particular the maps 
$\std{V_i \cup V_0} \to \std{U_0}$ are open.}
  It follows easily that
$\Upsilon_{bc} \cap \std{V_i}$ is of  pure dimension $n$. \eprf

\section{The tropical homotopy}\label{ss10.5}
In this rather technical section we construct a homotopy in $\G_\infty^w$ that we shall use in  \ref{ss10.6} in order to insure that the homotopy we build fixes pointwise its final image at every time.

By \thmref{G-embed-1c}, there exists an $A$-definable, continuous, injective map $\a: \Upsilon_{bc} \to \G_\infty^w$, with image $W \nsubseteq [0,\infty]^w$, where $w$ is a finite $A$-definable set.    
We may assume for some coordinate $x_i$ (resp. $x_j$), $x_i \circ \a$ (resp. $x_j \circ \a$) is the restriction of $x_h$ (resp. $x_v$).
Indeed, we may add two points $h,v$ to $w$ which we view as $A$-definable, i.e. fixed by  the action of the Galois group and
replace $\alpha $ by $x \mapsto (\alpha (x), x_h (x), x_v (x))$.
We shall denote by $\v$ and $\h$ the projections $\G_\infty^w \to \G_\infty$
on the $v$ and $h$ coordinate, respectively.


We write $[x_i=x_j]$ for   $\{a \in [0,\infty]^w: x_i(a)=x_j(a) \}$, and similarly $[x_i=0]$, etc.

Since $\Upsilon_{bc}$ is $\si$-compact via $(x_h,x_v)$, $W$ is $\si$-compact via $(\h,\v)$.
In particular, $W \m [\v=\infty]$ is $\si$-compact via $\v$,  and hence closed  in  $\G_\infty^w \m [\v=\infty]$;   so
$W \meet \G^w$ is closed in $\G^w$.  

We   let $G$ act on $W$, so that $\a:  \Upsilon_{bc} \to \G_\infty^w$ is equivariant.  By re-embedding $W$
in   $\G_\infty^{w \times G}$, via $w \mapsto (\si(w))_{\si \in G}$, we may assume $G$ acts on the coordinate set $w$,
and the induced action of $G$ on    $\G_\infty^w$ extends the action of $G$ on $W$.  
We still denote by $\xi_i$ the functions on $W$ that are the composition of the restriction of
$\xi_i$ to $ \Upsilon_{bc}$ with $\alpha^{-1}$.

\tcb{In \lemref{G-homotopy}, we shall show the existence, entirely within $\G_\infty^w$, 
of a definable deformation retraction from $(W \meet \G^w) \union [\h=\infty]$ to a definably compact
subset $W_0$. 
Furthermore we shall show that when $W$ has pure dimension $n$,   one can insure
$W_0 \cap W$ has also pure dimension $n$.
Then, in \lemref{G-homotopy2}, we shall extend this result to
$(W^o \m [\v=\infty])  \union [\h=\infty]$, for some
 z-dense and z-open definable  subset $W^o$ of $W$.
This will be used in an essential way in the final part of the proof given in  \ref{ss10.6}.
}

\begin{lem}\label{G-homotopy} Let 
\[W' = (W \meet \G^w) \union [\h=\infty].\]
There exists an $A$-definable deformation retraction 
\[H_\G: [0,\infty] \times W'  \to W'\] 
whose image is a definably compact subset 
$W_0$ of $W'$ and such that $H_\G$ leaves the $\xi_i$ invariant, fixes $[\h=\infty]$, and is $G$-equivariant.  
Moreover, one may require the following to hold:
\begin{enumerate}
\item There exists an $A$-definable open subset $W_o$ of $W$ containing
$W_0 \m [\h=\infty]$ and $m \in \Nn$, $c \in \G (A)$, 
such that $x_i \leq (m  + 1)x_h + c$ on $W_o$, for every $i \in w$.
\item If $W$ has pure dimension $n$,  then
$W_0 \cap W$ has also pure dimension $n$.
\end{enumerate}
\end{lem}  
 
In this lemma, we take $0$ to be the initial point, $\infty$ the final point.   On $\G_\infty$, we view   $\infty$ as the unique simple point.  In this sense the 
 flow is still  away from the simple points,
as for the other homotopies.  Moreover, starting  at any given point, the flow will terminate at a finite time.   The homotopy we obtain will in fact be
a semigroup action, i.e. $H_\G(s,H_\G(t,x))=H_\G(s+t,x)$,
in particular it will satisfy $(*)$ (in the form: $H_\G(\infty,H_\G(t,x))=H_\G(\infty,x)$).

\prf For the convenience of the reader we shall  divide the proof into 3 steps.

\medskip
\step1 \kern-0.5em  \tcb{\textit{Preliminaries.}} \tcb{We start by choosing} 
an $A$-definable cell decomposition $\mathcal{D}$ of $\G^w$, compatible with $W \meet \G^w$ and with $[x_a=0]$ and $[x_a = x_b]$
where $a,b \in w$,
and such that  each $\xi_i$ is linear on each cell of $\mathcal{D}$.   
  We also assume 
$\mathcal{D}$ is invariant under both the Galois action \tcb{of $\Aut(\acl(A) / A)$} and the $G$-action on $w$.  This can be achieved
as follows.  Begin with a finite set of pairs $(\a_j,c_j) \in \Qq^w \times \G^w$, such that each of the subsets of $\G^w$ referred to above is defined by inequalities of the form  $\a_j v - c_j \odot_j 0$, where $\odot_j$  is $<$ or $>$ or $=$.
Take the closure of  this set under the Galois action and the $G$-action.    A cell  of $\mathcal{D}$ is any nonempty set defined by conditions  $\a_j v - c_j \odot_j 0$,
where $\odot_j$ is any function from the set of indices to $\{<,>,=\}$.  Such a cell is an open convex subset of its affine span.   

Any bijection $b: w \to \{1,\ldots,|w| \}$ yields a bijection $b_{*}: \G^w \to \G^{|w|}$; the image of $c_j$
under these various bijections depends on the choice of $b$ only up to reordering.  Thus $b_{*}(c_j)$ 
gives a well-defined subset of $\G$, which belongs to $\G(A)$.  
Let $\bA$ be the convex subgroup of $\G=\G(\Uu)$ generated by $\G(A)$, and let $B=\G(\Uu)/\bA$.
For each cell $C$ of $\mathcal{D}$, let $\beta C$ be the image of $C$ in $B^w$.  Note that $\beta C$ may
have smaller dimension than $C$; notably, $\beta C= (0)$ iff $C$ is bounded.  At all events $\beta C$ is
a cell defined by homogeneous linear equalities and inequalities.   When $\G(A)\neq (0)$, $\beta C$ is always
a closed cell, i.e. defined by weak inequalities.  

 For any $C \in \mD$, let $C_\infty$ 
be the closure of $C$ in $\G_\infty^w$.  Let $\mD_0$ be the set of   cells $C \in \mD$
such that $C_\infty \m \G^w \nsubseteq [\h = \infty]$.  Equivalently, $C \in \mD_0$ if and only if for each $i \in w$, an inequality
of the form $x_i \leq m \h + c$ holds on $C$, for some $m \in \Nn$ and $c \in \G \tcb{(A)}$.  
Other equivalent conditions are that $x_i \leq m \h$ on $\b C$ for some $i$, or that there
exists no $e \in \beta C$ with $\h(e)=0$ but $x_i(e) \neq 0$.  
 Let
 \[W_0 = (W' \cap (\union_{C \in  \mD_0}C)) \union [\h = \infty].\] 
  It is clear that $W_0$ is a definably compact subset of $\G_\infty^w$, 
  contained in $W'=(W \cap\G^w) \union [\h = \infty]$.   
    
More generally, define a quasi-ordering $\leq_C$ on $w$  by:  $i \leq_C j$
if for some $m \in \Nn$, $x_i(c) \leq m x_j(c)$ for all $c \in \b C$.  Since the decomposition respects the hyperplanes $x_i=x_j$, we have $i \leq_C j$ or $j \leq_C i$ or both. Thus $\leq_C$ is
a linear quasi-order.  
     Let $\b' C = \b C \meet [\h = 0]$.  We have $\b ' C = 0$ iff $\h$ is $\leq_C$-maximal iff $C \in \mD_0$.           If $C \in \mD_0$, let $e_C=0$.  Otherwise, $\b' C$ is a nonzero rational linear cone, in the positive quadrant.  Let
 $e_C$ be the barycenter of $\b' C \meet [\sum x_i = 1]$ (here we view $\b' C$ as a cone in $\mathbb{Q}_+^w$).  
 Thus $e_C$ belongs to $\mathbb{Q}_+^w$ and is a nonzero element of $\b' C$
which is $G$ and Galois invariant.    
    
 For $t \in \G$, we have $te_C:=e_C t \in \G^w$.   If $e_C \neq 0$ then $\G e_C$ is unbounded in $\G^w$, so for any $x \in C$ there exists $t \in \G$   such that $x -  t e_C \notin C$.    Let $\tau(x)$ be the unique smallest such $t$.   Note that $\tau (x) >0$.
    

\medskip
 \step2 \kern-0.5em {\itshape \tcb{Construction of $H_\G$ and continuity.}} We will now define $H_\G: [0,\infty] \times C \to \G^w$  separately on each cell $C \in \mD$ by induction on the dimension of $C$, as follows.  If $C \in \mD_0$,
$H_\G(t,x)=x$.  Assume $C \in \mD \m \mD_0$.
 If $x \in C$ and $t \leq \tau(x)$, let  $H_\G(t,x)= x-te_C$.    So $H_\G(\tau(x),x)$ lies in a lower-dimensional cell $C'$.
For $t \geq \tau(x)$ let    $H_\G(t,x) =  H_\G(t-\tau(x),\tau(x))$.   
For fixed $a$, $H_\G(t,a)$ thus traverses finitely many
cells as $t \to \infty$, with strictly decreasing dimensions, thus ultimately reaching $W_0$.

We claim that $H_\G$ is continuous on $[0,\infty] \times \G^w$.   
To see this fix $a \in C \in \mD$ and let $(t',a')  \to (t,a)$.  We need to show that 
$H_\G(t',a') \to H_\G(t,a)$.  By curve selection it suffices to consider $( t',a')$
varying along some line $\lambda$ approaching $(t,a)$.   For some cell $C'$ we have $a' \in C'$ eventually along this line.  

If $a' \in W_0$ then $a \in W_0$ since $W_0$ is closed.  In this case we have $H_\G(a',t')=a', H_\G(a,t)=a$,
and $a' \to a$ tautologically.  Assume therefore that $a' \notin W_0$, so 
$C' \notin \mD_0$ and $e' \neq 0$, where     $e'=e_{C'}$.

Consider first the case: $t' \leq \tau(a')$ (cofinally along $\lambda$).   Then by
definition we have  $H_\G(t',a') = a' - t' e'$.
Now $C$ must be a boundary face of $C'$, cut out from the closure of $C'$ by certain
hyperplanes $\alpha_j v -c_j=0$ $(j \in J(C,C'))$.     We have $\alpha_j v = c_j$ for $v \in  C$,
and (we may assume) $\alpha_j v \geq  c_j$ for $v \in C'$. 

If $\gamma_j = \alpha_j e' > 0$ for some $j$,   fix such a $j$.    As $t' \leq \tau(a')$, we have $\alpha_j (a' - t' e') = \alpha_j a' - \gamma_j t' \geq c_j$, so $t' \leq \gamma_j \inv (\alpha_j a' - c_j) $.  Now $a' \to a$ so $\alpha_j a' - c_j \to 0$.  Thus $t' \to 0$, i.e. $t=0$.
So $H_\G(t,a) = a $, and 
$H_\G(t,a) - H_\G(t',a') = a-(a'-t'e') = (a-a')+t'e \to 0$ (as $(t',a') \to (t,a)$ along $\lambda$).   

The remaining possibility is that $\alpha_j e' = 0$ for each $j \in J(C,C')$.   So $\alpha_j v=0$
for each $v \in \b' C'$.     Hence $\beta' C' \nsubseteq \beta C$.  Since $\beta' C \nsubseteq \beta' C'$,
it follows that $\beta' C = \beta' C'$ and so $e_C =e_{C'}$.  Now $(t,x) \mapsto x -te'$ is continuous
on all of $\G \times \G^w$ so on $C \union C'$, and hence again $H_\G(t',a') \to H_\G(t,a)$.

This finishes the case $t' \leq \tau(a')$.
In particular, $\tau(a') \to t^*$ for some $t^*$, and  letting $a''=H_\G(\tcb{\tau(a'), a'})$,  $a'' \to \tcb{H_\G (t^*,a)}$.
Now by induction on the dimension of the cell $C'$, we have $H_\G(t'-\tau(a'),a'') \to H_\G(t-t^*,a)$; it follows
that $H_\G(t',a') \to H_\G(t,a)$.   This shows continuity on $[0,\infty] \times \G^w$.

Note that if $C \in \mD \m \mD_0$, then $\xi_i$ depends only on coordinates $x_i$ with $\tcb{i} \leq_C \h$.   This follows from the fact that $\xi_i$ is bounded on any part of $C$ where $\h$ is bounded
(by assumption $\xi_i^{-1}(\infty) \nsubseteq D_0$);
so $\xi_i \leq m \h$ for some $m$, up to an additive constant.  Since $x_i(e_C) =0$ for $i \leq_C \h$, it follows that $\xi_i$ is
left unchanged  by the homotopy on $C$.    So along a path in the homotopy, $\xi_i$ takes only finitely many values (one on each cell); being continuous, it must
be constant.  In other words the $\xi_i$ are preserved.  The closures of the cells are also preserved,
hence, as $W \meet \G^w$ is closed, $W \meet \G^w$ is preserved by the homotopy.

Extend $H_\G$ to $W'$ by letting $H_\G(t,x)=x$ for $x \in W' \m \G^w$.    Thus, $W_0$  will be the image of the homotopy and
by construction $H_\G$ fixes $[\h=\infty]$.
We still have to prove  that $H_\G$ is continuous
at $(t,a)$ for   $a\in W' \m \G^w$, i.e. $\h(a)=\infty$.  We have to show that for $a'$ close to $a$, for all $t$, $H_\G(t,a')$ is also close to $a$.  
   If $a' \notin \G^w$ we have $H_\G(t,a')=a'$.
Assume $a' \in \G^w$; so $a' \in C$ for some $C \in \mD$.  If    $C \in \mD_0$, again we have $H_\G(t,a')=a'$.
Otherwise, \tcb{set $a'' = H_\G(\tau(a'),a')$. Thus, 
$a'' \notin C$ and  belongs to a cell of smaller dimension.}
We will show   that 
$H_\G(t,a')$ remains close to $a$ for $t \leq \tau(a')$.  In particular,   $a''$ is close to $a$;  so (inductively) $H_\G(t,a'')=H_\G(\tcb{\tau (a')}+t,a)$ is close to $a$.  Thus it suffices to show  
  for each coordinate $i \in w$  
that $x_i(a')$ remains close to $x_i(a)$.   If $i \leq_{C'} \h$ then the homotopy does not change  $x_i(a')$ so (as $a$ is fixed) we have $x_i(H_\G(t,a')) =x_i(a') \to x_i(a) = x_i (H_\G(t,a))$.   So assume $h <_C i$.  Since $\h(a)= \infty$ we have $\h(a') \to \infty$ and hence $x_i(a') \to \infty$.  
So $x_i(a) = \infty = x_i (H_\G(t,a))$.  
For any   $c=H_\G(t,a')$, $t \leq \tau(a')$,
 we have  $x_i(c)  \geq \h(c)/m =  \h(a')/m$ up to an additive constant.  Since $a' \to a$, $\h(a')$ is large, so $x_i(c)$ is large, i.e. close to $x_i(a)$.  This proves the continuity of $H_\G$ on $W'$.   
 This ends the proof of \lemref{G-homotopy} except for the additional items.

 \medskip
\step3 \kern-0.5em {\itshape \tcb{End of the proof.}}
 For (1), \tcb{note that by construction}, for each $i \in w$ there exists some $m_i \in \Nn$
and $c_i \in \G(A)$ 
such that $x_i \leq m_i x_h + c_i$ on $W_0 \meet \G^w$.  Set $m = \max_i m_i$ and $c = \max_i c_i$.
Now the open subset
of $W \meet  \G^w$ defined by $ W_o = \{x \in W \meet  \G^w; x_i < (m + 1) x_h + c, \forall i \in w\}$ does the job.
Now let us prove that one can require (2).
Set $M  = \vert \tcb{w} \vert (m + 1)$, $K = \vert \tcb{w} \vert c$
and let $L$ be the hyperplane $\sum_i x_i = M x_h + K$. Note  that $L$ is both $G$ and Galois invariant.
\tcb{We now consider the cell decomposition $\mathcal{D}'$ generated  by $L$ and $\mathcal{D}$ and we denote by
$\mathcal{D}'_0$} the corresponding set of ``bounded'' cells. We claim that replacing $\mathcal{D}$
by $\mathcal{D}'$ does the job.
\tcb{Let $C$ be a cell in $\mathcal{D}'_0$ which is contained in $W \cap \G^w$. Thus
 $C$ lies in the closure of a cell $C'$ in $\mathcal{D}$  of dimension $n$ and  contained in $W \cap \G^w$. Let $U$ be the half space  defined by 
$\sum_i x_i < M x_h + K$.
Thus $C'' = U \cap C'$ is a cell in $\mathcal{D}'_0$  of dimension $n$ contained in $W'$
and  $C$ lies in the closure of $C''$.}
This shows that after replacing $\mathcal{D}$
by $\mathcal{D'}$,
$W_0 \m  [\h = \infty]$ is of dimension $n$ at every point.
We still have to take care of $W \cap [\h = \infty]$.
Let $x$ 
be a point in $W \cap [\h = \infty]$. If some neighborhood of $\tcb{x}$ in $W$ is contained
in
$[\h = \infty]$, there is nothing to prove.
Otherwise, $x$ is in the closure of $W'$,
hence also  in the closure of image of $W'$ under the retraction attached to 
 $\mathcal{D'}$,  \tcb{$x$ being invariant under the retraction}. Since that image has dimension $n$ at all points, we are done.
 Finally note that it is possible to achieve (1) and (2) simultaneously.
  \eprf

While the construction of the $\G$-homotopy is essentially carried out in \lemref{G-homotopy}, we need to extend it to a more general situation in which,
e.g. $W \meet \G^w = \varnothing$, i.e. $W$ lies entirely within the $\infty$-boundary of $\G_\infty^w$.

\begin{lem}\label{G-homotopy2}   There exists a z-dense and  z-open $A$-definable subset $W^o$ of $W$ such that with 
\[W' =  (W^o \m [\v=\infty])  \union [\h=\infty],\]
there exists an $A$-definable deformation retraction 
\[H_{\G}: [0,\infty] \times W'  \to W'\]
whose image is a definably compact set \nomenclature{$H_{\G}$}{$\G$-homotopy}
$W_0$ of $W'$ and such that $H_{\G}$ leaves the $\xi_i$ invariant, fixes $[\h=\infty]$, and is
$G$-equivariant. 
Moreover, one may require the following to hold:
\begin{enumerate}
\item There exists an $A$-definable open subset $W_o$ of $W$ containing
$W_0 \m [\h=\infty],$ and $m \in \Nn$, $c \in \G (A)$, for $i \in w$,
such that $x_i \leq (m + 1)x_h + c$ on $W_o$, for every $i \in w$.
\item Let $W=\union_\nu W_\nu$ be the decomposition of $W$ into z-components. 
For each $\nu$ such that $W_\nu$ has pure dimension $n_\nu$,  
$W_0 \cap W_\nu$ has also pure dimension $n_\nu$.
\end{enumerate}
\end{lem}

\prf  
First assume $W$ is z-irreducible.
  Let $w^o$ be the set of all $i \in w$ such that the $i$-th projection $\pi_i:W \to \G_\infty$
does not take the constant value $\infty$ on $W$;  the set $w^o$ is Galois invariant.  
Clearly $\pi^o  = \Pi_{i \in w^o} \pi_i$ is a homeomorphism between $W$ and its  image.  
Note that $\pi^o(W) \meet \G^{w^o}$ is z-open and z-dense in $\pi^o(W)$, and disjoint from $[\v=\infty]$.
Set $W^o = \pi^o{}^{-1} (\pi^o(W) \cap \G^{w^o})$. 
Thus, either $W^o \cap  [\v=\infty] = \varnothing$ or $W$ is contained in $[\v=\infty]$  \tcb{(hence in $[\h =\infty]$)}. Set
$W' =  (W^o \m [\v=\infty])  \union [\h=\infty]$.
In the first case, applying
\lemref{G-homotopy} to $\pi^o(W) \cap \G^{w^o}$ and pulling back by $\pi^o$ we obtain the required homotopy $H_{\G} = H_{\G, W}
: [0,\infty] \times W'  \to W'$.  
Furthermore one may require there exists an $A$-definable open subset $W_o$ of $W$ containing
$W_0 \m [\h=\infty],$ and $m \in \Nn$, $c \in \G (A)$, for $i \in w_o$,
such that $x_i \leq (m + 1)x_h + c$ on $W_o$, for every $i \in w_o$.
When $i \notin  w_o$, $x_i \leq (m  + 1)x_h + c$ on $W_o$. 
Also one can require (2) holds.
The second case is obvious \tcb{(the homotopy is then the identity at all times)}.

In general let $W=\union_\nu W_\nu$ be the decomposition of $W$ into z-components. 
Define $W^o_\nu$ as above that and note that $W^o_\nu \cap W^o_{\nu'} = \varnothing$ if $\nu \not= \nu'$.     
Set $W^o=\union_\nu W_\nu^o$. It is a z-dense, z-open subset of $W$.
For each $\nu$, 
let $H_{\G, W_\nu} : [0,\infty] \times W'_{\nu} \to W'_{\nu}$
as above, with 
$W'_\nu=  (W^o_\nu \m [\v=\infty])  \union [\h=\infty]$.
The subsets $W'_\nu$ form a finite cover of
$W'$ by closed subsets. Hence the mappings $H_{\G, W_\nu}$ 
glue to a continuous
mapping $H_{\G, W} : [0,\infty] \times W'  \to W'$, \tcb{because they all agree with the trivial retraction on $[\h=\infty]$ which is the intersection of the sets $W'_\nu$}.
The process in \lemref{G-homotopy} and in the first paragraph of the present lemma being entirely canonical, \tcb{once an $A$-definable and $G$-invariant cell decomposition is chosen}, 
the retraction $H_{\G, W}$ obtained this way is $A$-definable and $G$-invariant. 
By construction the final image $W_0$ is definably compact.
For the additional items, for each $\nu$ one has open subsets  $W_{o, \nu}$ with corresponding
$m_{\nu}$ and $c_{\nu}$. One sets $W_o = \cup_{i \in w} W_{o, \nu}$,
$m = \max m_{ \nu}$, and $c = \max c_{\nu}$, which gives (1). By the construction  in \lemref{G-homotopy}
it is clear one can require (2) at the same time.
\eprf

\begin{lem}  \label{z-dense}  
Let $\Upsilon$ be an iso-definable  $\G$-internal  subset of $\std{V}$.
Let $\beta_0:  \std{V} \to [0, \infty]^{w_0}$ be a continuous $A$-pro-definable map, injective on  
 $\Upsilon$
as provided by \textup{\thmref{G-embed-1c}}.
Assume $\Upsilon$ is Zariski dense in $\std{V}$ in the sense of \textup{\ref{zt}}. 
Then we may enlarge $w_0$ to a finite $A$-definable set $w$ such that 
 $\beta_0$ factors through a continuous $A$-pro-definable map $\b: \std{V} \to [0, \infty]^{w}$  \textup{(}injective on  
 $\Upsilon$\textup{)} such that:
 \begin{enumerate}
 \item If
 $O$ is  a z-open z-dense subset of $\b (\Upsilon)$,  
then $\b \inv(O) \meet \Upsilon$ is a Zariski  open dense subset of $\Upsilon$. 
\item For any irreducible component $V_i$ of $V$, $\beta (\Upsilon \meet \std{V_i})$
is a z-component of $\beta (\Upsilon)$.
\end{enumerate}
\end{lem}

\prf    
Let $V_1, \ldots, V_r$ be the irreducible components  of $V$.  For each $V_j$,
let $x_j : V  \to [0, \infty]$ be a schematic distance function to $V_j$.
 Set 
$\beta( x) = (\beta_0 (x), x_1,\ldots,x_r)$.      
It follows from \lemref{infty},  that if $W$ is a z-closed subset of $[0, \infty]^{w}$,
then
$\b \inv (W)$ is Zariski closed. Thus, if   $Z \nsubseteq Y$ is z-closed (resp. z-open)
in $Y = \b (\Upsilon)$,
$\b \inv (Z) \meet \Upsilon$ is Zariski closed (resp. open) in $\Upsilon$.
Let us prove (1).
If $Z \nsubseteq Y$ is z-closed in $Y$ and contains no z-component of $Y$, suppose $\b \inv (Z)$ 
 contains some $\std{V_{j_0}} \meet \Upsilon$.  Then $\b \inv (Z) \union \union_{j \neq j_0} \std{V_j} $
contains $\Upsilon$, so $Z \union \union_{j \neq j_0} [x_j=\infty]$ contains $Y$.
It follows that $\union_{j \neq j_0} [x_j=\infty]$ contains $Y$ already.  But then
as $\std{V_j} = \b \inv ( [x_j=\infty])$ we have $\Upsilon \nsubseteq \union_{j \neq j_0} \std{V_j} $, contradicting the hypothesis on $\Upsilon$.
\tcb{For (2), let $C_j$, $j \in J$, denote the z-components of $Y$.
We have $\Upsilon \meet \std{V_i} \nsubset \cup_{j \in J} (\Upsilon \meet \beta^{-1} (C_j))$.
Since $\Upsilon \meet \std{V_i}$ is Zariski dense in $\std{V_i}$ and $V_i$ is irreducible,
it follows that, for some $j_i$,
$\Upsilon \meet \std{V_i}$ is contained in the Zariski closed set  $\Upsilon \meet \beta^{-1} (C_{j_i})$.
Hence, $\beta (\Upsilon \meet \std{V_i})$ is contained in $C_{j_i}$.
Since each $\beta (\Upsilon \meet \std{V_i})$ is z-closed  in $Y$ and the sets
$\beta (\Upsilon \meet \std{V_i})$ are mutually not included one in another, it 
follows that $\beta (\Upsilon \meet \std{V_i}) = C_{j_i}$.}
\eprf

\section{End of the proof}\label{ss10.6}

In  \ref{ss10.4}, we have constructed 
a continuous $A$-pro-definable retraction $\beta_{bc}$   from $\std{\tcb{V_0} \cup D_0} \to \Upsilon_{bc}$, sending $v$ to the final value of 
$t \mapsto  H_{bc}(t,v)$. 
Furthermore, by  \lemref{Upzdense}, $\Upsilon_{bc}$ is Zariski dense in $\std{V}$,
and we may assume
that, for every irreducible component
$V_i$ of $V$, $\Upsilon_{bc} \cap \std{V_i}$ is of  pure dimension $n=\dim(V)$.
By \thmref{G-embed-1c}, there exists a 
continuous $A$-pro-definable map $\beta:  \std{V} \to [0, \infty]^{w}$ for some finite $A$-definable set $w$, injective on  
 $\Upsilon_{bc}$. One denotes by $\alpha$ its  restriction to $\Upsilon_{bc}$.
After enlarging $w$, we may assume we are in the setting of   \ref{ss10.5},
in particular that  with the notation therein, $\v=x_v,\h=x_h$ for some $h,v \in w$.
Also,  after adding schematic distance functions to the irreducible components of $V$, we may assume that the conclusions of
 \lemref{z-dense} hold for $\beta$ and $\Upsilon_{bc}$.
We set $W = \alpha (\Upsilon_{bc})$ and
we define $W^o$, $W'$, $H_\G$, $W_0$, $W_o$, $m$ and $c$ as in \lemref{G-homotopy2}.

Note that $V \m V_0$ contains no irreducible component of $V$ \tcb{(recall $V_0$ is the preimage of $U_0$ in $V$)}. Indeed, if $V_i$ is an irreducible component of $V$,
$D_0 \meet V_i$ is nonempty of dimension $n - 1$ and $u$ restricts to a finite morphism
$D_0 \meet V_i  \to U$, thus $u (D_0 \meet V_i )$ contains $U_0$. 
By  \lemref{inflation} 
there exists an $A$-definable homotopy 
$H_{inf} : [0, \infty] \times \std{V} \rightarrow \std{V}
$
respecting the functions $\xi_i$ and the group action $G$ and fixing pointwise $\std{D_0}$
with image contained in
$\std{\tcb{V_0} \cup D_0} $.  
(In fact, by \lemref{separate} the image is contained in
$\std{Z}$ with $Z$ a v+g-closed bounded definable subset of
$V$ with $Z \meet \tcb{(V \m V_0)} \nsubseteq D_0$.)  
For each $i \in w$, set
$\phi_i = \min(x_i, (m+1) x_h   + c)$.
Note that, outside $D_0$,
the functions $\phi_i \circ \b$
are v+g-continuous with values in $\G$.
Furthermore, the functions $\phi_i$ are definable over a finite Galois extension of $A$ and permuted by the Galois group.
Thus, by \lemref{inflation}, we may also
require that the functions
$\phi_i \circ \b$ are preserved by $H_{inf}$ \tcb{away from $\std{D_0}$, hence, 
since $H_{inf}$ fixes pointwise $\std{D_0}$, that the functions
$\phi_i \circ \b$ are preserved by $H_{inf}$ everywhere}.
Recall $W_o$ is an open subset of $W$ containing $W_0 \m [\h = \infty ]$, so
$\alpha^{-1} (W_o)$ is open in $\Upsilon_{bc}$. Thus,
$\alpha^{-1} (W_o)$ has pure dimension $n=\dim(V)$.
Since  the restriction of $\phi_i$ to $W_o$ is just the $i$-th coordinate function, it follows from
\propref{maxdimabh} (2) that $\alpha \inv (W_o)$ is  fixed pointwise by $H_{inf}$.
Hence so is $\alpha \inv(W_0 \m [\h = \infty])$, and thus also $\alpha \inv(W_0)$.
 By construction  
$H_{inf}$ satisfies $(*)$.

\medskip

We will define $\tcb{h}$ as the composition (or concatenation) of homotopies
\[\tcb{h} = H_\G^\alpha \circ ( (H_{\widetilde{base}} \circ H_{curves}) \circ H_{inf}): I'' \times \std{V} \longrightarrow \std{V}\]
where $H^{\alpha}_{\Gamma}$ is to be constructed, and $I''$ denotes the generalized interval obtained by gluing $[\infty,0]$, $I'$ and $[0, \infty]$.
Being the composition of homotopies satisfying $(*)$, $\tcb{h}$ satisfies $(*)$.

Since the image of $H_{inf}$ is contained in the domain of $H_{bc}$, the first composition 
makes sense.

The set $W^o$  is a z-dense, z-open subset of $W$. Hence, by \lemref{z-dense} (1),
$\a \inv(W^o)$ is a Zariski  open dense subset of $\Upsilon_{bc}$. 
Let $O$ be a Zariski dense open subset of $V$ such that
$\std{O} \meet \Upsilon_{bc} = \a \inv(W^o)$.
By construction of $H_{inf}$, the image $I_{inf}$ of $H_{inf}$ is contained in $\std{O \union D_0}$.  Thus $\beta_{bc}(I_{inf})$ is 
a definably compact subset of $ \beta \inv(W') \meet \Upsilon_{bc}$.  Note that $\beta$ restricts to a homeomorphism $\alpha_1$
between this set and a definably compact subset $W_1$ of $W$.    
One sets  $H_\G^\alpha(t,x) = \alpha_1 \inv H_\G(t,\alpha_1(x))$:  in short, 
$H_\G^\alpha$ is $H_\G$ conjugated by $\alpha$, restricted to an appropriate definably compact set.  So
$\tcb{h}$ is well-defined by the above quadruple composition.   
    
Since $H_{inf}$ fixes  $\alpha \inv (W_0)$, and $W_0$ is the image of $H_\G$,
$H_{inf}$ fixes the image of $\tcb{h}$.   On the other hand $H_{bc}$ fixes $\Upsilon_{bc}$ and hence the subset 
$\alpha \inv (W_0) \nsubseteq \Upsilon_{bc}$.  Thus $\tcb{h}$ fixes its own image $\Upsilon = \alpha \inv (W_0)$.  
\tcb{It follows from  \propref{compactimage} that}
$\Upsilon$
is definably compact and  $\a$ is a homeomorphism from $\Upsilon$ to the definably compact subset $W_0$ of $  \G_\infty^w$.

We have thus constructed a homotopy
$\tcb{h} : I'' \times \std{V} \to \std{V}$
satisfying the statement of the theorem together with conditions (1), (2) and (4).
We shall now check that (3), \tcb{(5), (6)} and (7) also hold.

Let us start by checking (3), that is, $\tcb{h}$ is  Zariski generalizing, i.e. 
for any Zariski open subset $U$ of $\tcb{V}$, $\std{U}$  is   invariant under $\tcb{h}$.  
This property clearly holds for
the first three homotopies in the concatenation; let us check it for $H_\G^\alpha$.
By \corref{z-iso} it is enough to prove 
$H_\G$ is Zariski generalizing.
Consider 
a definable continuous function
$\eta : W' \to \G_{\infty} $ such that  $W' \m \eta \inv (\infty) \not= \varnothing$.
Pick a point $x$ in $W'$ with $\eta (x)$ finite.
By construction of $H_{\G}$, for some finite $t_0$, $H_{\G}(t_0, x)$ lies in $W_0$.
Thus, 
the function $t \mapsto \eta (H_{\G}(t, x))$ 
can only take finite values for finite $t$,
since a definable continuous function
$[0, t_0] \to \G_{\infty}$ which is nonconstant can take only finite values.

Let us now check (6), that is,    $\Upsilon$ is   the image of the set of simple points. Set $e = e_{I''}$.
\tcb{Let $p$ be a point in $\std{V}$.
Since $\Upsilon$ is iso-definable $\G$-internal, 
by orthogonality to $\Gamma$ there exists a definable subset $D$ of $V$ containing $p$ such that $h(e,x)=h(e,p)$ for every (simple) point $x$ of $D$. 
}


We now prove (5).  By \lemref{abh-char} (5), integrating a function into $\stda{V}$
on an element of $\stda{V}$ gives an element of $\stda{V}$.  We will use this repeatedly below.    In particular by (6), it suffices for (5) to show that  the image of the simple points lies in $\stda{V}$.  
 Now (5) is clear  for the inflation homotopy, as this homotopy is a finite cover of   the standard
affine homotopy $ I \times  \Aa^n \to \std{\Aa^n}$ (the image of a simple point being a tensor power of the image of a point on
$\Aa^1$).  By the remark on integration, precomposing with the inflation homotopy will not spoil (5).  Composing with a homotopy
taking place purely on the skeleton obviously does not add to the image of $h(e,V)$, as it adds no new points to this image.  It remains
to consider the inductive step.  
 Inductively, we may assume (5) holds for the skeleton
of the base homotopy.    In relative dimension one, any element of $\std{V/U}$ is in fact in $\stda{V}$.  Hence again by
 transitivity every element of $V$ moves through $\stda{V}$ throughout the homotopy.  
 
 It remains to prove  (7), i.e. that
given a finite family
 of closed irreducible subvarieties $W_i$ of $V$, one can assume
$\Upsilon  \meet \std{W_i}$ has pure dimension $\dim(W_i)$.
We already proved one can achieve each $\Upsilon_{bc} \cap \std{V_i}$ is of  pure dimension $n$.
It follows that  each $\a (\Upsilon_{bc} \cap \std{V_i}) $ is of pure dimension $n$.
By the conclusion of \lemref{z-dense} (2) which holds for
$\beta$ and $\Upsilon_{bc}$, the sets $\a (\Upsilon_{bc} \cap \std{V_i}) $
are the z-components of $W$. It follows from
 \lemref{G-homotopy2} (2) that one can achieve that 
 $\a (\Upsilon_{bc} \cap \std{V_i}) \meet W_0$ is of pure 
 dimension $n$. Since $\a$ restricts to a homeomorphism between
 $W_0$ and $\Upsilon$, it follows that each
$ \Upsilon \cap \std{V_i}$ is of pure dimension $n$.
With these choices,  for any $W_i$ of dimension $n$, 
 $\Upsilon  \meet \std{W_i}$ has pure dimension $\dim(W_i)$.
 Let us now deal with the
 case where some
 $W_i$ are of dimension $m_i < n$.
 \tcb{We may require all such $W_i$ are} contained in the hypersurface $D_0$ considered in   \ref{ss10.2}.
All reductions go through and
 when at the end of   \ref{ss10.2} we replace $V$ by $V_1$, it is enough to replace 
 $W_i$ by its strict transform.
The restriction $u_{W_i}$ of $u$ to $W_i$ is a finite morphism. Set $W_i' = u (W_i)$.
By  construction, the homotopies
$H^{\alpha}_{\Gamma}$, $H_{curves}$ and $H_{inf}$ fix pointwise the intersection of $\std{W_i}$ with their domains.
Now note that the pseudo-Galois morphism $U' \to U$ considered in \thmref{omin-finite} may be chosen so to factor through any given finite surjective morphism
$U'' \to U$. Thus, we may assume $D_0 \times_U U' \to U'$ is a \tcb{generically trivial covering}.
Let $W'_i$ be an irreducible component of $W_i \times_U U'$ and denote by $C_i$ its image under the projection \tcb{to} $U'$.
By the induction hypothesis,  we may require the  base homotopy $h'$ at the beginning of  \ref{ss10.4}
satisfies (7) for all $C_i$ associated to some  $W_i$ of dimension $< n$.
Let $\Upsilon'_i$ be the final image of $C_i$ under the retraction 
$h'$. By hypothesis it 
has pure dimension $m_i$. 
Since $W'_i \to C_i$ is \tcb{generically} an isomorphism \tcb{and $\Upsilon'_i$ is Zariski dense in $C_i$ by induction},  the same holds for
the preimage $\Upsilon''_i$ of $\Upsilon'_i$
in $\std{W'_i}$. The morphism
$\std{W'_i} \to \std{W_i}$ being continuous with finite fibers, it follows that the image
$\Upsilon_i$ of $\Upsilon''_i$ in $ \std{W_i}$ also has pure dimension $m_i$. 
By construction the final image of $ \std{W_i}$
under $H_{\widetilde{base}}$ is equal to $\Upsilon_i$, which proves (7).


This ends the proof of \thmref{1}.
 \qed 

\medskip

\begin{rem}\label{numberofintervals} 
In the proof of \thmref{1} one uses the induction hypothesis for the base $U$, lifted to a certain o-minimal cover (using the same generalized interval).  The homotopy on $U$ is (in a certain order) lifted and  composed with  three additional homotopies: inflation, 
  the relative curve homotopy, and the homotopy internal to $\G$.  Each of these use   the standard interval from $\infty$ to $0$
  (in reverse order, in the case of the homotopy internal to $\G$).
 %
The number $h(n)$ of basic intervals needed for an $n$-dimensional variety  thus satisfies $h(1)=1$,  $h(n+1) \leq  h(n)+3$, so $h(n) \leq 3n-2$. 

For a homotopy whose interval cannot be contracted to a standard one consider $\Pp^1 \times \Pp^1$.  With the natural
choice of fibering in curves, the proof of \thmref{1} will work even without the inflation homotopy. 
 It will
lead to an iterated homotopy to a point:  first
collapse to $\{\pt\} \times \Pp^1$, then to $\{\pt\} \times \{\pt\}$. 
 \end{rem}



\section{Variation in families}   \label{ss10.7}
 
Consider  a commutative diagram 
\begin{equation*}
\xymatrix{ X \ar[rr]^{h} \ar[rd] && Y \ar[ld] \\ & T}
\end{equation*}
of pro-definable maps, with $T$ a definable set.
We shall refer  to the family of maps $h_{\tau}: X_{\tau} \to Y_{\tau}$ obtained by restriction to fibers above $\tau \in T$   
as {\em uniformly pro-definable}.   \index{uniformly pro-definable}

Consider   a situation where $(V,X)=(V_{\tau},X_{\tau})$ are given uniformly in a parameter ${\tau}$, varying in a definable set $T$.  For each ${\tau}$, \thmref{1} guarantees the existence of a strong deformation retraction $h_{\tau}: I \times \std{X_{\tau}} \to \std{X_{\tau}}$,
and a definable homeomorphism $j_{\tau}: W_{\tau} \to h_{\tau}(e_I,\std{X_{\tau}})$, with $W_{\tau}$ a definable subset of $\G_\infty^{w(\tau)}$.
Such statements are often automatically uniform in the parameter ${\tau}$.   For instance if $X_{\tau},Y_{\tau}$ are uniformly definable families
of definable sets, and for each $a$ there exists an $a$-definable bijection $X_{\tcb{a}} \to Y_{\tcb{a}}$, then automatically there must
exist a uniformly definable bijection $h_{\tau}: X_{\tau} \to Y_{\tau}$.  Indeed if $H$ is the collection of all $\varnothing$-definable subsets of $X \times_T Y$,
then for any $a \in T$, for some $h \in H$, $h_a: X_a \to Y_a$ is a bijection.  \tcb{By compactness} the family of all formulas
asserting that $h_{\tau}$ is not a bijection $X_{\tau} \to Y_{\tau}$ is inconsistent.  Hence a finite subset is inconsistent; i.e. there exists a finite
set $h^1,\ldots,h^r \in H$ such that for any $a \in T$, for some $i \leq r$, $h^i_a$ is a bijection $X_a \to Y_a$.  
\tcr{For any $\tau$, let $i_0$ be the smallest $i \leq r$ such that  $h^i_{\tau}$ is a bijection, and
let $h_{\tau}: X_{\tau} \to Y_{\tau}$  be
equal to  $h^{i_0}_{\tau}$.}
More generally, if each $h_a$ has some property $P$ which is ind-definable (i.e. the family of all definable maps for which it holds is an ind-definable family),
then one can find $h$ such that each $h_\tau$ has this property.  
(See a fuller explanation in \cite{HK}, introductory section on compactness and gluing.)  

Here the pro-definable \tcb{map} $h_{\tau}$ is given by an infinite collection of definable maps, so compactness does not directly apply.  
 Nevertheless the theorem 
is uniform in the parameter ${\tau}$.  The reason is that $h_{\tau}$ is determined by its restriction to the simple points,
and on these, the homotopy moves along  $\stda{V}$, \tcb{which  is endowed  with a canonical strict ind-definable structure by  \ref{ssbertini}.}
We state this as a separate proposition.

 \begin{prop}  \label{1u} Let $V_{\tau}$ be a quasi-projective variety, $X_{\tau}$ a
  definable subset of $V_{\tau} \times \G_\infty^{\ell}$, definable uniformly in $\tau \in T$ over some base set $A$.  
Then there exists a uniformly pro-definable family $h_{\tau}: I \times \std{X_{\tau}} \to \std{X_{\tau}}$, a finite set $w(\tcb{\tau})$,
  a definable set $W_{\tau} \nsubseteq \G_\infty^{w(\tcb{\tau})}$, and
 $j_{\tau}: W_{\tau} \to h_{\tau}(0,\std{X_{\tau}})$,   pro-definable uniformly in ${\tau}$, such that for each $\tcb{\tau} \in T$, 
 $h_{\tau}$ is a   deformation retraction, and $j_{\tau}: W_{\tau} \to h_{\tau}(0,\std{X_{\tau}})$ is a pro-definable homeomorphism. 
  Moreover, \tcb{we may require} \textup{(1)} and \textup{(2)}  of  \textup{\thmref{1}}  to hold if the $\xi_i$ and the group action are given uniformly, as can   \textup{(4)}, \textup{(5)}, \textup{(6)} and \textup{(7)}.
\end{prop}
   
\prf    For any $a \in T$, we have  $h_a, j_a$ with the stated properties, by \thmref{1}.  By   \thmref{1} (5), $h_a$ restricts to $h_a ^{\#}:    V_a \times I \to   \stda{V_a}$.
Note that in principle $I=I_a$ depends on $a$.  However,  as $\dim(V_a)$ is bounded by some $m$, \tcb{it follows from  \remref{numberofintervals} that} $I_a$ is a union of  at most $\tcb{3m - 2}$ copies of $[0,\infty]$; extending the homotopy trivially to be constant to the left, we may assume it is a gluing of exactly   $\tcb{3m-2}$
copies of this interval, so it does not depend on $a$.
  We have: \begin{enumerate} 
\item[(1)]  Given finitely many  $A$-definable functions $\xi_i: \coprod_{\tau \in T} V_\tau  \to \G_{\infty}$, one can choose $h_a$
to respect the $\xi_i$, i.e. $\xi_i (h_a^{\#}(t, x) )=\xi_i(x)$ for all ${\tau}$.  
\item[(2)] Assume given, in addition, a finite algebraic group action on $V_a$ 
given uniformly in $a$.  Then the homotopy retraction can be chosen to be equivariant.
 \item[(4)] Let $x \in X$ and let $c=h_a^{\#}(e_I,x)$ be the final image of $x$. Also let $t \in I$, and $p=h(t,x)$.  Then for generic $y \models p$, $h_a^{\#}(e_I,y)=c$; i.e.
 $\models (dp_y)h_a^{\#}(e_I,y)=c$. 
\item[(7)]  Each irreducible component $V'$ is left invariant by $h_a^{\#}$; and if $X\meet V'$ contains an open subset of $V'$, then $h_a^{\#}(0,V')$  has pure dimension equal to $\dim(V')$.
\item[(5$'$)] $h_a^{\#}$ extends to a homotopy $H_a: \std{X_{a}} \to \std{X_a}$.
\end{enumerate}

Now the validity of (5$'$) for $h_a^{\#}$ is an ind-definable property of $a$, by  \propref{def-vg-crit},  and  (1), (2) and (7) are obviously ind-definable (using the classical fact that the irreducible components of $V_a$ are $\ACF$-definable uniformly in $a$).
Property (4) is also stated in an ind-definable way.
 
Hence by the compactness \tcb{and} gluing argument mentioned above, one can find a uniformly definable family $h_\tau$ with the same properties.  Now
let $H_\tau(p) = \int_{x \models p} h_\tau(x)$.  By (5$'$), this is a homotopy $H_\tau: \std{X_{\tau}} \to \std{X_\tau}$.   Property (5) of \thmref{1} holds by definition.  Property (6) is proved
in the same way as in \thmref{1}.   
\eprf

\begin{rem} 
We proved above that irreducible components are preserved, but not the full  Zariski generalization property \thmref{1} (3), as it is not  an ind-definable property on the face of it.   It can still
be achieved uniformly; this can be seen in one of two ways:

-   either by   following the proof of (3), carrying the parameter $\tau$ along; 

- or else by proving that a stronger ind-definable 
property holds; namely that there is a uniformly definable family of varieties, such that the Zariski closure of $h(x,t)$ is an element of this family, and is increasing with $t$ along $I$.  In the case of a definably compact set $X$ contained in the smooth locus of $V$, the proof of \thmref{1simp} gives this in a very simple form:
the Zariski closure of  $h(x,t)$ is $\{x\}$ if $t=\infty$, and equals $V$ otherwise.  
\end{rem}

  
\chapter{The smooth case}\label{secns}

{\small \noindent \textbf{Summary.} 
In this chapter we examine the simplifications occuring in the proof of the main theorem in the smooth case.
We also note the birational character of the definable homotopy type in \remref{birational}.
\par\bigskip}

\section{Statement}
 For definable sets avoiding the singular locus it is possible to prove the following variant of
 \thmref{1}.  The proof
uses the same ingredients but is considerably simpler in that only birational versions of most parts of the construction are required.
 For clauses (1), (2) and (4), the homotopy internal to $\G$ is not required; and a global inflation homotopy is applied only once, rather than iterated
 at each dimension.    For clause (3), a single final use of the $\G$-homotopy is added.

Given  an algebraic variety $V$ over a field,  one denotes by  $V_{\mathrm{sing}}$ its singular locus, i.e. its nonsmooth locus.


\begin{thm} \label{1simp} Let $V$ be a quasi-projective variety  over a valued field $F$ and let $X$ be a v-\tcb{clopen $F$-definable subset of $V \m V_{\mathrm{sing}}$}.
  Then   there exists an $F$-definable homotopy
 $h : I \times \std{X}\to \std{X}$ between the identity and a continuous map to  an $F$-iso-definable subset  definably 
homeomorphic to a definable subset of $ w' \times \G^w$, for some finite $F$-definable sets $w$ and $w'$.

Moreover, one can require the following additional properties for $h$:
\begin{enumerate}
\item Given finitely many v-continuous $F$-definable functions $\xi_i: X \to \G$, one can choose $h$
to respect the $\xi_i$, i.e. $\xi_i (h(t, x) )=\xi_i(x)$ for all $\tcb{t}$.   
\item Assume given, in addition, a finite algebraic group action on $V$.  Then the homotopy   can be chosen to be equivariant.
\item If $\std{X}$ is definably compact,   $h$ can be taken to be a deformation retraction.
\item   Clauses \textup{(3)} to \textup{(6)} of \textup{\thmref{1}} hold.  Also, if $V$ has dimension $d$ at each point $x \in X$, then  each point $q$ of the image of $h$ is strongly stably dominated with $\dim(q) = d$.  \end{enumerate} 

In particular this holds for $X=V$ when $V$ is smooth. 
 \end{thm}

\tcb{By \remref{interv}}, when $\std{X}$ is definably compact the conclusion is stronger than \thmref{1} in that the interval is the standard interval $[0,\infty]$.  
If $\std{X}$ is not definably compact, the conclusion is also weaker in that we do not assert that the final image is fixed by the homotopy.

The finite set $w'$ can be dispensed with if $\G(F) \neq (0)$, or if $\std{X}$ is connected, but not otherwise, as can be seen by considering the case when $X$ is finite.  
 Indeed, when $\G(F) = (0)$
the only nonempty finite $F$-definable subset of $\G^n$ is $\{0\}$, but one can
have arbitrarily large finite $F$-definable subsets in $\G^n_{\infty}$ for $n$ large enough.

\section{Proof and remarks}
The proof depends on two lemmas.  The first recaps the proof of \thmref{1}, but on a Zariski dense open set $\tcb{V_0}$ only.  The second uses smoothness to
enable  a stronger form of inflation, serving to move into $\tcb{V_0}$.

While the theorem requires the functions $\xi_i$ to be v-continuous, this need not be assumed in \lemref{1simp1} since any definable function
 is v-continuous on some Zariski dense open set.   But then $X$ need not be explicitly mentioned, since one can add the \tcb{valuation of the} characteristic function of $X$
to the list of $\xi_i$.   The proof of this lemma uses only an iteration of the curves homotopy, without inflation or the $\G$-homotopy.

\begin{lem}\label{1simp1}  Let $V$ be a quasi-projective variety defined over $F$.     Then   there exists a 
Zariski open dense subset $V_0$ of $V$,  and 
an $F$-definable deformation retraction 
 $H : I \times \std{V_0}\to \std{V_0}$ whose  image is an $F$-iso-definable subset $S_0$,  definably 
homeomorphic to an $F$-definable subset of $w' \times  \G^w$, for some finite $F$-definable sets $w'$ and $w$.   

Moreover: 

\begin{enumerate}
\item Given finitely many $F$-definable functions $\xi_i: V \to \G$, one can choose $h$ 
to respect the $\xi_i$, i.e. $\xi_i (h(t, x) )=\xi_i(x)$ for all $\tcb{t}$.   
\item Assume given, in addition, a finite algebraic group action on $V$.
 Then $V_0$ and the \tcb{deformation} retraction can be chosen to be equivariant.
\end{enumerate}
 \end{lem}

\prf Find a Zariski open $V_1$ with  $\dim(V \m V_1)< \dim(V)$,
 and a morphism $\pi: V_1 \to U$, \tcb{with $U$ normal}, whose fibers are curves.  Let $h_{curves}$ be
the homotopy described in  \ref{ss10.3}.  It is \tcb{v+g-}continuous outside some  subvariety \tcb{$W$} of $U$
with $\dim(\tcb{W})<\dim(U)$;  replace $V_1$ by $V_1= \tcb{V \m \pi \inv(W)}$.  So $h_{curves}$ is \tcb{v+g-}continuous on $V_1$ and its image $S_1$ is  relatively $\G$-internal over $U$.  
By (a greatly simplified version of) the results of  \ref{ss6.4}, 
\tcb{there exists a finite pseudo-Galois covering $U' \to U$, such that, setting $V'_1 = U' \times_U V_1$ and $S'_1 = U' \times_U S_1$,
there exists a definable v+g-continuous morphism $g : V'_1 \to U'\times \Gamma^n_{\infty}$ such that $g$ induces a homeomorphism between
$\std{S'_1}$ and its image in $\std{U'}\times \Gamma^n_{\infty}$. 
}


\begin{claim} \tcb{After replacing $V_1$ by a Zariski dense open subset,  one may assume there exists a definable isomorphism
between $S'_1$ and an iso-definable  subset of $\std{U'} \times \{1,\ldots,N\} \times \G^n$ relatively $\G$-internal over $U'$, 
for some positive integers $N$, $n$,
inducing a homeomorphism between
$\std{S'_1}$ and its image in $\std{U'} \times \{1,\ldots,N\} \times \G^n$.}
\end{claim}
\begin{proof}[Proof of the claim]  After removing a nowhere dense subvariety, we may assume $V_1$ is a disjoint union of irreducible components, and work within each 
component separately.  \tcb{So we may assume $V_1$ is irreducible.
We may also assume $V'_1$ is irreducible.
The set of points of $S'_1$ which are mapped by  $g$ to $U' \times \G^n$ is Zariski open in $S_1$; thus, if it is nonempty it must be dense in $V'_1$, and after shrinking $V_1$ again we may assume
$S'_1$ maps to $U' \times \G^n$.
Otherwise $S'_1$ maps to $U' \times (\G_\infty^n \m \G^n)$.}
 In this case we can remove a proper subvariety and decompose the rest into finitely many algebraic pieces, each mapping into one   hyperplane at $\infty$ of $\G_\infty^n$.
 \tcb{Then one  concludes the proof by induction on $n$.}    \end{proof} 

We may thus assume  \tcb{there exists a definable isomorphism
between $S'_1$ and an iso-definable  subset of $\std{U'} \times \{1,\ldots,N\} \times \G^n$ relatively $\G$-internal over $U'$, 
for some positive integers $N$, $n$,
inducing a homeomorphism between
$\std{S'_1}$ and its image in $\std{U'} \times \{1,\ldots,N\} \times \G^n$.}
\tcb{By induction, there exists a Zariski dense open $U'_0$ of $U_0$ and an $F$-definable deformation retraction
$h' : I  \times \std{U'_0} \to \std{U'_0}$ satisfying the conclusions of the lemma. Furthermore we may assume the
pseudo-Galois covering $U' \to U$ restricts to a pseudo-Galois covering $U'_0 \to U_0$ for some dense open subset $U_0$ of $U$
and that $h'$ is the lifting of an $F$-definable deformation retraction
$h : I  \times \std{U_0} \to \std{U_0}$ satisfying the conclusions of the lemma. Set $V_0 = \pi^{-1} (U_0)$.
Using  \thmref{omin-finite} as in  \ref{ss10.4}, we may arrange that $h$  lifts to a homotopy $H_{\widetilde{base}} : I \times (\std{S_1} \meet \std{V_0}) \to \std{S_1} \meet \std{V_0}$.}
  The homotopies can be taken to meet conditions (1) and (2).  Composing, we obtain a deformation retraction of $\std{V_0}$ to a subset $S$, 
  and a homeomorphism $\a: S \to Z \nsubseteq  \{1,\ldots,M\} \times \G ^m$, defined over $\acl(A)$.  We may assume $M>1$.
  As in  \thmref{G-embed-1c} we can obtain an $A$-definable homeomorphism into $(\{1,\ldots,M\} \times \G^m)^w $.   
\eprf

\begin{lem} \label{1simp2} Let $V$ be a subvariety of $\Pp^n$, and let $a \in V$ be a smooth point.  Then the standard metric on $\Pp^n$ restricts to a good metric on some v-open neighborhood of $a$ in $V$.  
\end{lem}  

\prf  For sufficiently large $\alpha$, the set of points of distance $\geq \alpha$ from $a$ may be represented
as the $\Oo$-points of a scheme over $\Oo$ with
good reduction, whose special fiber is irreducible, in fact a linear variety. 

This can be done as follows.  We may assume $V \nsubseteq \Aa^n$, and $a = (0)$.  As $a$ is smooth, $V$ is a complete intersection near $0$,
and we may localize further and assume
 it is cut out by polynomials $f,\ldots, h$ in affine coordinates $x_1,\ldots,x_n$, whose number $\ell$ is the codimension of $V$.

We can write $f=f_1+f_2$, where $f_1$ is linear and $f_2$ consists of higher degree terms and similarly for $g,\ldots,h$.   
The vectors $f_1,\ldots,h_1$ generate an $\ell$-dimensional subspace of the space with basis $x_1,\ldots,x_n$.

By performing   row   operations, we may assume $f_1,\ldots,h_1$ have coefficients in $\Oo$, and further
that their coefficient vectors generate a lattice of rank $\ell$ in $\Oo^n$.  (In fact, permuting 
 the variables if necessary, and performing   row   operations,  we can arrange that    modulo $\Oo x_{\ell+1} + \ldots + \Oo x_n$ we have $f_1=x_1, \ldots, h_1= x_\ell$.)
 
Of course, the nonlinear coefficients of   $f,\ldots,h$ have coefficients in the field $K$, some having valuation as negative as $-\val(c)$ say,
where $c \in \Oo$.   Let $F(x) = c \inv  f(cx), \ldots, H(x)= c \inv h(cx)$.    The intersection of $V$ with  $c \Oo^n$   is isomorphic to the intersection of $(F,\ldots,H)$ with $\Oo^n$.  But it is clear that $F,\ldots,H$ have coefficients in $\Oo$, and that they cut out a smooth scheme $S_c$ over $\Oo$.

For this $c$ or for any $c'$ with $\val(c') \geq \val(c)$,  $S_c(K)$ clearly admits a unique generic type, dominated by the 
generic type of the linear variety $S_c(k)$, via the residue map.
  \eprf

\begin{proof}[Proof of \thmref{1simp}]    \tcb{Let $V_0$, $H$ and $S_0$ be as provided by  \lemref{1simp1}.}  
As in the first few lines of the proof of \thmref{1}, we may choose a projective embedding equivariant with respect to the finite group action of (2).
\tcb{By \lemref{1simp2}, for any $x \in X$,  the standard  metric $d$ on $\Pp^n$ restricts to a good metric on 
some v-open neighborhood of $x$. Thus, there exists a definable function $\rho : X \to [0; \infty)$   
which is locally bounded and such that for any 
$x \in X$ and any $t \geq \rho(x)$, 
$B(x;d,t)$ is affine and has a unique generic type which we shall denote by $p (x, t)$.
Since $X$ is v-open and the functions $\xi_i$ are v-continuous, we may assume that for 
$t \geq \rho(x)$, $p (x, t)$ lies in $\std{X}$ and $t \mapsto \xi_i (p (x, t))$ is constant.
Since $X$ is v-closed, by  \lemref{majorize}
there exists a v+g-continuous function $g : X \to [0; \infty)$ such that
for every $x \in X$, $\rho (x) \leq g (x)$.
For $t \in [0, \infty]$ and $x \in X$, set  $H_{inf} (x, t) = p (x, \max(t, g(x)))$. It is a v+g-continuous  definable function $[0, \infty] \times X \to \std{X}$
which extends to a homotopy
$H_{inf} : [0, \infty] \times \std{X} \to \std{X}$.}
Note that the image of $H_{inf}$ is contained in $\std{V_0}$.  Define $h$ as the composition of $H$ and $H_{inf}$.


For clause (3), to ensure that the composition is also a deformation retraction, we compose with an additional homotopy internal to $\G$ as 
in \thmref{1}. 

The verification that the image of closed points is strongly stably dominated is as in \thmref{1}; moreover the homotopies of \lemref{1simp2}
are Zariski generalizing, while the inflation homotopy \lemref{1simp1} has final image consisting of points of maximal dimension; this proves
(4).
\end{proof}


\begin{rem}\label{interv}  \tcb{Under the hypotheses of \thmref{1simp}}, if $\std{X}$ is definably compact,   the  interval $I$ can in fact be taken to be $[0,\infty]$.  We sketch the argument.
   The proof above yields a composition of homotopies $H_{\G} \circ H_m \circ \cdots \cdot H_1 \circ H_{inf}$, where the $H_i$ for $i =1,\ldots,m$ are relative curve
homotopies using intervals $[0,\infty]$ oriented from $\infty$ to $0$, $H_{inf}$ uses a similar interval $[0,\infty]$, and
$H_\G$, the homotopy internal to $\G$, uses the same interval oriented in the opposite direction. 

For $k=0,\ldots,m$, set $H^k = H_{k} \circ \cdots \circ H_{inf}$, with $H^0 = H_{inf}$, and denote by 
$S_k$ the final image of $H^k$.   We wish to show by   induction on $k$ that the
 interval of $H^k$ can be contracted to a standard interval $[0,\infty]$.  It suffices to replace $H_k$ by a homotopy
 whose time interval is a closed interval in $\G$, by showing that for some $\alpha_k$, for all $t > \alpha_k$ and all $x \in S_k$,  $H_k(t,x)=x$.  
 
If we write $X$ as a finite union of definably compact subsets $X_\nu$ of affine open subsets of $V$, and show that the statement holds
for each $X_\nu$, then it holds for $X$.  In this way we can reduce to an affine situation.  
 
  Each $a \in S_k$ is a    strongly stably dominated point.   It is possible to find an \'etale neighborhood $V'$ of $X$ in $V$ and morphisms
  $f:V' \to W$ and  $g: W \to U$ such that 
  $W \nsubset U \times \Aa^1$ and $g$ is the projection, $(g \circ f)_* (a) = a'$ is a generically stable type on $U$,
  $a = \int_{a'} h$  where $h$ is a definable map   $U \to \std{\Aa^1}$, 
  and $H_k$ is compatible with the standard homotopy on $\Aa^1$, relative to $U$.  
   The decomposition $f: V' \to W$ and $g:W \to U$ is part of the construction of the homotopy, 
  while the integral decomposition of $a$ over $a'$ follows from the strong stable domination of $a$ (cf.  \tcb{\propref{stda-char-2}}).
  Moreover, as $a$ lies in the final image
  of the inflation homotopy, the Zariski closure has dimension $\dim(V)$, and it follows that  one can take $h(u)$ to be the generic type
  of a closed ball which is not reduced to a point.  Moreover the radius of this ball is a continuous definable function on $U$.
  By definable compactness, it is bounded above on $S_k$, say $\leq \alpha_k$.    It follows that $H_k(t,x)=x$ for $t > \alpha_k, x \in S_k$.
  This allows us to collapse the interval of $H^m$ to a standard interval $[0,\infty]$.  
   
   Recall now the homotopy within $\G$.   The composed curve homotopies $H_m \circ \cdots \cdot H_1$ act on a certain affine $\tilde{V}$,
   with final image $\tilde{S} \cong \Omega$; $\Omega$ is a definable subset of $\G_\infty^w$.  The homotopy $H_\G$ takes
   $\Omega$ to a definably compact set $S_\G$.  At this point, $H_{inf}$ is chosen so as to fix $S_\G$.   The final image
   of the composition $H_m \circ \cdots \cdot H_1 \circ H_{inf}$ is the definably compact set $S_m$.  
   Now $H_\G$ is applied,
   with time interval $[0,\infty]$.  But $H_\G$ moves each point of $S_m$ into $S_\G$ in finite time.  Since $S_m$ is definably
   compact, there is some time $t_\G$ such that by time $t_\G$, each point of $S_m$ is moved by $H_\G$ into $S_\G$ (and then frozen).
   Thus if $H_\G'$ is the restriction of $H_\G$ to time interval $[0,t_\G]$, then the composition $H_\G \circ H_m \circ \cdots \cdot H_1 \circ H_{inf}$
   also has final image fixed by $H_{inf}$ and by each $H_i$ and $H_\G$.  This gives a homotopy whose time interval is the concatenation of
   $[0,\infty] $ with $[0,t_\G]$; this is again isomorphic to $[0,\infty]$.    
 \end{rem}

\begin{remark}[A birational invariant] \label{birational} It follows from the proof of \thmref{1simp} 
that the definable homotopy type of
$\std{V \m V_{\mathrm{sing}}}$ 
(or more generally of $\std{X \m V_{\mathrm{sing}}}$ when $X$ is a v-\tcb{clopen} definable subset of $V$) 
is a birational invariant of $V$  (of the pair $(V,X)$).  
This rather curiously complements a theorem of Thuillier \cite{thuillier}.

As a referee pointed out, this remark requires only the inflation homotopy.  Let us spell this out without $X$, to simplify notation. 
 It suffices to show that if $U$ is a smooth variety and $W$ a Zariski dense open subset, then $\std{U}$ and $ \std{W}$ are pro-definably homotopy equivalent.  
 \tcb{Indeed, let $H=H_{inf}$ be the inflation homotopy on $U$ as considered in the proof of \thmref{1simp} and denote by $Z$ its image. 
 Note that  if   $z  \in Z$, then $H(t,z)=z$ for all $t$. By density of simple points we may assume $z=H(t_0,z_0)$ with $z_0$ a simple point and $t_0$ the stopping time provided by the cut-off function.
 It is enough to prove that if $M$ is a base, $z=H(t_0,z_0)$, $x \models z|M$,
and $y$ realizes the generic type over $M(x)$, of the ball of valuative radius $t$ around $x$, then $y \models z|M$.  Indeed
$y$ still falls in the ball of valuative radius $t_0=\val(c_0)$ around $z_0$, and has the same image as $x$ under the dominating function $\res(c_0 \inv y)$.  
 Thus, $H$ provides a deformation retraction of $\std{U}$ to $Z$. Since 
 $Z \nsubseteq  \std{W}$, the restriction of $H$ to $\std{W}$ provides a deformation retraction of $\std{W}$ to $Z$. Thus,
$\std{U}$ and $\std{W}$ are both  definably homotopy equivalent to $Z$.}   
\end{remark}


\chapter{An equivalence of categories}\label{secec}

{\small \noindent \textbf{Summary.}
In this chapter we deduce from \thmref{1} an equivalence of categories between the homotopy category
 of definable subsets of quasi-projective varieties over a given valued field and the homotopy category
of definable subsets of some $\G_\infty^w$.
\par\bigskip}

\section{Statement of the equivalence of categories}
 \def\io{\iota}Let $F$ be a valued field.  \tcb{We fix a base set $A \subset \VF \cup \G$ with $F = \VF (A)$.}
Let $V$ be an algebraic variety over $F$; by a {\em semi-algebraic} subset of $\std{V}$ we mean  a subset of the form $\std{X}$, where $X$ is a definable \index{semi-algebraic subset of $\std{V}$}
subset of $V$. \tcb{If $X$ is \tcb{$A$-definable}, we say $\std{X}$
is $A$-{\em semi-algebraic}.}

Let $C_{\VF}$ be the category \nomenclature{$C_{\VF}$}{the category of semi-algebraic subsets}
 of semi-algebraic subsets of $\std{V}$, with $V$ a quasi-projective  variety over $F$; the morphisms are  pro-definable continuous maps. 
 We could also say that the objects are definable subsets of $V$, but the morphisms $U \to U'$ are still 
 pro-definable continuous maps $\std{U} \to \std{U'}$.

 Let $C_\G$ be the category of definable subsets  $X$ of $\G_\infty^w$ (for various definable finite sets $w$), \nomenclature{$C_{\G}$}{the category of $\G$-sets}
with definable continuous maps.    Any such map is piecewise given by an element of  $\GL_w(\Qq)$ composed with a translation, 
and with coordinate projections and inclusions $x \mapsto (x, \infty)$ and $x \mapsto (x,0)$.
Let $C_\G^{i}$ be the category of  topologically  $\G$-internal  subsets  $X$ of $\std{V}$, for various varieties $V$, \nomenclature{$C_{\G}^i$}{the category of topologically $\G$-internal sets}
with continuous definable  maps.

These categories can be viewed as ind-pro definable:  more precisely 
$\mathrm{Ob}_C$ is an ind-definable set, and for $X,Y \in \mathrm{Ob}_C$, $\mathrm{Mor}(X,Y)$ is a pro-ind-definable set.  
The three categories admit natural functors to the category  $\mathrm{TOP}$ of  topological spaces with \nomenclature{$\mathrm{TOP}$}{the category of topological spaces}
continuous maps. 
But usually we will be interested only in the subcategory consisting of $A$-definable objects and morphisms.
It can be defined in the same way as in the first place, only replacing \tcr{definability  by $A$-definability}. 
We shall denote these categories by
$C^A_{\VF}$, $C^A_\G$ and $C^{A, i}_\G$.


There is a natural functor $\iota: \tcb{C^A_\G} \to \tcb{C^{A,i}_\G}$, commuting with the natural  functors to $\mathrm{\mathrm{TOP}}$;  namely,
given \tcb{an $A$-definable subset} $X \nsubseteq \G_\infty^n$, let $\io(X) = \{p_\g: \g \in X \}$, where $p_\g$ is as defined above  \lemref{gtop0}.
By  \lemref{gtop0} and \lemref{prodtop}, the map $\g \mapsto p_\g$ induces a homeomorphism $X \to \io(X)$.   

\begin{lem} \label{ec1}The functor   $\iota: \tcb{C^A_\G} \to \tcb{C^{A,i}_\G}$ is an equivalence of categories.    \end{lem}

\prf It is clear that the functor is fully faithful.  Essential surjectivity follows from \thmref{G-embed-0} \tcb{and \remref{rem:basetopint}}.   
\eprf

We now consider the corresponding homotopy categories $HC_{\VF}^{\tcb{A}}$,   $H\tcb{C^A_\G}$ and $HC_\G^{\tcb{A}, i}$.  These categories have the same objects as the original ones, but the morphisms are
factored out by  (strong) homotopy equivalence.   Namely two morphisms $f$  and $g$ from $X$ to $Y$ are identified if there exists an \tcb{$A$-definable}  generalized interval
$I= \tcb{[i_I,e_I]}$ and a continuous 
\tcb{$A$-}pro-definable
map $h: X \times I \to Y$ with $\tcb{h_{i_I}}=f$ and $\tcb{h_{e_I}}=g$.  One may verify
that composition preserves equivalence; the image of $\mathrm{Id}_X$ is the identity morphism in the category.
\nomenclature{$HC^A_{\VF}$}{homotopy category of $C^A_{\VF}$} \nomenclature{$HC^A_{\G}$}{homotopy category of $C^A_{\G}$} \nomenclature{$HC_{\G}^{A,i}$}{homotopy category of $C_{\G} ^{A,i}$}

 The equivalence $\io$ above induces an equivalence
$H\tcb{C^A_\G} \to HC^{\tcb{A}, i}_\G$.   As a reader pointed out, the same retraction was considered by Berkovich in the setting of Berkovich spaces.

\begin{lem}  \label{ec2}  \tcb{Let $w$ be an $A$-definable finite set.} For \tcb{an $A$-}definable subset $X \nsubseteq \G_\infty^w$, let $C(X) =  \{x \in \Aa^w: \val(x) \in X \}$.  Then
the inclusion $\io(X) \nsubseteq \std{C(X)}$ is a homotopy equivalence.
\end{lem}

\prf For $t \in [0, \infty]$ one sets $H_0=G_m(\Oo)$, $H_{\infty}=\{1\}$,
and for $t > 0$, with $t = \val (a)$, 
$H_t$ denotes the subgroup  $1+a \Oo$ of $G_m(\Oo)$.
For $x$ in $C(X)$
one denotes by $p(H_t x)$  
the unique $H_t$-translation invariant stably dominated type on $H_t x$.
In this way one defines \tcb{an $A$-definable} homotopy
$[0,\infty]  \times C (X) \to \std{C(X)}$ by sending
$(x, t)$ to 
$p(H_tx)$,
whose canonical extension
$[0,\infty]  \times \std{C (X)} \to \std{C(X)}$
is a  deformation retraction with image 
$\io(X)$.
\eprf

\tcb{We shall prove the following statement in the next section:}
\begin{thm} \label{catequiv} The categories $HC^{\tcb{A}}_{\G}$ and $HC^{\tcb{A}}_{\VF}$ are  equivalent by an  equivalence
respecting the subcategories of definably compact objects.
\end{thm}

\section{Proof of the equivalence of categories}

To prove \thmref{catequiv}, we introduce a category $C_2^{\tcb{A}}$ 
defined as follows.
Objects of $C_2^{\tcb{A}}$ are pairs $(X,\pi)$, with $X$ an object of $C^{\tcb{A}}_{\VF}$ and $\pi: X \to X$
a continuous \tcb{$A$-}definable retraction with \tcb{topologically} $\G$-internal image, which is homotopic to the
identity
$\Id : X \to X$ via \tcb{an $A$-definable} homotopy $h :  I \times X  \to X$
with $h_{i_I} = \Id$, $h_{e_I} = \pi$, and $\pi \circ h_t = h_t \circ \pi = \pi$
for every $t$ in $I$.
A morphism $f: (X,\pi) \to (X',\pi')$ in $C_2^{\tcb{A}}$ is a continuous \tcb{$A$-}definable map $f: X \to X'$ such that $f \circ \pi= \pi' \circ f$.
  We define a homotopy equivalence
relation $\sim_2$ on $\mathrm{Mor}_{C_2^{\tcb{A}}} ((X,\pi),(X',\pi'))$ by  $f \sim_2 g$ if there exists a continuous \tcb{$A$-}definable
$h: I \times X   \to X'$, with $h_{i_I}=f$ and $h_{e_I}=g$, such that $h_t \circ \pi = \pi' \circ h_t$ for all $t$.    
Note that
$f \sim_2 f \circ \pi$ and
$f \sim_2 \pi' \circ f$. In particular,
$f \sim_2 \pi' \circ f \circ \pi$.
Again one
checks that this is a congruence and that one can define a quotient category denoted by  $HC_2^{\tcb{A}}$. 

There is an obvious functor $C_2^{\tcb{A}} \to C^{\tcb{A}}_{\VF}$ forgetting $\pi$, and also a functor $C_2^{\tcb{A}} \to C_\G^{\tcb{A, i}}$, mapping 
$(X, \pi)$ to $\pi(X)$.  One checks that the natural maps on morphisms are well-defined  and that they induce
functors $HC_2^{\tcb{A}} \to HC^{\tcb{A}}_{\VF}$ and $HC_2^{\tcb{A}} \to HC_\G^{\tcb{A, i}}$.
  To prove
the theorem, it suffices therefore to prove, keeping in mind \lemref{ec1}, that each of these two functors is essentially surjective and fully faithful, and to observe that they restrict to functors on the definably compact objects, essentially surjective
on  definably compact objects.

 (If the categories are viewed as ind-pro-definable,  these functors are  
 morphisms of ind-pro-definable objects, but we do not claim that a direct definable equivalence exists.) 

\begin{lem} \label{ec4}The functor $HC_2^{\tcb{A}} \to HC^{\tcb{A}}_{\VF}$ is surjective on objects, and fully faithful. \end{lem}

\prf  Surjectivity on objects is given by \thmref{1}.   
Consider
  $(X,\pi)$ and $ (X',\pi')$ in $\mathrm{Ob} HC_2^{\tcb{A}} = \mathrm{Ob} C_2^{\tcb{A}}$.  Let  $f: X \to X' $ be a morphism of $C^{\tcb{A}}_{\VF}$.  
  Then the composition $\pi' \circ f  \circ \pi$ is homotopy equivalent to $f$, since $\pi \sim \mathrm{Id}_X$ and $\pi' \sim \mathrm{Id}_{X'}$, and is 
a morphism of $C_2^{\tcb{A}}$.  This proves surjectivity of $\mathrm{Mor}_{HC_2^{\tcb{A}}}((X,\pi), (X',\pi')) \to \mathrm{Mor}_{HC^{\tcb{A}}_{\VF}}(X,X')$.
For injectivity,
let
$f, g : (X, \pi) \to (X', \pi')$ with $f \sim g$ in $C^{\tcb{A}}_{\VF}$.
Thus, $\pi' \circ f \circ \pi$ and
$\pi' \circ g \circ \pi$ are homotopic in $C_2^{\tcb{A}}$.
Since
$f \sim_2 \pi' \circ f \circ \pi$ and
$g\sim_2 \pi' \circ g \circ \pi$, it
follows that $f \sim_2 g$.
\eprf

\begin{lem} \label{ec5} The functor  $HC_2^{\tcb{A}} \to HC_{\G}^{\tcb{A, i}}$ is essentially surjective and fully faithful. \end{lem}

\prf   To prove essential surjectivity it suffices to consider objects of the form $\io(X)$,
with $X \in \mathrm{Ob} \tcb{C^A_\G}$.  For these, \lemref{ec2} does the job. 

Let    $(X,\pi), (X',\pi') \in \mathrm{Ob} HC_2^{\tcb{A}} = \mathrm{Ob} C_2^{\tcb{A}}$.   Let $g: \pi(X) \to \pi'(X')$ be a morphism of $C_\G^{\tcb{A, i}}$.  
Then $g \circ \pi: X \to X'$ is a morphism of $C_2^{\tcb{A}}$.
So even $\mathrm{Mor}_{C_2^{\tcb{A}}}((X,\pi), (X',\pi')) \to \mathrm{Mor}_{C_\G^{\tcb{A,i}}}(X,X')$ is surjective.

To prove injectivity, suppose $g_1$ and $ g_2 : X \to X'$ are $C_2^{\tcb{A}}$-morphisms, and
$h: I \times \pi (X) \to \pi' (X')$ is a homotopy between  $g_1 | \pi(X) $ and $g_2 | \pi(X)$.  
We wish to show that $g_1$ and $ g_2$ are $C_2^{\tcb{A}}$-homotopic.  Now for $i=1,2$,  
$g_i$ and $\pi' \circ g_i \circ \pi$ have the same
image in $\mathrm{Mor} (\pi (X), \pi' (X'))$, and there is a homotopy between $g_i$ and 
$\pi' \circ g_i \circ \pi$, $ i = 1, 2$,
as remarked before.  So we may
assume $g_i = \pi' \circ g_i \circ \pi$ for $i=1,2$.  
Define
$H: I \times X \to X'$ by $H(t, x) = \pi' \circ h(t, \pi(x))$.  This  is a $C_2^{\tcb{A}}$-homotopy
between $g_1$ and $g_2$
showing that $g_1$ and $g_2$ have the same class as  morphisms in $HC_2^{\tcb{A}}$.
\eprf

\begin{remark}\label{necinf}Note that in the definition of the category $HC^{\tcb{A}}_{\G}$ one cannot replace $\G_\infty$  by $\G$.
Indeed, consider the triangle $T$ in  $\G^2_\infty$ consisting in  those $(x, y)$ with $0 \leq x,y \leq \infty$ belonging to
one of the lines $y=0$, $x=y$, and $x=\infty$. There does not exist a homotopy equivalence
$g : T \to T'$ with $T'$ a definable subset of some $\G^n$ (or some $\G^w$ with finite definable $w$).
Indeed, assume such a $g$ exists and consider a homotopy inverse $f : T' \to T$. 
\tcb{Note that any  definable subset $X \not= T$ of $T$ which is definably connected retracts to a point. It follows that 
any 
homotopy equivalence $T \to T$ is surjective, so $f \circ g$ should be surjective.}
In particular, $f$ should be surjective. On the other hand, $T'$ should be definably path connected,
hence definably connected.
But a continuous surjective definable $f : T' \to T$ with $T'$ a definably connected subset of 
some $\G^n$ cannot exist, since $(y \circ f)^{-1} (\infty)$ would be a nontrivial clopen.
\end{remark}

 \section{Remarks on homotopies over imaginary base sets}

\tcb{Note that  \thmref{1} is valid over an arbitrary base set $A$, including imaginaries, when $X$ is a constructible
subset of $V$.  More generally, if $X$ and the $\xi_i$ are defined over $A \meet (\VF \union \G)$, the theorem follows,
simply by applying it over $A'= A \meet (\VF \union \G)$.}

Is \thmref{1} true in full generality over an arbitrary base?   Here is an indication that the answer may be positive, at least over a finite extension.  
Assume $(V,X)$ are given as in \thmref{1}, but over a base $A$ including  imaginary elements.
 A homotopy $h_c$ is definable over additional field parameters $c$, satisfying the conclusion
  of \thmref{1} over $A(c)$.  By the uniformity results of  \ref{ss10.7}, there exists an $A$-definable set $Q$ such that any parameter $c \in Q$
  will do.  One can find a definable type $q$ on $Q$, over a finite extension $A'$ of $A$ (i.e. $A'=A(a'), a' \in \acl(A)$).  We know that
  $q = \int_r f$, with $r$ an $A$-definable type on $\G^n$, and $f$ an $A$-definable $r$-germ of a function into $\std{Q}$.
  Define $h(t,v)  =  \lim_{u \in r} \int_{c \models f(u)} h_c(t,v) $.  Then $h(t,v)$ is an $A'$-definable homotopy. \tcb{The final image of $h$ is clearly $\G$-parameterized, and has property (5) of \thmref{1};   isotriviality, as well as the condition of being topologically $\G$-internal, should follow.}

\chapter{Applications to the topology of Berkovich spaces} \label{berkovich}

{\small \noindent \textbf{Summary.}
In this final chapter we   deduce from our main results general tameness statements about the topology of Berkovich spaces.
In \thmref{sc} we prove the existence of strong retractions to skeleta for analytifications of definable subsets
of quasi-projective varieties. \thmref{func-bi} is about functoriality and birationality statements for these retractions.
In \thmref{prosimplicial}, we show that, in the compact case, these analytifications are homeomorphic to
the projective limit of embedded finite simplicial complexes, under a compactness assumption.  
In \thmref{finitelymany} we prove finiteness of homotopy types in families in a strong sense.
We prove local contractibility in \thmref{lc} and
a result on homotopy equivalence of upper level sets of definable functions in \thmref{fpoi}.   
 \tcb{All these results are based on a certain surjection from the stable completion of a variety  over a maximal immediate extension of the algebraic closure of a field $F$,
 to the Berkovich space of that variety over $F$.    In the final section, we describe an injection in the opposite direction (over an algebraically
 closed field) which in general provides an identification between points of Berkovich analytifications and Galois orbits of stably dominated points.}
\par\bigskip}

\medskip

\section{Berkovich spaces}Set $\Rr_\infty = \Rr \union \{\infty\}$.  Let $F$ be a valued field with $\val (F) \subset \Rr_{\infty}$, and let $\bF= (F, \Rr)$ be viewed 
 as a substructure of a model of $\ACVF$ (in the $\VF$ and $\G$-sorts).   Here $\Rr=(\Rr,+)$ is viewed as an  ordered  abelian group. \nomenclature{$\bF$}{the structure $(F, \Rr)$}
 
 Let $V$ be an algebraic variety over $F$, and    let $X$ be an $\bF$-definable subset of
 the variety $V$; or more generally, of $V \times \G_\infty^n$.     We define the  Berkovich space $B_{\bF}(X)$ to be the space \nomenclature{$B_{\bF}(X)$}{Berkovich space of $X$}
 of types over $\bF$, in $X$, that are almost orthogonal to $\G$.   Thus for any $\bF$-definable function $f : X \to \G_\infty$
 and any $a \models p$, we have $f(a) \in \G_{\infty}(\bF) = \Rr_{\infty}$.  So $f(a)$ does not depend on $a$, and we denote it by $f(p)$.  
We endow $B_{\bF}(X)$ with a topology by 
defining a pre-basic open set to have
 the form $\{p \in X \meet U: \val (f) (p) \in W \}$, where $U$ is an affine open subset of $V$, $f$ is regular on $U$, 
 and $W$ is an open subset of $\Rr_\infty$.    A basic open set is a finite intersection of pre-basic ones. 
 This construction is functorial, thus, if 
$f : X \to X'$ is an $\bF$-definable morphism between
$\bF$-definable subsets of algebraic varieties over $F$,
one denotes by 
 $B_{\bF}(f) : B_{\bF}(X) \to B_{\bF}(X')$ the induced morphism.
 When we wish to consider $q \in B_{\bF}(X)$ as a type, rather than a point, we will write it as $q|\bF$. 
 
 When $V$ is an algebraic variety over $F$, $B_{\bF}(V)$ can be identified with the underlying topological space
of the Berkovich analytification $V^{an}$ of $V$. \nomenclature{$V^{an}$}{Berkovich analytification of $V$}
Recall that the underlying set of $V^{an}$ may be described
as the set of pairs $(x, u_x)$ with $x$ a point (in the schematic sense)
of $V$ and $u_x : F (x) \to \Rr_{\infty}$ a valuation extending $\val$ on the residual field $F (x)$, cf. \cite{duc-b}.
Such a pair $(x, u_x)$  determines a rational point $c_x \in V (F (x))$ whose type $p_x$ belongs to $B_{\bF}(V)$.
This correspondence is clearly bijective and a homeomorphism.
\tcr{It follows from Theorems 3.4.8 (i) and 3.5.1 (i) of \cite{berk} that $V^{an}$ is Hausdorff, since under our conventions an algebraic variety is always assumed to be separated.}
When $X$ is an $\bF$-definable subset of $V$, $B_{\bF}(X)$ is a semi-algebraic subset of $B_{\bF}(V)$ in the sense of \cite{duc}; 
conversely any semi-algebraic subset has this form.

 An element of $B_{\bF}(X)$ has the form $\tp(a/\bF)$, where $\bF(a)$ is an extension whose value group remains $\Rr$.  To see the relation to stably dominated types, note that if there
 exists an $\bF$-definable stably dominated type $p$ with $p|\bF = \tp(a/\bF)$, then $p$ is unique;
 in this case the Berkovich point can be directly identified with this element of $\std{X}$.  
 If there exists a stably dominated type $p$ defined over a finite Galois extension $F'$ of $F$,
 $\bF'=(F',\Rr)$, with $p|\bF = \tp(a/\bF)$, then the Galois orbit of $p$ is unique;
 in this case the relation between Berkovich points and points of $\std{X}$ is similar to the relation between closed points of $\spec(V)$ and points of $V(F^{\alg})$.  In general the Berkovich
 point of view relates to ours in rather the same way that Grothendieck's schematic
 points relate to Weil's points of the universal domain.   We proceed to make this more explicit.

  Let $K$ be a maximally
 complete algebraically closed field, containing $F$, with value group $\Rr$, and residue field equal to the algebraic closure of the residue field of $F$.   Such a   $K$ is unique up to isomorphism over $\bF$ by Kaplansky's theorem, and it will be convenient  to pick a copy of this field $K$ and denote it $F^{max}$. \nomenclature{$F^{max}$}{a certain maximally
 complete algebraically closed field containing $F$}
 
 We have a restriction
 map  from types  over $F^{max}$ to types over $\bF$.  On the other hand we have 
 an injective restriction map from stably dominated types defined over $F^{max}$, to types
 defined over $F^{max}$.  Composing these maps,
 we obtain a map   from the set of  stably dominated  types in $X$ defined over $F^{max}$ to the set of types over $\bF$ on $X$ whose image is contained
 in $B_{\bF}(X)$. Indeed, if
 $q$ lies in the image of this map,   then $q=\tp(c/\bF)$ for some $c$ with
$\tp(c/F^{max})$ orthogonal to $\G$, and it follows that $\G(\bF(c)) \nsubseteq \G(F^{max}(c)) = \G(F^{max}) = \G(\bF)$.
This defines a continuous map
\[ \pi_X: \std{X}(F^{max}) \to B_{\bF}(X).\] 
We shall sometimes omit the subscript when there is no ambiguity. 
 
\begin{lem} \label{br1}Let $X$ be an $\bF$-definable subset of an  algebraic variety over $F$. 
The mapping $ \pi: \std{X}(F^{max}) \to B_{\bF}(X)$ is surjective. In case $F= F^{max}$,
$\pi$ is a homeomorphism.\end{lem}

\prf  Suppose $q=\tp(c/\bF)$ is almost orthogonal to $\G$.  Let 
$L=F(c)^{max}$.  Then $\G(\bF) = \G(\bF(c)) = \G(L)$.  The field $F^{max}$ embeds 
into $L$ over $\bF$; taking it so embedded, let $p=\tp(c/F^{max})$.  Then
$p$ is almost orthogonal to $\G$, and $q=p|\bF$.  Since $F^{max}$ is maximally complete,
$p$ is orthogonal to $\G$, cf. \thmref{maxcomp}.  

In case $F= F^{max}$,  $\pi$ is also injective since $p|F$ determines $p$, for a stably dominated type based on $F$.
Thus $\pi$ is a continuous bijection;   since in this case the definitions of the topologies coincide on both sides, it is a homeomorphism.
\eprf



 %

Recall \ref{deftop}, and the remarks on definable topologies there.
 
\begin{prop} \label{B1.2}  Let $X$ be an $\bF$-definable subset of  an  algebraic variety $V$ over $F$.  Let  $\pi: \std{V}(F^{max}) \to B_{\bF}(V)$
be the natural map.     Then  $\pi \inv (B_{\bF}(X)) = \std{X}(F^{max})$, and $\pi: \std{X}(F^{max}) \to B_{\bF}(X)$ is a  closed map.  Moreover,
the following conditions are equivalent:   
\begin{enumerate}
\item   $\std{X}$ is definably compact;
\item     $X$ is bounded and v+g-closed;
\item  $\std{X}(F^{max})$ is compact;
\item   $B_{\bF}(X)$ is  compact;
\item $B_{\bF}(X)$ is closed in $B_{\bF}(V')$, where $V'$ is any complete $F$-variety
containing $V$.
\end{enumerate}
 
The natural map $B_{\bF'}(X) \to B_{\bF}(X)$ is also closed, if $F \leq F'$ and $\G(F') \leq \Rr$.  In particular, $B_{\bF}(X)$ is closed in $B_{\bF}(V)$
iff $B_{\bF'}(X)$ is closed in $B_{\bF'}(V)$. 
 \end{prop} 

\prf   The equality $\pi \inv (B_{\bF}(X)) = \std{X}(F^{max})$ is clear from the definitions.  
Let us consider the five conditions.  

The equivalence of (1) and (2) is  already known by \thmref{compactiffboundedclosed}.

 Assume (2).   We wish to prove (3) over $F^{max}$.  As $X$ is bounded, there exists a finite affine cover $V = \union V_i$,
closed immersions   $g_i: V_i \to \Aa^n$, and balls $B_i = \{x \in \Aa^n:  v(x_j) \geq b_i \}$, such that $X \nsubseteq \union_i g_i \inv (B_i)$.
 It suffices to prove (3) for $X \meet g_i \inv (B_i)$.  Thus we may assume  $X \nsubseteq B = \{x \in \Aa^n:  v(x_i) \geq b \}$.
 
 By \lemref{br1}, the natural map
 $\std{B} (F^{max}) \to B_{F^{max}}(B)$ is a homeomorphism.
 Let us first prove that this space is compact.
 Consider the polynomial ring $A={F^{max}}[X_1,\ldots,X_n]$.  Each element
 $p \in B_{F^{max}}(B)$ determines a map $v_p: A \to \Rr_\infty$.
 This provides an embedding $\Phi : B_{F^{max}}(B) \to \Fn (A, \Rr_\infty)$,
 with $\Fn (A, \Rr_\infty)$ the space of functions from $A$ to 
 $\Rr_\infty$. If one endows $\Fn (A, \Rr_\infty)$ with the Tychonoff topology,
 $\Phi$ induces a homeomorphism between 
 $B_{F^{max}}(B)$ and its image $\Phi (B_{F^{max}}(B))$. For $f$ in $A$,
 denote by $d_f$ the degree of $f$, by $a_f$ the smallest valuation of a coefficient of $f$,
 and set $b_f = b d_f + a_f$.
 Since $v_p(f) \geq b_f$ for any $p \in B_{F^{max}}(B)$,
 $\Phi (B_{F^{max}}(B))$ is contained in $\prod_{f \in A} [b_f, \infty]$, which is compact by Tychonoff's theorem.
 On the other hand, $\Phi (B_{F^{max}}(B))$ is clearly closed, being the set of  functions $u: A \to \Rr_{\infty}$  such that $u (fg) = u (f) + u (g)$,
 $u (f + g) \geq \min (u (f), u (g))$, 
 $u$ restricts to $\val$ on $F^{max}$, and $u(X_i) \geq b$ for every $i$.
It follows that $B_{F^{max}}(B)$ is compact.
 The definable set $X$, being  v+g-closed in $B$, is \tcb{a positive Boolean combination of}
 algebraic equalities $f_i = 0$ and weak inequalities $\val(g_i) \leq \val (h_i)$ \tcb{by \propref{vandg}}.
 Thus $\Phi (B_{F^{max}}(X))$ is the subset of
 $\Phi (B_{F^{max}}(B))$  similarly
 defined by the conditions
 $u (f_i)= \infty$ and  $u (g_i) \leq u(h_i)$, hence is closed.
It follows that  $\std{X}(F^{max}) = B_{F^{max}}(X)$ is compact.   This gives (3).

\tcb{If  (3) holds, then (4) also, since $\pi (\std{X}(F^{max})) = B_{\bF}(X)$.}
If $V'$ is any complete $F$-variety containing $V$,
 the inclusion $B_{\bF}(X) \to B_{\bF}(V')$ is continuous,
 and $B_{\bF}(V')$ is Hausdorff, so (4) implies (5). 
 
  On the other hand  if (1) fails, let $V'$ be  some complete
 variety containing $V$.   There exists an $F^{max}$-definable type  on $\std{X}$ with limit point $q$ in $\std{V'} \m \std{X}$.  So $\pi(q)$ is in $B_{\bF}(V')$ and in the closure of $B_{\bF}(X)$, but not in $B_{\bF}(X)$.  This finishes the proof of  the equivalence of (1-5).

Now the restriction of a closed map $\pi$ to a set of the form $\pi \inv(W)$
is always closed, as a map onto $W$.  So to prove the closedness property of $\pi$, we may take $X=V$, and moreover by embedding $V$
in a complete variety we may assume $V$ is complete.  In this case $X=V$ is v+g-closed and bounded, so $\std{X}(F^{max})$ is compact by condition (3).
As $B_{\bF}(X)$ is Hausdorff, $\pi$ is closed.  The proof that $B_{\bF'}(X) \to B_{\bF}(X)$  is also closed is identical, and taking $X=V$ we obtain
the statement on the base invariance of the closedness of $X$.  We could alternatively use the proof of \lemref{rattop}.
\eprf

 \begin{prop}\label{B1.3} Assume $X$ and $W$ are $\bF$-definable subsets of some  algebraic variety  over $F$.  
 \begin{enumerate} 
 \item Let $h_0:X \to \std{W}$ be an $\bF$-definable function.  Then $h_0$ 
 induces functorially a  function $\widetilde{h}: B_{\bF}(X) \to B_{\bF}(W)$ such that $\pi_W \circ h_0 = \widetilde{h} \circ \pi_X \circ i$, 
 with $i : X \to \std{X}$ the canonical inclusion. 
 \item
 Any continuous $\bF$-definable function $h: \std{X} \to \std{W}$ 
 induces a continuous  function $\widetilde{h}: B_{\bF}(X) \to B_{\bF}(W)$ such that $\pi_W \circ h = \widetilde{h} \circ \pi_X$. 
 \item The same applies if either $X$ or $W$ is a definable subset of $\G_\infty^n$ and we read 
 $B_{\bF}(X)=X(\bF)$, respectively $B_{\bF}(W) = W(\bF)$. 
 \end{enumerate} 
 \end{prop}

\prf     Define
$\widetilde{h}:  B_{\bF}(X) \to B_{\bF}(W)$ \tcr{similarly as in \ref{canon-ext}}.  Namely, let  $p \in  B_{\bF}(X)$.  
We view 
$p$ as a type over $\bF$, almost orthogonal to $\G$.  Say $p|\bF = \tp(c/\bF)$.   Let
 $d \models h_0(c) | \bF(c)$.  Since $h_0(c)$ is stably dominated, $\tp(d/\bF(c))$ is almost orthogonal
 to $\G$, hence so is  $\tp(cd/\bF)$, and thus also
$\tp(d/\bF)$.  Let 
$\widetilde{h}(c) =\tp(d/\bF) \in B_{\bF}(W)$.   Then   $\widetilde{h}(c)$ depends only on $\tp(c/\tcb{\bF})$, so we can let $\widetilde{h}(p)= \widetilde{h} (c)$.   

For the second part, 
 let $h_0 = h | X$ be the restriction of $h$ to the simple points. \tcb{It is v+g-continuous and}
by \lemref{hbasic}, $h$ is the unique continuous extension of $h_0$.   Define $\widetilde{h}$ as in (1).  
Let $\pi_X: \std{X}(F^{max}) \to B_{\bF}(X)$ and $\pi_W: \std{W}(F^{max}) \to B_{\bF}(W)$ be the restriction maps as above.  It is clear from  the definition that $\widetilde{h}(\pi_X(p)) = \pi_W( h(p))$.  
 (In case $F^{max}$ is nontrivially valued, this is also clear from the density of simple points, since $\widetilde{h} \circ \pi_X$ and $\pi_W\circ h$ agree on the simple points of $\std{X}(F^{max})$.) 

It remains to prove continuity.   By the discussion above, $\pi_X$ is a surjective and closed map. 
Let $Z$ be a closed subset of $B_{\bF}(X)$. By continuity of
$\pi_W \circ h$, $\pi_X \inv ( \widetilde{h} \inv(Z))
=h \inv ( \pi_W \inv(Z))$ is closed, hence
$ \pi_X (\pi_X \inv ( \widetilde{h} \inv(Z))) = \widetilde{h} \inv(Z)$ is closed.

(3) The proof goes through in both cases.  
 \eprf

 If $f : X \to Y$ is  an $\bF$-definable map and $b$ is a point in $Y$, we denote by
 $X_b$ the fiber $f^{-1} (b)$ over $b$. Similarly, if $q$ is a point of 
 $B_{\bF} (Y)$, $B_{\bF}(X)_q$ denotes the fiber
over $q$ of the induced mapping
$B_{\bF} (X) \to B_{\bF} (Y)$.

\begin{lem}  \label{B1.4}  Let $X$ be an $\bF$-definable subset of $V \times \G_\infty^n$ with  $V$ a variety over $F$. 
\begin{enumerate}
\item Let 
$f: X \to Y$ be an $\bF$-definable map, with $Y$  an $\bF$-definable subset of some variety over $F$. Let
$q \in B_{\bF}(Y)$, and assume $U$ is an $\bF$-definable subset of $X$, and  
$\std{U_b}$ is closed in $\std{X_b}$ for any $b \models q|\bF$.     Then
$B_{\bF}(U)_q$ is closed in $B_{\bF}(X)_q$. 
\item
Similarly if $g: X \to \mathbb{R}_\infty$ is an $\bF$-definable function, and 
$\std{g}|\std{X}_b$ is continuous  for any $b \models q|\bF$, then $B_{\bF}(g)$ induces a continuous map on $ B_{\bF}(X)_q \to \mathbb{R}_\infty$.
\item 
  More generally,
if  $g: X \to V'$ is an $\bF$-definable map into some variety $V'$, and $\tcb{g}|X_b$ is v+g-continuous for any $b \models q|\bF$, then 
$B_{\bF}(g)$ induces a continuous map 
$B_{\bF}(X)_q: B_{\bF}(X)_q \to B_{\bF}(Z)$. 
\end{enumerate}
 \end{lem}

\prf  
Indeed if $r \in B_{\bF}(X)_q \m B_{\bF}(U)_q$, let $c \models r|\bF$,
$b=f(c)$.    
We have $c \in X_b \m U_b$, so there exists a  definable function $\a_b: X_b \to \G_\infty$ and an open neighborhood
$E_c$ of $\a_b(c)$ such that $\a_b \inv(E_b) \nsubset X_b \m U_b$.  
By \lemref{rattop}, 
  $\a_b$ can be taken to be $\bF(b)$-definable, and in fact to be a continuous
  function of the valuations of some $F$-definable regular functions, and elements of $\G(\bF)$.    There exists an $\bF$-definable function $\a$
on $X$ with $\a_b = \a | X_b$.  
Now $\a$ separates $r$ from $B_{\bF}(U)_q$ on $B_{\bF}(X)_q$, showing that $U$ is closed in $B_{\bF}(X)_q$. 

The statement on continuity (2) follows immediately:  if $Z$ is a closed subset of $\G_\infty$, then $g \inv(Z) \meet \std{X_b}$ is closed in each $\std{X_b}$, hence $g \inv(Z) \meet B_{\bF}(U)_q$ is closed. 

The more general statement (3) follows since to show that a map into $B_{\bF}(Z)$ is continuous, it suffices to show that the composition 
with $B_{\bF}(s)$ is continuous for any definable, continuous $s: Z' \to \G_\infty$, with $Z'$ Zariski open in $Z$.
\eprf


The following lemma will be applied when $W$ is also over $Y$ and $h: X \to \std{W/Y}$; but a referee pointed out that
the more general statement is also valid, and simpler.
 %
 %

\begin{lem} \label{B1.5}  Let $X$, $Y$ and $W$ be $\bF$-definable subsets of some  algebraic variety  over $F$. 
Let $f: X \to Y$   be a  v+g-continuous,  $\bF$-definable map, 
and $h: X \to \std{W}$   an $\bF$-definable map inducing   $H: \std{X/Y} \to \std{W}$.   Assume $H | \std{X}_b$ is continuous for every $b \in Y$. 
Then for any $q \in B_{\bF}(Y)$,  
$h$ induces a continuous function $\widetilde{h}_q: B_{\bF}(X)_q \to B_{\bF}(W)$.  \end{lem}

\prf
 The topology on  $B_{\bF}(W)$  is the coarsest one such that
$B_{\bF}(g)$ is continuous for any v+g-continuous definable $g: W \to \G_\infty$.  Composing with $B_{\bF}(g)$, we see that we may 
assume $W=  \G_\infty$.   We have $h: X \to \G_\infty$, inducing $H: \std{X/Y} \to \G_\infty$, and we assume 
$H | \std{X}_b$ is continuous for $b \in Y$.  We have to show that a continuous $\widetilde{h}_q: B_{\bF}(X)_q \to \Gamma_\infty$ is induced.  
 
 In case the map $\std{X} \to \G_\infty$ induced from $h$ is continuous, by \propref{B1.3} $\widetilde{h}$ is continuous, and hence the restriction to each fiber $B_{\bF}(X)_q$ 
is continuous.

In general, let $X'$ be the graph of $h$.  
The projection $X' \to \Gamma_{\infty}$ being v+g-continuous, 
a natural, continuous function $B_{\bF}(X')_q \to \mathbb{R}_\infty$ is induced, by the above special case.  It remains to prove that the projection
map $B_{\bF}(X')_q \to B_{\bF}(X)_q$ is a homeomorphism with inverse induced by $x \mapsto (x, \tcb{h}(x))$.  When $q=b \in Y$ is a simple point,
this follows from the continuity of  $H | \std{X}_b$.  Hence by \lemref{B1.4}, it is true in general.
\eprf

  In the Berkovich category, as in  \ref{ph} and throughout the paper, by deformation retraction we mean a strong deformation retraction.
We continue to write   $\pi: \std{V}(F^{max}) \to B_{\bF}(V)$
for the natural map, defined above  \lemref{br1}.

 \begin{cor}\label{B1.6}    \leavevmode
 \begin{enumerate}
\item  Let $X$ be an $\bF$-definable subset of some  algebraic variety  over $F$. 
Let $h: I \times \std{X} \to \std{X}$ be an $\bF$-definable  deformation retraction, with 
 image $h(e_I, \std{X})=Z$.  Let  $\bI = I(\Rr_\infty) $ and $ \bZ = \pi (Z(F^{max}))$.
 Then $h$ induces a   deformation retraction  
 $\widetilde{h}: \bI \times B_{\bF}(X) \to B_{\bF}(X)$ with image $\bZ$.
 \item   Let $X \to Y$ be an $\bF$-definable morphism between $\bF$-definable subsets of some  algebraic variety  over $F$. 
  Let $h: I \times \std{X/Y} \to \std{X/Y}$ be an $\bF$-definable deformation retraction satisfying $(\ast)$, with fibers $h_y$ having image $Z_y$. 
Let $q \in B_{\bF}(Y)$.  Then $h$ induces a deformation retraction $\widetilde{h}_q: \bI \times B_{\bF}(X)_q \to B_{\bF}(X)_q$, with image $\bZ_q$.  
\item  Assume in addition there exists a definable $\Ups \nsubseteq \G_\infty^n$ and definable homeomorphisms $\a_y: Z_y \to \Ups$,
given uniformly in $y$.  Then $\bZ_q \cong \Ups$.    More generally if $\Ups \nsubseteq \G_\infty^w$ with $w$ a finite, Galois invariant subset of a finite field extension $F'$ of $F$, 
$\a_y: Z_y \to \Ups$, then $\bZ_q \cong \Ups/G$, where $G= \mathrm{Gal}(F'/F)$ is acting naturally on $w$.
\end{enumerate}
 \end{cor}
  
\prf  (1) follows from \propref{B1.3}; the statement on the image is easy to verify.   (2) follows similarly from \lemref{B1.5}.     
For (3), define $\beta:  X \to \Ups$ by $\beta(x) = \a_y( h(e_I,x))$  for $x \in X_y$, $e_I$ being the final point of $I$.    Then
$\alpha_y\inv \circ \beta (x)= h(e_I,x)$, $\beta(h(t,x))=\beta(x)$, $\beta(\alpha_y \inv(x))=x$.  Applying $B_{\bF}$ and restricting to the fiber over $q$
we obtain continuous maps $\beta,\alpha_y \inv$ by \lemref{B1.4}; the identities survive, and give the result.     
 \eprf  
  
 %

\section{Retractions to skeleta}
Let $V$ be an algebraic variety over 
a valued field $F$ with
$\val (F) \subset \Rr_{\infty}$ and let
 $S$ be an $\bF$-iso-definable $\G$-internal subset of $\std{V}$.
According to \thmref{G-embed-0}, there exists an $\bF$-definable 
embedding $S \to \G_\infty^w$, where $w$ is a finite set.  Let $F'$ be a finite Galois extension of $F$,
such that $\Aut(F^{\alg}/F')$ fixes each point of $w$. 
We shall say $S$ {\em splits} over $F'$. \index{splitting of a $\G$-internal set}
Then there exists an $\bF'$-definable embedding
$S \to \G_\infty^n$, $n = |w|$.  It follows that $S(\bF'')=S(\bF')$ whenever   $F'' \geq F'$ is a valued field extension with  
    $\G(F'') \nsubseteq \Rr$.  The image $S_{\bF}$ of $S$ in $B_{\bF}(X)$ is thus homeomorphic to $S(\bF')/\mathrm{Gal}(F'/F)$.
    The image $S_{\bF''}$ of $S$ in $B_{\bF''}(X)$ is homeomorphic to $S(\tcb{{\bF'}})$.  
  Note that the canonical map $\std{V}(F^{max}) \to B_{\bF'}(V)$ 
restricts to an injective map on $S$, since $S(F^{max}) \nsubseteq S(\bF')$.   

\medskip

For our purposes, a $\Qq$-{\em tropical structure} \index{$\Qq$-tropical structure} on a topological space $X$ is a homeomorphism of $X$
with a subspace $\bS$ of $[0,\infty]^n$ defined as a finite Boolean combination of equalities or   inequalities   between terms 
$\sum \a_i x_i + c$ with $\a_i \in \Qq,  \alpha_i \geq 0, c \in \Rr$.  Since $\bS$ is definable in $(\Rr,+,\cdot)$, 
 $X$ is  homeomorphic to a   finite simplicial complex. \tcb{Recall that a valued field extension $L$ of a valued field $F$ is called an {\em  Abhyankar extension} \index{Abhyankar extension} if the transcendence degree of $L / F$ is
 equal to the sum of the transcendence degree of the  residue field  extension  and the $\Qq$-rank
of $\G(L)/\G(F)$.}
  
From \thmref{1} and \corref{B1.6} we obtain:

\begin{thm} \label{sc}  Let $X$ be an $\bF$-definable subset of a quasi-projective algebraic variety $V$ over 
a valued field $F$ with
$\val (F) \subset \Rr_{\infty}$.  There exists 
a  \textup{(}strong\textup{)} deformation retraction 
 $H: \bI \times B_{\bF}(X) \to B_{\bF}(X)$, whose image $\bZ$ is of the form $S_{\bF}$ with $S$ an $\bF$-iso-definable $\G$-internal subset of $\std{V}$.
 Thus,  $\bZ$ has a $\Qq$-tropical structure, in particular it is  homeomorphic to   a finite simplicial complex.    
 Furthermore each 
 \tcb{point $q$ of $\bZ$, as a type over $\bF$, extends to a unique
 stably dominated type $p$ and this type is strongly stably dominated. Restricted to $F$, $q$ determines 
 an Abhyankar extension of
 the valued field $F$.}
 \end{thm}

 \prf
   Let $S$ be the final image provided by  \thmref{1} assuming (5) holds. Thus $S$ consists of strongly stably dominated types and we have
 an $\bF$-definable homeomorphism $h: W \to S$, where $W$ is a subset of $\G_\infty^w$, with $w$ a finite $\bF$-definable set.  So for $a \in W(\bF) = W(\Rr)$, 
 $p=h(a)$ is strongly stably dominated over $\bF$, and extends the restriction to $\bF$, which is the image in $\bZ$ of $h(a)$.
 For the last point,  $p$ is defined over $F \union A$ where $A$ is a finitely generated $\Qq$-subspace of $\Rr$.
 Let $F'$ be an Abhyankar extension of $F$, with value group equal to $\val(F) +A$.  Then $F'(p)$ is Abhyankar over $F'$, 
 and hence over $F$.   
\eprf


\begin{example}\label{ellipticber}Let us revisit the  elliptic curve example of \exref{elliptic} in the Berkovich setting.
Assume for instance $F = \Qq_3$ and set $\lambda = 3$.
So $C_3$ is the projective model of the curve $y^2 = x (x-1) (x - 3)$.
We have seen in 
\exref{elliptic} that its skeleton $K'$ in $\std{C_3}$  is  a combinatorial circle. 
This circle admits a $ \Qq_3$-definable embedding in $\G^{\{i, -i\}}$, it
splits over $\Qq_3 (i)$ and conjugation acts on it by exchanging the points in the fibers
of $K' \to K$. Thus, \tcb{for $F= \Qq_3 (i)$},
$\tcb{B_{\bF} (C_3)}$ has the homotopy type of a circle, while \tcb{for $F'= \Qq_3$},
$\tcb{B_{\bF'} (C_3)}$ retracts to a segment, thus is contractible.
\end{example}

We now state some functorial properties of the deformation retraction above.   Like \thmref{sc}, 
these were  proved by Berkovich  assuming the base field $F$ is nontrivially valued, and that $U$ and $V$ can be embedded in proper varieties which admit a pluri-stable model over the ring of integers of $F$.  We thank Vladimir Berkovich for suggesting these statements to us.

\medskip

Whenever we write $B_{\bF}(V)$, we assume the valuation on $F$ is real valued,  allowing the case that the valuation is trivial.   
If $\bF'$ is an extension of $\bF$, we write $B_{\bF'}(U)$ for $B_{\bF'}(U\otimes F')$.
\begin{thm}\label{func-bi}  Let $U$ and $V$ be quasi-projective algebraic varieties over a valued field $F$ with value group contained in $\Rr$.  
Let $X$ and $Y$ be $\bF$-definable subsets of $U$ and $V$, respectively.

\begin{enumerate}
\item  There exists a finite separable extension $F'$ of $F$ such that, for any  non-archimedean field $F''$ over $F'$, the canonical map $B_{\bF''}(X)  \to B_{\bF'}(X)$ is a homotopy 
equivalence.  In fact, 
there exists a deformation retraction of  $B_{\bF'}(X)$ to $\bZ'$  as in \textup{\thmref{sc}}   that  lifts to a deformation retraction of  
$B_{\bF''}(X)$ to $\bZ''$,  for which the canonical map $\bZ'' \to \bZ'$ is a homeomorphism. 
 \item   There exists a finite separable extension $F'$ of $F$ such that, for any  non-archimedean field extension $F''$ of $F'$, 
 the canonical map $B_{\bF''}(X\times Y)  \to B_{\bF''}(X)  \times B_{\bF''}(Y) $ is a homotopy equivalence. 
\item  Let $U$ be smooth and $U'$ be a dense open subset of $U$. Then the canonical embedding $B_{\bF}(U') \to B_{\bF}(U)$  
is a homotopy equivalence.  
\end{enumerate}
\end{thm}

\prf   Let us prove (1).   The homotopy of \thmref{1} is $F$-definable, and so functorial on $F''$-points for any $F'' \geq F$.  Denote by $S$ its final image.
Choose a finite Galois extension $F'$ of $F$ that splits $S$.
For any $F \leq F' \leq F''$, the homotopy of  $B_{\bF''}(X)$ is compatible with the homotopy of  
$B_{\bF'}(X)$ via the natural map $B_{\bF''}(X) \to B_{\bF'}(X)$ (restriction of types).  The final   
image of the homotopies is respectively $S_{F''}$ and $S_{F'}$; we noted that these are homeomorphic images of $S$ as soon as $F'$ splits $S$
and hence homeomorphic via the natural map.

(2)   follows similarly from   \corref{product-homotopy-eq} (which was devised  precisely with the present motivation) and its proof. 
Indeed, as in the proof of 
 \corref{product-homotopy-eq}, let us consider definable deformations retractions for $X$ and $Y$ with final images $S$ and $T$.
Recall the homotopy equivalence $  \std{X \times Y} \to \std{X} \times \std{Y}$
 in \corref{product-homotopy-eq} was part of a commutative diagram
  \begin{equation*}\xymatrix{
     \std{X \times Y}
	\ar[d]_{\pi_X \times \pi_Y} \ar[r]^{} & 
	S \tensor T  \ar[d]^{{\pi_S \times \pi_T}}\\
	\std{X} \times \std{Y} \ar[r]_{} &S \times T,
	}
	\end{equation*}
	whose horizontal morphisms are definable retractions and that $\pi_S \times \pi_T$ was proven to be a homeomorphism.
	 Choose a finite Galois extension $F'$ which splits both $S$ and $T$ (in fact it would be enough to require $F'$ to  split one of $S$ and $T$).
	 It is then clear that for any $F''\geq F'$,
	the homotopy equivalence $  \std{X \times Y} \to \std{X} \times \std{Y}$
	induces
	a  homotopy equivalence
	$B_{\bF''}(X \times Y)  \to B_{\bF''}(X) \times B_{\bF''}(Y)$.

(3)  follows directly from \remref{birational}.
\eprf

The following result  was previously known when $X$ is a smooth  projective curve \cite{berk}.  

\begin{thm}  \label{prosimplicial}  Let $X$ be  an $\bF$-definable subset of a quasi-projective  algebraic variety $V$ over 
a valued field $F$ with
$\val (F) \subset \Rr_{\infty}$ and assume $B_{\bF}(X)$ is  compact.
Then there exists a family  $(X_i)_{i\in I}$ of
finite simplicial complexes of dimension $\leq \dim V$,
embedded in $B_{\bF}(X)$, 
where $I$ is a directed partially ordered set, such that $X_i$ is a subcomplex of $X_j$ for $i < j$,
with   deformation retractions $\pi_{i, j}: X_j \to X_i$ for $i<j$, and  deformation retractions $\pi_i :  B_{\bF}(X) \to X_i$
 for $i \in I$, satisfying $\pi_{i, j} \circ \pi_j = \pi_i$ for $i < j$, such that the canonical map from
 $B_{\bF}(X)$ to the  projective limit of the spaces $X_i$ is a homeomorphism.
\end{thm}

\prf    Let the index
set $I$ consist of all $\bF$-definable continuous maps $j: \std{X} \to \std{X}$, such that 
there exists an $\bF$-definable deformation retraction $H_j : I \times \std{V} \to \std{V}$ as in \thmref{1},
restricting to a  deformation retraction $H_j^X : I \times \std{X} \to \std{X}$ such that
$j(x)=H_j^X(e_I,x)$.
Here we insist that $H_j$ satisfies  \tcb{clause} (7) of \thmref{1} \tcb{for the irreducible components of $V$}. Let us denote by $T_j$ the final image of
$H_j$ and by $S_j$ that of $H_j^X$. Thus $S_j = j(\std{X}) = \std{X} \meet T_j$. Let $X_j$ denote the image of $S_j (F^{max})$ in
$B_{\bF} (X)$. Thus $X_j$ is homeomorphic to $S_j (\acl(\bF))/ \mathrm{Gal}(F^{\alg}/F)$.   
Say that $j_1 \leq j_2$ if $S_{j_1} \nsubseteq S_{j_2}$. 
In this case, $j_1 | S_{j_2} :  S_{j_2} \to S_{j_1}$ is a deformation retraction through the homotopy $j_2 \circ H_{j_1}(t, \cdot)$.

Let $\pi_{j_1,j_2} $ be the induced map $X_{j_2} \to X_{j_1}$.  It is a deformation retraction.
Let us prove the system is directed, i.e. given $j_1$ and $j_2$ there exists $j_3$ with $j_1,j_2 \leq j_3$.   To see this, for $j = j_1, j_2$,
let $\a_j: T_j \to \G_\infty^{\tcb{w_j}}$ be a definable injective map, \tcb{with $w_j$ a finite $F$-definable set,} and let $j_3$ belong to a homotopy $H_{j_3}$ respecting
the functions $x \mapsto \a_{j_1} (H_{j_1} (e_I, x))$,
$x \mapsto \a_{j_2} (H_{j_2} (e_I, x))$ and preserving the irreducible components of $V$. Then by \propref{maxdimabh} (2), since $H_{j_3}$ satisfies \tcb{clause} (7)
of \thmref{1} \tcb{for the irreducible components of $V$},
 $H_{j_3}$ fixes $T_{j_1}$ and $T_{j_2}$ pointwise, thus
  $H^X_{j_3}$ fixes $S_{j_1}$ and $S_{j_2}$ pointwise and 
the image of $j_3$ includes them both.

We have a natural surjective map $\pi_j : B_{\bF}(X)\to X_j$ for each $j$, induced by the mapping $j$; 
it satisfies $\pi_{i, j} \circ \pi_j = \pi_i$ for $i < j$ and it is a deformation retraction.

This yields a continuous   map from $\theta : B_{\bF}(X) \to \limproj_j X_j$.  The image is dense since each $\pi_j$ is surjective; as
 $B_{\bF}(X)$ is compact the image is closed, so $\theta$ is surjective.
We now show that $\theta$ is injective.    
Let $p \neq q \in B_{\bF}(X)$; view them as types almost orthogonal to $\G$.    For any open affine $U$
and   regular $f$ on $U$, for some $\a$, either $x \notin U$ is in $p$ or $\val (f) = \alpha$ is in $p$;
this is because $p$ is almost orthogonal to $\G$.  Thus as $p \neq q$, for some open affine $U$
and some regular $f$ on $U$, either $p \in U$ and $q \notin U$, or vice versa, or $p,q \in U$ and for some regular 
$f$ on $U$, 
$f(x)=\alpha \in p$, $f(x)=\beta \in q$, with $\alpha \neq \beta$.  Let $H$ be as in \thmref{1} respecting $U$ and $\val (f)$, and let $j$ be a corresponding \tcb{retraction}.  Then clearly $\pi_j (p) \neq\pi_j(q)$.  
Thus, $\theta$ is a continuous bijection and
by compactness it is a homeomorphism.
\eprf 

\begin{remark}  
Let $\Sigma$ be (image of) the direct limit of the $X_i$ in $B_{\bF}(X)$.
   Note that $\Sigma$ contains all rigid  points of $B_{\bF}(X)$ (that is, images of simple points under the mapping $\pi$ in \lemref{br1}):
   this follows from \thmref{1}, by finding a homotopy to a skeleton $S_x$ fixing a given simple point $x$ of $\std{X}$.
  We are not certain whether $\Sigma$ can be taken to be the whole of $B_{\bF}(X)$.
    But given a stably dominated type $p$ on $X$, letting $S_p=S_x$
  for $x \models p$ and  averaging the 
  homotopies with image $S_x$ over $x \models p$, we obtain a definable homotopy whose final image is a 
  continuous, definable image of $S_p$.     In this way we can express $B_{\bF}(X)$ as a direct limit of a system of finite simplicial complexes, with continuous transition maps.
\end{remark}


 \section{Finitely many homotopy types}\label{finiteht}

 We will now prove that a uniform family of Berkovich spaces runs through only finitely many homotopy types.
 
 In the definable setting, for \tcb{stable completions}, the situation is different.  Consider a family of triangles in $\G^2$; they may be the skeleta of elliptic curves, and so homotopy equivalent 
 to them.  Two triangles are  definably  homotopy equivalent iff they are definable homeomorphic.  But there are many definable homeomorphism types of triangles, or even of segments;
 indeed $[0,\a]$ and $[0,\beta]$ are definably homeomorphic iff $\beta$ is a rational multiple of $\a$.
 
 On the other hand, if we expand $\G$ to be a model of the theory $\RCF$ of real closed fields, then it is known that only finitely many homeomorphism types appear in a given definable family.
 Using the uniform version of \thmref{1}, this extends to uniformly definable families of \tcb{stable completions}.
 
\tcb{For applications to  Berkovich spaces in terms of the usual topological homotopy type, or even homeomorphism type of skeleta,
the expansion to $\RCF$ is harmless. In the setting of stable completions, we explain in \remref{omin-comb} how it can be avoided.}

\medskip

Part (1) of the following theorem is a special case of part (2); we single it out as we will prove it first.  \tcb{We consider a uniformly definable family  of 
definable subset of $\Pp^m$.}

 \begin{thm}\label{finitelymany}
 \tcb{Let $V$ be a variety defined over a valued field $F$.
 Let $Y$ be an $\bF$-definable subset of $V \times \G^r$, for some $r$, and let $X$ be an 
 $\bF$-definable subset of $Y \times \Pp^m$ for some $m$.
 Denote by $f : X \to Y$ the projection on the first factor.}
 \begin{enumerate}\item
   For $b \in Y$, 
let $X_{b} = f \inv({b})$.  Then there are finitely many possibilities for the homotopy type of $B_{\bF (b)} (X_{b})$, as ${b}$ runs through $Y$.   
More generally, let $\tcb{U_1 \nsubset  \dots \nsubset U_{\ell} =  X}$ be a chain of $\bF$-definable sets.  Then as ${b}$ runs through $Y$ there are finitely many possibilities for the homotopy type of the \tcb{tuple $(B_{\bF (b)}(X_b \meet U_i))$}.  
\item For any valued field extension $F \leq F'$ with $\G(F') \leq \Rr$ and $q \in B_{\bF'}(Y)$,
let  $B_{\bF'}(X)_q$ denote the fiber over $q$ of the 
canonical map $B_{\bF'}(X) \to B_{\bF'}(Y)$. Then there are only finitely many possibilities for the homotopy type of
$B_{\bF'}(X)_q$  as $q$ runs over  $B_{\bF'}(Y)$ and $F'$ over extensions of $F$.  
More generally, let $\tcb{U_1 \nsubset  \dots \nsubset U_{\ell} = X}$ be a chain of $\bF$-definable sets.  Then as $q$ runs over  $B_{\bF'}(Y)$ 	and $F'$ over extensions of $F$ there are finitely many possibilities for the homotopy type of the \tcb{tuple $(B_{\bF' (b)}(X_b \meet U_i))$}.  
\end{enumerate}
 \end{thm}

\prf \tcb{Let us start by proving the first statement in (1) under the assumption that for any $b \in Y$, 
 $X_{b}$ is Zariski closed in $\Pp^m$.}

According to the uniform version of \thmref{1}, \propref{1u}, there exists an $\bF$-definable map $W \to Y$ with finite
fibers $W({b})$ over ${b} \in Y$,  and uniformly in ${b} \in Y$ an $\bF (b)$-definable  homotopy retraction $h_b$ on $X_{b}$ preserving the given data, with final image $Z_{b}$, and an $\bF (b)$-definable 
homeomorphism $\phi_{b}: Z_{b} \to S_{b} \nsubseteq \G_\infty^{W({b})}$.

\begin{claim} We may find, definably uniformly in $b$,
an $\bF (b)$-definable subset $T_b \nsubseteq \G_\infty^n$, a \tcb{finite} $\bF (b)$-definable set $W_!(b)$,   and for $w \in W_!(b)$, a definable homeomorphism $\psi_w: Z_b \to T_b$,
such that $H_b = \{\psi_{w'} \inv \circ \psi_w: w,w' \in W_!(b) \}$ is a group of homeomorphisms of $Z_b$,
and $H'_b =  \{\psi_{w'}   \circ \psi_w \inv: w,w' \in W_!(b) \}$ is a group of homeomorphisms of $T_b$. 
\end{claim}
\begin{proof}[Proof of the claim]  In fact for a fixed
$b$, one can  pick some $W(b)$-definable homeomorphism $\psi_b$ of $Z_b$ onto a definable subspace of $\G_\infty^n$;
let $\Xi_b =\{\psi_w: w \in W_!(b) \}$ be the set of automorphic conjugates of $\psi_b$ over $\bF(b)$; and verify that 
$H_b$ is a finite group, $\Xi_b$ is a principal torsor for $H_b$, and so $H'_b$ is also a finite group (isomorphic to $H_b$).
Thus, for a fixed $b$, one can do the construction as stated, obtaining the stated properties. 
Now the fact that the $\psi_w$ are conjugates of $\psi_b$ is not an ind-definable property of $b$.  
But the consequences mentioned in the claim\tcr{\textemdash}that $\psi_w$ is a definable homeomorphism, and the compositional properties\tcr{\textemdash}are clearly ind-definable, and in fact definable,
 properties of $b$.  Hence   by the compactness \tcb{and} gluing argument we may find $W_!(b)$ and $ \Xi_b$ uniformly in $b$, with
the required  properties.    In particular, there exists an $\bF$-definable map $W_! \to Y$ with 
fibers $W_!({b})$ over ${b} \in Y$. \end{proof}


By stable embeddedness of $\G$, and elimination of imaginaries in $\Gamma$, we may write $T_b = T_{\rho(b)}$ where $\rho: Y \to \G^m$ is a definable function.
Let $\G^*$ be an expansion of $\G$ to $\RCF$.  Then by  \remref{omin-tri1}, $T_b$ runs through finitely many $\G^*$-definable
homeomorphism types as $b$ runs through $Y$.   Similarly, the pair $(T_b,H'_b)$ runs through finitely many $\G^*$-definable
equivariant homeomorphism types (e.g. we may find an $H'_b$-invariant cell decomposition of $T_b$ and 
describe the action combinatorially). In particular, for $b \in Y$, $(T_{b}(\Rr),H'_b)$ runs through finitely
many homeomorphism types (i.e. isomorphism types of pairs $(U,H)$ with $U$ a topological space, $H$ a finite 
group acting on $U$ by auto-homeomorphisms).

By  \corref{B1.6} we have, for ${b} \in Y$, a deformation retraction
of $B_{\bF (b)}(X_{b})$ to $B_{\bF (b)}(Z_b)$.  Pick $w \in W_!(b)$, and let $W^*(b)$ be the set of realizations of $\tp(w/\bF(b))$.
If $w,w' \in W^*(b)$ then $w' = \si(w)$ for some automorphism $\si$ fixing $\bF(b)$; we may take it to fix $\G$ too.
It follows that $\psi_w \inv  \circ \psi_{w'} = \si | Z_b$.  Conversely, if $\si$ is any automorphism of $W_!(b)$, it may be
extended by the identity on $\G$, and it follows that $\psi_{\si(w)} = \psi_w \circ \si$; so $W^*(b)$ is a torsor
of $H^*(b)=\{ \psi_w \inv \circ \psi_{w'}: w, w' \in W^*(b) \}$, which is a group.  Let $H_*(b) = \{\psi_w \circ \psi_{w'} \inv: w, w' \in W^*(b) \}$.
It follows that $H_*(b)$ is a group, and for any $w \in W^*(b)$, $\psi_w$ induces a bijection 
$Z_b / H^*(b) \to T_b/H_*(b)$; moreover it is the same bijection, i.e. it does not depend on the choice of $w \in W^*(b)$.

We are interested in the case $\G(\bF(b)) = \G(\bF)=\Rr$.  In this case,   
since $H^*(b)$ acts by automorphisms over $\bF(b)$, two $H^*(b)$-conjugate elements of $Z_b$ have the same 
image in $B_{\bF (b)}(X_b)$. On the other hand two non-conjugate elements have distinct images in $T_b/H_*(b)$,
and so cannot have the same image in $B_{\bF (b)}(X_b)$.  It follows that $B_{\bF (b)}(Z_b)$, $Z_b(\bF (b))/ H^*(b)$ and $T_b (\Rr)/H_*(b) $
are canonically isomorphic.
By compactness and definable compactness considerations one deduces that these isomorphisms between  $B_{\bF (b)}(Z_b)$, $ Z_b(\bF (b)) / H^*(b)  $ and $T_b(\Rr) /H_*(b) $ are in fact homeomorphisms.  
It is only for this reason  that we \tcb{required $X_b$ to be Zariski closed in the beginning of the proof}.

 The number of possibilities for $H_*(b)$ is finite and bounded, since $H'_b$ is a group of finite size, bounded independently of $b$, and $H_*(b)$ is a subgroup of $H'_b$.      Since the number of equivariant homeomorphism types of $(T_b(\Rr),H'_b)$ is bounded, 
we are done with the first statement in (1).

With the help of \corref{B1.6}, this proof goes through for non-simple Berkovich points too.
Let $q \in B_{\bF}(Y)$, and view it as a type over $\bF$.   By \corref{B1.6} (2), $B_{\bF}(X)_q$ has the homotopy type
of $\bZ_q$.  
Let $b \models q$, pick  $w \in W_!(b)$ and let notation be as above.  Let $b'=(b,w)$ and let $q' = \tp(b,w / \bF)$.  
Let $X' = X \times_Y W_!$.   By \corref{B1.6} (2) applied to the pullback of the retraction $I \times \std{X/Y} \to \std{X/Y}$ to $\std{X'/ W_!}$,
   $B_{\bF}(X')_{q'}$ retracts to a space $\bZ_{q'}$
  which is  homeomorphic to $T_b (\Rr)$.  By the same reasoning as above, it follows that $\bZ_q$ is homeomorphic to 
    $\bZ_{q'}$ {\em modulo} a certain subgroup $H^*(b)$ of $H_b$, and also homeomorphic to $T_b$ modulo $H_*(b)$
    for a certain subgroup of $H'_b$, so again the number of possibilities is bounded. This holds uniformly when $F$ is replaced by any valued field extension, and the first statement in (2) follows.
    
\tcb{The proof goes through directly to provide the generalization to chains. 
In particular  we can now remove the hypothesis that $X_{b}$ is Zariski closed in $\Pp^m$, after replacing
$U_1 \nsubset  \dots \nsubset U_{\ell} =  X$ by $U_1 \nsubset  \dots \nsubset U_m  \nsubset U_{\ell +1} =  Y \times \Pp^m$.}
\eprf

 \begin{rem} \label{omin-tri1} In the expansion of $\G$ to a real closed field, definable subsets of $\G_\infty^n$ are locally contractible and definably compact subsets of $\G_\infty^n$ admit a  definable triangulation, compatible with
 any given definable partition into finitely many  subsets.    By taking the closure in case the sets are not compact, it follows that given a definable family of semi-algebraic subsets of $\Rr_\infty^n$,
 there exist a finite number of rational polytopes (with some faces missing), such that each member of the family is homeomorphic  to at least one such  polytope.  In particular the number of definable homotopy types is finite.
 In fact it is known that the number of definable homeomorphism types is finite.  See 
\cite{coste}, \cite {vddries-tame}.
\end{rem}

\begin{rem}\label{omin-tri2} Eleftheriou has shown \cite{elef} that there exist abelian groups interpretable in
$\Th(\Qq,+,<)$ that cannot be definably and homeomorphically embedded in affine space within $\DOAG$.   By  \thmref{G-embed-1c},
the skeleta of abelian varieties can be so embedded.  It would be good to bring out the additional structure
they have that ensures this embedding.  
\end{rem}

\begin{rem}\label{omin-comb}\tcb{Let us explain how to avoid the use of the expansion to $\RCF$ in the 
setting of stable completions.
It is shown
in the thesis of Eleftheriou, see also \cite{elef-st} p.~1115, that a definable subset $X$ of $\G^n$ may be partitioned into finitely 
 many linear cells.  This decomposition is defined by some formulas $\phi_j(x,\a_X)$ requiring some parameters $\a_X$ (the coefficients in the linear equations are from $\Qq$;
 the parameters refer to the inhomogeneous part of the equations).   One can easily determine a $0$-definable set $P$
 such that $\a_X \in P$, and such that  for $\a \in P$, the formulas  $\phi_j(x,\a)$ determine a cell complex with the same adjacencies.
 Choosing $\b \in P(\Rr,+,<)$, and letting $X_\beta$ be defined by  the formulas $\phi_j(x,\beta)$, we obtain a topological  space
 whose homeomorphism type clearly does not depend on the choice of $\beta$.  By refining one can see that it also does not depend on the choice of the formulas
 $\phi_j$ (though strictly speaking, that is not needed for our finiteness statements).   We call this homeomorphism type the {\em combinatorial
 homeomorphism type} of $X$.  \index{combinatorial homeomorphism type}
 For instance, all triangles have the same combinatorial homeomorphism type, though as
 explained above they have distinct definable homotopy type.  Now any definable family of definable subsets of $\G^n$ runs through 
 a finite number of combinatorial homeomorphism types. It follows that for any definable family of quasi-projective varieties, there exists a finite set $\Omega$ of combinatorial homeomorphism types such that
   the stable completion of any variety in the family admits a skeleton with combinatorial homeomorphism type in $\Omega$. Similar considerations and  finiteness statements apply to $\G_\infty^n$, and filtered definable spaces $(X,X_1,\ldots,X_n)$ where
 for $X \nsubset \G_\infty^n$, we let $X_k$ be the subset of points exactly $k$ of whose coordinates are $\infty$.
}
\end{rem}

\section{More tame topological properties}

\begin{thm}[Local contractibility]\label{lc}Let $X$ be an $\bF$-definable subset of an algebraic variety $V$ over 
a valued field $F$ with
$\val (F) \subset \Rr_{\infty}$. The space $B_{\bF}(X)$  is locally contractible.  
\end{thm}

\prf We may assume $V$ is affine.
Since the topology of $B_{\bF}(X)$ is generated by open  subsets of the form
$B_{\bF} (X')$ with $X'$ definable in $X$, it is enough to prove
that every point $x$ of $B_{\bF}(X)$ admits a contractible neighborhood.
By \thmref{1} and \corref{B1.6}, there exists a strong deformation retraction
$H : I \times B_{\bF} (X) \to B_{\bF} (X)$
with image a subset $\Upsilon$ which is homeomorphic to a semi-algebraic subset of some
$\Rr^n$. Denote by $\varrho$ the retraction $B_{\bF} (X) \to \Upsilon$.
By (4) in \thmref{1}  one may assume
that $\varrho (H (t, x)) = \varrho (x)$ for every $t$ and $x$.
Recall that any semi-algebraic subset $Z$ of 
$\Rr^n$
is locally contractible:  one may assume $Z$ is bounded, then its closure
$\overline{Z}$ is compact and semi-algebraic and the statement follows from
the existence of triangulations of $\overline{Z}$ compatible with the inclusion
$Z \hookrightarrow \overline{Z}$ and having any given point of $Z$ as vertex.
I\tcb{t} is thus possible to pick a contractible
neighborhood $U$ of $\varrho (x)$ in $\Upsilon$.
Since the set $\varrho \inv (U)$ is invariant by the homotopy $H$,
it retracts to $U$, hence is contractible.
\eprf

\begin{rem}\label{lcab}\tcb{As noted in \remref{lastrem}, if $x$ is an Abhyankar point in the sense of  \defref{abdef},
it follows from \thmref{sgc} and \propref{ncom13}, together with the proof of
\thmref{lc},
that $x$ admits a basis of neighborhoods that  strongly retracts to $x$.}
\end{rem}

\begin{rem}Berkovich proved in \cite{berk1} and  \cite{berk2} local contractibility of smooth 
non-archimedean analytic spaces, and raised the question of the singular case.   His proof uses de Jong's results on alterations.  
\end{rem}

Let us give another application of our results, in the spirit of results of  Abbes and Saito \cite{AS} 5.1 and Poineau \cite{poi} Th\'eor\`eme 2.

\begin{thm}\label{fpoi}Let $X$ be an $\bF$-definable subset of a quasi-projective algebraic variety over 
a valued field $F$ with 
$\val (F) \subset \Rr_{\infty}$ and let 
 $G : X \to \G_{\infty}$ be an  $\bF$-definable map. 
 Consider the corresponding
 map
 $\bG: B_{\bF}(X) \to \Rr_\infty$. Then there is a finite partition of
 $\Rr_\infty$ into intervals such that
 the fibers of
 $\bG$ over each interval have the same homotopy type.
 Also, if one sets 
\tcb{$ B_{\bF}(X)_{\leq \varepsilon}$ 
to be the preimage of 
$(-\infty, \varepsilon]$},
there exists 
a finite partition of
 $\Rr_\infty$ into intervals such that
for each interval $I$
the inclusion
\tcb{$B_{\bF}(X)_{\leq \varepsilon} \to B_{\bF}(X)_{\leq \varepsilon'}$},
for $\varepsilon \tcb{<} \varepsilon'$ both in $I$,
is a homotopy equivalence.
 \end{thm}

\prf
Consider a strong deformation retraction
of $\std{X}$ leaving the fibers of $G$ invariant, as provided by \thmref{1}.
By \corref{B1.6} it induces a retraction $\varrho$ 
of $B_{\bF} (X)$ onto a subset $\Upsilon$ 
such that there exists a homeomorphism $h : \Upsilon \to S$ with 
$S$  a semi-algebraic subset  of some
$\Rr^n$. 
By construction $\bG$ factors as 
$\bG = g \circ \varrho$ with $g$ a function $S \to  \Rr_\infty$. Furthermore, we may assume that
$g' :=h \inv \circ g$ is a semi-algebraic function $S$.
Thus, it is enough to  prove that there is a finite partition of
 $\Rr_\infty$ into intervals such that
 the fibers of
 $g'$ over each interval have the same homotopy type
 and that if $ S_{\tcb{\leq} \varepsilon}$
is the locus of $g' \tcb{\leq} \varepsilon$,
there exists 
a finite partition of
 $\Rr_\infty$ into intervals such that
for each interval $I$
the inclusion
$ S_{\tcb{\leq} \varepsilon} \to  S_{\tcb{\leq} \varepsilon'}$,
for $\varepsilon \tcb{<} \varepsilon'$ both in $I$,
is a homotopy equivalence.
But such statements are well-known 
in o-minimal geometry, cf., e.g., \cite{coste} Theorem 5.22.
\eprf

\section{\tcb{The lattice completion}}

\tcb{The previous constructions depended on a canonical map from the stable completion to the Berkovich space.
In this section we will introduce a different and more direct connection between the Berkovich space and 
the stable completion.    Our construction involves a preliminary 
base change to a canonical {\em completion} of the given base, one involving   imaginary
elements from the sort $S_n$, as well as the field sort and $S_1=\G$.}

\tcb{Let $F$ be a  valued field.  The usual completion of $F$ as a valued field can be viewed as a subfield of a maximal immediate extension $F^{max}$
of $F^{\alg}$, consisting of ``rigid'' points, i.e. points invariant under $\Aut(F^{max}/F)$.  
The completion is well-defined
up to a {\em unique} $F$-isomorphism; in particular there is no dependence on the choice of $F^{max}$.   The completion 
 is functorial in   extensions that do not augment the value group at $\infty$.}

\tcb{We now take the sorts $S_n$  into consideration.    Recall the linear topology of \ref{modules}.
Define the {\em lattice completion} \index{lattice completion} $\uF$ of $F$ to consist of the
completion $F^c$ of $F$ in the field sort, and the 
 closure of $S_n(F^c)$ in $S_n(F^{max})$ in the $S_n$-sorts.  \nomenclature{$\uF$}{lattice completion of $F$}
As the linear topology
is Hausdorff, it is clear that points of $S_n(\uF)$ are fixed by $\Aut(F^{max}/F)$.   Thus, up to a unique isomorphism
over $F$, the lattice completion $\uF$ is well-defined and independent of the choice of $F^{max}$.  In fact
it is functorial for extensions that do not augment the value group at $0$ or at $\infty$.}

\tcb{Let $L$ be a valued field.  If $\Lam$ and  $\Lam'$ are two lattices in $L^n$, there exists $M \in \GL_n(L)$ with $M \Lam = \Lam'$;
$\val (\det(M))$ does not depend on the choice of $M$, we call it the {\em relative volume} and denote it by $\vol(\Lam',\Lam)$.  \index{relative volume}
 Thus, if one sets $\vol(\Lam') = \vol(\Lam',\Oo^n)$, we may also write the relative volume as $\vol(\Lam')-\vol(\Lam)$.
We say a family of lattices is {\em directed} (respectively, reverse directed) if any two lattices in the family \index{directed family of lattices}
is contained in (respectively, contain) a third.}

\begin{lem}  \label{lcom1}  \tcb{Let $L$ be a valued field and let $\Lam$ be a rank $n$ sublattice of $L^n$. Consider a 
directed family $(\Lam_i)$ of rank $n$ sublattices of  $\Lam$.  Assume $\vol(\Lam_i)-\vol(\Lam) \to 0$ in $\G(L)$.  Then $\Lam_i \to \Lam$ in $S_n(L)$.}  \end{lem}

\prf   \tcb{Let $w_i$ be the seminorm corresponding to $\Lam_i$,
and $w$ to $\Lam$.  We have to show that for any $v$ we have $w_i(v) \to w(v)$; equivalently
if $w(v)=0$ we have to show that $w_i(v)\to 0$.  
Let $\a_i= \vol(\Lam_i)-\vol(\Lam)$; so $\a_i \to 0$.   We claim that  
$0 \geq w_i(v) \geq   - \alpha_i$ from which the statement follows.  To see this, fix $i$ and set $\Lam_i = \Lam'$, $w_i = w'$, $\alpha_i =\alpha$.
Since $\Lam$ and $\Lam'$ can be simultaneously diagonalized in some basis,  we may
 assume that $\Lam=\Oo^n$ and $\Lam'= \oplus \Oo d_k$. Since $\val(d_k) \geq 0$
 and $\sum_k \val (d_k) = \alpha$, we have $0 \leq \val (d_k) \leq \alpha$. It follows that 
 for any $v$ with $w(v)=0$,  $0 \geq w ' (v) \geq - \alpha$.}   \eprf

\begin{lem}  \label{lcom3} \tcb{Assume $L$ is maximally complete. Consider a family of rank $n$ 
lattices $\Lam_i$ in $L^n$, directed  under inclusion or reverse inclusion.  Assume, for any subspace $U$ of $L^n$, 
that 
$\vol(\Lam_i \meet U) \to \gamma_U \in \G(L)$.  Then there exists a unique rank $n$ lattice 
$\Lam$ in $L^n$ with $\Lam_i \to \Lam$.   Moreover, $\vol(\Lam_i)-\vol(\Lam) \to 0$.} \end{lem}

\prf  \tcb{Uniqueness is clear since the linear topology is Hausdorff.   To prove the remaining assertions,
it suffices, by \lemref{lcom1}, to find a rank $n$ lattice $\Lam$ in $L^n$ such that 
$\vol(\Lam_i,\Lam) \to 0$
 in $\G(L)$, and $\Lam$ contains the $\Lam_i$ in the inclusion case (respectively, in the reverse inclusion case, is contained in the $\Lam_i$).}

  \tcb{Consider first the one-dimensional case.   Then $\Lam_k= \{x: \val(x) \geq \alpha_k\}$ with $\alpha_k \to \gamma$.  
  We set $\Lam = \{x:  \val (x) \geq \gamma\}$, and the statement  is clear.  Note that $\Lam=\meet_k \Lam_k$
  in case the $\Lam_k$ form a descending chain.}
  
  \tcb{Recall that when $V$ is an $n$-dimensional vector space, with dual $V^*$,
  the dual of a lattice $\Lam$ in $V$ is $\Lam^*:=  \{x \in V^*: (\forall y \in \Lam) \val(x \cdot y) \geq 0\}$.
  Beginning with $\Oo^n$ as the standard lattice of $K^n$,
 we take the standard lattice of the dual space  to be the dual lattice of the standard lattice of a given space,
 and the standard lattice of a subspace $U$ to be the intersection with $U$ of the standard lattice, 
 and of a quotient $V/U$ to be the image of the standard lattice.
  Duality reverses inclusion and volume, i.e. $\vol(\Lam_1^*,\Lam_2^*)=- \vol(\Lam_1,\Lam_2)$.  
  Also, $\vol(\Lam^* \meet U^{\perp} ) = -  \vol(\Lam + U) = - \vol(\Lam) +  \vol(\Lam \meet U)$;
  so the convergence assumption goes through to the dual.  
  Thus by passing
  to duals if necessary, 
  it suffices to prove the statement in the case   that the  $\Lam_i$ are reverse directed.}
  
\tcb{In this case, let $\Lam$ be the intersection of all $\Lam_i$.  
We argue first that 
  $\Lam$ spans $L^n$.  Thus fix  a subspace $U$ of $L^n$ of dimension $n-1$; we have to show 
  it does not contain $\Lam$.     By induction, $\Lam_i \meet U $
converges to a lattice $\Lam'$ of $U$.  Modulo $U$, the lattices $\Lam_i + U$ have volume
$\vol(\Lam_i) - \vol(\Lam_i \meet U)$ which   converges to $\g - \vol(\Lam')$, and it follows
that they contain some nonzero element $c+U$.  Now viewing $c+U$ as a coset of $U$ in $L^n$,   maximal completeness implies that $(\meet_i \Lam_i ) \meet (c+U) \neq \varnothing$, 
so as $c+U$ is disjoint from $U$, we see that $\Lam$ is not contained in $U$.    Thus indeed $\Lam$ spans $L$, i.e. $L \Lam = L^n$.    In particular, $\Lam \meet W \neq (0)$ for any one-dimensional $W$ (if $0 \neq  w \in W$, 
then $ w \in c \Lam$ for some $c \in K$ so $c \inv w \in W \meet \Lam$).} 

\tcb{To see that $\Lam$ is a lattice, since $L$ is maximally complete, it suffices to show that
$\Lam \meet W$ is a lattice for any one-dimensional $W$.    This is clear by the one-dimensional case treated
above.}      \eprf

\begin{rem} \label{lcom4} \tcb{The convergence assumption in \lemref{lcom3} holds automatically when  the value group is $\Rr$, provided $|\vol(\Lam_i)|$ is bounded in $\G(L)$.  
Indeed for any subspace $U$,
$\vol(\Lam_i \meet U)$ and $\vol (\Lam_i /U)$ are both   monotone in $\Lam_i$ (e.g. both increasing if
the $\Lam_i$ are increasing) and their sum $\vol(\Lam_i)$  is bounded (in absolute value), hence they both tend towards a real limit value.}
\end{rem}


\section{\tcb{Berkovich points as Galois orbits}}

\tcb{In this section we fix a valued field $F$ with $\val (F) \subset \Rr_{\infty}$.
Let $V$ be an algebraic variety over $F$. For any  
base set $A$ containing $F$ with $\G(A)=\Rr$, we define $B_A(V)$ to be the space of types on $V$ over $A$
that are almost orthogonal to $\G$. We shall be concerned with the case when $A = \uF$.}
   
%
%

%
%
%

 \begin{lem} \label{lcom5} \tcb{Let $V$ be an algebraic variety over $F$. The restriction map $B_{\uF} (V) \to B_{\bF}(V)$ is bijective.}  \end{lem}
 
 \prf  \tcb{Since the mapping $\std{V} (F^{max}) \to B_{\bF}(V)$ factors through the map $\std{V} (F^{max}) \to B_{\uF} (V) \to B_{\bF}(V)$,
surjectivity follows from  \lemref{br1}.
}
 
     \tcb{For injectivity, consider
 $q'$ and $q''$  in $B_{\uF}(V)$ with the same restriction
  to $B_{\bF} (V)$.  
 Let $c' \models q'$ and $c'' \models q''$.   We can find models $M'$ and $M''$  containing respectively  $\uF(c')$ and $\uF(c'')$, and with  value group $\Rr$.  
 Let $K'$  and $K''$ be  maximally complete algebraically closed
 valued fields, with value group $\Rr$,
 containing respectively $\uF(c')$ and $\uF(c'')$. By enlarging one of them, we may assume their residue fields have the same transcendence degree over the residue field of $F$, so that they are isomorphic over $\bF$.  Since
 $\tp(c/\bF)=\tp(c'/\bF)$, there exists an $\bF$-isomorphism $\beta: K' \to K''$ with $\beta(c)=c'$. But 
 $\beta | \uF$
 must be the identity, since each point of $\uF$ is the unique limit point of some sequence (or net) of
 elements of $\bF$ or $S_n(\bF)$.  Thus $q'=\tp(c'/\uF)=\tp(c''/\uF)=q''$.}  
 \eprf

\begin{lem} \label{lcom9}\tcb{Let $V$ be an algebraic variety over  $F$.  The natural map $\std{V}(\uF) \to B_{\uF} (V)$ sending $q$ to $q|\uF$ is
surjective.  Hence when $F$ is algebraically closed, it is bijective.}
\end{lem}

\prf   \tcb{We may assume $F$ is complete as a valued field. Let $p \in B_{\uF} (V)$. We may assume that $V$ is affine
and that $p$ is Zariski dense in $V$.  We first give the argument assuming $\val(F)$ is dense in $\Rr$; the general case differs
only notationally, and will be explained below.}

 \tcb{Let $H$ be the affine coordinate ring of $V$, $H=\union_{d \geq 0} H_d$ where 
  $H_d$ is the space of polynomials of degree at most $d$, modulo those vanishing on $V$.  
 Let $M_d = \{f \in H_d(F):  (\val (f) \geq 0) \in p\}$ and let $\mathcal{F}_d$ 
   be the family of all lattices of $H_d$, generated by a finite subset of $M_d$.  
    We view it as ordered
   by inclusion.}  
   
  \tcb{We wish to show that $\mathcal{F}_d$ admits a limit lattice containing it.  For this purpose we may 
   replace  $\mathcal{F}_d$ by the subfamily of elements  of  $\mathcal{F}_d$ containing
   some fixed such lattice $\Lam_0$; then this family has a lower bound.  Let us prove the existence of an upper bound.
The category of complete $\Rr$-valued fields  admits algebraically independent 
   amalgamation; the one-dimensional case is easy, and the general case follows inductively using \cite{frob} Lemma 6.18.
   Hence,     using the
    Zariski density of $p$, one can find realizations $c_1,\ldots,c_N$ of $p$ in some $\Rr$-valued field, 
 where $N=\dim(H_d)$, 
   such that no nonzero element of $H_d$ vanishes on all the $c_i$.  Then $(c_1,\ldots,c_N)$
   generate a lattice in $H_d^*$, whose dual lattice contains $M_d$ and hence all elements of $\mathcal{F}_d$.
   In particular, $\vol (\Lam)$ is bounded above and below, for $\Lam \in \mathcal{F}_d$.  
By \lemref{lcom3}  and \remref{lcom4},  $\mathcal{F}_d$
 has a unique limit lattice $\Lam_d$ in $F^{max}$.    A code for $\Lam_d$ lies in $\uF$ by
 definition of the latter.}  
  
\tcb{Let $w_d$ be the seminorm on $H_{d}$ corresponding to $\Lam_d$. 
We now prove that  these seminorms satisfy the condition in
 \remref{rem:5.1.5}, namely that
for any $f_i \in  H_{d_1}, f_2 \in  H_{d_2} $ we have $w_{d_1+d_2}(f_1f_2) = w_{d_1}(f_1)+w_{d_2}(f_2)$, for every $d_1$ and $d_2$.
Indeed, find a sequence (or net) of lattices $\Omega_j \in \mathcal{F}_{d_1+d_2}$ 
 such that $\Omega_j \to \Lambda_{d_1+d_2}$ while $\Omega_j \meet H_{d_i} \to \Lambda_{d_i}$.  Let $w'_j$ be the seminorm 
 corresponding to $\Omega_j$; then $w'_j(f_1f_2) \to w_{d_1+d_2}(f_1f_2)$ while $w'_j(f_i) \to w_{d_i}(f_i)$, and the condition follows
 by continuity.  
One deduces from \remref{rem:5.1.5} that there exists a unique $q$ in $\std{V}(\uF)$ such that $J_d(q)=\Lam_d$, for every $d$,
in the notation of \ref{modules}.} 

\tcb{Note that $\Lam_d \meet H_d(F)=M_d$.  Indeed by definition we have $M_d \subset \Lam_d$.  On the other hand
if $f \in H_d(F)$ and   $f \notin M_d$ then  $p(x) \vdash \val f(x) \leq   \alpha $ for some negative $\alpha \in   \Qq$,
and it follows by continuity that $w_d(f) \leq \alpha$ so $f \notin \Lam_d$.  Thus $q|B_F = p|B_F$.  In fact  since $\val(F)$ is dense in $\Rr$
it follows that 
  $q|B_{\bF}= p$, and the statement follows.}
  
  \tcb{When $\val(F)$ is not necessarily dense in $\Rr$, we must define $\mathcal{F}_d$ as the family of all lattices
  $\oplus_{i=1}^N \alpha_i \Oo f_i$, where $f_1,\ldots,f_N$ is a basis for $H_d(F)$, $p(x) \vdash \val(f_i(x))=- \alpha_i$,
  and $\gamma \Oo := \{x : \val(x) \geq \gamma\}$.  These lattices are still defined over $\bF$ (though they may not have a basis
  in $F$), and the proof goes through as before.}
   \eprf

\begin{rem} \label{lcom11} \tcb{At this point we   recover  the functorial base change of \cite{poi-ang}.  Namely when $F$ is algebraically closed, it follows from 
\lemref{lcom5} and \lemref{lcom9} that we have a canonical bijection between
$B_{\bF} (V)$ and $ \std{V}(\uF)$  and thus, by canonical extension of stably dominated types, 
we get a canonical map 
from $B_{\bF} (V)$ to $\std{V}(\Uu)$.
   In fact we discovered the approach of this section 
upon  reflecting upon  Poineau's theorem.  
Another proof was independently given 
in \cite{by}.}   \end{rem}

\medskip

\tcb{Let 
$V$ be a   variety over $F$. 
Let $G=\Aut(F^{\mathrm{\alg}}/F^{\mathrm h})$ be the absolute Galois group of the Henselization $F^{\mathrm h}$  of $F$; 
so $G$ acts also on $\overline{F^{\alg}}$.  
Recall that  $G$ is the group of 
valued field automorphisms of $F^{\mathrm{\alg}}$ over $F$, 
and $B_{\bF}(V) = B_{\mathbf F^{\mathbf{\alg}}}(V) /G$, 
with $\mathbf F^{\mathbf{\alg}}$ the structure $(F^{\mathrm{\alg}}, \Rr)$.
Composing the maps of \lemref{lcom5} and \lemref{lcom9}
one gets a $G$-equivariant bijection
\begin{equation*}\std{V}(\overline{F^{\mathrm{\alg}}})\longrightarrow  
   B_{\mathbf F^{\mathrm{\alg}}}(V)
\end{equation*} 
which one checks to be a homeomorphism and
whose inverse induces a 
natural   homeomorphism}
\begin{equation*}\tcb{\varrho:
B_{\bF}(V) \longrightarrow \std{V}(\overline{F^{\mathrm{\alg}}})/G.}
\end{equation*}

\begin{defn}\label{abdef}\tcb{We call a point of 
$B_{\bF}(V)$  
which restricted to $F$  determines
 an Abhyankar extension of
 the valued field $F$ an {\em Abhyankar point} of $B_{\bF}(V)$.} \index{Abhyankar point}
\end{defn}

\begin{prop}\label{ncom13}  \tcb{The homeomorphism $\varrho$ induces a bijection between Abhyankar points of $B_{\bF}(V)$  
and  $G$-orbits of strongly stably dominated points of 
$\std{V}(\overline{F^{\alg}})$.}
\end{prop}

\prf  \tcb{Let $p$ be an Abhyankar point in $B_{\bF}(V)$. Let
$c \models p| F$ and  $d=\mathrm{trdeg}_F F(c)$. Then there exist $F$-definable functions $f_1,\ldots,f_k$ and $g_1,\ldots,g_\ell$, with $k+\ell=d$,
such that $f_1(c),\ldots,f_k(c)$ are elements of the residue field, algebraically independent over the residue
field of $F$ and  $g_1(c),\ldots,g_\ell (c)$ are elements of the value group, linearly independent over $\val(F)$.
Now over $\bF$, $f_1(c),\ldots,f_k(c),g_1(c),\ldots,g_\ell(c)$ are algebraically independent
elements of $\RES_{\bF}$, hence $\tp(c/\bF)$ is strongly stably dominated. For the other direction see the proof of \thmref{sc}.} 
\eprf

\begin{rem}\label{lastrem}\tcb{The direct connection between Berkovich points and stably dominated points over the lattice completion provided by the homeomorphism $\varrho$
immediately yields another proof of \thmref{sc}.   Simply, given $V$ over $F$ and data over $\bF$, 
the homotopy of \thmref{1} is defined over $\bF$; and being $\bF$-definable, in particular $\uF$-definable,
it takes points of $\std{V}(\uF)$ to points of $\std{V}(\uF)$, so it restricts to a homotopy on $\std{V}(\uF)$;
the  isomorphism $\varrho$ translates this to a homotopy on $B_\bF(V)$.} 

\tcb{Moreover, by \thmref{sgc} and \propref{ncom13}, 
 one may ask that the homotopy fix any given Abhyankar point $p$ of $B_\bF(V)$.
Thus the proof of \thmref{lc}, 
along with the analogous fact in the o-minimal case, shows that $p$ admits a basis of neighborhoods  that
 strongly retract to $p$.}
\end{rem}



\backmatter

\cleardoublepage

 \printindex
 
 \cleardoublepage
 
\makeatletter
\immediate\closeout\@nomenclaturefile
\makeatother
\def\chopit#1@[#2]\begingroup #3\nomeqref#4|nompageref{\item {#2}, #3}
\def\nomenclatureentry#1#2{\chopit#1, #2}
 \begin{thenomenclature}
\rightskip0pt plus 4em\itemindent-2em\leftskip2em
\input \jobname.nlo
\end{thenomenclature}

 \end{document}